

 
%
\input amssym.def
%
%
\magnification=\magstephalf 
\hsize=5.8 true in        
\vsize=8.9 true in        
\hoffset=.6 true in       
\baselineskip=24pt        
\parindent=2.3em          
\def\beginchapter{\vfill\eject\topglue 1in \mark{}\mark{2}} 
\headline={\if\firstmark 2 {\tenrm\hfil\folio} \fi}
\footline={\if\firstmark 2 {\hfil} \else {\tenrm\hfil\folio\hfil} \fi}
\def\chapter#1{\beginchapter\centerline{{\bf\uppercase{CHAPTER {#1}}}}}
\def\headspace{\bigskip}  
\def\preheadspace{\headspace\penalty-200}                                                  
\def\postheadspace{\nobreak\headspace}                                            
\def\chaptertitle#1{\centerline{{\bf#1}}}
\def\firsthead#1{\centerline{{\bf#1}}}
\def\secondhead#1{\leftline{{\bf#1}}}

%
%
%
\def\sc#1{{\sevenrm\uppercase{#1}}} 
%
\def\xproclaim#1. #2\par{\medbreak{\bf#1.\enspace}{\sl#2}\par
                         \ifdim\lastskip<\medskipamount\removelastskip\penalty55\medskip\fi}
\def\definition#1\par{\medbreak D\sc{efinition}.\enspace{\sl#1}\par
                         \ifdim\lastskip<\medskipamount\removelastskip\penalty55\medskip\fi}
\def\proof{{\bf Proof.}\enspace}  
\def\remark{R\sc{emark}.\enspace}
\def\xremark#1{R\sc{emark}\enspace#1.\enspace}
\def\QED{\vbox{\hrule\hbox{\vrule height 6pt \kern 6pt \vrule}\hrule}}
%
\def\xcite#1{[{\bf{#1}}]} 
\def\xciteplus#1#2{[{\bf{#1}},{ #2}]} 
%
\def\xdim{\mathop{\sl dim\/}\nolimits} 
\def\xlim{\mathop{\sl lim\/}\nolimits} 
\def\xmax{\mathop{\sl max\/}\nolimits} 
\def\xmin{\mathop{\sl min\/}\nolimits} 
\def\xsup{\mathop{\sl sup\/}\nolimits} 
\def\xinf{\mathop{\sl inf\/}\nolimits} 
\def\supp{\mathop{\rm supp}\nolimits} 
\def\xsupp{\mathop{\sl supp\/}\nolimits} 
\def\dist{\mathop{\rm dist}\nolimits} 
 
%

%
\def\TAG#1#2{(\vphantom{#1}#2)}
%
%
%
%
\def\pz{\phantom0}
\def\tint{\textstyle\int}
\def\dint{\displaystyle\int}
\def\tsum{\textstyle\sum}
\def\tsuml{\textstyle\sum\limits}
\def\NN{{\Bbb N}}

\def\RR{{\Bbb R}}

\def\FF{{\Bbb F}}
\def\TT{{\Bbb T}}
\def\l#1{\ell_{\vphantom{N}}^{#1}}                                     
\def\L#1{L_{\vphantom{N}}^{#1}}                                       
\def\Lz#1{L_0^{#1}}                                       
\def\SL#1{{\cal L}_{#1}}                              
\def\Xpw{X_{p,w}}                                         
\def\finitel#1#2{\ell^{#1}_{#2}}
\def\paren#1{\left({#1}\right)}
\def\vector#1{\left\langle{#1}\right\rangle}
\def\set#1{\left\{{#1}\right\}}                        
\def\setwlimits#1#2#3{{\set{#1}}_{#2}^{#3}}           
\def\span#1{{\left[#1\right]}}                       
\def\spansub#1#2{{\span{#1}}_{#2}}                    
\def\class#1#2{{\span{#1}}_{#2}}                    
\def\abs#1{\left|{#1}\right|}
\def\abspower#1#2{{\abs{#1}}^{#2}}            
\def\norm#1{\left\|{#1}\right\|}                     
\def\normsub#1#2{{\norm{#1}}_{#2}}                 
\def\normsubpower#1#2#3{{\norm{#1}}_{#2}^{#3}}      
\def\xnorm#1{\left| \mkern-1.5mu \norm{#1} \mkern-1.5mu \right|}
\def\cequiv{\mathrel{\mathop{\equiv}\nolimits_c}}
\def\injects{\hookrightarrow}

\def\cinjects{\buildrel\rm c\over\hookrightarrow}

\def\Tmapsto#1{\buildrel#1\over\mapsto}
\def\within#1#2{\mathrel{\mathop{\approx}\limits^{#1}_{#2}}}
\def\finitespan{\mathop{{\rm span}}\nolimits}
\def\xfinitespan{\mathop{{\sl span\/}}\nolimits}
\def\seqsum#1#2{\paren{{#1}\oplus{#1}\oplus\cdots}_{#2}}    
\def\sumsum#1#2{\paren{\textstyle\sum\nolimits^{\oplus} #1}_{#2}}
\def\xsubpower#1#2{\!{{}\atop{}_{#1}^{#2}}} 
\def\normalize#1#2#3{{#1\over\norm{#1}}\xsubpower{#2}{#3}} 
\def\sumprime{\sum\nolimits^{\prime}} 
                            
\def\ratio{{\normsub{f_n}{2}}\big/{\normsub{f_n}{p}}}       
\def\wdef{w=\set{w_n}=\set{\ratio}}                       
\def\wstopower{{w_n}^{2p \over p-2}}                         
%
%
\def\Rpa{R^p_{\alpha}}
\def\Rpb{R^p_{\beta}}
\def\Had{H_{\alpha}^{\delta}}
\def\psum#1{\paren{#1\oplus#1}_p}
\def\brackpsum#1{\span{#1\oplus#1}_p}
\def\Ipsum#1#2{\paren{\textstyle\sum_{#2}^{\oplus} #1}_{{\rm Ind},p}}
\def\con{\mathchoice{\kern 1.0 pt\raise 2.1 pt\hbox{$\cdot$}}
                    {\kern 1.0 pt\raise 2.1 pt\hbox{$\cdot$}}
                    {\kern 0.9 pt\raise 1.1 pt\hbox{$\scriptstyle\cdot$}\kern -.2 pt}
                    {\kern 0.4 pt\raise 1.0 pt\hbox{$\scriptscriptstyle\cdot$}\kern -.1 pt}}
\def\dotprec{\mathrel{\dot{\prec}}}
\def\dotpreceq{\mathrel{\dot{\preceq}}}
\def\triless{\mathrel{\triangleleft}}
\def\distiso{\buildrel{\rm dist}\over\sim}
\def\S#1{{\cal#1}}
\def\E{\S{E}}
\def\T{\S{T}}
\def\B{\S{B}}

\def\SD{{\cal D}}

%
\def\ra{\rightarrow}

\def\usa{\nearrow}
\def\dsa{\searrow}
\def\ua{\uparrow}
\def\da{\downarrow}
\def\cra{\buildrel\rm c\over\ra}

\def\cusa{\buildrel\rm c\over\usa}
\def\cdsa{\buildrel\rm c\over\dsa}
\def\cua{\mathrel{\ua\scriptstyle{{\rm c}}}}
\def\cda{\mathrel{\da\scriptstyle{{\rm c}}}}
\def\Tra#1{\buildrel#1\over\ra}

\def\Tdsa#1{\buildrel#1\over\dsa}

\def\TRua#1{\ua{\scriptstyle#1}}
\def\TLda#1{{\scriptstyle#1}\da}
\def\TRda#1{\da{\scriptstyle#1}}
\def\vequiv{\scriptstyle{|||}}
\def\lp{\ell^p}
\def\ltwo{\ell^2}
\def\lpsumltwo{\sumsum{\ltwo}{\lp}}
%
%
\begingroup               
\nopagenumbers            
\baselineskip=12pt        
\def\ctr#1{\centerline{#1}}
\def\bsk{\bigskip}
%
\topglue 1 true in
{
\ctr{{\bf CONSTRUCTIONS OF ${\cal L}_p$ SPACES}}
\bsk
\ctr{{\bf FOR $p\in(1,\infty)\setminus\{2\}$}}
}
\vfil 
{
\ctr{By}
\bsk
\ctr{GREGORY MICHAEL FORCE}
\bsk
\ctr{Bachelor of Arts}
\ctr{Coe College}
\ctr{Cedar Rapids, Iowa}
\ctr{1977}
\bsk   
\ctr{Master of Science}
\ctr{Purdue University}
\ctr{West Lafayette, Indiana}
\ctr{1981}
}
\vfil 
{
\ctr{Submitted to the Faculty of the Graduate College}
\ctr{of the Oklahoma State University}
\ctr{in partial fulfillment of}
\ctr{the requirements for}
\ctr{the Degree of}
\ctr{DOCTOR OF EDUCATION}
\ctr{December, 1995}
}
\vglue 1 true in
\pageno=-1
\eject
\endgroup
%
%
\begingroup               
\baselineskip=12pt        
\def\signitureline#1{\vskip 12pt\centerline{\hbox to 4 true in {\hrulefill}}\smallskip\centerline{{#1}\vphantom{Sign}}}
%
%
\topglue 1 true in
{
\centerline{{\bf CONSTRUCTIONS OF ${\cal L}_p$ SPACES}}
\bigskip
\centerline{{\bf FOR $p\in(1,\infty)\setminus\{2\}$}}
}
\vfil
{
\centerline{\hbox to 4 true in {Thesis Approved:\hss}}
\vglue 24pt 
\signitureline{Thesis Advisor}
\signitureline{}
\signitureline{}
\signitureline{}
\signitureline{}
\signitureline{Dean of the Graduate College}
}
\vglue 2 true in
\eject
\endgroup 
%
%
\begingroup               
\raggedright
\topglue 1 true in        
\centerline{{\bf ACKNOWLEDGEMENTS}}
\bigskip
 I would like to thank the members of my committee: 
 John Wolfe, who has been a friend,
 Alan Noell and David Ullrich, who were influential especially in the early part of my graduate career, 
 David Webster, who helped make the education component of my degree more agreeable, 
 and of course my thesis advisor, Dale Alspach. 
 Dale has shown quite a lot of patience.
 I could not have taken the topic of my thesis as far as\break 
 I did without his considerable help and insight. 

 Finally, I would like to thank my parents, Ron and Winnie Force, who have\break 
 always been supportive of my academic pursuits.

\vfill
\eject
\endgroup
%
%
\begingroup               
\parindent=3.2em          
\baselineskip=12pt        
%
\def\xspace{\item{}{}}                      
%
\def\adotfill{\leaders \hbox to .5em {\hss.\hss} \hfill} 
%
\def\pagebox#1{\hbox to 2em {\hss#1}}
%
\def\chap#1. #2...#3:{\xspace\item{#1.}#2 \adotfill\pagebox{#3}\xspace} 
%
\def\head#1...#2:{\item{}\phantom{THE }#1 \adotfill\pagebox{#2}} 
%
\def\sub#1...#2:{\item{}\phantom{THE The }#1 \adotfill\pagebox{#2}} 
%
\def\otherline#1...#2:{\line{#1 \adotfill\pagebox{#2}}}
%
\topglue 1 true in        
\centerline{{\bf TABLE OF CONTENTS}}
\vglue 27pt               
\line{Chapter \hfil Page}
\vglue 5pt                
\chap I. INTRODUCTION...1:
         \head Preliminaries for $\SL{p}$ Spaces...1:  
               \sub The $\SL{p}$ Spaces...1:  
               \sub The Relations $\injects$ and $\cinjects$...3:  
               \sub The Classical $\SL{p}$ Spaces...4:  
               \sub Elementary Constructions...5:  
         \head Preliminaries for Banach Spaces...6:
         \head Overview of Chapters...10:
\chap II. THE NONCLASSICAL $\SL{p}$ SPACES OF ROSENTHAL...11: 
         \head The Space $X_p$...11:  
               \sub The Space $\Xpw$...11:
               \sub Rosenthal's Inequality...14:
               \sub The Complementation of $\Xpw$ in $\L{p}$...15:
               \sub The Mutual Isomorphism of the Spaces $\Xpw$...20:
               \sub The Isomorphism Type of $X_p$...27:
               \sub Complementation and Imbedding Relations for $X_p$...40:
         \head The Space $B_p$...43:
               \sub The Space $X_{p,v^{(N)}}$...43:
               \sub The Space $B_p$...44:
               \sub Sums of $B_p$...51:
         \head Sums Involving $X_p$ or $B_p$...54:
         \head Concluding Remarks...66:
\chap III. THE TENSOR PRODUCT CONSTRUCTION OF SCHECHTMAN...68:
         \head The Tensor Product Construction...68:  
         \head The Isomorphic Distinctness of $X_p^{\otimes m}$ and $X_p^{\otimes n}$...73:  
         \head The Sequence Space Realization of $X_p^{\otimes n}$...80:  
\chap IV. THE INDEPENDENT SUM CONSTRUCTION OF ALSPACH...87:
         \head The Independent Sum $\sumsum{X_i}{I,w}$...88:
         \head The Complementation of $\sumsum{X_i}{I,w}$ in $\L{p}(\Omega)$...89:
         \head Independent Sums with Basis...96:
         \head The Independent Sum $\sumsum{X}{I}$...99:
         \head The Space $D_p$...108:
         \head Sums Involving $D_p$...115:
\chap V. THE CONSTRUCTION AND ORDINAL INDEX OF BOURGAIN, \hfil\break ROSENTHAL, AND SCHECHTMAN...124:
        \head Preliminaries...124:
        \head The Ordinal Index...126:
               \sub A General Ordinal Index $h$...126:
               \sub Motivation from $\L{p}$...127:
               \sub The Space $\paren{\overline{B}^{\delta},\prec}$...129:
               \sub A Characterization of $\L{p}\injects B$...130:
               \sub The Ordinal Index $h_p(\delta,\pz)$...131:
               \sub The Ordinal Index $h_p$...136:
         \head The Disjoint and Independent Sum Constructions...137:
               \sub The Interaction of the Constructions and the Ordinal Index...141:
         \head The Complementation of $\Rpa$ in $\L{p}$...146:
               \sub Preliminaries...146:
               \sub The Isomorphism of $Z_{\NN}^p$ and $\L{p}$...148:
               \sub The Complementation of $\Rpa$ in $Z_{\NN}^p$...155:
         \head Concluding Remarks...158:
\xspace
\bigskip
\otherline BIBLIOGRAPHY...160:
\bigskip
\vfill
\eject
\endgroup
%
%
%
\pageno=1
\begingroup 
\raggedright
%
\chapter{I}
\headspace
\chaptertitle{INTRODUCTION}
\headspace

 The $\SL{p}$ spaces are defined in terms of their finite-dimensional subspaces. 
 However, in the category of separable infinite-dimensional Banach spaces, 
 the $\SL{p}$ spaces for $1<p<\infty$ with $p\ne2$ are those spaces
 which are isomorphic to complemented subspaces of $\L{p}$, 
 but not isomorphic to the Hilbert space $\l{2}$. 

 Rosenthal \xcite{RI}, Schechtman \xcite{S}, Alspach \xcite{A}, and Bourgain \xcite{B-R-S} have \hfil\break 
 developed methods of constructing $\SL{p}$ spaces for $1<p<\infty$ with $p\ne2$
 which have a probabilistic aspect. 
 These methods have enlarged the set of known $\SL{p}$ spaces from the classical examples
 [$\l{p}$, $\l{2}\oplus\l{p}$, $\seqsum{\l{2}}{\l{p}}$, and $\L{p}$]
 to a family indexed by the countable ordinals. 
 We will examine these constructions, 
 provide some details,
 clarify a few points, 
 and to some extent interrelate the constructed spaces with respect to the relation $\cinjects$. 

\preheadspace
\firsthead{Preliminaries for $\SL{p}$ Spaces}
\headspace
\secondhead{The $\SL{p}$ Spaces}
\postheadspace
 The $\SL{p}$ spaces were introduced by Lindenstrauss and Pe\char32lczy\'nski in \xcite{L-P},
 and were studied further by Lindenstrauss and Rosenthal in \xcite{L-R}.   
 The definition and some basic results are presented below.    

\definition Let $1\le p\le\infty$ and $1\le\lambda<\infty$.  
            A Banach space $X$ is an $\SL{p,\lambda}$ space if 
            for each finite-dimensional subspace $Z$ of $X$, 
            there is a finite-dimensional subspace $Y$ of $X$ 
            containing $Z$ such that $d(Y,\finitel{p}{n})\le\lambda$,  
            where $n=\dim(Y)$ and $d(Y,\finitel{p}{n})$ is the \hfil\break 
            Banach-Mazur distance between $Y$ and $\finitel{p}{n}$.  
            Finally, a Banach space is an $\SL{p}$ space if \hfil\break 
            it is an $\SL{p,\gamma}$ space for some $1\le\gamma<\infty$.   

 Let $1<p<\infty$ where $p\ne2$.   
 In \xciteplus{L-P}{Example 8.2}, it is shown that   
 $\l{p}$, $\l{2}\oplus\l{p}$, $\seqsum{\l{2}}{\l{p}}$, and $\L{p}$   
 are mutually nonisomorphic $\SL{p}$ spaces,   
 although this is more easily seen in light of the subsequent results of \xcite{L-R}.
 These spaces are the classical $\SL{p}$ spaces.   

 Let $X$ be a Banach space. 
 A bounded linear mapping $P:X\to X$ is called a \hfil\break 
 projection if $P^2=P$. 
 Let $Y$ be a closed subspace of $X$. 
 Then $Y$ is called a \hfil\break 
 complemented subspace of $X$
 if there is a (bounded linear) projection $P:X\to X$ mapping $X$ onto $Y$. 
 If $Y$ is a complemented subspace of $X$, 
 $P:X\to X$ is the \hfil\break 
 (bounded linear) projection mapping $X$ onto $Y$, 
 and $Z$ is the null space of $P$, 
 then $X=Y\oplus Z$. 
 Conversely, if $X=Y\oplus Z$ for some closed subspace $Z$ of $X$, 
 then $Y$ is a complemented subspace of $X$ (as is $Z$). 

 We will restrict our attention to separable infinite-dimensional $\SL{p}$ spaces for \hfil\break 
 $1<p<\infty$ with $p\ne2$.    
 For these spaces, \xcite{L-P} and \xcite{L-R} each contribute one \hfil\break 
 implication in the following characterization, but in greater generality.    

\xproclaim {Theorem 1.1}. Let $1<p<\infty$ where $p\ne2$,   
                          and let $X$ be a separable infinite-dimensional Banach space.    
                          Then $X$ is an $\SL{p}$ space if and only if    
                          $X$ is isomorphic to a complemented subspace of $\L{p}$ but $X$ is not isomorphic to $\l{2}$.  

 The essence of the forward implication \xciteplus{L-P}{Theorem 7.1} is the following.

\xproclaim {Proposition 1.2}. Let $1<p<\infty$ and let $X$ be an $\SL{p}$ space.
                              Then $X$ is \hfil\break 
                              isomorphic to a complemented subspace of $\L{p}(\mu)$ for some measure $\mu$.   

\remark In the above proposition, analogous statements for $p=1$ and $p=\infty$ are false.  
        For $p=1$, \xciteplus{L-P}{Example 8.1} provides a counterexample. 
        For $p=\infty$, any separable infinite-dimensional $C(K)$ space provides a counterexample, as noted in \xcite{L-P}.
        However, by \xciteplus{L-P}{Corollary 2 of Theorem 7.2},
        if $X$ is an $\SL{1}$ space, then $X$ is isomorphic to a subspace of $\L{p}(\mu)$ for some measure $\mu$. 

 The essence of the reverse implication \xciteplus{L-R}{Theorem 2.1} is the following.    

\xproclaim {Proposition 1.3}. Let $1<p<\infty$ and let $X$ be (isomorphic to)
                              a complemented subspace of $\L{p}(\mu)$ for some measure $\mu$.  
                              Then either $X$ is an $\SL{p}$ space or $X$ is \hfil\break 
                              isomorphic to a Hilbert space.   

\remark In the above proposition,
        modified versions hold for $p=1$ and $p=\infty$ \xciteplus{L-R}{Theorem 3.2}.   
        If $X$ is (isomorphic to) a complemented subspace of $\L{1}(\mu)$ for some measure $\mu$,  
        then $X$ is an $\SL{1}$ space. 
        If $X$ is (isomorphic to) a complemented \hfil\break 
        subspace of a $C(K)$ space,  
        then $X$ is an $\SL{\infty}$ space. 

 Let us assume the hypotheses of Theorem 1.1.   
 The hypothesis that $X$ is infinite-dimensional excludes a class of spaces which are trivially $\SL{p}$.   
 The hypothesis that $X$ is separable allows us to replace the $\L{p}(\mu)$ of Proposition 1.2 by $\L{p}=\L{p}(0,1)$.  
 As noted in \xcite{L-P} and \xcite{L-R},   
 the $\SL{2}$ spaces are precisely the spaces which are isomorphic to Hilbert spaces.   
 However, the only separable infinite-dimensional Hilbert space \hfil\break 
 (up to isometry) is $\l{2}$.   
 Thus we may replace the Hilbert space of Proposition 1.3 by $\l{2}$.   
 The conclusion of Theorem 1.1 now follows.   

\preheadspace
\secondhead{The Relations $\injects$ and $\cinjects$}
\postheadspace

 Let $X$ and $Y$ be Banach spaces.
 We write $X\injects Y$ if $X$ is isomorphic to a closed subspace of $Y$.  
 We write $X\cinjects Y$ if $X$ is isomorphic to a complemented subspace of $Y$.  
 Of course if $X\cinjects Y$, then $X\injects Y$. 
 If $X\cinjects Y$, then $X^*\cinjects Y^*$. 
 However if $X\injects Y$, it does not follow that $X^*\injects Y^*$. 
 If $X$ is a closed subspace of $Y$ with $X\cinjects Y$, 
 it does not follow that $X$ itself is a complemented subspace of $Y$. 
 The relations $\injects$ and $\cinjects$ are reflexive and transitive, but not antisymmetric.   

 We write $X\equiv Y$ if $X\injects Y$ and $Y\injects X$. 
 We write $X\cequiv Y$ if $X\cinjects Y$ and \hfil\break 
 $Y\cinjects X$.    
 We write $X\sim Y$ if $X$ is isomorphic to $Y$. 
 The relations $\equiv$, $\cequiv$, and $\sim$ \hfil\break 
 are equivalence relations. 
 Let $\class{\phantom{0}}{\sim}$, $\class{\phantom{0}}{\cequiv}$, and $\class{\phantom{0}}{\equiv}$ 
 denote equivalence classes under \hfil\break 
 $\sim$, $\cequiv$, and $\equiv$, respectively.  
 Then $\class{X}{\sim}\subset\class{X}{\cequiv}\subset\class{X}{\equiv}$.  

 If $X\equiv X'$ and $Y\equiv Y'$, then $X\injects Y$ if and only if $X'\injects Y'$.  
 Similarly, 
 if $X\cequiv X'$ and $Y\cequiv Y'$, then $X\cinjects Y$ if and only if $X'\cinjects Y'$.  
 Thus $\injects$ and $\cinjects$ induce partial orderings on equivalence classes under $\equiv$ and $\cequiv$, respectively. 

\preheadspace
\secondhead{The Classical $\SL{p}$ Spaces}
\postheadspace

 Let $2<p<\infty$.  
 Then $\l{2}$ and the classical separable infinite-dimensional $\SL{p}$ spaces  
 are related by $\injects$ as in diagram \TAG{1.1}{1.1} below, where  
 $X\to Y$ denotes $X\injects Y$ but $Y\not\injects X$,
 $X\equiv Y$ denotes $X\injects Y$ and $Y\injects X$,
 and the absence of a relation symbol between $X$ and $Y$ implies 
 $X\not\injects Y$ and $Y\not\injects X$,
 unless some relation is implied by the transitivity of $\injects$. 
 The same conventions will apply in future diagrams relating spaces by $\injects$. 

$$\matrix{\ltwo&&&&&&\cr
          &\dsa&&&&&\cr
          &&\ltwo\oplus\lp&\ra&\seqsum{\l{2}}{\l{p}}&\ra&\L{p}.\cr
          &\usa&&&&&\cr
          \lp&&&&&&\cr}\eqno{\TAG{1.1}{1.1}}$$     

 Let $1<p<\infty$ where $p\ne2$.  
 Then $\l{2}$ and the classical separable infinite-{dimensional} $\SL{p}$ spaces  
 are related by $\cinjects$ as in diagram \TAG{1.2}{1.2} below.   
 Conventions \hfil\break 
 analogous to those described above will apply in this and in future diagrams relating spaces by $\cinjects$      
 (with $\cinjects$, $\cra$, and $\cequiv$ replacing $\injects$, $\ra$, and $\equiv$, respectively).     

$$\matrix{\ltwo&&&&&&\cr
          &\cdsa&&&&&\cr
          &&\ltwo\oplus\lp&\cra&\seqsum{\l{2}}{\l{p}}&\cra&\L{p}.\cr
          &\cusa&&&&&\cr
          \lp&&&&&&\cr}\eqno{\TAG{1.2}{1.2}}$$     

 The positive relations asserted to exist above follow routinely from well-known results.  
 Of course $\l{2}\cinjects\l{2}\oplus\l{p}$ and $\l{p}\cinjects\l{2}\oplus\l{p}$.   
 Letting $\Bbb F$ denote the scalar field,   
 $$\l{2}\oplus\l{p}\sim\l{2}\oplus\seqsum{\Bbb F}{\l{p}}\cinjects\l{2}\oplus\seqsum{\l{2}}{\l{p}}\sim\seqsum{\l{2}}{\l{p}}.$$  
 Khintchine's inequality \xciteplus{W}{I.B.8} for the  Rademacher functions $\set{r_n}$ shows that \hfil\break 
 $\spansub{r_n}{\L{p}}\sim\l{2}$.   
 Moreover, for $2<p<\infty$, the orthogonal projection of $\L{p}$ onto $\spansub{r_n}{\L{p}}$ is bounded.      
 Hence for $2<p<\infty$, and for $1<p<2$ by duality, $\l{2}\cinjects\L{p}$.   
 It follows that      
 
 $$\seqsum{\l{2}}{\l{p}}\cinjects\seqsum{\L{p}}{\l{p}}\sim\L{p}.$$

 Some of the the negative results are another matter, although  
 $\l{2}\not\injects\l{p}$,   
 $\l{p}\not\injects\l{2}$,   
 $\l{2}\oplus\l{p}\not\injects\l{2}$, and    
 $\l{2}\oplus\l{p}\not\injects\l{p}$,   
 all follow from the fact that
 $\l{r}\not\injects\l{s}$
 for $r,s\in[1,\infty)$ with $r\ne s$. 
 The fact that $\seqsum{\l{2}}{\l{p}}\not\injects\l{2}\oplus\l{p}$ for $2<p<\infty$
 is \xciteplus{RI}{Lemma for Corollary 14}, presented below as Lemma 2.23.  
 The fact that $\L{p}\not\injects\seqsum{\l{2}}{\l{p}}$ for $2<p<\infty$ is \xciteplus{L-P 2}{Theorem 6.1}.  

\preheadspace
\secondhead{Elementary Constructions}
\postheadspace

 Fix $1<p<\infty$ where $p\ne2$.  

 Let $X$ and $Y$ be separable infinite-dimensional Banach spaces such that  
 $X\cinjects\L{p}$ and $Y\cinjects\L{p}$.   
 Then $X\oplus Y\cinjects\L{p}\oplus\L{p}\sim\L{p}$.   
 Note that since $\l{2}$ is prime, if $X\not\sim\l{2}$ and $Y\not\sim\l{2}$,    
 then $X \oplus Y \not\sim \l{2}$.   
 Hence if $X$ and $Y$ are $\SL{p}$ spaces,  
 then $X \oplus Y$ is an $\SL{p}$ space.  

 A result of Pe\char32lczy\'nski \xciteplus{P}{Proposition $(*)$}, presented below as Lemma 2.8,   
 states that for Banach spaces $V$ and $W$ which are isomorphic to their squares in the sense that   
 $V \oplus V \sim V$ and $W \oplus W \sim W$,   
 if $V \cinjects W$ and $W \cinjects V$,  
 then $V \sim W$.   

 Suppose $X$ and $Y$ are as above and are isomorphic to their squares.   
 If $X\cinjects Y$, then $X\oplus Y\sim Y$    
 [since $X\oplus Y$ and $Y$ are isomorphic to their squares, \hfil\break 
 $X\oplus Y\cinjects Y\oplus Y\sim Y$, and   
 $Y\cinjects X\oplus Y$].       
 If $X$ and $Y$ are incomparable in the sense that 
 $X\not\cinjects Y$ and $Y\not\cinjects X$,  
 then $X\oplus Y$ is isomorphically distinct from both $X$ and $Y$  
 [since $X\oplus Y\sim X$ would imply that $Y\cinjects X$, and    
 $X\oplus Y\sim Y$ would imply that $X\cinjects Y$].   
 Hence if $X$ and $Y$ are $\SL{p}$ spaces   
 which are isomorphic to their squares,    
 then the $\SL{p}$ space $X\oplus Y$  
 is isomorphically distinct from both $X$ and $Y$   
 if and only if $X$ and $Y$ are incomparable in the sense mentioned above.  

 From the list $\l{2}$, $\l{p}$, $\l{2}\oplus\l{p}$, $\seqsum{\l{2}}{\l{p}}$, $\L{p}$ of five spaces,     
 the only \hfil\break 
 incomparable pair of spaces is $\set{\l{2},\l{p}}$.   
 However, $\l{2}\oplus\l{p}$ has already been included in the list.   

 Let $Z$ be a separable infinite-dimensional Banach space such that $Z\cinjects\L{p}$.   
 Then $\seqsum{Z}{\l{p}}\cinjects\seqsum{\L{p}}{\l{p}}\sim\L{p}$.  
 Note that $\l{p}\cinjects\seqsum{Z}{\l{p}}$,  
 whence $\seqsum{Z}{\l{p}}\not\sim\l{2}$    
 and $\seqsum{Z}{\l{p}}$ is an $\SL{p}$ space. 
 The space $\seqsum{\l{2}}{\l{p}}$ is an example.   
 However, from the list $\l{2}$, $\l{p}$, $\l{2}\oplus\l{p}$, $\seqsum{\l{2}}{\l{p}}$, $\L{p}$ of five spaces,     
 no space arises from this method of construction which has not already been included in the list.    

\preheadspace
\firsthead{Preliminaries for Banach Spaces}
\postheadspace

 We now introduce some terminology used in the study of Banach spaces.
 The presentation is unavoidably terse and a bit disjointed. 
 General references for this \hfil\break 
 material include \xcite{L-T} and \xcite{W}. 
 Throughout the following discussion, $X$ and $Y$ will denote Banach spaces.  
 
 A Banach space is a complete normed vector space.
 Classical examples include the space $\L{p}(0,1)$ for $1\le p\le\infty$,
 with $\norm{f}_p=\paren{\int_{(0,1)} \abs{f}^p}^{{1\over p}}$ for $1\le p<\infty$
 and $\norm{f}_{\infty}=\mathop{{\rm ess}}\sup\abs{f}$ for $p=\infty$, 
 and the space $\l{p}$ for $1\le p\le\infty$, with \hfil\break 
 $\norm{\set{a_i}}_{\l{p}}=\paren{\sum\abs{a_i}^p}^{{1\over p}}$ for $1\le p<\infty$
 and $\norm{\set{a_i}}_{\l{\infty}}=\sup\abs{a_i}$ for $p=\infty$. 
 Here $\int$ denotes Lebesgue integration.
 Functions $f,g\in\L{p}(0,1)$ are identical as elements of $\L{p}(0,1)$
 if they agree except on a set of measure zero, 
 which is to say that strictly speaking,
 the elements of $\L{p}(0,1)$ are equivalence classes of functions.

 Given Banach spaces $X_1,X_2,\ldots$ and $1\le p<\infty$, 
 $\paren{X_1\oplus X_2\oplus\cdots}_{\l{p}}$ is the set of all sequences $\set{x_i}$ with $x_i\in X_i$ such that 
 $\norm{\set{x_i}}=\paren{\sum\norm{x_i}_{X_i}^p}^{{1\over p}}<\infty$. 
 The sum $\paren{X_1\oplus X_2\oplus\cdots}_{\l{p}}$ is a Banach space, 
 and will also be denoted $\sumsum{X_i}{\l{p}}$. 

 Suppose $T:X\to Y$ is a linear operator. 
 Then $T$ is said to be bounded if \hfil\break 
 $\norm{T}=\sup_{x\in X\setminus\set{0}} {\norm{T(x)}\over\norm{x}} < \infty$.
 A linear operator is bounded if and only if it is \hfil\break 
 continuous. 

 Suppose $T:X\to Y$ is a bounded linear operator. 
 Then $T$ is said to be an \hfil\break 
 isomorphism if $T$ has an inverse $T^{-1}:Y\to X$ which is a bounded linear operator. 
 If $T$ is a bijection, then $T$ is an isomorphism by the open mapping theorem. 
 If there is an isomorphism $S:X\to Y$, then $X$ and $Y$ are said to be isomorphic, and we write $X\sim Y$. 
 If $X\sim Y$, the Banach-Mazur distance between $X$ and $Y$ is \hfil\break 
 $d(X,Y)=\inf_S\set{\norm{S}\|S^{-1}\|}$, 
 where the infimum is taken over all isomorphisms \hfil\break 
 $S:X\to Y$. 

 Suppose $T:X\to Y$ is a bounded linear operator. 
 Then $T$ is called an isomorphic imbedding of $X$ into $Y$ if $T$ is an injection onto a closed subspace $Y'$ of $Y$. 
 If there is an isomorphic imbedding $S:X\to Y$, we write $X\injects Y$.

 Suppose $P:X\to X$ is a bounded linear operator. 
 Then $P$ is called a projection if $P^2=P$.
 Suppose $P:X\to X$ is a projection.
 Then $P(X)$ is a closed subspace of $X$,
 and each $x\in X$ has a unique representation as $x=y+z$ where $y\in P(X)$ and $P(z)=0$. 
 Moreover, $I-P:X\to X$ is a projection as well, where $I:X\to X$ is the identity mapping. 
 The range $R=P(X)$ and null space $N=(I-P)(X)$ of $P$ are said to be complemented subspaces of $X$, 
 and $X=R\oplus N$. 
 We write $R\cinjects X$ and $N\cinjects X$. 
 More generally, we write $Y\cinjects X$ if $Y$ is isomorphic to a complemented subspace of $X$. 
 
 The Rademacher functions $r_k:[0,1]\to\set{-1,1}$ for $k\in\NN$ are defined by \hfil\break 
 $r_k(t)=\mathop{{\rm sgn}}\sin(2^k\pi t)$. 

 For expressions $A$ and $B$ and constants $K_1$ and $K_2$,
 we write $A\within{K_1}{K_2}B$ to signify that $A\le K_1 B$ and $B\le K_2 A$. 
 We also write $A\approx B$ if $K_1$ and $K_2$ exist but are not specified. 
 If so indicated, $A\approx B$ will refer to an approximation rather than to a pair of inequalities.

 Khintchine's inequality states that for $1\le p<\infty$,
 there is a constant $K_p$ such that 
 for all scalars $a_1,a_2,\ldots$,
 for the Rademacher functions $r_1,r_2,\ldots$,
 and for all \hfil\break 
 $N\in\NN$, 
 ${1/K_p} \paren{\sum_{i=1}^N \abs{a_i}^2}^{{1\over2}} \le \norm{\sum_{i=1}^N a_i r_i}_p
                                                       \le K_p \paren{\sum_{i=1}^N \abs{a_i}^2}^{{1\over2}}$. 
 This inequality could also be expressed as 
 $\norm{\sum_{i=1}^N a_i r_i}_p \within{K_p}{K_p} \paren{\sum_{i=1}^N \abs{a_i}^2}^{{1\over2}}$. 

 A sequence $\set{x_i}$ in $X$ is said to be a (Schauder) basis for $X$ if for each $x\in X$, 
 there is a unique sequence $\set{a_i}$ of scalars such that $x=\sum a_i x_i$, 
 with convergence in the norm of $X$.  

 Given a sequence $\set{x_i}$ in $X$,
 the closed linear span of $\set{x_i}$ in $X$ will be denoted $\span{x_i}_X$, or simply $\span{x_i}$ if the context is clear. 
 Such a sequence is called a basic sequence if $\set{x_i}$ is a basis for $\span{x_i}_X$. 
 
 Given a sequence $\set{x_i}$ in $X$,
 the series $\sum x_i$ is said to converge unconditionally if any of the following equivalent conditions hold:
 (a) $\sum \epsilon_i x_i$ converges for all $\set{-1,1}$-valued sequences $\set{\epsilon_i}$, 
 (b) $\sum x_{\sigma(i)}$ converges for all permutations $\sigma$ of $\NN$, or
 (c) $\sum x_{n(i)}$ converges for all increasing sequences $\set{n(i)}$ in $\NN$. 

 A basis $\set{x_i}$ for $X$ is said to be unconditional if for each sequence of scalars for which $\sum a_i x_i$ converges, 
 the convergence is unconditional. 
 If $\set{x_i}$ is an unconditional basis for $X$,
 then for $P_E:\span{x_i}\to\span{x_i}$ defined by $P_E\paren{\sum_{i=1}^{\infty} a_i x_i}=\sum_{i\in E} a_i x_i$, 
 we have $\sup_{E\subset\NN} \norm{P_E} < \infty$. 

 Suppose $\set{x_i}$ is a basic sequence in $X$.
 A sequence $\set{y_j}$ in $X$ is called a block basic sequence (with respect to $\set{x_i}$)
 if $y_j\ne0$ for all $j\in\NN$
 and there are disjoint nonempty finite $E_1,E_2,\ldots\subset\NN$ with $\max E_j<\min E_{j'}$ for $j<j'$
 and scalars $a_1,a_2,\ldots$ such that 
 $y_j=\sum_{i\in E_j} a_i x_i$ for all $j\in\NN$. 
 Suppose $\set{y_j}$ is a block basic sequence \hfil\break 
 (with respect to $\set{x_i}$).
 Then $\set{y_j}$ is a basic sequence. 
 If $\set{x_i}$ is unconditional, then $\set{y_j}$ is unconditional as well. 

 Suppose $\set{x_i}$ and $\set{y_i}$ are bases for $X$ and $Y$, respectively.
 Then $\set{x_i}$ and $\set{y_i}$ are said to be equivalent if for all sequences $\set{a_i}$ of scalars, 
 $\sum a_i x_i$ converges if and only if $\sum a_i y_i$ converges. 
 If $\set{x_i}$ and $\set{y_i}$ are equivalent, 
 then there is a natural \hfil\break 
 isomorphism between $X$ and $Y$ by the closed graph theorem. 

 Suppose $\set{x_i}$ and $\set{y_i}$ are normalized bases for $X$ and $Y$,
 respectively, which are equivalent.
 Let $K$ be a positive constant. 
 Then $\set{x_i}$ and $\set{y_i}$ are said to be $K$-equivalent
 if for all sequences $\set{a_i}$ of scalars such that $\sum a_i x_i$ and $\sum a_i y_i$ converge, 
 $\norm{\sum a_i x_i}\within{K}{K}\norm{\sum a_i y_i}$. 

 A random variable is a measurable function on a probability space $(\Omega,\mu)$. 
 For $N\in\NN$, random variables $X_1,X_2,\ldots,X_N$ on $\Omega$ are said to be independent
 if for all Borel sets $B_1,B_2,\ldots,B_N$,
 $\mu\paren{\bigcap_{i=1}^N \set{t:X_i(t)\in B_i}} = \prod_{i=1}^N \mu\paren{\set{t:X_i(t)\in B_i}}$. \hfil\break 
 Random variables $X_1,X_2,\ldots$ on $\Omega$ are said to be independent
 if $X_1,X_2,\ldots,X_N$ are independent for each $N\in\NN$. 

\preheadspace
\firsthead{Overview of Chapters}
\postheadspace

 We briefly discuss the content of the succeeding chapters.

 Chapter II reviews the construction of Rosenthal \xcite{RI}. 
 Rosenthal's work is based on the study of the span in $\L{p}$ for $2<p<\infty$
 of sequences of independent mean zero random variables. 
 A few nonclassical $\SL{p}$ spaces were found by Rosenthal, 
 principal among them the space $X_p$.
 Chapter II includes a complete ordering of these spaces with respect to the (partial order) relation $\cinjects$. 
 
 Chapter III reviews the construction of Schechtman \xcite{S}. 
 Schechtman takes \hfil\break 
 Rosenthal's space $X_p$ and iterates a tensor product operation  
 to produce a sequence of $\SL{p}$ spaces. 
 Chapter III includes a section on the sequence space realization of \hfil\break 
 Schechtman's spaces, 
 expanding on a remark found in \xcite{S}.

 Chapter IV reviews the construction of Alspach \xcite{A}.
 Alspach's work generalizes the construction of Rosenthal,
 and generates spaces by means of a notion of \hfil\break 
 independent sum, 
 but has only been available in manuscript form.  
 A few nonclassical $\SL{p}$ spaces were found by Alspach,  
 principal among them a space denoted $D_p$.
 Chapter IV includes a complete ordering of these and Rosenthal's spaces with respect to $\cinjects$. 

 Chapter V reviews the construction of Bourgain, Rosenthal, and Schechtman \hfil\break 
 \xcite{B-R-S}. 
 These authors iterate and intertwine a notion of disjoint sum and a notion of independent sum 
 to generate a family of $\SL{p}$ spaces indexed by the countable ordinals,
 and distinguish these spaces isomorphically by means of an isomorphic invariant, \hfil\break 
 introduced in \xcite{B-R-S}, 
 which assigns an ordinal number to each separable Banach space. 

 Each chapter has a diagram relating the spaces under discussion with respect to $\cinjects$. 
 These diagrams are \TAG{1.2}{1.2}, \TAG{2.27}{2.27}, \TAG{3.2}{3.2}, \TAG{4.10}{4.10}, and \TAG{5.5}{5.5}.

%
%
%
\chapter{II}
\headspace
\chaptertitle{THE NONCLASSICAL $\SL{p}$ SPACES OF ROSENTHAL}
\headspace
 Let $1<p<\infty$ where $p\ne2$.  Rosenthal \xcite{RI} was the first to extend the list of 
 separable infinite-dimensional $\SL{p}$ spaces beyond the four previously known 
 isomorphism types:  $\L{p}$, $\l{p}$, $\l{2}\oplus\l{p}$, and $\seqsum{\l{2}}{\l{p}}$.  
 The principal $\SL{p}$ spaces which \hfil\break 
 Rosenthal constructed are $X_p$ and $B_p$,
 to be discussed presently. Using the newly revised list of six $\SL{p}$ spaces,
 Rosenthal constructed a few more such spaces by forming direct sums 
 (pairwise and in the sense of $\l{p}$ for sequences) of these six.
\preheadspace
\firsthead{The Space $X_p$}
\postheadspace
 In contrast to most classical Banach spaces, $X_p$ does not have a preferred
 standard realization. Let $2<p<\infty$. One realization of $X_p$ is as the 
 closed linear span in $\L{p}$ of a sequence $\set{f_n}$ of independent
 symmetric three-valued random variables such that the ratios 
 $\ratio$ approach zero slowly 
 (in a sense to be made precise). On the other hand, given positive weights 
 $w_n$ approaching zero slowly in the same sense, another realization of 
 $X_p$ is as the set of all sequences $\set{x_n}$ in $\l{p}$ for which the 
 weighted $\l{2}$ norm $\left(\sum|w_n x_n|^2\right)^{1\over2}$ is finite.  For the 
 conjugate index $q$, $X_q$ is defined to be the dual of $X_p$.     
\preheadspace
\secondhead{The Space $\Xpw$}
\postheadspace
 We first examine the sequence space realization of $X_p$.    


\definition Let $2<p<\infty$ and let $w=\set{w_n}$ be a sequence of positive scalars.   
            Define $\Xpw$ to be the set of all sequences $x=\set{x_n}$ of scalars for which both
            $\sum\abspower{x_n}{p}$ and $\sum\abspower{w_n x_n}{2}$ are finite.  
            For $x\in\Xpw$, define the norm $\normsub{x}{\Xpw}$ to be the maximum of
            $\left(\sum\abspower{x_n}{p\vphantom{2}}\right)^{1\over p}$ and 
            $\left(\sum\abspower{w_n x_n}{2}\right)^{1\over2}$. 

 Thus $\normsub{x}{\Xpw}$ is the maximum of the $\l{p}$ norm of $x$ and the weighted $\l{2}$ norm of $x$.  
 Under this norm, it is a routine matter to show that $\Xpw$ is a Banach space with unconditional standard basis. 
 The isomorphism type of $\Xpw$ depends on the sequence $w=\set{w_n}$ of weights,
 as partially outlined in the following proposition \xcite{RI}.

\xproclaim {Proposition 2.1}. Let $2<p<\infty$ and let $w=\set{w_n}$ be a sequence of positive scalars.
                              \item{(a)}If $\inf{w_n}>0$, then $\Xpw$ is isomorphic to $\l{2}$.
                              \item{(b)}If $\sum{\wstopower}<\infty$, then $\Xpw$ is isomorphic to $\l{p}$.
                              \item{(c)}If there is some $\epsilon>0$ for which $\set{n\colon w_n\ge\epsilon}$ and 
                              $\set{n\colon w_n<\epsilon}$ are both infinite and for which 
                              $\sum_{w_n<\epsilon} {\wstopower} < \infty$,
                              then $\Xpw$ is isomorphic to $\l{2}\oplus\l{p}$.
                              \item{(d)}Otherwise, $w$ satisfies condition $(*)$: 
     $$\hbox{{\sl for each\/} $\epsilon>0$, } \sum_{w_n < \epsilon} {\wstopower} = \infty.\eqno{(*)}$$  

\proof 

\item{(a)}Suppose $\xinf{w_n}=C>0$ and let $x=\set{x_n}\in\Xpw$. 
          Then $$\textstyle\normsub{x}{\l{p}} \le \normsub{x}{\l{2}}=\left(\sum\abspower{x_n}{2}\right)^{1\over2} \le %
                 {1\over C}\left(\sum\abspower{w_n x_n}{2}\right)^{1\over2}.$$
          Hence $$\textstyle\left(\sum\abspower{w_n x_n}{2}\right)^{1\over2} \le 
                  \normsub{x}{\Xpw} \le \max\set{{1\over C},1}\left(\sum\abspower{w_n x_n}{2}\right)^{1\over2},$$ 
          so $\Xpw$ is isomorphic to $\l{2}$ via the mapping 
          $\set{x_n}\mapsto\set{w_n x_n}$.
\item{(b)} Suppose $\sum\wstopower<\infty$ and let $x=\set{x_n}\in\Xpw$.  
          Then by H\"older's inequality with conjugate indices $p'={p\over2}>1$ and $q'={p\over p-2}$, we have 
          $$\textstyle\sum\abspower{w_n x_n}{2}=\sum|{w_n}^2 {x_n}^2| \le %
                \left(\sum {w_n}^{2{p\over{p-2}}}\right)^{{{p-2}\over p}}
                \left(\sum\abspower{x_n}{2{p\over2}}\right)^{{2\over p}}.$$ 
          Let $K=\left(\sum {w_n}^{2{p\over{p-2}}}\right)^{{{p-2}\over2p}}$. 
          Then $\left(\sum\abspower{w_n x_n}{2}\right)^{{1\over2}} \le 
                K\left(\sum\abspower{x_n}{p\vphantom{2}}\right)^{{1\over p}}$. 
          Hence $$\textstyle\normsub{x}{\l{p}} \le \normsub{x}{\Xpw} \le \max\set{1,K}\normsub{x}{\l{p}},$$ 
          so $\Xpw$ is isomorphic to $\l{p}$ via the formal identity mapping.
\item{(c)}The hypothesis of part (c) is equivalent to the hypothesis that $\NN$ is the  
          disjoint union of two infinite sets $N_1$ and $N_2$ for which  
          $\xinf_{n\in{N_1}} w_n >0$ and \hfil\break 
          $\sum_{n\in{N_2}} \wstopower<\infty$. 
          Thus part (c) follows from parts (a) and (b) and the unconditionality
          of the standard basis of $\Xpw$.  
\item{(d)}Condition $(*)$ is equivalent to the conjunction of the negations
          of the hypotheses of parts (a), (b), and (c). \qquad\QED  

\medskip
\xremark{1} We will show later that for fixed $2<p<\infty$, all spaces $\Xpw$   
        for $w$ satisfying condition $(*)$ are mutually isomorphic, but isomorphically distinct 
        from $\l{2}$, $\l{p}$, and $\l{2}\oplus\l{p}$ (as well as $\seqsum{\l{2}}{\l{p}}$ and $\L{p}$). 
        Thus part (d) is indeed a different case, and part (d) does not split into subcases.    

\xremark{2} Let $2<p<\infty$. If $\xinf{w_n}=0$ (as occurs in parts (b), (c), and (d)), then $\Xpw$ contains a   
        complemented subspace isomorphic to $\l{p}$, since some sub\-sequence of $w$ satisfies the 
        hypothesis of part (b).     
        Hence in parts (b), (c), and (d), $\Xpw$ is not isomorphic to $\l{2}$. 
        We will show later that the spaces $\Xpw$ are 
        isomorphic to complemented subspaces of $\L{p}$.
        Thus only part (a) does not yield an $\SL{p}$ space,
        while parts (b) and (c) yield known $\SL{p}$ spaces,
        and part (d) yields a previously unknown $\SL{p}$ space. 
        The spaces $\Xpw$ for $w$ satisfying 
        condition $(*)$ will be our sequence space realizations of $X_p$. 


\preheadspace
\secondhead{Rosenthal's Inequality}
\postheadspace
 Rosenthal proved the following fundamental probabilistic inequality \hfil\break 
 \xciteplus{RI}{Theorem 3}, which (in its corollary) relates $\Xpw$ with the 
 closed linear span of a sequence of independent mean zero random variables in $\L{p}$ ($2<p<\infty$). 
 
\xproclaim {Theorem 2.2}. Let $2<p<\infty$. 
                          There is a constant $K_p$, depending only on $p$, such that if 
                          $f_1,\ldots,f_N$ are independent mean zero random variables in $\L{p}$, then
                          \item{(a)} $\normsub{\sum_{n=1}^{N}{f_n}}{p} \le 
                                     K_p \max\set{\paren{\tsum_{n=1}^N \norm{f_n}_p^p }^{{1\over p}},
                                                  \paren{\tsum_{n=1}^N \norm{f_n}_2^2 }^{{1\over 2}} }$, and
                          \item{(b)} $\normsub{\sum_{n=1}^{N}{f_n}}{p} \ge 
                                     {1\over2} \max\set{\paren{\tsum_{n=1}^N \norm{f_n}_p^p }^{{1\over p}},
                                                        \paren{\tsum_{n=1}^N \norm{f_n}_2^2 }^{{1\over 2}} }$. 
                          \item{}    If in addition $f_1,\ldots,f_N$ are assumed to be symmetric,
                                     then the constant $1\over2$ can be replaced by $1$.

\remark It is shown in \xcite{J-S-Z} that $K_p$ is of order $p/{\log{p}}$. 

 The proof of Rosenthal's inequality will not be presented, but we deduce its \hfil\break 
 corollary \xcite{RI}.  

\xproclaim {Corollary 2.3}. Let $2<p<\infty$, 
                            let $\set{f_n}$ be a sequence of independent mean zero random variables in $\L{p}$,
                            and let $\wdef$.
                            Then $\spansub{f_n}{\L{p}}$ is isomorphic to $\Xpw$,
                            and $\set{f_n}$ in $\L{p}$ is equivalent to the standard basis of $\Xpw$.  

\proof Without loss of generality,   
       suppose each $f_n$ is of norm one in $\L{p}$, so that $w_n=\normsub{f_n}{2}$. 
       Let $f\in\finitespan{\set{f_n}}$ and express $f$ as $\sum_{n=1}^{N}{c_n f_n}$. 
       Then by Theorem 2.2, we have 
       $$\textstyle\normsub{\sum_{n=1}^{N}{c_n f_n}}{p} \within{K_p}{2} %
         \max\set{\left(\sum_{n=1}^{N}\abspower{c_n}{p}\right)^{1\over p}, 
         \left(\sum_{n=1}^{N}\abspower{c_n w_n}{2}\right)^{1\over2}}.$$ 
       Hence $\spansub{f_n}{\L{p}}$ is isomorphic to $\Xpw$ via the mapping 
       $\sum{c_n f_n} \mapsto \set{c_n}$, and $\set{f_n}$ in $\L{p}$ is 
       equivalent to the standard basis of $\Xpw$. \qquad\QED  

  
\xremark{1} Let $2<p<\infty$. Given a sequence $w=\set{w_n}$ of positive   
        scalars for which $\xsup{w_n}\le1$, $\set{w_n}$ can be realized as 
        $\set{\ratio}$ for $\set{f_n}$ satisfying the hypotheses of Corollary 2.3. 
        If $\xsup w_n>1$, then $\Xpw\sim X_{p,w'}$ for some sequence \hfil\break 
        $w'=\set{w'_n}$ satisfying $\xsup w'_n\le1$.  
        Thus there is a complete correspondence between the sequence spaces $\Xpw$
        and the function spaces $\spansub{f_n}{\L{p}}$ 
        for $\set{f_n}$ satisfying the hypotheses of Corollary 2.3.  

\xremark{2} For fixed $2<p<\infty$,
        the spaces $\spansub{f_n}{\L{p}}$ for $\set{f_n}$ satisfying the hypotheses of Corollary 2.3 and 
        $\wdef$ satisfying condition $(*)$ of Proposition 2.1 will be our function space realizations of $X_p$. 
\preheadspace
\secondhead{The Complementation of $\Xpw$ in $\L{p}$} 
\postheadspace
 Let $2<p<\infty$. In its sequence space realizations, it is not so clear that $X_p$ is an $\SL{p}$ space. 
 However, we will soon show that in its function space realizations, the complementation of $\spansub{f_n}{\L{p}}$
 in $\L{p}$ follows if the sequence $\set{f_n}$ satisfies certain additional hypotheses. On the other hand, 
 in its function space realizations, the isomorphic structure of $X_p$ is not so clear.  
 We will go back and forth between realizations, depending on their relative advantages at the time. 

 Suppose $f_n$ is a symmetric three-valued random variable. Let $\alpha_n$ be the 
 positive value attained by $|f_n|$ and let $\mu_n$ be the measure of the set 
 on which $f_n$ is nonzero. Then for $1\le r<\infty$, we have  
 $$\normsub{f_n}{r}=\left({\alpha_n}^r \mu_n\right)^{{1\over r}}=\alpha_n {\mu_n}^{{1\over r}}.$$
 Let $2<p<\infty$. Then $w_n=\ratio={\mu_n}^{{1\over 2}-{1\over p}}={\mu_n}^{{{p-2}\over 2p}}$. 
 Hence $$\wstopower=\mu_n.$$
 This provides an interpretation for condition $(*)$ of Proposition 2.1 in terms of properties of  
 a sequence $\set{f_n}$ of independent symmetric three-valued random variables, namely    
 $$\hbox{for each $\epsilon>0$, } \sum_{\mu_n<\epsilon}{\mu_n}=\infty.$$
 Let $q$ be the conjugate index of $p$. Then     
 $$\normsub{f_n}{p}\normsub{f_n}{q}={{\alpha_n}^2}{\mu_n}^{{1\over p}+{1\over q}}={{\alpha_n}^2}\mu_n
                                   ={\left({\alpha_n}{\mu_n}^{{1\over2}}\right)}^2=\normsubpower{f_n}{2}{2}.$$ 
 This provides a way to interrelate the $\L{p}$, $\L{q}$, and $\L{2}$ norms of a  
 symmetric three-valued random variable. We will find this useful in the proof of the next theorem, 
 where we show that a certain projection is bounded in both $\L{2}$ and $\L{p}$ norms. 
 We will make explicit use of the fact that if $f_n$ is a symmetric three-valued random variable of 
 norm one in $\L{p}$, then 
 $$\normsub{{f_n\over\normsubpower{f_n}{2}{2}}}{q}=
   {\normsub{f_n}{q}\over\normsubpower{f_n}{2}{2}}={1\over\normsub{f_n}{p}} = 1.\eqno{\TAG{2.1}{2.1}}$$

\remark If the scalars are complex,
        the hypothesis that $f_n$ is a symmetric three-valued random variable  
        can be replaced by the hypothesis that $f_n$ is a mean zero \hfil\break 
        random variable for which $\abs{f_n}$ is $\set{0,\alpha_n}$-valued for $\alpha_n \ne 0$.  
  
 Rosenthal proved the following theorem \xciteplus{RI}{Theorem 4}, which (in its corollary) establishes that for $2<p<\infty$, 
 the spaces $\Xpw$ are isomorphic to complemented subspaces of $\L{p}$. To prove the theorem, we use the  
 following probabilistic inequality \xciteplus{RI}{Lemma 2b}, which we state without proof.   

\xproclaim {Lemma 2.4}. Let $1\le q<2$ and let $f_1,\ldots,f_N$ be independent mean zero random variables in $\L{q}$. Then
                        $$\normsub{\tsum_{n=1}^{N}{f_n}}{q}
                          \le 2 \paren{\tsum_{n=1}^N \norm{f_n}_q^q }^{{1\over q}}.$$
                        If in addition $f_1,\ldots,f_N$ are assumed to be symmetric, 
                        then the constant $2$ can be \hfil\break 
                        replaced by $1$.

\xproclaim {Theorem 2.5}. Let $1<p<\infty$ and let $\set{f_n}$ be a sequence of independent symmetric
                          three-valued random variables in $\L{p}$. Then there is a projection
                          $P\colon\L{p}\rightarrow\L{p}$ onto $\spansub{f_n}{\L{p}}$ with $\norm{P}\le C_p$,
                          where $C_2=1$, $C_p=K_p$ (the constant in Theorem 2.2) for $2<p<\infty$, and   
                          $C_p=C_q$ for conjugate indices $p$ and $q$.  

\proof If $p=2$, the orthogonal projection $\pi\colon\L{2}\rightarrow\L{2}$ onto $\spansub{f_n}{\L{2}}$  
       satisfies the requirements. We will presently show that for    
       $2<p<\infty$, the set-theoretic restriction of $\pi$ to $\L{p}$ 
       yields a bounded projection $P\colon\L{p}\rightarrow\L{p}$ onto $\spansub{f_n}{\L{p}}$   
       with $\norm{P}\le K_p$. This will suffice to prove the theorem in the general case,  
       since the adjoint then induces a projection $Q\colon\L{q}\rightarrow\L{q}$ 
       onto $\spansub{f_n}{\L{q}}$ with $\norm{Q}=\norm{P}$.  

       Let $2<p<\infty$, so that $\L{p}\subset\L{2}$. Let $\wdef$.
       Without loss of generality,  
       suppose $f_n$ is real-valued with $\normsub{f_n}{p}=1$. Then $w_n=\normsub{f_n}{2}$. 
       Let $\pi\colon\L{2}\rightarrow\spansub{f_n}{\L{2}}$ be the orthogonal projection 
       defined by  
       $$\pi(g)=\sum\left(\int_0^1 g(t){f_n\over\normsub{f_n}{2}}(t)\,dt\right){f_n\over\normsub{f_n}{2}}.$$
       Then $\normsub{\pi(g)}{2}\le\normsub{g}{2}$. We will show that if $g\in\L{p}$, then 
       $\pi(g)\in\L{p}$ and \hfil\break 
       $\normsub{\pi(g)}{p}\le K_p\normsub{g}{p}$. 
       Thus
       $$P(g)=\sum\left(\int_0^1 g(t){f_n\over\normsubpower{f_n}{2}{2}}(t)\,dt\right)f_n$$  
       defines a mapping $P\colon\L{p}\rightarrow\spansub{f_n}{\L{p}}$.
       Set-theoretically, $P$ is the restriction of $\pi$ to $\L{p}$.   
       It will follow that $P$ is a projection and $\norm{P}\le K_p$.    

       Fix $g\in\L{p}$ and let  
       $$x_n = \int_0^1 g(t){f_n\over\normsubpower{f_n}{2}{2}}(t)\,dt,$$   
       so that $\pi(g)=\sum x_n f_n$.  
       We will show that $\set{x_n}\in\Xpw$ and $\normsub{\set{x_n}}{\Xpw} \le \normsub{g}{p}$. 
       Corollary 2.3 will then yield
       $\normsub{\pi(g)}{p}=\normsub{\sum x_n f_n}{p} \le K_p\normsub{\set{x_n}}{\Xpw} \le K_p\normsub{g}{p}$.  

       First we examine the weighted $\l{2}$ norm of $\set{x_n}$.            
       Let        
       $$y_n = \int_0^1 g(t){f_n\over\normsub{f_n}{2}}(t)\,dt=x_n\normsub{f_n}{2}=x_n w_n.$$   
       Then            
       $$\left(\textstyle\sum\abspower{w_n x_n}{2}\right)^{{1\over2}} = 
         \normsub{\set{y_n}}{\l{2}} = 
         \normsub{\textstyle\sum y_n {f_n\over\normsub{f_n}{2}}}{\L{2}} = 
         \normsub{\pi(g)}{2} \le \normsub{g}{2} \le \normsub{g}{p}.\eqno{\TAG{2.2}{2.2}}$$

       Next we examine the $\l{p}$ norm of $\set{x_n}$. We verify that $\set{x_n}\in\l{p}$ 
       by testing against $\l{q}$. Let $\set{c_n}\in\l{q}$. Using Lemma 2.4 and equation \TAG{2.1}{2.1},    
       for each $N\in\NN$ 
       $$\eqalign{\normsub{\sum_{n=1}^{N} c_n {f_n\over\normsubpower{f_n}{2}{2}}}{q}&\le 
                  \left(\sum_{n=1}^{N} \normsubpower{c_n {f_n\over\normsubpower{f_n}{2}{2}}}{q}{q}
                  \right)^{{1\over q}}\cr&= 
                  \left(\sum_{n=1}^{N} \abspower{c_n}{q}\right)^{{1\over q}}\cr&\le 
                  \normsub{\set{c_n}}{\l{q}}.}$$ 
       Now by H\"older's inequality and the observation above, for each $N\in\NN$   
       $$\eqalign{\left|\sum_{n=1}^{N} c_n x_n \right|&=  
                  \left|\sum_{n=1}^{N} c_n\int_0^1 g(t){f_n\over\normsubpower{f_n}{2}{2}}(t)\,dt\right|\cr&=  
                  \left|\int_0^1 g(t)\sum_{n=1}^{N} c_n{f_n\over\normsubpower{f_n}{2}{2}}(t)\,dt\right|\cr&\le  
                  \normsub{g}{p}\normsub{\sum_{n=1}^{N} c_n{f_n\over\normsubpower{f_n}{2}{2}}}{q}\cr&\le  
                  \normsub{g}{p}\normsub{\set{c_n}}{\l{q}}.}$$  
       Hence $\set{x_n}\in\l{p}$ and $$\normsub{\set{x_n}}{\l{p}} \le \normsub{g}{p}.\eqno{\TAG{2.3}{2.3}}$$ 
                     
       Combining \TAG{2.2}{2.2} and \TAG{2.3}{2.3}, we see that     
       $\set{x_n}$ is indeed in $\Xpw$ and \hfil\break 
       $\normsub{\set{x_n}}{\Xpw} \le \normsub{g}{p}$. 


       Now by Corollary 2.3 (and the inequality appearing in its proof), we have \hfil\break 
       $\normsub{\sum x_n f_n}{p} \within{K_p}{1} \normsub{\set{x_n}}{\Xpw}$, so that 
       $$\normsub{\pi(g)}{p}=\normsub{\tsum x_n f_n}{p} 
         \le K_p\normsub{\set{x_n}}{\Xpw} \le K_p\normsub{g}{p}.$$ 
       Hence $P(g)=\pi(g)\in\spansub{f_n}{\L{p}}$ and $P$ is a projection from $\L{p}$ onto $\spansub{f_n}{\L{p}}$ 
       with $\norm{P} \le K_p$. \qquad\QED   

\remark If the scalars are complex,
        the hypothesis that each $f_n$ is symmetric and three-valued can be replaced by the hypothesis that each $f_n$ 
        is mean zero and $|f_n|$ is $\set{0,\alpha_n}$-valued for $\alpha_n \ne 0$, 
        but without the hypothesis of symmetry we have \hfil\break 
        $\norm{P}\le2C_p$.  

 We deduce the following corollary \xcite{RI}. 

\xproclaim {Corollary 2.6}. Let $2<p<\infty$ and let $w=\set{w_n}$ be a sequence of positive scalars. Then 
                            $\Xpw$ is isomorphic to a complemented subspace of $\L{p}$.   
                            If\/ $\inf{w_n}=0$, then $\Xpw$ is an $\SL{p}$ space. 
                            In particular, if $w$ satisfies condition $(*)$ of Proposition 2.1, 
                            then $\Xpw$ is an $\SL{p}$ space.   

\proof First suppose that $\xsup{w_n} \le 1$. Then $\set{w_n}$ can be realized as \hfil\break 
       $\set{\ratio}$ for a sequence 
       $\set{f_n}$ of independent symmetric (whence mean zero) three-valued random variables in $\L{p}$. 
       Hence $\Xpw$ is isomorphic to $\spansub{f_n}{\L{p}}$ by Corollary 2.3, and $\spansub{f_n}{\L{p}}$ is 
       complemented in $\L{p}$ by Theorem 2.5.    
                    
       Now suppose that $\xsup w_n >1$.   
       Let $N_0=\set{n\colon w_n\le1}$ and $N_1=\set{n\colon w_n>1}$.   
       Let $w_{[0]}=\set{w_n}_{n\in N_0}$ and $w_{[1]}=\set{w_n}_{n\in N_1}$,  
       and let $\set{1}=\set{1}_{n\in N_1}$ be the sequence with constant value one.    
       Let $w'=\set{w'_n}_{n=1}^{\infty}=\set{\min\set{w_n,1}}_{n=1}^{\infty}$,   
       whence $\xsup w'_n \le 1$ and $X_{p,w'}\cinjects\L{p}$.   
       Then
       $$X_{p,w} \sim X_{p,w_{[0]}} \oplus X_{p,w_{[1]}}
                 \sim X_{p,w_{[0]}} \oplus X_{p,\set{1}} \sim X_{p,w'} \cinjects \L{p},$$   
       where for an $N$-tuple $v=\set{v_1,\ldots,v_N}$ of positive scalars,   
       $X_{p,v}$ is defined in the \hfil\break 
       obvious way, and $X_{p,\emptyset}=\set{0}$.   

       If $\xinf{w_n}=0$, then $\Xpw$ contains a complemented subspace isomorphic to $\l{p}$, 
       whence $\Xpw$ is not isomorphic to $\l{2}$. Hence if $\xinf{w_n}=0$, then $\Xpw$ is an $\SL{p}$ space by Theorem 1.1.   
       Finally, note that if $w=\set{w_n}$ satisfies condition $(*)$ of Proposition 2.1, then 
       $\xinf w_n = 0$.    
       \qquad\QED 
\preheadspace
\secondhead{The Mutual Isomorphism of the Spaces $\Xpw$}
\postheadspace
 We will show that for fixed $2<p<\infty$,    
 all spaces $\Xpw$ for $w=\set{w_n}$ satisfying condition $(*)$ of Proposition 2.1 are mutually isomorphic,   
 and isomorphically distinct from the previously known $\SL{p}$ spaces.   
 These two results are our next major concerns.    
 The following proposition \xciteplus{RI}{Lemma 7} will be used in the proofs of both of these results.  

\xproclaim {Proposition 2.7}. Let $2<p<\infty$ and let $w=\set{w_n}$ be a sequence of positive scalars. 
          Suppose that $\set{E_j}$ is a sequence of disjoint nonempty finite subsets of $\NN$.
          Let \hfil\break 
          $b_j=\sum_{n\in E_j} {w_n}^{{2\over p-2}} e_n$ and $\tilde b_j=b_j\big/\normsub{b_j}{\l{p}}$,   
          where $\set{e_n}$ is the standard basis of $\Xpw$.    
          Let $v_j=\left(\sum_{n\in E_j} {w_n}^{{2p\over p-2}}\right)^{{p-2\over2p}}$ and $v=\set{v_j}$. Then   
          \item{(a)} $\set{\tilde b_j}$ is an unconditional basis for $\spansub{\tilde b_j}{\Xpw}$ which      
                     is isometrically equivalent to the standard basis of $X_{p,v}$, and
          \item{(b)} there is a projection $P\colon\Xpw \to \spansub{\tilde b_j}{\Xpw}$ with $\norm{P}=1$. 

\proof First we establish some notation. 
       Let $\ell_{2,w}$ be the Hilbert space of all sequences $x=\set{x_n}$ of scalars for which 
       $\normsub{x}{\ell_{2,w}}=\left(\sum\abspower{w_n x_n}{2}\right)^{{1\over2}}<\infty$, where the inner product  
       in $\ell_{2,w}$ is defined by $\vector{x,y}=\sum x_n \bar y_n {w_n}^2$
       (where $x=\set{x_n}$, $y=\set{y_n}$, and bar
       is complex conjugation).  
       Motivating the choice of the $b_j$ is the fact that 
       $$\normsubpower{b_j}{\l{p}}{p}=\sum\nolimits_{n\in E_j} {w_n}^{{2p\over p-2}} =
         \sum\nolimits_{n\in E_j} {w_n}^{{4\over p-2}} {w_n}^2 = 
         \normsubpower{b_j}{\ell_{2,w}}{2}.$$ 
       Let $\sigma_j$ denote the common value of 
       $\normsubpower{b_j}{\l{p}}{p}$, $\normsubpower{b_j}{\ell_{2,w}}{2}$, and $\sum_{n\in E_j} {w_n}^{{2p\over p-2}}$. 
       Note that $v_j={\sigma_j}^{{p-2\over2p}}$ by our definitions.   

       \item{(a)} The unconditionality of $\set{\tilde b_j}$ follows from the unconditionality of $\set{e_n}$ in $\Xpw$. 
                  We now examine the isometric equivalence of the bases. 
                  Let $J\in\NN$ and let \hfil\break 
                  $\lambda_1,\ldots,\lambda_J$ be scalars.
                  Then  
                $$\eqalignno{\normsubpower{\tsuml_{j=1}^{J} \lambda_j b_j}{\l{p}}{p}&= 
                             \normsubpower{\tsuml_{j=1}^{J} \lambda_j \tsuml_{n\in E_j} {w_n}^{{2\over p-2}} e_n}{\l{p}}{p}\cr&= 
                             \tsuml_{j=1}^{J} \abspower{\lambda_j}{p}\tsuml_{n\in E_j} {w_n}^{{2p\over p-2}}\cr&= 
                             \tsuml_{j=1}^{J} \abspower{\lambda_j}{p} \sigma_j &\TAG{2.4}{2.4}}$$ 
                  and 
                $$\eqalign{\normsubpower{\tsuml_{j=1}^{J} \lambda_j b_j}{\ell_{2,w}}{2}&= 
                           \normsubpower{\tsuml_{j=1}^{J} \lambda_j \tsuml_{n\in E_j} {w_n}^{{2\over p-2}} e_n}{\ell_{2,w}}{2}\cr&=
                           \tsuml_{j=1}^{J} \abspower{\lambda_j}{2}\tsuml_{n\in E_j} {w_n}^{{4\over p-2}} {w_n}^2 \cr&= 
                           \tsuml_{j=1}^{J} \abspower{\lambda_j}{2}\tsuml_{n\in E_j} {w_n}^{{2p\over p-2}} \cr&= 
                           \tsuml_{j=1}^{J} \abspower{\lambda_j}{2} \sigma_j.}$$ 
                  Normalizing each $b_j$ in $\l{p}$ and noting that $\normsub{b_j}{\l{p}}={\sigma_j}^{{1\over p}}$, 
                  we have 
                  $$\normsubpower{\tsuml_{j=1}^{J} \lambda_j \tilde b_j}{\l{p}}{p} =
                    \tsuml_{j=1}^{J} \abspower{\lambda_j}{p}\eqno{\TAG{2.5}{2.5}}$$ 
                  and 
                  $$\normsubpower{\tsuml_{j=1}^{J} \lambda_j \tilde b_j}{\ell_{2,w}}{2} =  
                    \tsuml_{j=1}^{J} \abspower{\lambda_j}{2} {\sigma_j\over{\sigma_j}^{{2\over p}}} = 
                    \tsuml_{j=1}^{J} \abspower{\lambda_j}{2} {\sigma_j}^{{p-2\over p}} =  
                    \tsuml_{j=1}^{J} \abspower{\lambda_j}{2} {v_j}^2.\eqno{\TAG{2.6}{2.6}}$$   
                  Thus  
                  $$\eqalign{\normsub{\tsuml_{j=1}^{J} \lambda_j \tilde b_j}{\Xpw}&=  
                             \max\set{\normsub{\tsuml_{j=1}^{J} \lambda_j \tilde b_j}{\l{p}}, 
                                      \normsub{\tsuml_{j=1}^{J} \lambda_j \tilde b_j}{\ell_{2,w}}}\cr&=  
                             \max\set{\left(\tsuml_{j=1}^{J} \abspower{\lambda_j}{p}\right)^{{1\over p}},  
                                      \left(\tsuml_{j=1}^{J} \abspower{v_j \lambda_j}{2}\right)^{{1\over2}}}.}$$  
                  Hence $\set{\tilde b_j}$ in $\Xpw$ is isometrically equivalent to the standard basis of $X_{p,v}$.  

       \item{(b)} We wish to define a projection $P\colon\Xpw \to \spansub{b_j}{\Xpw}$ with $\norm{P} = 1$. 
                  Recalling the inner product $\vector{\phantom0,\phantom0\!}$ previously introduced on $\ell_{2,w}$, let 
                  $\pi\colon\ell_{2,w} \to \spansub{b_j}{\ell_{2,w}}$ be the orthogonal projection defined by 
                  $$\pi(x)=\sum_{j=1}^{\infty}
                    \vector{x,{b_j\over\norm{b_j}}\xsubpower{\ell_{2,w}}{}}{b_j\over\norm{b_j}}\xsubpower{\ell_{2,w}}{}.$$
                  Then $\normsub{\pi(x)}{\ell_{2,w}} \le \normsub{x}{\ell_{2,w}}$. We will show that if 
                  $x\in\l{p}\cap\ell_{2,w}$, then $\pi(x)\in\l{p}$ and $\normsub{\pi(x)}{\l{p}}\le\normsub{x}{\l{p}}$.
                  Thus
                  $$P(x)=\sum_{j=1}^{\infty} \vector{x,{b_j\over\norm{b_j}}\xsubpower{\ell_{2,w}}{2}} {b_j}$$ 
                  defines a mapping $P\colon\l{p}\cap\ell_{2,w}\to\spansub{b_j}{\l{p}\cap\ell_{2,w}}$. 
                  Set-theoretically, $P$ is the restriction of $\pi$ to $\l{p}\cap\ell_{2,w}$.  
                  It will follow that if $x\in\l{p}\cap\ell_{2,w}=\Xpw$, then  
                  $$\normsub{P(x)}{\Xpw} = \max\set{\normsub{P(x)}{\ell_{2,w}},\normsub{P(x)}{\l{p}}} 
                                        \le\max\set{\normsub{x}{\ell_{2,w}},\normsub{x}{\l{p}}} 
                                         = \normsub{x}{\Xpw}.$$
                  Fix $x=\set{x_n}\in\l{p}\cap\ell_{2,w}$ and let
                  $$\lambda_j=\vector{x,{b_j\over\norm{b_j}}\xsubpower{\ell_{2,w}}{2}},$$ 
                  so that $\sum_{j=1}^J \lambda_j b_j$ is a partial sum of $\pi(x)$. 
                  We now show that $\pi(x)\in\l{p}$ and 
                  $\normsub{\pi(x)}{\l{p}}\le\normsub{x}{\l{p}}$.
                  As in equation \TAG{2.4}{2.4}, we have 
                  $$\normsubpower{\tsuml_{j=1}^J \lambda_j b_j}{\l{p}}{p} = 
                    \tsuml_{j=1}^J \abspower{\lambda_j}{p} \sigma_j,$$
                  where 
                  $$\eqalign{\lambda_j=\vector{x,{b_j\over\norm{b_j}}\xsubpower{\ell_{2,w}}{2}}&=
                                       {1\over\sigma_j}\vector{x,b_j}\cr&=
                                       {1\over\sigma_j}\tsum_{n\in E_j} x_n {w_n}^{{2\over p-2}} {w_n}^2\cr&=  
                                       {1\over\sigma_j}\tsum_{n\in E_j} x_n {w_n}^{{2(p-1)\over p-2}}.}$$ 
                  Now by H\"older's inequality, for $q={p\over p-1}$ we have 
                  $$\eqalign{\abs{\lambda_j}&=
                             {1\over\sigma_j}\abs{\tsum_{n\in E_j} x_n {w_n}^{{2(p-1)\over p-2}}}\cr&\le  
                             {1\over\sigma_j} 
                             \left(\tsum_{n\in E_j} \abspower{x_n}{p}\right)^{{1\over p}} 
                             \left(\tsum_{n\in E_j} {w_n}^{{2(p-1)\over p-2}q}\right)^{{1\over q}}\cr&= 
                             {1\over\sigma_j} 
                             \left(\tsum_{n\in E_j} \abspower{x_n}{p}\right)^{{1\over p}} 
                             \left(\tsum_{n\in E_j} {w_n}^{{2p\over p-2}}\right)^{{p-1\over p}}\cr&= 
                             {1\over\sigma_j} 
                             \left(\tsum_{n\in E_j} \abspower{x_n}{p}\right)^{{1\over p}} 
                             {\sigma_j}^{{p-1\over p}}\cr&=
                             {1\over{\sigma_j}^{{1\over p}}}
                             \left(\tsum_{n\in E_j} \abspower{x_n}{p}\right)^{{1\over p}}.}$$ 
                  Hence $\abspower{\lambda_j}{p}\sigma_j \le \sum_{n\in E_j} \abspower{x_n}{p}$.  
                  Referring again to equation \TAG{2.4}{2.4}, for each $J\in\NN$  
                  $$\normsubpower{\tsuml_{j=1}^J \lambda_j b_j}{\l{p}}{p}
                    =\tsuml_{j=1}^J \abspower{\lambda_j}{p} \sigma_j
                    \le \tsuml_{j=1}^J \tsuml_{n\in E_j} \abspower{x_n}{p}
                    \le \normsubpower{x}{\l{p}}{p}.\eqno{\TAG{2.7}{2.7}}$$
                  Hence $\pi(x)=\sum_{j=1}^{\infty} \lambda_j b_j \in \l{p}$ 
                  and $\normsub{\pi(x)}{\l{p}} \le \normsub{x}{\l{p}}$. 
                  \qquad\QED 

\medskip
 We continue with results leading to the conclusion that for fixed $2<p<\infty$, 
 all spaces $\Xpw$ for $w=\set{w_n}$ satisfying condition $(*)$ of Proposition 2.1 are mutually isomorphic.  
 The following result of Pe\char32lczy\'nski \xciteplus{P}{Proposition $(*)$} indicates the approach to be taken.  

\xproclaim {Lemma 2.8}. Let $X$ and $Y$ be Banach spaces. Suppose $X \cinjects Y$ and $Y \cinjects X$, 
                        where $X \sim X \oplus X$ and $Y \sim Y \oplus Y$. Then $X \sim Y$. 

\proof Let $X'$ be a closed subspace of $X$ such that $X\sim Y\oplus X'$. Then \hfil\break 
       $X\sim Y\oplus X'\sim Y\oplus Y\oplus X'\sim Y\oplus X$. Similarly, $Y\sim X\oplus Y$. Hence \hfil\break 
       $X\sim Y\oplus X\sim X\oplus Y\sim Y$. \qquad\QED 

 First we examine the matter of mutual complementation \xciteplus{RI}{Theorem 13}. 


\xproclaim {Proposition 2.9}. Let $2<p<\infty$ and let $w=\set{w_n}$ and $w'=\set{w'_n}$ be sequences of positive scalars 
                              satisfying condition $(*)$ of Proposition 2.1. Then 
                              $X_{p,w'} \cinjects \Xpw$. 

\proof By condition $(*)$, we may choose a sequence $\set{E_j}$ of disjoint nonempty finite subsets of $\NN$
       such that for each $j\in\NN$, 
       $$\paren{w_j^{\prime}}^{{2p\over p-2}} \le \sum_{n\in E_j} {w_n}^{{2p\over p-2}}
                                              \le \paren{2w_j^{\prime}}^{{2p\over p-2}}.$$
       Then for $v_j=\paren{\sum_{n\in E_j} {w_n}^{{2p\over p-2}}}^{{p-2\over2p}}$,
       $w_j^{\prime} \le v_j \le 2w_j^{\prime}$.
       Hence for $v=\set{v_j}$ and \hfil\break 
       $x\in X_{p,w'}$, $\normsub{x}{X_{p,w'}}\le \normsub{x}{X_{p,v}}\le 2\normsub{x}{X_{p,w'}}$. 
       Thus $X_{p,w'}\sim X_{p,v}$ via the formal identity mapping. 
       For $\tilde b_j$ as in Proposition 2.7, $X_{p,v} \sim \spansub{\tilde b_j}{\Xpw} \cinjects \Xpw$. 
       Hence $X_{p,w'} \cinjects \Xpw$. \qquad\QED 

 Next we examine the matter of $\Xpw$ being isomorphic to its square. As a pre\-{liminary},
 we show that a certain symmetric sum of $\Xpw$ is complemented in $\Xpw$ \hfil\break 
 \xciteplus{RI}{Proposition 12}. 
 This symmetric sum is a special case of a more general sum which we now define. 

 Let $2<p<\infty$. For each sequence $v=\set{v_j}$ of positive scalars, define a space $\ell_{2,v}$ as in the proof of 
 Proposition 2.7. For each $k\in\NN$, let $v^{(k)}=\setwlimits{v_j^{(k)}}{j=1}{\infty}$ be a sequence of positive 
 scalars, and let $X_k$ be a closed subspace of $X_{p,v^{(k)}}$.  
 Let \hfil\break 
 $\paren{X_1\oplus X_2\oplus\cdots}_{p,2,\{v^{(k)}\}}$ be the Banach space of all sequences $\set{x_k}$ with 
 $x_k\in X_k$ such that 
 $\norm{\set{x_k}}=\max\set{\paren{\sum\normsubpower{x_k}{\l{p}}{p\vphantom{2}}}^{{1\over p}},
                   \paren{\sum\normsubpower{x_k}{\ell_{2,v^{(k)}}}{2}}^{{1\over 2}}}<\infty$. 
 If each $v^{(k)}$ is identical to a fixed sequence $v$, we will denote   
 $\paren{X_1\oplus X_2\oplus\cdots}_{p,2,\{v^{(k)}\}}$ by 
 $\paren{X_1\oplus X_2\oplus\cdots}_{p,2,v}$. 

\xproclaim {Proposition 2.10}. Let $2<p<\infty$ and let $w=\set{w_n}$ be a sequence of positive scalars satisfying 
                               condition $(*)$ of Proposition 2.1. Let \hfil\break 
                               $\tilde X_{p,w} = \seqsum{\Xpw}{p,2,w}$.  
                               Then $\tilde X_{p,w} \cinjects \Xpw$. 

\proof By condition $(*)$, we may choose a sequence $\set{N_k}$ of disjoint infinite \hfil\break 
       subsets of $\NN$ such that for each $\epsilon>0$ and for each $k$, 
       $$\sum_{{\scriptstyle w_n<\epsilon \atop \scriptstyle n\in N_k}} {w_n}^{{2p\over p-2}} = \infty.$$
       Hence for each $k$, we may choose a sequence 
       $\setwlimits{E_j^{(k)}}{j=1}{\infty}$ of disjoint nonempty finite subsets of $N_k$ such that 
       $${w_j}^{{2p\over p-2}} \le \sum_{n\in E_j^{(k)}} {w_n}^{{2p\over p-2}} \le \paren{2w_j}^{{2p\over p-2}}.$$  
       Then for $v_j^{(k)} = \paren{\sum_{n\in E_j^{(k)}} {w_n}^{{2p\over p-2}}}^{{p-2\over 2p}}$,   
       $w_j \le v_j^{(k)} \le 2w_j$.  
       Hence for $v^{(k)}=\setwlimits{v_j^{(k)}}{j=1}{\infty}$ and $x_k\in\Xpw$, 
       $\normsub{x_k}{\ell_{2,w}}\le\normsub{x_k}{\ell_{2,v^{(k)}}}\le 2\normsub{x_k}{\ell_{2,w}}$. 
       Hence
       $$\paren{X_{p,w} \oplus X_{p,w} \oplus \cdots}_{p,2,w} \sim  
         \paren{X_{p,v^{(1)}} \oplus X_{p,v^{(2)}} \oplus \cdots}_{p,2,\{v^{(k)}\}}   
         \eqno{\TAG{2.8}{2.8}}$$   
       via the formal identity mapping. 

       Let $b_j^{(k)}=\sum_{n\in E_j^{(k)}} {w_n}^{{2\over p-2}} e_n$ (where $\set{e_n}$ is the standard basis of $\Xpw$).  
       Let \hfil\break 
       $\tilde b_j^{(k)} = b_j^{(k)}\Big/\normsub{b_j^{(k)}}{\l{p}}$. Then by part (a) of Proposition 2.7,
       and equations \TAG{2.5}{2.5} and \TAG{2.6}{2.6}, for each $k$ there is an isometry
       $T_k\colon X_{p,v^{(k)}} \to \spansub{\tilde b_j^{(k)}:j\in\NN}{\Xpw}$ with \hfil\break 
       $\normsub{T_k(y_k)}{\l{p}} = \normsub{y_k}{\l{p}}$ and 
       $\normsub{T_k(y_k)}{\ell_{2,w}} = \normsub{y_k}{\ell_{2,v^{(k)}}}$ for $y_k\in X_{p,v^{(k)}}$.
       Hence  
       $$\paren{X_{p,v^{(1)}} \oplus X_{p,v^{(2)}} \oplus \cdots}_{p,2,\{v^{(k)}\}} \sim   
         \paren{\spansub{\tilde b_j^{(1)}}{\Xpw}\oplus\spansub{\tilde b_j^{(2)}}{\Xpw}\oplus\cdots}_{p,2,w}
         \eqno{\TAG{2.9}{2.9}}$$   
       via the isometry $\set{y_k}\mapsto\set{T_k(y_k)}$. 

       The direct sum on the right side of \TAG{2.9}{2.9} should be thought of as an internal \hfil\break 
       direct sum of subspaces of $\Xpw$. We next show that 
       $$\paren{\spansub{\tilde b_j^{(1)}}{\Xpw}\oplus\spansub{\tilde b_j^{(2)}}{\Xpw}\oplus\cdots}_{p,2,w} \sim
         \spansub{\tilde b_j^{(k)}:j,k\in\NN}{\Xpw}
         \eqno{\TAG{2.10}{2.10}}$$   
       via the mapping $\set{z_k}\mapsto\sum z_k$.
       For each $k$, let $z_k=\sum_{j=1}^{\infty} \lambda_j^{(k)} \tilde b_j^{(k)} \in
                              \spansub{\tilde b_j^{(k)}:j\in\NN}{\Xpw}.$ 
       Then by equations \TAG{2.5}{2.5} and \TAG{2.6}{2.6}, and part (a) of Proposition 2.7, we have 
       $$\eqalign{
         \norm{\set{z_k}}
            &=\max\set{\paren{\tsuml_{k=1}^{\infty} \normsubpower{z_k}{\l{p}}{p\vphantom{2}}}^{{1\over p}},
                       \paren{\tsuml_{k=1}^{\infty} \normsubpower{z_k}{\ell_{2,w}}{2}}^{{1\over 2}}} 
         \cr&=\max\set{\paren{\tsuml_{k=1}^{\infty}
                  \normsubpower{\tsuml_{j=1}^{\infty} \lambda_j^{(k)} \tilde b_j^{(k)}}{\l{p}}{p\vphantom{2}}}^{{1\over p}},
                       \paren{\tsuml_{k=1}^{\infty}
                  \normsubpower{\tsuml_{j=1}^{\infty} \lambda_j^{(k)} \tilde b_j^{(k)}}{\ell_{2,w}}{2}}^{{1\over 2}}} 
         \cr&=\max\set{\paren{\tsuml_{k=1}^{\infty} \tsuml_{j=1}^{\infty}
                  \abspower{\lambda_j^{(k)}}{p}}^{{1\over p}},
                       \paren{\tsuml_{k=1}^{\infty} \tsuml_{j=1}^{\infty}
                  \abspower{v_j^{(k)}\lambda_j^{(k)}}{2}}^{{1\over 2}}}
         \cr&=\normsub{\tsuml_{k=1}^{\infty} \tsuml_{j=1}^{\infty} \lambda_j^{(k)} \tilde b_j^{(k)}}{\Xpw}
         \cr&=\normsub{\tsuml_{k=1}^{\infty} z_k}{\Xpw}.}$$
       Hence the mapping $\set{z_k}\mapsto\sum z_k$ is an isometry.

       By part (b) of Proposition 2.7, we have 
       $$\spansub{\tilde b_j^{(k)}:j,k\in\NN}{\Xpw} \cinjects \Xpw.
         \eqno{\TAG{2.11}{2.11}}$$   
       Combining \TAG{2.8}{2.8}, \TAG{2.9}{2.9}, \TAG{2.10}{2.10}, and \TAG{2.11}{2.11} yields 
       $$\seqsum{\Xpw}{p,2,w} \cinjects \Xpw.$$
       \QED 
       
 The complementation of $\tilde X_{p,w}$ in $\Xpw$ is the key to showing that $\Xpw$ is iso\-{morphic} to its square 
 \xciteplus{RI}{Proposition 11}. 

\xproclaim {Proposition 2.11}. Let $2<p<\infty$ and let $w=\set{w_n}$ be a sequence of positive scalars satisfying 
                               condition $(*)$ of Proposition 2.1. Then $\Xpw\sim\Xpw\oplus\Xpw$. 

\proof Let $\tilde X_{p,w}$ be as in Proposition 2.10. Then $\tilde X_{p,w} \cinjects \Xpw$. Let $Y$ be a closed subspace 
       of $\Xpw$ such that $\Xpw\sim\tilde X_{p,w}\oplus Y$. Note that $\tilde X_{p,w}\sim\Xpw\oplus\tilde X_{p,w}$. 
       Hence 
       $$\Xpw\oplus\Xpw \sim \Xpw\oplus\tilde X_{p,w}\oplus Y \sim \tilde X_{p,w}\oplus Y \sim \Xpw.$$ 
       \QED

\remark After noting that $\tilde X_{p,w} \sim \tilde X_{p,w}\oplus\tilde X_{p,w}$, we now see by Lemma 2.8 that 
        $\Xpw \sim \tilde X_{p,w}$. 

 The above results immediately yield the following theorem \xciteplus{RI}{Theorem 13}. 

\xproclaim {Theorem 2.12}. Let $2<p<\infty$ and let $w=\set{w_n}$ and $w'=\set{w'_n}$ be sequences of positive 
                           scalars satisfying condition $(*)$ of Proposition 2.1. Then $\Xpw \sim X_{p,w'}$. 

\proof The spaces $\Xpw$ and $X_{p,w'}$ satisfy the hypotheses of Lemma 2.8. \qquad\QED         

\remark For $p$, $w$, and $w'$ as above, 
        there is a constant $C_p$, depending only on $p$, such that $d\paren{\Xpw,X_{p,w'}}\le C_p$,
        where $d\paren{\Xpw,X_{p,w'}}$ is the Banach-Mazur distance between $\Xpw$ and $X_{p,w'}$

\definition Let $2<p<\infty$.    
            Define $X_p$ to be (the isomorphism type of) $\Xpw$    
            for any sequence $w=\set{w_n}$ of positive scalars satisfying condition $(*)$ of Proposition 2.1.    
            For the conjugate index $q$, define $X_q$ to be the dual of $X_p$.  

By Theorem 2.12, $X_p$ is well-defined.   

\preheadspace
\secondhead{The Isomorphism Type of $X_p$} 
\postheadspace
 We now present results leading to the conclusion that for    
 $2<p<\infty$ and for $w=\set{w_n}$ satisfying condition $(*)$ of Proposition 2.1,
 $\Xpw$ is isomorphically distinct from the previously known $\SL{p}$ spaces. 
 The first result    
 \xciteplus{RI}{Corollary 8} establishes an unusual property of $\Xpw$. 

\xproclaim {Proposition 2.13}. Let $2<p<\infty$ and let $w=\set{w_n}$
                               be a sequence of positive scalars satisfying condition $(*)$ of Proposition 2.1.
                               Then for each $N\in\NN$,  
           \item{(a)} there is a basic sequence $\set{\tilde b_j}$ in $\Xpw$, $2N$-equivalent to the
                      standard basis of $\l{2}$, such that for all distinct $j_1,\ldots,j_N \in \NN$,
                      $\set{\tilde b_{j_1},\ldots,\tilde b_{j_N}}$ is isometrically equivalent
                      to the standard basis of $\finitel{p}{N}$, and   
           \item{(b)} there is a basic sequence $\set{d_j}$ in $X_{p,w}^{*}$, $2N$-equivalent to the
                      standard basis of $\l{2}$, \hfil\break 
                      \indent such that for all distinct $j_1,\ldots,j_N \in \NN$,
                      $\set{d_{j_1},\ldots,d_{j_N}}$ is isometrically equivalent \hfil\break 
                      \indent to the standard basis of $\finitel{q}{N}$,   
                      where $q$ is the conjugate index of $p$.  
          
\proof Fix $N\in\NN$. 
       By condition $(*)$, we may choose a sequence $\set{E_j}$ of disjoint nonempty finite subsets of $\NN$ such that 
       $$\paren{{1\over 2N}}^{{2p\over p-2}}\le\sum\nolimits_{n\in E_j} {w_n}^{{2p\over p-2}} \le{1\over N}.$$
       Define $b_j$, $\tilde b_j$, $v_j$, and $v$ as in Proposition 2.7. 
       Recalling that \hfil\break 
       $v_j=\left(\sum_{n\in E_j} {w_n}^{{2p\over p-2}} \right)^{{p-2\over2p}}$, we have   
       $${1\over 2N}\le v_j \le \paren{{1\over N}}^{{p-2\over2p}} \le 1.$$
       Hence $\xinf{v_j}\ge{1\over2N}>0$, $\xsup{v_j}\le 1$, and $\xsup{v_j}^{{2p\over p-2}}\le{1\over N}$. 

       \item{(a)} By part (a) of Proposition 2.7, $\set{\tilde b_j}$ is a basic sequence in $\Xpw$ which is
                  isometrically equivalent to the standard basis of $X_{p,v}$.
                  Since $\xinf{v_j}>0$ and $\xsup{v_j}\le 1$, the proof of part (a) of Proposition 2.1 shows that 
                  the standard basis of $X_{p,v}$ is equivalent to the standard basis of $\l{2}$, with      
                  $\normsub{x}{X_{p,v}}\within{1}{2N}\normsub{x}{\l{2}}$ for every sequence $x=\set{x_n}$ of scalars.
                  Hence $\set{\tilde b_j}$ in $\Xpw$ is $2N$-equivalent to the standard basis of $\l{2}$.
                  \hfil\break 
                  Let $j_1,\ldots,j_N \in \NN$ be distinct and let $x_1,\ldots,x_N$ be scalars. 
                  Then by H\"older's inequality with conjugate indices $P={p\over2}$ and $Q={p\over p-2}$,  
                  and the fact that \hfil\break 
                  $\xsup{v_j}^{{2p\over p-2}} \le {1\over N}$, we have 
                  $$\eqalign{\tsum_{n=1}^{N} \abspower{v_{j_n} x_n}{2}=   
                             \tsum_{n=1}^{N} \abs{{x_n}^2 {v_{j_n}}^2}&\le    
                             \left(\tsum_{n=1}^{N} \abspower{x_n}{2{p\over2}}\right)^{2\over p} 
                             \left(\tsum_{n=1}^{N} {v_{j_n}}^{2{p\over p-2}}\right)^{p-2\over p}\cr&\le 
                             \left(\tsum_{n=1}^{N} \abspower{x_n}{p}\right)^{2\over p} 
                             \left(\tsum_{n=1}^{N} {1\over N} \right)^{p-2\over p}\cr&=  
                             \left(\tsum_{n=1}^{N} \abspower{x_n}{p}\right)^{2\over p}.}$$ 
                  Thus by part (a) of Proposition 2.7 and the above observation, we have   
                  $$\eqalign{\normsub{\tsum_{n=1}^{N} x_n \tilde b_{j_n}}{\Xpw}&= 
                             \max\set{\left(\tsum_{n=1}^{N} \abspower{x_n}{p} \right)^{{1\over p}},  
                                      \left(\tsum_{n=1}^{N} \abspower{v_{j_n} x_n}{2} \right)^{{1\over2}}}\cr&=   
                             \left(\tsum_{n=1}^{N} \abspower{x_n}{p} \right)^{{1\over p}}.}$$ 
                  Hence $\set{\tilde b_{j_1},\ldots,\tilde b_{j_N}}$ is isometrically equivalent to the 
                  standard basis of $\finitel{p}{N}$. 

       \item{(b)} Define $\ell_{2,w}$ and its inner product $\vector{\phantom{0},\!\phantom{0}}$ as in Proposition 2.7.  
                  Let $d_j=b_j\big/\normsubpower{b_j}{\l{p}}{p-1}$ and consider $d_j$ as an element of $X_{p,w}^*$ 
                  with action $\vector{\phantom{0},d_j}$. Then $\vector{\tilde b_j,d_{j^{\prime}}}=0$ for $j\ne j^{\prime}$, 
                  and 
                      $$\vector{\tilde b_j,d_j}={1\over\norm{b_j}}\xsubpower{\l{p}}{p}\vector{b_j,b_j}=
                        {\normsub{b_j}{\vphantom{1_2}}\over\norm{b_j}^{\vphantom{p}}}
                        {\!{{}_{\ell_{2,w}}^2\atop\!\!\!\!{}_{\l{p}}^p}}=1.$$           
                  Let $\set{\alpha_n}$ be a sequence of scalars and let $j_1,\ldots,j_N\in\NN$ be distinct. 
                  We are trying to prove that 
                      $$\normsub{\sum\nolimits_{n=1}^{\infty} \alpha_n d_n}{X_{p,w}^*} \within{2N}{1}  
                        \paren{\sum\nolimits_{n=1}^{\infty} \abspower{\alpha_n}{2}}^{{1\over2}}$$ 
                  and 
                      $$\normsub{\sum\nolimits_{n=1}^N \alpha_n d_{j_n}}{X_{p,w}^*} = 
                        \paren{\sum\nolimits_{n=1}^N \abspower{\alpha_n}{q}}^{{1\over q}}.$$ 
                  The proofs of these two relationships are quite similar. We introduce a shorthand to allow us to 
                  handle them simultaneously.   
                  Let $\sumprime$ denote $\sum_{n=1}^{\infty}$ in the first setting and $\sum_{n=1}^N$ in the second 
                  setting. Let $\tau_n$ denote $n$ in the first setting and $j_n$ in the second setting. 
                  Then for sequences $\set{\gamma_n}$ of scalars, we have  
                      $$\eqalignno{\normsub{\textstyle\sumprime \alpha_n d_{\tau_n}}{X_{p,w}^*}
                                &= \sup\Bigl\{
                                       \abs{\vector{x,\textstyle\sumprime \alpha_n d_{\tau_n}}} : 
                                       \normsub{x}{\Xpw}=1 
                                       \Bigr\}
                                       & \cr  
                              &\ge \sup\Bigl\{
                                       \big|\bigl\langle\textstyle\sumprime \gamma_n \tilde b_{\tau_n},
                                                    \textstyle\sumprime \alpha_n d_{\tau_n}\bigr\rangle\big| :
                                       \big\|\sumprime \gamma_n\tilde b_{\tau_n}\big\|_{\Xpw} = 1
                                       \Bigr\} \hbox{\qquad}
                                       &\TAG{2.12}{2.12} \cr
                                &= \sup\Bigl\{
                                       \abs{\textstyle\sumprime \gamma_n \bar\alpha_n} : 
                                       \big\|\textstyle\sumprime \gamma_n\tilde b_{\tau_n}\big\|_{\Xpw} = 1
                                       \Bigr\}.
                                       & \cr}$$ 
                  We will show that equality holds at \TAG{2.12}{2.12}. It will then follow by part (a) that  
                      $$\eqalign{\normsub{\textstyle\sum\nolimits_{n=1}^{\infty} \alpha_n d_n}{X_{p,w}^*}
                                    &=
                                 \sup\Bigl\{
                                     \abs{\textstyle\sum\nolimits_{n=1}^{\infty} \gamma_n \bar\alpha_n} : 
                                     \big\|\textstyle\sum\nolimits_{n=1}^{\infty} \gamma_n\tilde b_n \big\|_{\Xpw} = 1 
                                     \Bigr\} 
                                 \cr&{\within{2N}{1}\,}
                                 \sup\Bigl\{
                                     \abs{\textstyle\sum\nolimits_{n=1}^{\infty} \gamma_n \bar\alpha_n} : 
                                     \bigl(\textstyle\sum\nolimits_{n=1}^{\infty} \abspower{\gamma_n}{2}\bigr)^{{1\over2}}=1
                                     \Bigr\}
                                 \cr&=
                                 \bigl(\textstyle\sum\nolimits_{n=1}^{\infty} \abspower{\alpha_n}{2}\bigr)^{{1\over2}} }$$
                  and  
                      $$\eqalign{\big\|\textstyle\sum\nolimits_{n=1}^N \alpha_n d_{j_n} \big\|_{X_{p,w}^*}
                                    &=
                                 \sup\Bigl\{
                                     \big|\textstyle\sum\nolimits_{n=1}^N \gamma_n \bar\alpha_n \big|: 
                                     \big\|\textstyle\sum\nolimits_{n=1}^N \gamma_n\tilde b_{j_n} \big\|_{\Xpw} =1 
                                     \Bigr\} 
                                 \cr&=
                                 \sup\Bigl\{
                                     \big|\textstyle\sum\nolimits_{n=1}^N \gamma_n \bar\alpha_n \big| : 
                                     \bigl(\textstyle\sum\nolimits_{n=1}^N \abspower{\gamma_n}{p}\bigr)^{{1\over p}}=1
                                     \Bigr\}
                                 \cr&=
                                 \bigl(\textstyle\sum\nolimits_{n=1}^N \abspower{\alpha_n}{q}\bigr)^{{1\over q}},}$$
                  which is what we are trying to prove. \hfil\break 
                  We now show that equality holds at \TAG{2.12}{2.12}. It suffices to find a projection \hfil\break 
                  $P^{\prime}\colon\Xpw\to\Xpw$ of norm one which is the set-theoretic restriction to \hfil\break 
                  $\Xpw=\l{p}\cap\ell_{2,w}$
                  of the orthogonal projection $\pi^{\prime}\colon\ell_{2,w}\to\ell_{2,w}$ onto $\spansub{\tilde b_n}{\ell_{2,w}}$ 
                  in the first setting and onto $\xfinitespan \setwlimits{\tilde b_{j_n}}{n=1}{N}$ 
                  in the second setting. For then we will have 
                      $$\eqalign{\phantom{0}&{\phantom{=}}\;\,
                                 \sup\Bigl\{
                                     \abs{\vector{x,\textstyle\sumprime \alpha_n d_{\tau_n}}} : 
                                     \normsub{x}{\Xpw}=1 
                                     \Bigr\}
                                     \cr&= 
                                 \sup\Bigl\{
                                     \abs{\vector{x,\paren{P^{\prime}}^*\paren{\textstyle\sumprime \alpha_n d_{\tau_n}}}} : 
                                     \normsub{x}{\Xpw}=1 
                                     \Bigr\}
                                     \cr&= 
                                 \sup\Bigl\{
                                     \abs{\vector{P^{\prime}(x),\textstyle\sumprime \alpha_n d_{\tau_n}}} :
                                     \normsub{x}{\Xpw}=1
                                     \Bigr\}
                                     \cr&\le 
                                 \sup\Bigl\{
                                     \abs{\vector{P^{\prime}(x),\textstyle\sumprime \alpha_n d_{\tau_n}}} : 
                                     \normsub{P^{\prime}(x)}{\Xpw}=1 
                                     \Bigr\}
                                     \cr&=
                                 \sup\Bigl\{
                                     \big|\bigl\langle\textstyle\sumprime \gamma_n \tilde b_{\tau_n},
                                                  \textstyle\sumprime \alpha_n d_{\tau_n}\bigr\rangle\big| : 
                                     \big\|\textstyle\sumprime \gamma_n\tilde b_{\tau_n} \big\|_{\Xpw} = 1 
                                     \Bigr\},}$$
                  whence equality will hold at \TAG{2.12}{2.12}. Let $P^{\prime}\colon\Xpw\to\Xpw$ be defined by 
                      $$P^{\prime}(x)=\sumprime\vector{x,\normalize{b_{\tau_n}}{\ell_{2,w}}{2}} b_{\tau_n}.$$  
                  In either setting, $P^{\prime}$ is essentially the projection $P$ of part (b) of Proposition 2.7,
                  the only difference between the settings being the choice of $\set{E_j}$ on which the projection is based. 
                  In either setting, $\norm{P^{\prime}}=1$, as can be seen by \TAG{2.7}{2.7}. 
                  Thus equality indeed holds at \TAG{2.12}{2.12}. 
                  \qquad\QED
                             
\medskip
 Following Rosenthal \xcite{RI}, we say that a Banach space $X$ satisfies ${\cal P}_2$ if for each $\epsilon>0$ and each 
 sequence $\set{f_n}$ in $X$ equivalent to the standard basis $\set{e_n}$ of $\l{2}$, there is a subsequence $\set{g_n}$ 
 of $\set{f_n}$ such that $\set{g_n}$ is $(1+\epsilon)$-equivalent to $\set{e_n}$. 

 The following result \xcite{RI} restates part (b) of Proposition 2.13 in terms of ${\cal P}_2$. 


\xproclaim {Corollary 2.14}. Let $2<p<\infty$ and let $w=\set{w_n}$ be a sequence of positive scalars satisfying condition $(*)$ 
                             of Proposition 2.1.  Then $X_{p,w}^*$ is not isomorphic to any Banach space satisfying ${\cal P}_2$. 

\proof Suppose $X_{p,w}^*$ is isomorphic to a Banach space $Y$ satisfying ${\cal P}_2$. 
       Let \hfil\break 
       $K=d\paren{X_{p,w}^*,Y}$, the Banach-Mazur distance between $X_{p,w}^*$ and $Y$.  
       Let $\epsilon>0$.
       Choose $N\in\NN$ such that $\paren{1+\epsilon}\paren{K+\epsilon} < d\paren{\finitel{2}{N},\finitel{q}{N}}$, 
       the Banach-Mazur distance between $\finitel{2}{N}$ and $\finitel{q}{N}$,
       where $q$ is the conjugate index of $p$. 

       Choose a basic sequence $\set{d_j}$ in $X_{p,w}^*$ as in part (b) of Proposition 2.13. 
       Then $\set{d_j}$ is equivalent to the standard basis of $\l{2}$, but for all distinct $j_1,\ldots,j_N\in\NN$, 
       $\set{d_{j_1},\ldots,d_{j_N}}$ is isometrically equivalent to the standard basis of $\finitel{q}{N}$.    

       Choose an isomorphism $T\colon X_{p,w}^* \to Y$ such that $\norm{T\vphantom{T^{-1}}}\norm{T^{-1}}<K+\epsilon$. 
       Let \hfil\break 
       $\set{y_j}=\set{T(d_j)}$. 
       Then $\set{y_j}$ is equivalent to the standard basis of $\l{2}$. 

       Suppose $\set{y_{j_n}}$ is a subsequence of $\set{y_j}$ such that $\set{y_{j_n}}$ is $(1+\epsilon)$-equivalent 
       to the standard basis of $\l{2}$.   
       Then the standard basis of $\finitel{2}{N}$ is $(1+\epsilon)$-equivalent to $\set{y_{j_1},\ldots,y_{j_N}}$, 
       $\set{y_{j_1},\ldots,y_{j_N}}$ is $(K+\epsilon)$-equivalent to $\set{d_{j_1},\ldots,d_{j_N}}$, and 
       $\set{d_{j_1},\ldots,d_{j_N}}$ is isometrically equivalent to the standard basis of $\finitel{q}{N}$. 
       Hence the standard basis of $\finitel{2}{N}$ is $(1+\epsilon)(K+\epsilon)$-equivalent
       to the standard basis of $\finitel{q}{N}$, contrary to the choice of $N$. \qquad\QED 

 It is a fairly routine matter to show that for $2<p<\infty$,
 $\ell_2^*$, $\ell_p^*$, and $\paren{\l{2}\oplus\l{p}}^*$ satisfy ${\cal P}_2$. 
 We will show that for $2<p<\infty$, $\seqsum{\l{2}}{\l{p}}^*$ satisfies ${\cal P}_2$ as well.  
 Thus for $2<p<\infty$, the duals of the classical sequence space $\SL{p}$ spaces satisfy ${\cal P}_2$. 
 It follows that for $2<p<\infty$ and $w$ satisfying condition $(*)$ of Proposition 2.1, 
 $\Xpw$ is isomorphically distinct from the classical sequence space $\SL{p}$ spaces. 
 Rather than take this approach, however, we will show that 
 $\seqsum{\l{2}}{\l{p}}^*$ satisfies ${\cal P}_2$ for $2<p<\infty$
 as a lemma for a somewhat stronger result.   

 The following example \xciteplus{RI}{Sublemma 1} motivates the argument. 

\xproclaim {Example 2.15}. The space $\l{2}$ satisfies ${\cal P}_2$. 

\proof Let $\set{e_n}$ be the standard basis of $\l{2}$. 
       Suppose $\set{f_n}$ is a basic sequence in $\l{2}$ equivalent to $\set{e_n}$.  
       Then $\set{f_n}$ is weakly null, $\xinf\normsub{f_n}{\l{2}}>0$, and $\xsup\normsub{f_n}{\l{2}}<\infty$. 
       Let $\epsilon>0$ and choose $\delta>0$ and $\gamma>0$ such that $(1+\delta)^2<1+\epsilon$ and $(1+\gamma)^2<1+\delta$. 
       By the method of Bessaga and Pe\char32lczy\'nski \xciteplus{B-P}{Theorem 3}, 
       choose a subsequence $\set{g_n}$ of $\set{f_n}$ such that 
       $\set{g_n}$ is $(1+\delta)$-equivalent to a block basic sequence $\set{b_n}$ of $\set{e_n}$.  
       It remains to show that $\set{b_n}$ has a subsequence which is $(1+\delta)$-equivalent to $\set{e_n}$. 

       Note that $\set{b_n}$ is equivalent to $\set{e_n}$, whence $\xinf\normsub{b_n}{\l{2}}>0$ and   
       $\xsup\normsub{b_n}{\l{2}}<\infty$. 
       Choose a subsequence $\set{b_{\alpha(n)}}$ of $\set{b_n}$ such that $0<L=\lim\normsub{b_{\alpha(n)}}{\l{2}}$ exists, with
       $$L{1\over1+\gamma}<\normsub{b_{\alpha(n)}}{\l{2}}<L(1+\gamma)$$ 
       for all $n$. Then for scalars $\lambda_1,\lambda_2,\ldots$, we have 
       $$\normsub{\sum\nolimits_{n=1}^{\infty} \lambda_n b_{\alpha(n)}}{\l{2}}
         =\paren{\sum\nolimits_{n=1}^{\infty} \abspower{\lambda_n}{2} \normsubpower{b_{\alpha(n)}}{\l{2}}{2} }^{{1\over2}}
         \within{1+\gamma}{1+\gamma}L\paren{\sum\nolimits_{n=1}^{\infty} \abspower{\lambda_n}{2} }^{{1\over2}}.$$
       Hence $\set{b_{\alpha(n)}}$ is $(1+\delta)$-equivalent to $\set{e_n}$,  
       but $\set{g_{\alpha(n)}}$ is $(1+\delta)$-equivalent to $\set{b_{\alpha(n)}}$,  
       so $\set{g_{\alpha(n)}}$ is $(1+\epsilon)$-equivalent to $\set{e_n}$. \qquad\QED   

 The following result \xciteplus{RI}{Sublemma 1} is similar, but is more technical than \hfil\break 
 motivational. In our first application, $r=2$.      

\xproclaim {Lemma 2.16}. Let $1\le r<\infty$ and let $X$ be isomorphic to $\l{r}$. 
                         Suppose $\set{f_n}$ \hfil\break 
                         is a sequence in $X$ which is weakly null but not norm null. 
                         Then $\set{f_n}$ has a basic \hfil\break 
                         subsequence equivalent to the standard basis $\set{e_n}$ of $\l{r}$. 

\proof Note that $M=\sup\normsub{f_n}{X}<\infty$ since $\set{f_n}$ is weakly bounded. 
       Let $\set{g_n}$
       be a subsequence of $\set{f_n}$ such that $\xinf\normsub{g_n}{X}>0$. 
       Choose $0<\delta<1$ such that $\delta\le\inf\normsub{g_n}{X}$.  
       Fix an isomorphism $T\colon\l{r}\to X$ and its inverse $S\colon X\to\l{r}$. 

       By the method of Bessaga and Pe\char32lczy\'nski \xciteplus{B-P}{Theorem 3}, 
       choose a basic subsequence $\set{h_n}$ of $\set{g_n}$ such that
       $\set{h_n}$ is equivalent to a block basic sequence $\set{b_n}$ of $\set{T(e_n)}$, with 
       $\normsub{h_n-b_n}{X}<{\delta\over2}$ for each $n$. 
       Then for each $n$,  
       $$\normsub{b_n}{X}\ge\normsub{h_n}{X}-\normsub{h_n-b_n}{X}>\delta-{\delta\over2}={\delta\over2}$$
       and 
       $$\normsub{b_n}{X}\le\normsub{h_n}{X}+\normsub{b_n-h_n}{X}<M+{\delta\over2}.$$

       Hence $\set{S(b_n)}$ is a block basic sequence of $\set{e_n}$, 
       $\xinf\normsub{S(b_n)}{\l{r}}>0$, and \hfil\break 
       $\xsup\normsub{S(b_n)}{\l{r}}<\infty$. 
       Hence $\set{S(b_n)}$ is equivalent to $\set{e_n}$, 
       so $\set{b_n}$ is equivalent to $\set{e_n}$.
       Since $\set{h_n}$ is equivalent to $\set{b_n}$, $\set{h_n}$ is equivalent to $\set{e_n}$. \qquad\QED 

 Let $1\le q<\infty$ and let $N\in\NN$. 
 Let $\Gamma$ be an index set, either $\set{1,\ldots,N}$ or $\NN$. 
 We now introduce some notation for $X=\paren{\sum_{j\in\Gamma}^{\oplus} \l{2} }_{\l{q}(\Gamma)}$, 
 that is, $X=\paren{\l{2}\oplus\cdots\oplus\l{2}}_{\finitel{q}{N}}$ \hfil\break 
          ($N$ summands)
          or $X=\seqsum{\l{2}}{\l{q}}$. 
 Denote a generic $x\in X$ by $\set{x^{(j)}}_{j\in\Gamma}$,   
 with each $x^{(j)}\in\l{2}$. 
 For each $J\in\Gamma$, define $\pi_J\colon X\to\l{2}$ by $\pi_J\paren{\set{x^{(j)}}_{j\in\Gamma} }=x^{(J)}$.  
 Let $\set{e_k}$ be the standard basis of $\l{2}$.
 Let $\set{e_{i,j}}$ be the standard basis of $X$, 
 with $\pi_j(e_{i,j})=e_i$ and $\pi_{j'}(e_{i,j})=0_{\l{2}}$ for $j,j'\in\Gamma$ such that $j\neq j'$. 

 The following somewhat idealized example provides a model to be approximated.

\xproclaim {Example 2.17}. Let $1\le q<\infty$ and let $\Gamma$,  
                           $X=\paren{\sum_{j\in\Gamma}^{\oplus} \l{2} }_{\l{q}(\Gamma)}$,  
                           $\pi_j\colon X\to\l{2}$, and $\set{e_{i,j}}$ be as above.
                           Let $\set{\alpha_j}_{j\in\Gamma}$ be a sequence of nonnegative real numbers such that 
                           $\alpha=\paren{\sum_{j\in\Gamma} {\alpha_j}^q }^{{1\over q}}>0$. 
                           Suppose $\set{b_{[k]}}$ is a basic sequence in $X$ which is disjointly supported 
                           with respect to $\set{e_{i,j}}$, 
                           such that for each $j\in\Gamma$, 
                           $\normsub{\pi_j\paren{b_{[k]}}}{\l{2}}=\alpha_j$ for all $k$. 
                           Then $\set{b_{[k]}}$ is 1-equivalent to the standard basis of $\l{2}$. 


\proof For scalars $\lambda_1,\lambda_2,\ldots$, we have 
       $$\eqalignno{\normsub{\tsuml_{k=1}^{\infty} \lambda_k b_{[k]}}{X}
                 &= \left[\tsuml_{j\in\Gamma}
                          \normsubpower{\pi_j\paren{\tsuml_{k=1}^{\infty} \lambda_k b_{[k]}}}{\l{2}}{q}\right]^{{1\over q}}
              \cr&= \left[\tsuml_{j\in\Gamma} \paren{\tsuml_{k=1}^{\infty}
                          \abspower{\lambda_k}{2}\normsubpower{\pi_j\paren{b_{[k]}}}{\l{2}}{2} }^{{1\over2}q} \right]^{{1\over q}}
                 &\TAG{2.13}{2.13}   
              \cr&= \left[\tsuml_{j\in\Gamma} \paren{\tsuml_{k=1}^{\infty}
                         \abspower{\lambda_k}{2} {\alpha_j}^2 }^{{1\over2}q} \right]^{{1\over q}}
              \cr&= \left[\tsuml_{j\in\Gamma} {\alpha_j}^q\right]^{{1\over q}}
                    \left[\tsuml_{k=1\vphantom{j}}^{\infty\vphantom{N}} \abspower{\lambda_k}{2}\right]^{{1\over2}}
              \cr&= \alpha
                    \left[\tsuml_{k=1}^{\infty} \abspower{\lambda_k}{2}\right]^{{1\over2}}.}$$
       Hence $\set{b_{[k]}}$ is 1-equivalent to the standard basis of $\l{2}$. \qquad\QED

 The following lemma \xciteplus{RI}{Sublemma 3} shows the relevance of Example 2.17 \hfil\break 
 for $\Gamma=\set{1,\ldots,N}$ to the space $\seqsum{\l{2}}{\l{q}}$ for $1 \le q < 2$. 

\xproclaim {Lemma 2.18}. Let $1\le q<2$ and let $X=\seqsum{\l{2}}{\l{q}}$.
                         Denote a generic $x\in X$ by $\set{x^{(1)},x^{(2)},\ldots}$, with $x^{(1)},x^{(2)},\ldots\in\l{2}$. 
                         For each $n\in\NN$, 
                         define $P_n\colon X\to X$ by
                         $P_n\paren{\set{x^{(1)},x^{(2)},\ldots}}=\set{x^{(1)},\ldots,x^{(n)},0,0,\ldots}$
                         and define $Q_n\colon X\to X$ by \hfil\break 
                         $Q_n(x)=x-P_n(x)$. 
                         Suppose $Y$ is a subspace of $X$ isomorphic to $\l{2}$. 
                         Then \hfil\break 
                         $\lim_{n\to\infty}\norm{Q_n|_Y}=0$ and $\lim_{n\to\infty}\norm{P_n|_Y}=1$. 

\proof For each $n\in\NN$, $1-\norm{Q_n|_Y}\le\norm{P_n|_Y}\le1$. 
       Hence it suffices to show that $\lim_{n\to\infty}\norm{Q_n|_Y}=0$. 
       Fix an ordering of the standard basis $\set{e_{i,j}}$ of $X$.   

       Suppose the conclusion is false. 
       Then we may choose $0<\delta<1$ and \hfil\break 
       $y_1,y_2,\ldots\in Y$ of norm one such that
       $\normsub{Q_n\paren{y_n}}{X}\ge\delta$ for each $n$, 
       and (by the \hfil\break 
       reflexivity of $Y$) such that $\set{y_n}$ is weakly convergent.  
       Choose positive integers \hfil\break 
       $n_1<n_2<\ldots$ such that for $k<k'$, 
       $\normsub{Q_{n_{k'}}\paren{y_{n_k}}}{X}<{\delta\over8}$.

       Let $d_k=y_{n_{2k}}-y_{n_{2k-1}}$ and let $T_k=Q_{n_{2k}}$. 
       Then $\set{d_k}$ is weakly null, 
       $$\normsub{T_k\paren{d_k}}{X}\ge\normsub{Q_{n_{2k}}\paren{y_{n_{2k}}}}{X}-\normsub{Q_{n_{2k}}\paren{y_{n_{2k-1}}}}{X}
         >\delta-{\delta\over8}={7\over8}\delta,$$ 
       and for $k<k'$,
       $$\normsub{T_{k'}\paren{d_k}}{X}\le\normsub{Q_{n_{2k'}}\paren{y_{n_{2k}}}}{X}+\normsub{Q_{n_{2k'}}\paren{y_{n_{2k-1}}}}{X}
         <{\delta\over8}+{\delta\over8}={\delta\over4}.$$ 
       Note that $\normsub{d_k}{X}\ge\normsub{T_k(d_k)}{X}>{7\over8}\delta$, whence $\set{d_k}$ is not norm null.  
       Hence by the method of Bessaga and Pe\char32lczy\'nski \xciteplus{B-P}{Theorem 3} and Lemma 2.16,  
       we may choose a subsequence $\set{d_{\alpha(k)}}$ of $\set{d_k}$ such that
       $\set{d_{\alpha(k)}}$ is equivalent to a block basic sequence
       $\set{\tilde d_{\alpha(k)}}$ of the standard basis $\set{e_{i,j}}$ of $X$,
       and such that $\set{d_{\alpha(k)}}$ and $\set{\tilde d_{\alpha(k)}}$
       are equivalent to the standard basis $\set{e_k}$ of $\l{2}$,  
       where $\tilde d_{\alpha(k)} = d_{\alpha(k)} \cdot 1_{\supp \tilde d_{\alpha(k)} }$, \hfil\break 
       $\normsub{d_{\alpha(k)} - \tilde d_{\alpha(k)}}{X} < {\delta\over8}$,  
       and there is a $C>0$ such that for each $K\in\NN$,   
       $$\normsub{\sum\nolimits_{k=1}^K \tilde d_{\alpha(k)}}{X} \le C\normsub{\sum\nolimits_{k=1}^K e_k}{\l{2}} = CK^{1\over2}.$$  
       Hence 
       $$\normsub{T_{\alpha(k)}\paren{\tilde d_{\alpha(k)}}}{X} \ge 
         \normsub{T_{\alpha(k)}\paren{d_{\alpha(k)}\vphantom{\tilde d}}}{X} -
         \normsub{T_{\alpha(k)}\paren{d_{\alpha(k)}-\tilde d_{\alpha(k)}}}{X} >
         {7\over8}\delta-{\delta\over8}={3\over4}\delta,$$
       and for $k<k'$,
       $$\normsub{T_{\alpha(k')}\paren{\tilde d_{\alpha(k)}}}{X} \le
         \normsub{T_{\alpha(k')}\paren{d_{\alpha(k)}\vphantom{\tilde d}}}{X} < {\delta\over4}.$$

       Let $b_{\alpha(k)}=\paren{T_{\alpha(k)}-T_{\alpha(k+1)}\vphantom{\tilde d}}\paren{\tilde d_{\alpha(k)}}$. Then  
       $$\normsub{b_{\alpha(k)}\vphantom{\tilde d}}{X} \ge
         \normsub{T_{\alpha(k)}\paren{\tilde d_{\alpha(k)}}}{X} - \normsub{T_{\alpha(k+1)}\paren{\tilde d_{\alpha(k)}}}{X}  
         >{3\over4}\delta-{\delta\over4}={\delta\over2}.$$ 
       Hence for each $K\in\NN$,
       $$\normsub{\sum\nolimits_{k=1}^K \tilde d_{\alpha(k)}}{X} \ge
         \normsub{\sum\nolimits_{k=1}^K b_{\alpha(k)}}{X} =
         \paren{\sum\nolimits_{k=1}^K \normsubpower{b_{\alpha(k)}}{X}{q} }^{{1\over q}} >
         {\delta\over2}K^{{1\over q}}.$$
       Thus for each $K\in\NN$, ${\delta\over2}K^{{1\over q}} < CK^{{1\over2}}$,  
       which is impossible for sufficiently large $K$. 
       \qquad\QED

 We have laid the groundwork for the following result \xciteplus{RI}{Lemma 10}. 
 

\xproclaim {Proposition 2.19}. Let $1\le q<2$. Then $X=\seqsum{\l{2}}{\l{q}}$ satisfies ${\cal P}_2$.      

\proof Define $\pi_j\colon X\to\l{2}$ and the standard basis $\set{e_{i,j}}$ of $X$
       as in the dis\-{cussion} preceding Example 2.17. 
       Let $\set{e_k}$ be the standard basis of $\l{2}$.
       Fix an ordering of $\set{e_{i,j}}$. 

       Suppose $\set{f_{[k]}}$ is a basic sequence in $X$ equivalent to $\set{e_k}$. 
       Then $\set{f_{[k]}}$ is weakly null, $\xinf\normsub{f_{[k]}}{X}>0$, and $\xsup\normsub{f_{[k]}}{X}<\infty$. 
       Let $\epsilon>0$. Choose $\delta>0$, $\gamma>0$, and $\eta>0$ such that 
       $(1+\delta)^2<1+\epsilon$, $(1+\gamma)^2<1+\delta$, and $\eta={\gamma\over2}$, so that 
       $1+2\eta=1+\gamma$ and ${1+\eta}<1+\gamma$. 

       By the method of Bessaga and Pe\char32lczy\'nski \xciteplus{B-P}{Theorem 3}, 
       choose a subsequence $\set{g_{[k]}}$ of $\set{f_{[k]}}$ such that 
       $\set{g_{[k]}}$ is $(1+\delta)$-equivalent to a block basic sequence $\set{b_{[k]}}$ of $\set{e_{i,j}}$. 
       It remains to show that $\set{b_{[k]}}$ has a subsequence which is $(1+\delta)$-equivalent to $\set{e_k}$. 

       We will choose a subsequence of $\set{b_{[k]}}$ in such a way as to approximate the \hfil\break 
       situation of Example 2.17 
       for $\Gamma=\set{1,\ldots,N}$, after the application of the projection $P_N$ of Lemma 2.18 
       for sufficiently large $N$.  

       Note that $\set{b_{[k]}}$ is equivalent to $\set{e_k}$, whence 
       $\xinf\normsub{b_{[k]}}{X}>0$, $\xsup\normsub{b_{[k]}}{X}<\infty$, and $\spansub{b_{[k]}}{X}\sim\l{2}$. 
       By Lemma 2.18, we may choose $N\in\NN$ such that for all $x\in\spansub{b_{[k]}}{X}$, 
          $${1\over1+\gamma}\normsub{x}{X}\le\normsub{P_N(x)}{X}\le\normsub{x}{X},$$ 
       where $P_N$ is as in Lemma 2.18.
       Choose a subsequence $\set{b_{[\alpha(k)]}}$
       of $\set{b_{[k]}}$ such that for each $j\in\set{1,\ldots,N}$, 
       $L_j=\lim_{k\to\infty} \normsub{\pi_j\paren{b_{[\alpha(k)]}}}{\l{2}}$ exists. 
       Let 
          $$L=\lim_{k\to\infty} \normsub{P_N\paren{b_{[\alpha(k)]}}}{X}= 
              \lim_{k\to\infty} \paren{\tsuml_{j=1}^N \normsubpower{\pi_j\paren{b_{[\alpha(k)]}}}{\l{2}}{q}}^{{1\over q}}= 
              \paren{\tsuml_{j=1}^N {L_j}^q}^{{1\over q}}.$$ 
       Then $L\ge{1\over1+\gamma}\inf\normsub{b_{[\alpha(k)]}}{X}>0$ and some $L_j$ is nonzero. 
       Let $J_1=\set{1\!\le \!j\!\le \!N:L_j>0}$ and $J_0=\set{1\!\le \!j\!\le \!N:L_j=0}$.   
       Choose a subsequence $\set{b_{[\beta(k)]}}$ of $\set{b_{[\alpha(k)]}}$ such that 
       for each $j\in J_1$, 
          $$L_j{1\over1+\eta}<\normsub{\pi_j\paren{b_{[\beta(k)]}}}{\l{2}}<L_j(1+\eta)$$ 
       for all $k$, and for each $j\in J_0$, 
          $$L_j{1\over1+\eta}=0\le\normsub{\pi_j\paren{b_{[\beta(k)]}}}{\l{2}}<{L\eta\over N}$$ 
       for all $k$. 
       Then for scalars $\lambda_1,\lambda_2,\ldots$, we have  
          $$\eqalign{\left[\tsuml_{j=1}^N \paren{\tsuml_{k=1}^{\infty} 
                     \abspower{\lambda_k}{2} \normsubpower{\pi_j\paren{b_{[\beta(k)]}}}{\l{2}}{2} }^{{1\over2}q}\right]^{{1\over q}} 
                &\ge    
                     \left[\tsuml_{j=1}^N \paren{\tsuml_{k=1}^{\infty} 
                     \abspower{\lambda_k}{2} \paren{L_j{1\over1+\eta}}^2}^{{1\over2}q}\right]^{{1\over q}} 
               \cr&=   
                     {1\over1+\eta}\paren{\tsuml_{j=1}^N {L_j}^q}^{{1\over q}}
                              \paren{\tsuml_{k=1}^{\infty} \abspower{\lambda_k}{2}}^{{1\over2}} 
               \cr&=   
                     {1\over1+\eta}L\paren{\tsuml_{k=1}^{\infty} \abspower{\lambda_k}{2}}^{{1\over2}} }$$ 
       and 
          $$\eqalign{\left[\tsuml_{j=1}^N \paren{\tsuml_{k=1}^{\infty} 
                     \abspower{\lambda_k}{2} \normsubpower{\pi_j\paren{b_{[\beta(k)]}}}{\l{2}}{2} }^{{1\over2}q}\right]^{{1\over q}} 
            &\le 
                     \left[\tsuml_{j\in J_1} \paren{\tsuml_{k=1}^{\infty} 
                     \abspower{\lambda_k}{2} \normsubpower{\pi_j\paren{b_{[\beta(k)]}}}{\l{2}}{2} }^{{1\over2}q}\right]^{{1\over q}} 
            \cr&\phantom{=}\,\phantom{1}+
                     \left[\tsuml_{j\in J_0} \paren{\tsuml_{k=1}^{\infty} 
                     \abspower{\lambda_k}{2} \normsubpower{\pi_j\paren{b_{[\beta(k)]}}}{\l{2}}{2} }^{{1\over2}q}\right]^{{1\over q}} 
            \cr&\le
                     \left[\tsuml_{j\in J_1} \paren{\tsuml_{k=1}^{\infty} 
                     \abspower{\lambda_k}{2} \paren{L_j(1+\eta)}^2}^{{1\over2}q}\right]^{{1\over q}} 
            \cr&\phantom{=}\,\phantom{1}+
                     \left[\tsuml_{j\in J_0} \paren{\tsuml_{k=1}^{\infty} 
                     \abspower{\lambda_k}{2} \paren{{L\eta\over N}}^2}^{{1\over2}q}\right]^{{1\over q}} 
            \cr&\le
                     (1+\eta)\paren{\tsuml_{j\in J_1} {L_j}^q}^{{1\over q}}
                              \paren{\tsuml_{k=1}^{\infty} \abspower{\lambda_k}{2}}^{{1\over2}} 
            \cr&\phantom{=}\,\phantom{1}+
                     L\eta   \paren{\tsuml_{k=1}^{\infty} \abspower{\lambda_k}{2}}^{{1\over2}} 
            \cr&\le
                     (1+2\eta)L\paren{\tsuml_{k=1}^{\infty} \abspower{\lambda_k}{2}}^{{1\over2}} }$$ 
       (compare with equation \TAG{2.13}{2.13} and its consequents).
       Noting that 
          $${1\over1+\gamma}\normsub{\tsuml_{k=1}^{\infty} \lambda_k b_{[\beta(k)]}}{X} \le  
            \normsub{P_N\paren{\tsuml_{k=1}^{\infty} \lambda_k b_{[\beta(k)]}}}{X} \le 
            \normsub{\tsuml_{k=1}^{\infty} \lambda_k b_{[\beta(k)]}}{X}$$ 
       by the choice of $N$, and
          $$\eqalign{\normsub{P_N\paren{\tsuml_{k=1}^{\infty} \lambda_k b_{[\beta(k)]}}}{X}
                   &=\left[\tsuml_{j=1}^N
                     \normsubpower{\pi_j\paren{\tsuml_{k=1}^{\infty} \lambda_k b_{[\beta(k)]}}}{\l{2}}{q}\right]^{{1\over q}} 
                \cr&=
                     \left[\tsuml_{j=1}^N \paren{\tsuml_{k=1}^{\infty} 
                     \abspower{\lambda_k}{2}
                     \normsubpower{\pi_j\paren{b_{[\beta(k)]}}}{\l{2}}{2} }^{{1\over2}q}\right]^{{1\over q}} }$$
       (compare with equation \TAG{2.13}{2.13} and its antecedents), 
       we have 
          $${1\over1+\gamma}\normsub{\tsuml_{k=1}^{\infty} \lambda_k b_{[\beta(k)]}}{X} \le  
            (1+2\eta)L\paren{\tsuml_{k=1}^{\infty} \abspower{\lambda_k}{2}}^{{1\over2}}$$ 
       and 
          $${1\over1+\eta}L\paren{\tsuml_{k=1}^{\infty} \abspower{\lambda_k}{2}}^{{1\over2}} \le 
            \normsub{\tsuml_{k=1}^{\infty} \lambda_k b_{[\beta(k)]}}{X}.$$ 
       Hence
          $${1\over(1+\gamma)^2}\normsub{\tsuml_{k=1}^{\infty} \lambda_k b_{[\beta(k)]}}{X} \le  
            L\paren{\tsuml_{k=1}^{\infty} \abspower{\lambda_k}{2}}^{{1\over2}} \le 
            (1+\gamma)^2\normsub{\tsuml_{k=1}^{\infty} \lambda_k b_{[\beta(k)]}}{X}.$$ 
       Thus $\set{b_{[\beta(k)]}}$ is $(1+\delta)$-equivalent to $\set{e_k}$, 
       but $\set{g_{[\beta(k)]}}$ is $(1+\delta)$-equivalent to $\set{b_{[\beta(k)]}}$, 
       so $\set{g_{[\beta(k)]}}$ is $(1+\epsilon)$-equivalent to $\set{e_k}$. \qquad\QED 

 The preceding proposition, together with the following lemma \xcite{RI}, will lead to the main result  
 concerning the isomorphic distinctness of $\Xpw$. 

\xproclaim {Lemma 2.20}. Let $1<q<2$.
                         Suppose $X$ is a Banach space satisfying ${\cal P}_2$. \hfil\break 
                         Suppose $Y$ is isomorphic to a quotient space of $\l{q}$. 
                         Then $Z = X \oplus Y$ satisfies ${\cal P}_2$. 

\proof Let $\set{e_n}$ be the standard basis of $\l{2}$. 
       Suppose $\set{z_n}$ is a basic sequence in $Z$ equivalent to $\set{e_n}$. 
       Let $\epsilon>0$ and choose $\delta>0$ such that $(1+\delta)^2<1+\epsilon$. 

       Express each $z_n$ as $x_n \oplus y_n$ with $x_n\in X$ and $y_n\in Y$.  
       Then there is a bounded linear operator $T\colon\l{2}\to Y$ such that $T(e_n)=y_n$ for all $n$  
       [$e_n\mapsto z_n=x_n\oplus y_n\mapsto y_n$]. 
       The adjoint $T^*$ induces a bounded linear operator from a closed subspace of $\l{p}$ to $\l{2}$,
       where $p$ is the conjugate index of $q$. 
       Hence $T^*$ is compact since $2<p<\infty$ \hfil\break 
       \xciteplus{R}{Theorem A2}.  
       Thus $T$ is compact as well.
       Moreover, $\set{e_n}$ is weakly null. 
       Hence \hfil\break 
       $\lim_{n\to\infty} \normsub{y_n}{Y}=\lim_{n\to\infty} \normsub{T(e_n)}{Y}=0$. 
 
       Choose a subsequence $\set{y_{\alpha(n)}}$ of $\set{y_n}$ such that 
       $\set{z_{\alpha(n)}} = \set{x_{\alpha(n)} \oplus y_{\alpha(n)} }$ is \hfil\break 
       $(1+\delta)$-equivalent to $\set{x_{\alpha(n)}}$. 
       Choose a subsequence $\set{x_{\beta(n)}}$ of $\set{x_{\alpha(n)}}$ such that 
       $\set{x_{\beta(n)}}$ is $(1+\delta)$-equivalent to $\set{e_n}$, 
       as we may since $X$ satisfies ${\cal P}_2$. Then \hfil\break 
       $\set{z_{\beta(n)}} = \set{x_{\beta(n)} \oplus y_{\beta(n)} }$ is $(1+\epsilon)$-equivalent to $\set{e_n}$. 
       \qquad\QED
       
 Finally we present the theorem \xciteplus{RI}{Theorem 9} which (in its corollary) \hfil\break 
 establishes that 
 for $2<p<\infty$ and $w$ satisfying condition $(*)$ of Proposition 2.1, 
 $\Xpw$ is iso\-{morphically} distinct from the classical sequence space $\SL{p}$ spaces.   
 
\xproclaim {Theorem 2.21}. Let $2<p<\infty$ and let $w=\set{w_n}$ be a sequence of positive scalars satisfying 
                           condition $(*)$ of Proposition 2.1. 
                           Let $V$ be a closed subspace of $\l{p}$. 
                           Then $\Xpw$ is not a continuous linear image of $\seqsum{\l{2}}{\l{p}}\oplus V$. 

\proof Equivalently, we show that for $Y$ isometric to a quotient space of $\l{q}$, \hfil\break 
       where $q$ is the conjugate index of $p$,
       $X_{p,w}^*$ is not isomorphic to a closed subspace of \hfil\break 
       $\seqsum{\l{2}}{\l{q}}\oplus Y$.

       Let $Y$ be isometric to a quotient space of $\l{q}$. 
       By Corollary 2.14, $X_{p,w}^*$ is not isomorphic to any Banach space satisfying ${\cal P}_2$. 
       However, $\seqsum{\l{2}}{\l{q}}\oplus Y$ satisfies ${\cal P}_2$ (as do all of its closed subspaces)
       by Proposition 2.19 and Lemma 2.20.
       \qquad\QED

 The following corollary \xciteplus{RI}{Corollary 14} extracts only part of the information available from the preceding theorem.   

\xproclaim {Corollary 2.22}. Let $2<p<\infty$ and let $w=\set{w_n}$ be a sequence of positive scalars satisfying 
                             condition $(*)$ of Proposition 2.1. 
                             Then $\Xpw$ is isomorphically distinct from 
                             $\l{2}$, $\l{p}$, $\l{2}\oplus\l{p}$, and $\seqsum{\l{2}}{\l{p}}$. 

\proof Each of the spaces $\l{2}$, $\l{p}$, $\l{2}\oplus\l{p}$, and $\seqsum{\l{2}}{\l{p}}$
       is a continuous \hfil\break 
       linear image of $\seqsum{\l{2}}{\l{p}}\oplus\l{p}$. 
       However, $\Xpw$ is not such an image, as \hfil\break 
       established by Theorem 2.21. \qquad\QED 

\preheadspace
\secondhead{Complementation and Imbedding Relations for $X_p$}
\postheadspace

 The following lemma \xciteplus{RI}{Corollary 14} distinguishes the isomorphism types  
 of two classical sequence space $\SL{p}$ spaces,   
 and is used in the proof that \hfil\break 
 $\seqsum{\l{2}}{\l{p}} \not\injects X_p$ for $2<p<\infty$.   

\xproclaim {Lemma 2.23}. Let $2<p<\infty$. Then $\seqsum{\l{2}}{\l{p}}\not\hookrightarrow\l{2}\oplus\l{p}$. 

\proof Suppose $T\colon\seqsum{\l{2}}{\l{p}}\to\l{2}\oplus\l{p}$ is an isomorphic imbedding.   
       Let $P\colon\l{2}\oplus\l{p}\to\l{2}\oplus\set{0_{\l{p}}}$ and $Q\colon\l{2}\oplus\l{p}\to\set{0_{\l{2}}}\oplus\l{p}$
       be the obvious projections, with $P+Q=I$, the identity operator on $\l{2}\oplus\l{p}$.     

       For each $N\in\NN$, let $X_N$ be the set of all $s^{(1)}\oplus s^{(2)}\oplus \cdots \in \seqsum{\l{2}}{\l{p}}$
       with $s^{(n)}=0_{\l{2}}$ if $n\le N$.  
       Then each $X_N$ is a subspace of $\seqsum{\l{2}}{\l{p}}$ isometric to $\seqsum{\l{2}}{\l{p}}$.  
       
       We will show that $\xlim_{N\to\infty} \norm{PT|_{X_N}} = 0$.  
       Assuming this for now, \hfil\break 
       $\xlim_{N\to\infty} \norm{P|_{T\paren{X_N}} } = 0$ as well, so we may choose $N\in\NN$ such that \hfil\break 
       $\norm{I|_{T\paren{X_N}}-Q|_{T\paren{X_N}} } = \norm{P|_{T\paren{X_N}} } < 1.$   
       Hence $Q|_{T\paren{X_N}}\colon T\paren{X_N}\to\set{0_{\l{2}}}\oplus\l{p}$ is an isomorphic imbedding,  
       and for an isomorphic imbedding $R\colon\l{2}\to\seqsum{\l{2}}{\l{p}}$,      
       the operator $QTR\colon\l{2}\to\set{0_{\l{2}}}\oplus\l{p}$ is an isomorphic imbedding as well.  
       However, no such imbedding exists, and the lemma will follow.   

       It remains to show that $\xlim_{N\to\infty} \norm{PT|_{X_N}}$ is indeed zero.  
       Suppose \hfil\break 
       $\xlim_{N\to\infty} \norm{PT|_{X_N}} \ne 0$. 
       Then we may choose $\epsilon>0$ and a normalized sequence $\set{x_N}$ with $x_N\in X_N$ such that  
       $\normsub{PT\paren{x_N}}{\l{2}\oplus\set{0}}\ge\epsilon$ for each $N$.  
       Let \hfil\break 
       $\tau_N\colon\seqsum{\l{2}}{\l{p}}\to\seqsum{\l{2}}{\l{p}}$
       be the truncation operator defined by \hfil\break 
       $\tau_N\paren{s^{(1)}\oplus s^{(2)}\oplus \cdots} = 
                     s^{(1)}\oplus\cdots\oplus s^{(N)}\oplus 0_{\l{2}}\oplus 0_{\l{2}}\oplus \cdots $.    
       Choose positive integers \hfil\break 
       $N_1<N_2<\cdots $ such that for $\tilde x_{N_k} = \tau_{N_{k+1}}\paren{x_{N_k}}$,  
       ${1\over2} \le \normsub{\tilde x_{N_k}}{\seqsum{\l{2}}{\l{p}} } \le 1$ and \hfil\break 
       $\normsub{PT\paren{\tilde x_{N_k}} }{\l{2}\oplus\set{0}} \ge {\epsilon\over2}$.  
       Then $\set{\tilde x_{N_k}}$ is equivalent to the standard basis of $\l{p}$.  
       Hence \hfil\break 
       $PT|_{\spansub{\tilde x_{N_k} }{\seqsum{\l{2}}{\l{p}} } }$ 
       induces a bounded linear operator from $\l{p}$ into $\l{2}$,
       so \hfil\break 
       $PT|_{\spansub{\tilde x_{N_k} }{\seqsum{\l{2}}{\l{p}} } }$ 
       must be compact.  
       Hence some subsequence $\set{PT\paren{\tilde x_{N_{k(\alpha)}}}}$ of \hfil\break 
       $\set{PT\paren{\tilde x_{N_k}}}$ converges in norm.
       Since $\set{\tilde x_{N_k}}$ is weakly null,    
       $\set{PT\paren{\tilde x_{N_k}}}$ is weakly null as well.    
       Hence $\set{PT\paren{\tilde x_{N_{k(\alpha)}}}}$ must converge to $0_{\l{2}\oplus\l{p}}$ in norm,
       contrary to \hfil\break 
       $\normsub{PT\paren{\tilde x_{N_k}} }{\l{2}\oplus\set{0}} \ge {\epsilon\over2}$ for all $k$.  
       \qquad\QED 

 We are now ready to see how $X_p$ is related to the classical $\SL{p}$ spaces under the relations $\injects$ and $\cinjects$.   
 Recall that $X \equiv Y$ means $X \injects Y$ and $Y \injects X$.   

\xproclaim {Proposition 2.24}. Let $2<p<\infty$. Then 
                               \item{(a)} $X_p\injects\l{2}\oplus\l{p}$,  
                               \item{(b)} $\l{2}\oplus\l{p}\cinjects X_p$,  
                               \item{(c)} $X_p\equiv\l{2}\oplus\l{p}$,   
                               \item{(d)} $X_p\not\cinjects\l{2}\oplus\l{p}$,   
                               \item{(e)} $\sumsum{\l{2}}{\l{p}}\not\injects X_p$,   
                               \item{(f)} $X_p\not\cinjects\sumsum{\l{2}}{\l{p}}$,   
                               \item{(g)} $\L{p}\not\injects X_p$, and  
                               \item{(h)} parts (b), (d), and (f) hold for $1<p<2$ by duality. 

\proof 

\item{(a)}  
       We norm $\l{2}\oplus\l{p}$ by $\normsub{a\oplus b}{\l{2}\oplus\l{p}} = \max\set{\normsub{a}{\l{2}},\normsub{b}{\l{p}} }$. 
       Let $w=\set{w_n}$ be a sequence of positive scalars satisfying condition $(*)$ of Proposition 2.1.  
       Then $\Xpw\sim X_p$. 
       Define $T\colon\Xpw\to\l{2}\oplus\l{p}$ by $T\paren{\set{x_n}}=\set{w_n x_n}\oplus\set{x_n}$.  
       Then $T$ is an isometry.   
       It follows that $X_p\injects\l{2}\oplus\l{p}$.  


\item{(b)}
       Let $w=\set{w_n}$ be a sequence of positive scalars such that 
       $w_{[1]}=\set{w_{3n-2}}$ satisfies $\xinf w_{3n-2} > 0$,     
       $w_{[2]}=\set{w_{3n-1}}$ satisfies $\sum\paren{w_{3n-1}}^{{2p\over p-2}}<\infty$, and       
       $w_{[3]}=\set{w_{3n}}$ satisfies condition $(*)$ of Proposition 2.1.  
       Then $w$ satisfies condition $(*)$ as well.  
       Hence
       $$X_p\sim\Xpw\sim X_{p,w_{[1]}}\oplus X_{p,w_{[2]}}\oplus X_{p,w_{[3]}} \sim\l{2}\oplus\l{p}\oplus X_p.$$
       It follows that $\l{2}\oplus\l{p}\cinjects X_p$.  

\item{(c)}
       The fact that $X_p\equiv\l{2}\oplus\l{p}$ is an immediate consequence of parts (a) and (b).  

\item{(d)}
       Suppose $X_p\cinjects\l{2}\oplus\l{p}$.  
       Then $X_p$ is a continuous linear image of $\l{2}\oplus\l{p}$, contrary to Theorem 2.21.  
       It follows that $X_p\not\cinjects\l{2}\oplus\l{p}$.  

\item{(e)}
       Suppose $\sumsum{\l{2}}{\l{p}} \injects X_p$. Then   
       $\sumsum{\l{2}}{\l{p}} \injects X_p \injects \l{2}\oplus\l{p}$ by part (a), so
       $\sumsum{\l{2}}{\l{p}} \injects \l{2}\oplus\l{p}$, contrary to Lemma 2.23.  
       It follows that $\sumsum{\l{2}}{\l{p}} \not\injects X_p$.    

\item{(f)}
       Suppose $X_p\cinjects\sumsum{\l{2}}{\l{p}}$.   
       Then $X_p$ is a continuous linear image of $\sumsum{\l{2}}{\l{p}}$, contrary to Theorem 2.21.
       It follows that $X_p\not\cinjects\sumsum{\l{2}}{\l{p}}$.   

\item{(g)}
       Suppose $\L{p}\injects X_p$. 
       Then $\sumsum{\l{2}}{\l{p}} \injects \L{p} \injects X_p$,   
       so $\sumsum{\l{2}}{\l{p}} \injects X_p$, contrary to part (e).    
       It follows that $\L{p}\not\injects X_p$.   

\item{(h)}
       Parts (b), (d), and (f) are the parts involving $\cinjects$.   
       \qquad\QED

\medskip
 Building on diagrams \TAG{1.1}{1.1} and \TAG{1.2}{1.2}, for $2<p<\infty$ we have  
$$\matrix{\ltwo&&&&&&\cr
          &\dsa&&&&&\cr
          &&\ltwo\oplus\lp&\ra&\sumsum{\l{2}}{\l{p}}&\ra&\L{p},\cr
          &\usa&\vequiv&&&&\cr
          \lp&&X_p&&&&\cr}\eqno{\TAG{2.14}{2.14}}$$     
and for $1<p<\infty$ where $p\ne2$, we have   
$$\matrix{\ltwo&&&&\sumsum{\l{2}}{\l{p}}&&\cr
          &\cdsa&&\cusa&&\cdsa&\cr
          &&\ltwo\oplus\lp&&&&\L{p}.\cr
          &\cusa&&\cdsa&&\cusa&\cr
          \lp&&&&X_p&&\cr}\eqno{\TAG{2.15}{2.15}}$$     


%
\preheadspace
\firsthead{The Space $B_p$} 
\postheadspace
 Let $2<p<\infty$. 
 The Banach space $B_p$ is of the form $\paren{X_1 \oplus X_2 \oplus \cdots}_{\l{p}}$, 
 where each space $X_N$ is isomorphic to $\l{2}$, 
 but $\setwlimits{X_N}{N=1}{\infty}$ is chosen so that \hfil\break 
 $\xsup_{N\in\NN} d\paren{X_N,\l{2}}=\infty$,     
 where $d\paren{X_N,\l{2}}$ is the Banach-Mazur distance between $X_N$ and $\l{2}$.   
 Each space $X_N$ is of the form $X_{p,v^{(N)}}$ where $v^{(N)}$ is an appropriately chosen constant sequence.   
 The specifics are outlined below.   
 For the conjugate index $q$, $B_q$ is defined to be the dual of $B_p$.   

\preheadspace
\secondhead{The Space $X_{p,v^{(N)}}$} 
\postheadspace
 
 Let $2<p<\infty$ and fix $N\in\NN$.   
 Let $v_j^{(N)}=\paren{{1\over N}}^{{p-2\over2p}}$ for each $j\in\NN$,  
 and let $v^{(N)}$ be the constant sequence $\setwlimits{v_j^{(N)}}{j=1}{\infty}$.  
 Then $X_{p,v^{(N)}}$ is isomorphic to $\l{2}$ by part (a) of Proposition 2.1.  

 The following observation \xcite{RI} concerning $X_{p,v^{(N)}}$ is analogous to Propositions 2.7 and 2.13,   
 but starts with $v^{(N)}$ and produces $w^{(N)}$ rather that the reverse.     
 The lemma eventually leads to information about $B_p$.    

\xproclaim {Lemma 2.25}. Let $2<p<\infty$ and fix $N\in\NN$.   
                         Let $v^{(N)}=\setwlimits{v_j^{(N)}}{j=1}{\infty}$ where   
                         $v_j^{(N)}=\paren{{1\over N}}^{{p-2\over2p}}$ as above.   
                         Then there is a sequence $w^{(N)}=\set{w^{(N)}_n}_{n=1}^{\infty}$ of \hfil\break 
                         positive scalars satisfying condition $(*)$ of Proposition 2.1, 
                         a basic sequence $\setwlimits{\tilde b_j^{(N)}}{j=1}{\infty}$ in $X_{p,w^{(N)}}$, 
                         and a basic sequence $\setwlimits{d_j^{(N)}}{j=1}{\infty}$ in $X_{p,w^{(N)}}^*$ such that 
              \item{(a)} $\setwlimits{\tilde b_j^{(N)}}{j=1}{\infty}$
                         is isometrically equivalent to the standard basis of $X_{p,v^{(N)}}$, 
              \item{(b)} there is a projection
                         $P_N \colon X_{p,w^{(N)}} \to \spansub{\tilde b_j^{(N)}:j\in\NN}{X_{p,w^{(N)}}}$ of norm one,   
              \item{(c)} $\setwlimits{\tilde b_j^{(N)}}{j=1}{\infty}$ is $2N$-equivalent to the standard basis of $\l{2}$,  
                         but for all distinct \hfil\break 
                         $j_1,\ldots,j_N\in\NN$,   
                         $\set{\tilde b_{j_1}^{(N)},\ldots,\tilde b_{j_N}^{(N)}}$  
                         is isometrically equivalent to the standard basis of $\finitel{p}{N}$, and 
              \item{(d)} $\setwlimits{d_j^{(N)}}{j=1}{\infty}$ is $2N$-equivalent to the standard basis of $\l{2}$,  
                         but for all distinct 
              \item{}    $j_1,\ldots,j_N\in\NN$,   
                         $\set{d_{j_1}^{(N)},\ldots,d_{j_N}^{(N)}}$  
                         is isometrically equivalent to the standard basis
              \item{}    of $\finitel{q}{N}$,
                         where $q$ is the conjugate index of $p$.  
                 
\proof Choose a sequence $\setwlimits{E_j^{(N)}}{j=1}{\infty}$ of disjoint nonempty finite subsets of $\NN$   
       and a sequence $w^{(N)}=\set{w^{(N)}_n}$ of positive scalars satisfying condition $(*)$ of \hfil\break 
       Proposition 2.1 such that for each $j\in\NN$,
       $\sum_{n\in E_j^{(N)}} \paren{w^{(N)}_n}^{{2p \over p-2}} = {1\over N}$.  
       [We may take $E_j^{(N)}$ of cardinality $j$ and $\paren{w^{(N)}_n}^{{2p \over p-2}} = {1 \over jN}$ for $n\in E_j^{(N)}$.]   
       Then for each $j\in\NN$, \hfil\break 
       $v_j^{(N)}=\paren{\sum_{n\in E_j^{(N)}} \paren{w^{(N)}_n}^{{2p \over p-2}} }^{{p-2\over2p}}$.   
       Let $b_j^{(N)}=\sum_{n\in E_j^{(N)}} \paren{w^{(N)}_n}^{{2\over p-2}} e_n$ and \hfil\break 
       $\tilde b_j^{(N)}=b_j^{(N)}\Big/\normsub{b_j^{(N)}}{\l{p}}$  
       (analogous to $b_j$ and $\tilde b_j$ in Proposition 2.7),    
       where $\set{e_n}$ is the standard basis of $X_{p,w^{(N)}}$.     
       Then parts (a) and (b) follow from Proposition 2.7. 

       Note that $\setwlimits{E_j^{(N)}}{j=1}{\infty}$ satisfies the condition in the proof of Proposition 2.13.   
       Let $b_j^{(N)}$ and $\tilde b_j^{(N)}$ be as above (analogous to $b_j$ and $\tilde b_j$ in Proposition 2.13),    
       and let \hfil\break 
       $d_j^{(N)}=b_j^{(N)}\Big/\normsubpower{b_j^{(N)}}{\l{p}}{p-1}$ 
       (analogous to $d_j$ in Proposition 2.13, and considered as an element of $X_{p,w^{(N)}}^*$).   
       Then parts (c) and (d) follow from Proposition 2.13.
       \qquad\QED

\preheadspace
\secondhead{The Space $B_p$} 
\postheadspace

 The following definition was suggested above, but we now present it formally.  

\definition Let $2<p<\infty$. For each $N\in\NN$, let $v^{(N)}=\setwlimits{v_j^{(N)}}{j=1}{\infty}$ where   
            $v_j^{(N)}=\paren{{1\over N}}^{{p-2\over2p}}$ as above.   
            Define $B_p$ to be $\paren{X_{p,v^{(1)}} \oplus X_{p,v^{(2)}} \oplus \cdots}_{\l{p}}$.    
            For the \hfil\break 
            conjugate index $q$, define $B_q$ to be the dual of $B_p$.   

 The following proposition \xcite{RI} is the first step in showing that $B_p$ is an $\SL{p}$ space.    
 The subsequent proposition \xcite{RI} is somewhat stronger. 

\xproclaim {Proposition 2.26}. Let $1<p<\infty$ where $p\ne2$.   
                                Then $B_p\cinjects\L{p}$.   


\proof First suppose $2<p<\infty$.  
       For each $N\in\NN$, let $v^{(N)}$ be as above. 
       Then as in the first part of the proof of Corollary 2.6,
       for each $N\in\NN$ 
       there is a sequence $\set{f_j^{(N)}}_{j=1}^{\infty}$
       of independent symmetric three-valued random variables in $\L{p}$ 
       such that $X_{p,v^{(N)}}\sim\span{f_j^{(N)}:j\in\NN}_{\L{p}}\cinjects\L{p}$, 
       where the isomorphism is uniform in $N$ by the proof of Corollary 2.3,
       and the complementation is uniform in $N$ by Theorem 2.5. 
       Hence 
       $$B_p=\paren{X_{p,v^{(1)}} \oplus X_{p,v^{(2)}} \oplus \cdots}_{\l{p}} \cinjects
             \seqsum{\L{p}}{\l{p}}\sim\L{p},$$  
       and $B_p\cinjects\L{p}$.   
       The result now holds for $1<p<2$ by duality.   
       \qquad\QED

\xproclaim {Proposition 2.27}. Let $1<p<\infty$ where $p\ne2$.   
                                Then $B_p\cinjects\sumsum{X_p}{\l{p}}$.   

\proof First suppose $2<p<\infty$.  
       For each $N\in\NN$ \hglue4pt let $v^{(N)}$, $w^{(N)}$, and \hfil\break 
       $\setwlimits{\tilde b_j^{(N)}}{j=1}{\infty}$ be as in Lemma 2.25. 
       Then by parts (a) and (b) of Lemma 2.25,
       there is a projection
       $P_N \colon X_{p,w^{(N)}} \to \spansub{\tilde b_j^{(N)}:j\in\NN}{X_{p,w^{(N)}}}$ of norm one,   
       and there is an isometry
       $T_N \colon \spansub{\tilde b_j^{(N)}:j\in\NN}{X_{p,w^{(N)}}} \to X_{p,v^{(N)}}$. 
       Thus by the remark following Theorem 2.12,
       for any sequence $w$ satisfying condition $(*)$ of Proposition 2.1, 
       $$\eqalign{B_p =          \paren{X_{p,v^{(1)}} \oplus X_{p,v^{(2)}} \oplus \cdots}_{\l{p}}
                      &\cinjects \paren{X_{p,w^{(1)}} \oplus X_{p,w^{(2)}} \oplus \cdots}_{\l{p}}
                  \cr &\sim      \seqsum{\Xpw}{\l{p}}.}$$
       Hence $B_p\cinjects\seqsum{\Xpw}{\l{p}}$.
       The result now holds for $1<p<2$ by duality.   
       \qquad\QED
       
\remark Alternatively, the proof of parts (a) and (b) of Lemma 2.25 could be slightly modified
        to produce a sequence $w=\set{w_n}$ of positive scalars satisfying \hfil\break 
        condition $(*)$ of Proposition 2.1 
        such that $B_p \cinjects \seqsum{\Xpw}{\l{p}}$,
        without the passage through $\paren{X_{p,w^{(1)}} \oplus X_{p,w^{(2)}} \oplus \cdots}_{\l{p}}$. 

 Let $2<p<\infty$. 
 We will show that $B_p^*$ is not isomorphic to any Banach space \hfil\break 
 satisfying ${\cal P}_2$.   
 This will distinguish $B_p$ isomorphically from  
 $\l{2}$, $\l{p}$, $\l{2}\oplus\l{p}$, and \hfil\break 
 $\seqsum{\l{2}}{\l{p}}$.    
 The proof follows the same pattern as the proof that $X_p^*$ is not \hfil\break 
 isomorphic to any Banach space satisfying ${\cal P}_2$.   
 The following proposition \xcite{RI} is \hfil\break 
 analogous to Proposition 2.13.

\xproclaim {Proposition 2.28}. Let $2<p<\infty$. Then for each $N\in\NN$,  
           \item{(a)} there is a basic sequence $\setwlimits{\dot b_j^{(N)}}{j=1}{\infty}$ in $B_p$,   
                      $2N$-equivalent to the standard basis of $\l{2}$,   
                      such that for all distinct $j_1,\ldots,j_N \in {\NN}$,   
                      $\set{\dot b_{j_1}^{(N)},\ldots,\dot b_{j_N}^{(N)} }$ is 
                      isometrically equivalent to the standard basis of $\finitel{p}{N}$, and   
           \item{(b)} there is a basic sequence $\setwlimits{\dot d_j^{(N)}}{j=1}{\infty}$ in $B_p^*$,
                      $2N$-equivalent to the standard basis
                      \hfil\break\indent 
                      of $\l{2}$, such that for all distinct $j_1,\ldots,j_N \in {\NN}$,
                      $\set{\dot d_{j_1}^{(N)},\ldots,\dot{\tilde d}_{j_N}^{(N)} }$ is isometrically
                      \hfil\break\indent 
                      equivalent to the standard basis of $\finitel{q}{N}$, where $q$ is the conjugate index of $p$.  
          
\proof Fix $N\in\NN$.
       Let $v^{(N)}$, $w^{(N)}$, $\tilde b_j^{(N)}$, and $d_j^{(N)}$ be as in Lemma 2.25.   
       Let $T_N \colon \spansub{\tilde b_j^{(N)}:j\in\NN}{X_{p,w^{(N)}}} \to X_{p,v^{(N)}}$ 
       be the isometry cited in the proof of Proposition 2.27, 
       and let $S_N\colon\span{\tilde b_j^{(N)}:j\in\NN}_{X_{p,w^{(N)}}}^* \to X_{p,v^{(N)}}^*$   
       be the isometry $S_N=\paren{T_N^{-1} }^*$. Let \hfil\break 
       $\iota_N: X_{p,v^{(N)}} \to B_p$ 
       and $\kappa_N: X_{p,v^{(N)}}^* \to B_p^*$
       be the obvious isometric injections.

       Now $\setwlimits{\tilde b_j^{(N)}}{j=1}{\infty}$
       and $\setwlimits{d_j^{(N)}}{j=1}{\infty}$
       have the properties asserted in parts (c) and (d) of Lemma 2.25. 
       Let $\dot b_j^{(N)}=\iota_N\paren{T_N\paren{\tilde b_j^{(N)} } }$.
       Then the sequence $\setwlimits{\dot b_j^{(N)}}{j=1}{\infty}$ in $B_p$
       is isometrically equivalent to $\setwlimits{\tilde b_j^{(N)}}{j=1}{\infty}$,  
       and part (a) follows.  

       Let $\tilde d_j^{(N)}$ be the restriction of $d_j^{(N)}$ to
       $\spansub{\tilde b_j^{(N)}:j\in\NN}{X_{p,w^{(N)}}}$.    
       Then $\setwlimits{\tilde d_j^{(N)}}{j=1}{\infty}$ 
       is isometrically equivalent to $\setwlimits{d_j^{(N)}}{j=1}{\infty}$
       by the argument in the proof of part (b) of Proposition 2.13,  
       where it is shown that equality holds at \TAG{2.12}{2.12}. Let \hfil\break 
       $\dot d_j^{(N)}=\kappa_N\paren{S_N\paren{\tilde d_j^{(N)} } }$.   
       Then the sequence $\setwlimits{\dot d_j^{(N)}}{j=1}{\infty}$ in $B_p^*$ is isometrically
       \hfil\break 
       equivalent to $\setwlimits{\tilde d_j^{(N)}}{j=1}{\infty}$ and $\setwlimits{d_j^{(N)}}{j=1}{\infty}$,   
       and part (b) follows.  
       \qquad\QED

 The proof of the following corollary \xcite{RI} is virtually identical to the proof of Corollary 2.14,  
 with $B_p^*$ replacing $X_{p,w}^*$, $\dot d_j^{(N)}$ replacing $d_j$, and Proposition 2.28 \hfil\break 
 replacing Proposition 2.13.  

\xproclaim {Corollary 2.29}. Let $2<p<\infty$. Then $B_p^*$ is not isomorphic to any Banach space satisfying ${\cal P}_2$.  

 The following theorem \xcite{RI} now follows as in the proof of Theorem 2.21,       
 with $B_p^*$ replacing $X_{p,w}^*$ and Corollary 2.29 replacing Corollary 2.14.  

\xproclaim {Theorem 2.30}. Let $2<p<\infty$ and let $V$ be a closed subspace of $\l{p}$.  
                           Then $B_p$ is not a continuous linear image of $\seqsum{\l{2}}{\l{p}} \oplus V$. 

 The following corollary \xciteplus{RI}{Corollary 14} is analogous to Corollary 2.22.  

\xproclaim {Corollary 2.31}. Let $1<p<\infty$ where $p\ne2$.
                             Then $B_p$ is isomorphically distinct from  
                             $\l{2}$, $\l{p}$, $\l{2}\oplus\l{p}$, and $\seqsum{\l{2}}{\l{p}}$.  
                             In particular, $B_p$ is an $\SL{p}$ space.   

\proof First suppose $2<p<\infty$.  
       Then each of the spaces $\l{2}$, $\l{p}$, $\l{2}\oplus\l{p}$, and $\seqsum{\l{2}}{\l{p}}$   
       is a continuous linear image of $\seqsum{\l{2}}{\l{p}} \oplus \l{p}$,      
       but by Theorem 2.30, $B_p$ is not such an image.   
       Finally, $B_p\cinjects\L{p}$ by Proposition 2.26,   
       but the fact that $B_p\not\sim\l{2}$ has just been established.  
       Hence $B_p$ is an $\SL{p}$ space.  
       The result now holds for $1<p<2$ by duality.  
       \qquad\QED 


 We now know that $B_p$ is isomorphically distinct from the classical sequence space $\SL{p}$ spaces.   
 We present next some results to distinguish $B_p$ isomorphically from $X_p$ and $\L{p}$.   
 The first result \xcite{RI} will distinguish $B_p$ from $X_p$, and the three subsequent results will refine the distinction.    

\xproclaim {Proposition 2.32}. Let $1<p<\infty$ where $p\ne2$. Then $\seqsum{\l{2}}{\l{p}} \cinjects B_p$. 

\proof First suppose $2<p<\infty$.  
       Let $v^{(N)}=\setwlimits{v_j^{(N)}}{j=1}{\infty}$ where $v_j^{(N)}=\paren{{1\over N}}^{p-2\over2p}$ as above.   
       Choose a doubly indexed sequence $\set{E_J^{(N)}}_{J,N\in\NN}$ of disjoint nonempty finite
       subsets of $\NN$ such that for each $J,N\in\NN$, 
       $$\sum\nolimits_{j\in E_J^{(N)}} \paren{v_j^{(N)}}^{{2p \over p-2}} =
         \sum\nolimits_{j\in E_J^{(N)}} {1 \over N} \ge 1.$$
       [We may take $E_J^{(N)}$ of cardinality $N$.]
       Let $u_J^{(N)} = \paren{ \sum\nolimits_{j\in E_J^{(N)}} \paren{v_j^{(N)}}^{{2p \over p-2}} }^{{p-2\over2p}}$ and 
       let $u^{(N)}=\setwlimits{u_J^{(N)}}{J=1}{\infty}$.  
       Then $\xinf u_J^{(N)} \ge 1$.
       Hence by part (a) of Proposition 2.1, and the inequality appearing in its proof, $X_{p,u^{(N)}}$ is isometric to $\l{2}$. 
       Moreover, by Proposition 2.7, $X_{p,u^{(N)}} \cinjects X_{p,v^{(N)}}$,  
       and the implied projection is of norm one.
       Hence  
       $$\seqsum{\l{2}}{\l{p}} \sim \paren{X_{p,u^{(1)}} \oplus X_{p,u^{(2)}} \oplus \cdots}_{\l{p}} \cinjects  
                                    \paren{X_{p,v^{(1)}} \oplus X_{p,v^{(2)}} \oplus \cdots}_{\l{p}} = B_p.$$ 
       The result now holds for $1<p<2$ by duality.
       \qquad\QED

 The following lemma \xcite{RI} is a modification of Lemma 2.18.     
 The proof is virtually identical,  
 with $\l{r}$ replacing $\l{2}$ and $K^{{1\over r}}$ replacing $K^{{1\over2}}$.   

\xproclaim {Lemma 2.33}. Let $1<q<r\le2$ and let $X=\paren{X_{p,v^{(1)}}^*\oplus X_{p,v^{(2)}}^*\oplus\cdots}_{\l{q}}$,    
                         where $p$ is the conjugate index of $q$.
                         Denote a generic $x\in X$ by $\set{x^{(1)},x^{(2)},\ldots}$,  
                         with each $x^{(k)}\in X_{p,v^{(k)}}^*$.     
                         For each $n\in\NN$, define $P_n\colon X\to X$ by \hfil\break 
                         $P_n\paren{\set{x^{(1)},x^{(2)},\ldots}}=\set{x^{(1)},\ldots,x^{(n)},0,0,\ldots}$  
                         and define $Q_n\colon X\to X$ by \hfil\break 
                         $Q_n(x)=x-P_n(x)$.    
                         Suppose $Y$ is a subspace of $X$ isomorphic to $\l{r}$.  
                         Then \hfil\break 
                         $\lim_{n\to\infty} \norm{Q_n|_Y} = 0$ and $\lim_{n\to\infty} \norm{P_n|_Y} = 1$.     

 As a corollary, we have the following \xcite{RI}.   

\xproclaim {Lemma 2.34}. Let $1<q<r<2$. Then $\l{r}\not\injects B_q$.   

\proof Suppose $\l{r}\injects B_q$.   
       Then $\l{r}\injects X=\paren{X_{p,v^{(1)}}^*\oplus X_{p,v^{(2)}}^*\oplus\cdots}_{\l{q}}$,    
       where $p$ is the conjugate index of $q$,    
       since $B_q = B_p^*\sim\paren{X_{p,v^{(1)}}^*\oplus X_{p,v^{(2)}}^*\oplus\cdots}_{\l{q}}$.        
       Let $T\colon\l{r}\to X$ be an isomorphic imbedding and let $Y=T(\l{r})$.  
       For each $n\in\NN$, let $P_n\colon X\to X$ and $Q_n\colon X\to X$ be as in Lemma 2.33,   
       with $P_n+Q_n=I$, the identity operator on
       $X$.  
       By Lemma 2.33, we may choose $N\in\NN$ such that   
       $\norm{I|_Y-P_N|_Y}=\norm{Q_N|_Y}<1$.  
       Hence $P_N|_Y\colon Y\to P_N(Y)$ is an isomorphism.   
       Now $Y\sim\l{r}$ and $P_N(Y)\sim\l{2}$,  
       so $P_N|_Y$ induces an isomorphism between $\l{r}$ and $\l{2}$.   
       However, no such isomorphism exists, and the lemma follows.   
       \qquad\QED

 We state without proof \xciteplus{RII}{Corollary 4.2}.  

\xproclaim {Lemma 2.35}. Let $1<q\le r\le2$. Then $\l{r}\injects X_q$.  

 The following observation \xcite{RI} will distinguish $B_p$ from $\L{p}$.  

\xproclaim {Lemma 2.36}. Let $2<p<\infty$. Then $\sumsum{X_p}{\l{p}}\injects\sumsum{\l{2}}{\l{p}}$.  

\proof By part (a) of Proposition 2.24, $X_p\injects\l{2}\oplus\l{p}$.   
       Hence, letting $\FF$ denote the scalar field,   
       $$\eqalign{\sumsum{X_p}{\l{p}}&\injects \sumsum{\paren{\l{2}\oplus\l{p}}}{\l{p}}           \cr
                                     &\sim     \sumsum{\l{2}}{\l{p}}\oplus\sumsum{\l{p}}{\l{p}}   \cr
                                     &\sim     \sumsum{\l{2}}{\l{p}}\oplus\l{p}                   \cr
                                     &\sim     \sumsum{\l{2}}{\l{p}}\oplus\sumsum{\FF}{\l{p}}  \cr
                                     &\sim     \sumsum{\paren{\l{2}\oplus\FF}}{\l{p}}          \cr
                                     &\sim     \sumsum{\l{2}}{\l{p}}.}$$  
       \QED

 Collecting our results and deducing simple consequences yields the following.  

\xproclaim {Proposition 2.37}. Let $2<p<\infty$. Then 
                               \item{(a)} $B_p\injects\sumsum{\l{2}}{\l{p}}$,  
                               \item{(b)} $\sumsum{\l{2}}{\l{p}}\cinjects B_p$,  
                               \item{(c)} $B_p\equiv\sumsum{\l{2}}{\l{p}}$,  
                               \item{(d)} $B_p\not\cinjects\sumsum{\l{2}}{\l{p}}$,  
                               \item{(e)} $X_p\injects B_p$,   
                               \item{(f)} $B_p\not\injects X_p$,   
                               \item{(g)} $X_p\not\cinjects B_p$,   
                               \item{(h)} $\L{p}\not\injects B_p$, and    
                               \item{(i)} parts (b), (d), and (g) hold for $1<p<2$ by duality. 

\proof 
       \item{(a)} We know $B_p\cinjects\sumsum{X_p}{\l{p}}\injects\sumsum{\l{2}}{\l{p}}$  
                  by Proposition 2.27 and Lemma 2.36.
                  It follows that $B_p\injects\sumsum{\l{2}}{\l{p}}$.  

       \item{(b)} Part (b) is a restatement of Proposition 2.32.  

       \item{(c)} The fact that $B_p\equiv\sumsum{\l{2}}{\l{p}}$ is an immediate consequence of parts (a) and (b).   

       \item{(d)} Suppose $B_p\cinjects\sumsum{\l{2}}{\l{p}}$.  
                  Then $B_p$ is a continuous linear image of $\sumsum{\l{2}}{\l{p}}$, contrary to Theorem 2.30.    
                  It follows that $B_p\not\cinjects\sumsum{\l{2}}{\l{p}}$.  

       \item{(e)} We know $X_p\injects\l{2}\oplus\l{p}\cinjects\sumsum{\l{2}}{\l{p}}\cinjects B_p$
                  by part (a) of Proposition 2.24 and part (b) above.     
                  It follows that $X_p\injects B_p$.   

       \item{(f)} Suppose $B_p\injects X_p$.   
                  Then $\sumsum{\l{2}}{\l{p}}\cinjects B_p\injects X_p\injects\l{2}\oplus\l{p}$ 
                  by part (b) above and part (a) of Proposition 2.24, 
                  so $\sumsum{\l{2}}{\l{p}}\injects\l{2}\oplus\l{p}$, contrary to Lemma 2.23.   
                  It follows that $B_p\not\injects X_p$.   

       \item{(g)}  Suppose $X_p\cinjects B_p$.  
                   Then $X_q\cinjects B_q$, where $q$ is the conjugate index of $p$.  
                   Hence for $1<q<r<2$, $\l{r}\injects X_q\cinjects B_q$ by Lemma 2.35, 
                   so $\l{r}\injects B_q$, contrary to Lemma 2.34.   
                   It follows that $X_p\not\cinjects B_p$.    

       \item{(h)} Suppose $\L{p}\injects B_p$.  
                  Then $\L{p}\injects B_p\injects\sumsum{\l{2}}{\l{p}}$ by part (a) above,
                  so \hfil\break 
                  $\L{p}\injects\sumsum{\l{2}}{\l{p}}$, contrary to \xciteplus{L-P 2}{Theorem 6.1}.    
                  It follows that $\L{p}\not\injects B_p$.  

       \item{(i)} Parts (b), (d), and (g) are the parts involving $\cinjects$. \qquad\QED 


\medskip
 Building on diagrams \TAG{2.14}{2.14} and \TAG{2.15}{2.15}, for $2<p<\infty$ we have  
 $$\matrix{\ltwo&&&&B_p&&\cr
          &\dsa&&&\vequiv&&\cr
          &&\ltwo\oplus\lp&\ra&\sumsum{\l{2}}{\l{p}}&\ra&\L{p},\cr
          &\usa&\vequiv&&&&\cr
          \lp&&X_p&&&&\cr}\eqno{\TAG{2.16}{2.16}}$$     
 and for $1<p<\infty$ where $p\ne2$, we have   
 $$\matrix{\ltwo&&&&\sumsum{\l{2}}{\l{p}}\cra B_p&&\cr
          &\cdsa&&\cusa&&\cdsa&\cr
          &&\ltwo\oplus\lp&&&&\L{p}.\cr
          &\cusa&&\cdsa&&\cusa&\cr
          \lp&&&&X_p&&\cr}\eqno{\TAG{2.17}{2.17}}$$     

\preheadspace
\secondhead{Sums of $B_p$}
\postheadspace

 We now present results leading to the conclusion that $B_p\sim B_p\oplus B_p$ and \hfil\break 
 $\sumsum{B_p}{\l{p}}\sim B_p$.   
 Along the way, we will show that the sequence used in the \hfil\break 
 definition of $B_p$ can be modified to some extent     
 without changing the isomorphism type of the space.     

\xproclaim {Lemma 2.38}. Let $2<p<\infty$. 
                         Let $r=\set{r_n}$ and $s=\set{s_n}$ be sequences of \hfil\break 
                         positive scalars, and suppose that $\inf_{n\in\NN} s_n = 0$.  
                         For each $n\in\NN$, let $r^{(n)}$ be the \hfil\break 
                                                                     constant sequence $\set{r_n,r_n,\ldots}$ and  
                                                let $s^{(n)}$ be the constant sequence $\set{s_n,s_n,\ldots}$.      
                         Let \hfil\break 
                         $B_{p,r}=\paren{X_{p,r^{(1)}}\oplus X_{p,r^{(2)}}\oplus \cdots}_{\l{p}}$ and 
                         $B_{p,s}=\paren{X_{p,s^{(1)}}\oplus X_{p,s^{(2)}}\oplus \cdots}_{\l{p}}$. 
                         Then \hfil\break 
                         $B_{p,r} \cinjects B_{p,s}$.  

\proof Fix a subsequence $\set{s_{\alpha(n)}}$ of $\set{s_n}$ such that for each $n\in\NN$, 
       $s_{\alpha(n)} \le r_n$.
       Let $S_{\alpha(n)} = s_{\alpha(n)}^{\pz2p\over p-2}$ and $R_n = r_n^{\pz2p\over p-2}$.
       Then $S_{\alpha(n)} \le R_n$ for each $n$.  
       Let $\set{K_n}$ be the sequence of positive integers such that for each $n\in\NN$, 
       $$K_n S_{\alpha(n)} \le R_n < (K_n+1) S_{\alpha(n)} \le 2K_n S_{\alpha(n)} \le 2^{2p\over p-2} K_n S_{\alpha(n)}.$$ 


       Fix $n\in\NN$.  
       Let $\setwlimits{E_j^{(n)}}{j=1}{\infty}$ be a sequence of disjoint subsets of $\NN$
       such that each $E_j^{(n)}$ has cardinality $K_n$.  
       Then for each $j\in\NN$, 
       $$\sum_{k\in E_j^{(n)}} S_{\alpha(n)} \le R_n < 2^{2p\over p-2} \sum_{k\in E_j^{(n)}} S_{\alpha(n)}.$$
       Let $t_n=\paren{\sum_{k\in E_j^{(n)}} S_{\alpha(n)} }^{p-2\over2p}$ [which does not depend on $j$]. 
       Then $t_n \le r_n < 2t_n$.  
       Hence for $t^{(n)}=\set{t_n,t_n,\ldots}$ and $x \in X_{p,t^{(n)}}$, 
       $\norm{x}_{X_{p,t^{(n)}}} \le \norm{x}_{X_{p,r^{(n)}}} \le 2\norm{x}_{X_{p,t^{(n)}}}$.  
       Thus $X_{p,r^{(n)}}\sim X_{p,t^{(n)}}$ via the formal identity mapping.  
       Moreover, 
       $X_{p,t^{(n)}} \cinjects X_{p,s^{(\alpha(n))}}$ by Proposition 2.7, 
       where the implied projection is of norm one.  

       Release $n$ as a free variable.  
       Then for each $n\in\NN$,  
       $X_{p,r^{(n)}} \sim X_{p,t^{(n)}} \cinjects X_{p,s^{(\alpha(n))}}$,  
       where the isomorphism $X_{p,r^{(n)}} \sim X_{p,t^{(n)}}$ is uniform in $n$.   
       It follows that 
       $$\eqalign{B_{p,r}=\paren{X_{p,r^{(1)}}\oplus X_{p,r^{(2)}}\oplus \cdots}_{\l{p}}
                    &\sim \paren{X_{p,t^{(1)}}\oplus X_{p,t^{(2)}}\oplus \cdots}_{\l{p}}
           \cr &\cinjects \paren{X_{p,s^{(\alpha(1))}}\oplus X_{p,s^{(\alpha(2))}}\oplus \cdots}_{\l{p}}
           \cr &\cinjects \paren{X_{p,s^{(1)}}\oplus X_{p,s^{(2)}}\oplus \cdots}_{\l{p}}
                    \cr &=B_{p,s}.}$$
       \QED

\remark For $2<p<\infty$, the space $B_p$ is of the form $B_{p,s}$
        where $s=\set{s_n}$ and $B_{p,s}$ are as above, with $\xinf_{n\in\NN} s_n =0$.     

\xproclaim {Lemma 2.39}.  Let $2<p<\infty$.  
                          Let $r=\set{r_n}$, $r^{(n)}$, and $B_{p,r}$ be as in Lemma 2.38.  
                          Then $B_{p,r} \sim B_{p,r} \oplus B_{p,r}$.  

\proof 
       Recall that $B_{p,r}=\paren{X_{p,r^{(1)}}\oplus X_{p,r^{(2)}}\oplus \cdots}_{\l{p}}$.
       For each $n\in\NN$, let $\setwlimits{z_k^{(n)}}{k=1}{\infty}$ represent an element of $X_{p,r^{(n)}}$.      
       Define a projection $P\colon B_{p,r} \to B_{p,r}$ by
       $P\paren{\set{z_k^{(1)}}\oplus\set{z_k^{(2)}}\oplus\cdots}=\paren{\set{x_k^{(1)}}\oplus\set{x_k^{(2)}}\oplus\cdots}$, 
       where for $k,n\in\NN$, 
       $x_k^{(n)}=z_k^{(n)}$ if $k$ is even and $x_k^{(n)}=0$ if $k$ is odd.  
       Then the image of $B_{p,r}$ under $P$ is isomorphic to $B_{p,r}$, 
       as is the kernel of $P$.  
       Hence $B_{p,r} \sim B_{p,r} \oplus B_{p,r}$.  
       \qquad\QED


 By the remark above, we have the following corollary (true for $1<p<2$ by duality) of Lemma 2.39.   

\xproclaim {Corollary 2.40}. Let $1<p<\infty$ where $p\ne2$.  
                              Then $B_p \sim B_p \oplus B_p$.    

 We also have the following corollary of Lemmas 2.38 and 2.39.  

\xproclaim {Corollary 2.41}. Let $2<p<\infty$.  
                               Let $r=\set{r_n}$ and $s=\set{s_n}$ be sequences of \hfil\break 
                               positive scalars such that
                               $\inf_{n\in\NN} r_n = 0$ and 
                               $\inf_{n\in\NN} s_n = 0$.  
                               Let $r^{(n)}$, $s^{(n)}$, $B_{p,r}$, and $B_{p,s}$ be as in Lemma 2.38.    
                               Then $B_{p,r} \sim B_{p,s}$.  

\proof The spaces $B_{p,r}$ and $B_{p,s}$ satisfy the hypotheses of Lemma 2.8.  
       \qquad\QED

\xremark{1} Recalling the remark above,  
            one consequence of Corollary 2.41 is that for $2<p<\infty$, and for $1<p<2$ by duality,
            the isomorphism type of $B_p$ does not depend on the specific sequence
            $\set{\paren{1\over N}^{{p-2\over2p}}}_{N=1}^{\infty}$ used in its definition,   
            but only on the fact that the infimum of the sequence is zero.   

\xremark{2} Let $2<p<\infty$.  
            Then $B_p$ is of the form 
            $\paren{X_{p,w^{(1)}} \oplus X_{p,w^{(2)}} \oplus\cdots}_{\l{p}}$ 
            where for each $N\in\NN$, 
            $w^{(N)}$ is a sequence $\setwlimits{w_k^{(N)}}{k=1}{\infty}$ of positive scalars.  
            The above remark gives a sufficient condition for 
            $B_p \sim \paren{X_{p,w^{(1)}} \oplus X_{p,w^{(2)}} \oplus\cdots}_{\l{p}}$ 
            in the case where each $w^{(N)}$ is a constant sequence.  
            Although the details will not be given, \hfil\break 
            $B_p \sim \paren{X_{p,w^{(1)}} \oplus X_{p,w^{(2)}} \oplus\cdots}_{\l{p}}$ 
            if and only if the following two conditions hold: \hfil\break 
            (a) for each $N\in\NN$, $w^{(N)}$ fails condition $(*)$ of Proposition 2.1, and 
            (b) there is an increasing sequence $\setwlimits{\alpha(N)}{N=1}{\infty}$ of positive integers 
                and a sequence $\setwlimits{S_N}{N=1}{\infty}$ of \hfil\break 
                infinite subsets of $\NN$ 
                such that for each $N\in\NN$,  
                $c_N = \lim\inf_{k\in S_N} w_k^{(\alpha(N))} > 0$, \hfil\break 
                but $\xlim_{N\to\infty} c_N = 0$.  

 Just as $B_p \oplus B_p \sim B_p$,    
 $\seqsum{B_p}{\l{p}} \sim B_p$, as shown below.   

\xproclaim {Corollary 2.42}. Let $1<p<\infty$ where $p\ne2$.  
                             Then $\seqsum{B_p}{\l{p}} \sim B_p$. 


\proof First suppose that $2<p<\infty$.  
       Then $B_p$ is of the form $B_{p,s}$ where \hfil\break 
       $s=\set{s_n}$ satisfies $\xinf_{n\in\NN} s_n = 0$, 
       and $s^{(n)}$ and $B_{p,s}$ are as in Lemma 2.38.  
       Let $S$ be the sequence $\set{s_1;s_1,s_2;s_1,s_2,s_3;\ldots} = \setwlimits{\setwlimits{s_n}{n=1}{T}}{T=1}{\infty}$.  
       Then $S$ has infimum zero as well.  
       Hence $B_{p,S} \sim B_{p,s}$ by Corollary 2.41.    
       It follows that 
       $$\eqalign{\seqsum{B_{p,s}}{\l{p}}
               &= \seqsum{ \paren{X_{p,s^{(1)}} \oplus X_{p,s^{(2)}} \oplus \cdots}_{\l{p}} }{\l{p}}
        \cr &\sim \paren{\textstyle\sum\nolimits_{T\in\NN}^{\oplus}
                         \textstyle\sum\nolimits_{1\le n\le T}^{\oplus} X_{p,s^{(n)}} }_{\l{p}}
        \cr &\sim B_{p,S}
        \cr &\sim B_{p,s}.}$$
       The result now holds for $1<p<2$ by duality.  
       \qquad\QED

\preheadspace
\firsthead{Sums Involving $X_p$ or $B_p$}
\postheadspace
 As observed by Rosenthal \xcite{RI}, a few more $\SL{p}$ spaces can be constructed by \hfil\break 
 forming sums involving $X_p$ or $B_p$.    
 The resulting spaces are $\sumsum{\l{2}}{\l{p}}\oplus X_p$, $B_p\oplus X_p$, and $\sumsum{X_p}{\l{p}}$.     
 The following proposition \xcite{RI} shows that these spaces cannot be \hfil\break 
 distinguished by the relation $\injects$.  

\xproclaim {Proposition 2.43}. Let $2<p<\infty$. Then  
                               \item{(a)} $B_p\oplus X_p\cinjects\sumsum{X_p}{\l{p}}$  
                                          (whence the same is true for $1<p<2$ by duality),
                               \item{(b)} $\sumsum{\l{2}}{\l{p}}$, $B_p$,
                                          $\sumsum{\l{2}}{\l{p}}\oplus X_p$, $B_p\oplus X_p$,
                                          and $\sumsum{X_p}{\l{p}}$ are equivalent under $\equiv$, and 
                               \item{(c)} letting $Y$ denote any of the five spaces of part (b) and 
                                          letting $X$ denote either
                               \item{}    $\l{2}\oplus\l{p}$ or $X_p$, we have  
                                          $X\injects Y\injects\L{p}$ but $\L{p}\not\injects Y\not\injects X$. 

\proof
       \item{(a)}
             By Proposition 2.27,
             we have $B_p\oplus X_p \cinjects \sumsum{X_p}{\l{p}}\oplus X_p \sim \sumsum{X_p}{\l{p}}$.


       \item{(b)}
             Consider the chains 
             $$\sumsum{\l{2}}{\l{p}}\cinjects B_p\cinjects B_p\oplus X_p\cinjects\sumsum{X_p}{\l{p}}$$
             and 
             $$\sumsum{\l{2}}{\l{p}}\cinjects\sumsum{\l{2}}{\l{p}}\oplus X_p\cinjects B_p\oplus X_p\cinjects\sumsum{X_p}{\l{p}}$$
             established by part (b) of Proposition 2.37 and part (a) above.  
             Now \hfil\break 
             $\sumsum{X_p}{\l{p}}\injects\sumsum{\l{2}}{\l{p}}$ by Lemma 2.36, which completes each of the two cycles.  
             It follows that the listed spaces are equivalent under $\equiv$.    

       \item{(c)}
             We know $\l{2}\oplus\l{p}\cinjects\sumsum{\l{2}}{\l{p}}\cinjects\L{p}$ 
             but $\L{p}\not\injects\sumsum{\l{2}}{\l{p}}\not\injects\l{2}\oplus\l{p}$
             as in the discussion of diagrams \TAG{1.1}{1.1} and \TAG{1.2}{1.2}.  
             The result now follows from the fact that 
             $X\equiv\l{2}\oplus\l{p}$ by part (c) of Proposition 2.24 and 
             $Y\equiv\sumsum{\l{2}}{\l{p}}$ by part (b) above.  
       \qquad\QED

\medskip
 Building on diagram \TAG{2.16}{2.16}, for $2<p<\infty$ we have  
$$\matrix{\ltwo&&&&B_p&&&&&&\cr
          &\dsa&&&\vequiv&&&&&&\cr
          &&\ltwo\oplus\lp&\ra&\sumsum{\ltwo}{\lp}&\equiv&B_p\oplus X_p&\equiv&\sumsum{X_p}{\lp}&\ra&\L{p}.\cr
          &\usa&\vequiv&&\vequiv&&&&&&\cr
          \lp&&X_p&&\sumsum{\ltwo}{\lp}\oplus X_p&&&&&&\cr}\eqno{\TAG{2.18}{2.18}}$$

 As we have seen, the relation $\injects$ is inadequate to distinguish
 $\sumsum{\l{2}}{\l{p}} \oplus X_p$,  
 $B_p \oplus X_p$, and 
 $\sumsum{X_p}{\l{p}}$ isomorphically.   
 We will distinguish these three spaces via the relation $\cinjects$.  
 The next three results will distinguish $B_p \oplus X_p$ and $\sumsum{X_p}{\l{p}}$.   
 The first result is a corollary of Lemma 2.34.   

\xproclaim {Lemma 2.44}. Let $1<q<r<2$. 
                         Suppose $S \colon \l{r} \to B_q$ is a bounded linear operator. 
                         Then given a sequence $\set{\epsilon_n}$ of positive scalars,   
                         there is a normalized block basic sequence $\set{x_n}$ 
                         of the standard basis $\set{e_k}$ of $\l{r}$ such that
                         $\normsub{S(x_n)}{B_q}<\epsilon_n$ for each $n\in\NN$.    


\proof It suffices to show that 
       there is a normalized block basic sequence $\set{x_n}$ 
       of the standard basis $\set{e_k}$ of $\l{r}$ such that
       $\normsub{S(x_n)}{B_q}\le{\norm{S}\over n}$ for each $n\in\NN$,    
       for the result will then follow upon passing to an appropriately chosen subsequence of $\set{x_n}$.  

       We define $\set{x_n}$ by induction, 
       where each $x_n$ is of the form $\sum_{k\in E_n} \lambda_k e_k$, 
       each $E_n$ is a finite subset of $\NN$, 
       each $\set{\lambda_k \colon k \in E_n}$ is a set of nonzero scalars, and \hfil\break 
       $\xmax E_i < \xmin E_j$ for $1\le i<j$.   

       Let $x_1 = \sum_{k\in E_1} \lambda_k e_k$ be a normalized block of $\set{e_k}$. 
       Then $\normsub{S(x_1)}{B_q}\le{\norm{S}\over1}$. \hfil\break 
       Suppose normalized disjointly supported blocks $x_1,\ldots,x_N$ have been chosen, 
       where $x_n=\sum_{k\in E_n} \lambda_k e_k$ and $\normsub{S(x_n)}{B_q}\le{\norm{S}\over n}$ for each $1\le n\le N$,    
       and $\xmax E_i < \xmin E_j$ for $1\le i<j\le N$.   
       Let $M=\max E_N$. 
       Then as we verify below, 
       we may choose $x_{N+1} \in \finitespan\set{e_k\colon k\ge M+1}$ of norm one 
       such that $\normsub{S\paren{x_{N+1}}}{B_q} \le {\norm{S} \over N+1}$.

       Suppose for a moment that no such $x_{N+1}$ exists.  
       Let $X_{M+1}=\spansub{e_k\colon k\ge M+1}{\l{r}}$, which is isometric to $\l{r}$.  
       Then for each normalized $x\in X_{M+1}$, 
       $\normsub{S(x)}{B_q} > {1\over2}{\norm{S} \over N+1}$.
       Hence $S|_{X_{M+1}}$ induces an isomorphic imbedding of $\l{r}$ into $B_q$. 
       However, by Lemma 2.34, no such imbedding exist. 
       Thus $x_{N+1}$ can be chosen as claimed, and the result follows. 
       \qquad\QED 

\xproclaim {Lemma 2.45}. Let $1<q<r<2$. Then $\sumsum{\l{r}}{\l{q}} \not\injects B_q \oplus X_q$.  

\proof Suppose $\sumsum{\l{r}}{\l{q}} \injects B_q \oplus X_q$.  
       Let $T\colon\sumsum{\l{r}}{\l{q}} \to B_q \oplus X_q$ be an isomorphic imbedding.      
       Let $Q\colon B_q \oplus X_q \to B_q \oplus \set{0_{X_q}}$ be the obvious projection.  
       Then $QT\colon\sumsum{\l{r}}{\l{q}} \to B_q \oplus \set{0_{X_q}}$ is a bounded linear operator.    

       We will show that there is a subspace $X$ of $\sumsum{\l{r}}{\l{q}}$, isometric to $\sumsum{\l{r}}{\l{q}}$,   
       such that $\norm{Q|_{T(X)}}<1$,    
       whence $(I-Q)|_{T(X)}$ induces an isomorphic imbedding of $\sumsum{\l{r}}{\l{q}}$ into $X_q$.   
       However by \xciteplus{S}{Proposition 2}, presented below as Lemma 3.7,
       no such imbedding exists, and the lemma will follow.     


       Let $\set{e_{m,n}}$ be the standard basis of $\sumsum{\l{r}}{\l{q}}$,   
       where for each $n\in\NN$,
       $\setwlimits{e_{m,n}}{m=1}{\infty}$ is isometrically equivalent to the standard basis of $\l{r}$.    
       By Lemma 2.44, for each $n\in\NN$ we may choose a normalized block basic 
       sequence $\set{x_k^{(n)}}_{k=1}^{\infty}$ of $\set{e_{m,n}}_{m=1}^{\infty}$ such that  
       $\normsub{QT\paren{x_k^{(n)}}}{B_q} < {1\over\norm{T^{-1}}2^{k+n}}$.   
       Let 
       $X=\span{x_k^{(n)}\colon k,n\in\NN}$.  
       Then $X$ is isometric to $\sumsum{\l{r}}{\l{q}}$.  
       Let 
       $\set{\lambda_k^{(1)}} \oplus \set{\lambda_k^{(2)}} \oplus \cdots \in \seqsum{\l{r}}{\l{q}}$ 
       be of norm one.  
       Then  
       $$\eqalign{\normsub{QT\paren{\tsuml_{n=1}^{\infty} \tsuml_{k=1}^{\infty} \lambda_k^{(n)} x_k^{(n)} }}{B_q} &=
                  \normsub{\tsuml_{n=1}^{\infty} \tsuml_{k=1}^{\infty} \lambda_k^{(n)} QT\paren{x_k^{(n)}} }{B_q} \cr &\le
                  \tsuml_{n=1}^{\infty} \tsuml_{k=1}^{\infty} \normsub{QT\paren{x_k^{(n)}} }{B_q} \cr &<
                  \tsuml_{n=1}^{\infty} \tsuml_{k=1}^{\infty} {1\over\norm{T^{-1}}2^{k+n}} \cr &=   
                  {1\over\norm{T^{-1}}}.}$$
       Hence $\big\|QT|_X\big\|<{1\over\norm{T^{-1}}}$,  
       so $\norm{Q|_{T(X)}} \le \norm{T^{-1}} \big\|QT|_X\big\| < 1$.    
       Thus $(I-Q)|_{T(X)}$ induces an isomorphic imbedding of $\sumsum{\l{r}}{\l{q}}$ into $X_q$,   
       where $I$ is the formal identity mapping, but no such imbedding exists. 
       \qquad\QED

\xproclaim {Proposition 2.46}. Let $1<p<\infty$ where $p\ne2$. Then $\sumsum{X_p}{\l{p}}\not\cinjects B_p\oplus X_p$. 

\proof First let $1<q<2$ and suppose $\sumsum{X_q}{\l{q}}\cinjects B_q\oplus X_q$. 
       For \hfil\break 
       $1<q<r<2$, $\l{r}\injects X_q$ by Lemma 2.35,
       so $\sumsum{\l{r}}{\l{q}}\injects\sumsum{X_q}{\l{q}}\cinjects B_q\oplus X_q$. \hfil\break 
       Hence $\sumsum{\l{r}}{\l{q}}\injects B_q\oplus X_q$, contrary to Lemma 2.45. 
       It follows that \hfil\break 
       $\sumsum{X_q}{\l{q}}\not\cinjects B_q\oplus X_q$. 
       The result now holds for $2<p<\infty$ by duality.   
       \qquad\QED       

 The next two results will distinguish $\sumsum{\l{2}}{\l{p}} \oplus X_p$ and $B_p \oplus X_p$ isomorphically.   
 The lemma isolates some preliminary calculations.  

\xproclaim {Lemma 2.47}. Let $2<p<\infty$ with conjugate index $q$, and let $n\in\NN$. 
                         Let $X_{p,v^{(n)}}$ be as in the definition of $B_p$,  
                         and let $v_n$ denote $\paren{{1\over n}}^{{p-2\over2p}}$,
                         the value taken by the constant sequence $v^{(n)}$.   
                         Let ${\cal B}_n$ be the closed unit ball of $X_{p,v^{(n)}}$. 
                         Then for $M_n\in\NN$
                         such that $M_n\le v_n^{-{2p\over p-2}} = n$,   
                         $$\sup_{\set{d_m}\in{\cal B}_n} \abs{\tsuml_{m=1}^{M_n} d_m} = M_n^{{1\over q}}.$$   
                         Moreover, for $K\in\NN$ and $\set{\lambda_k}\in\l{2}$,   
                         $$\sup_{\set{d_{k,\ell}}\in{\cal B}_n} \abs{\tsuml_{k=1}^K \tsuml_{\ell=1}^n \lambda_k d_{k,\ell}}
                           = n^{1\over q} \paren{\tsuml_{k=1}^K \abs{\lambda_k}^2}^{1\over2}.$$   

\proof Let $M\in\NN$ and let $\set{d_m}$ be a sequence of scalars.  
       Then by H\"older's inequality, 
       $$\eqalign{\abs{\tsuml_{m=1}^M d_m}
                = \abs{\tsuml_{m=1}^M 1 d_m}
             &\le \paren{\tsuml_{m=1}^M 1^q}^{{1\over q}} \paren{\tsuml_{m=1}^M \abs{d_m}^p}^{{1\over p}}
           \cr &= M^{{1\over q}} \paren{\tsuml_{m=1}^M \abs{d_m}^p}^{{1\over p}} }$$
       and 
       $$\eqalign{\abs{\tsuml_{m=1}^M d_m}
                = {1\over v_n} \abs{\tsuml_{m=1}^M 1 d_m v_n} 
             &\le {1\over v_n} \paren{\tsuml_{m=1}^M 1^2}^{{1\over2}} \paren{\tsuml_{m=1}^M \abs{d_m v_n}^2}^{{1\over2}}
           \cr &= {1\over v_n} M^{{1\over 2}} \paren{\tsuml_{m=1}^M \abs{d_m v_n}^2}^{{1\over2}}.}$$
       Suppose $\set{d_m}\in{\cal B}_n$.  
       Then
       $\paren{\sum_{m=1}^M \abs{d_m}^p}^{{1\over p}} \le 1$ and 
       $\paren{\sum_{m=1}^M \abs{d_m v_n}^2}^{{1\over2}} \le 1$.
       Hence 
       $\abs{\sum_{m=1}^M d_m} \le M^{1\over q}$ and  
       $\abs{\sum_{m=1}^M d_m} \le {1\over v_n} M^{1\over 2}$.     
       It follows that  
       $$\sup_{\set{d_m}\in{\cal B}_n} \abs{\tsuml_{m=1}^M d_m} \le \min\set{M^{{1\over q}}, {1\over v_n} M^{{1\over2}} }.$$
         
       Let $M_n\in\NN$ such that $M_n\le v_n^{-{2p\over p-2}}$. 
       Then $M_n^{{1\over q}-{1\over2}} = M_n^{{p-1\over p}-{1\over2}} = M_n^{{p-2\over2p}} \le {1\over v_n}$, 
       so $M_n^{1\over q} \le {1\over v_n}M_n^{{1\over2}}$.   
       Hence with no loss of sharpness,    
       $$\sup_{\set{d_m}\in{\cal B}_n} \abs{\tsuml_{m=1}^{M_n} d_m} \le M_n^{{1\over q}}.$$   
       Let $\tilde d_m = {1\over M_n} M_n^{{1\over q}} = M_n^{{1\over q}-1} = M_n^{-{1\over p}}$   
       for $1\le m\le M_n$, and $\tilde d_m = 0$ otherwise.  
       Then  
       $\sum_{m=1}^{M_n} \abs{\tilde d_m}^p = 1$ and
       $\sum_{m=1}^{M_n} \abs{\tilde d_m v_n}^2 = v_n^2 M_n^{{2\over q}-1} = \paren{v_n M_n^{{1\over q}-{1\over2}}}^2 \le 1$,   
       whence $\set{\tilde d_m}\in{\cal B}_n$.       
       Moreover, $\abs{\sum_{m=1}^{M_n} \tilde d_m} = M_n^{{1\over q}}$.   
       Hence   
       $$\sup_{\set{d_m}\in{\cal B}_n} \abs{\tsuml_{m=1}^{M_n} d_m} \ge M_n^{{1\over q}}.$$   
       It follows that  
       $$\sup_{\set{d_m}\in{\cal B}_n} \abs{\tsuml_{m=1}^{M_n} d_m} = M_n^{{1\over q}}. \eqno{\TAG{2.19}{2.19}}$$   

       Let $K\in\NN$, let$\set{\lambda_k}\in\l{2}$, and let $\set{d_{k,\ell}}$ be a sequence of scalars.  
       Note that ${1\over v_n} n^{1\over2} = n^{p-2\over2p} n^{1\over2} = n^{p-1\over p} = n^{1\over q}$.  
       Then by H\"older's inequality,  
       $$\eqalign{\abs{\tsuml_{k=1}^K \tsuml_{\ell=1}^n \lambda_k d_{k,\ell}}
                  &\le  
                  \paren{\tsuml_{k=1}^K \tsuml_{\ell=1}^n \abs{\lambda_k}^q}^{1\over q} 
                  \paren{\tsuml_{k=1}^K \tsuml_{\ell=1}^n \abs{d_{k,\ell}}^p}^{1\over p} 
              \cr &=
                  n^{1\over q}
                  \paren{\tsuml_{k=1}^K \abs{\lambda_k}^q}^{1\over q} 
                  \paren{\tsuml_{k=1}^K \tsuml_{\ell=1}^n \abs{d_{k,\ell}}^p}^{1\over p}}$$ 
       and    
       $$\eqalign{\abs{
                  \tsuml_{k=1}^K \tsuml_{\ell=1}^n 
                  \lambda_k 
                  d_{k,\ell}
                      }
                  &=
                  {1\over v_n} \abs{\tsuml_{k=1}^K \tsuml_{\ell=1}^n \lambda_k d_{k,\ell} v_n}
              \cr &\le  
                  {1\over v_n} \paren{\tsuml_{k=1}^K \tsuml_{\ell=1}^n \abs{\lambda_k}^2}^{1\over2} 
                  \paren{\tsuml_{k=1}^K \tsuml_{\ell=1}^n \abs{d_{k,\ell} v_n}^2}^{1\over2} 
              \cr &=
                  {1\over v_n} n^{1\over2} \paren{\tsuml_{k=1}^K \abs{\lambda_k}^2}^{1\over2} 
                  \paren{\tsuml_{k=1}^K \tsuml_{\ell=1}^n \abs{d_{k,\ell} v_n}^2}^{1\over2} 
              \cr &=
                  n^{1\over q} \paren{\tsuml_{k=1}^K \abs{\lambda_k}^2}^{1\over2} 
                  \paren{\tsuml_{k=1}^K \tsuml_{\ell=1}^n \abs{d_{k,\ell} v_n}^2}^{1\over2}.}$$ 
       Suppose $\set{d_{k,\ell}}\in{\cal B}_n$.    
       Then $\paren{\sum_{k=1}^K \sum_{\ell=1}^n \abs{d_{k,\ell}}^p}^{1\over p} \le 1$ and \hfil\break 
            $\paren{\sum_{k=1}^K \sum_{\ell=1}^n \abs{d_{k,\ell} v_n}^2}^{1\over2} \le 1$.    
       Hence  
       $\abs{\sum_{k=1}^K\sum_{\ell=1}^n\lambda_k d_{k,\ell}} \le n^{1\over q}\paren{\sum_{k=1}^K\abs{\lambda_k}^q}^{1\over q}$    
       and \hfil\break 
       $\abs{\sum_{k=1}^K\sum_{\ell=1}^n\lambda_k d_{k,\ell}} \le n^{1\over q}\paren{\sum_{k=1}^K\abs{\lambda_k}^2}^{1\over2}$.   
       It follows that  
       $$\eqalign{\sup_{\set{d_{k,\ell}}\in{\cal B}_n} \abs{\tsuml_{k=1}^K \tsuml_{\ell=1}^n \lambda_k d_{k,\ell}}
                  &\le n^{1\over q} \min\set{\paren{\tsuml_{k=1}^K \abs{\lambda_k}^q}^{1\over q},   
                                             \paren{\tsuml_{k=1}^K \abs{\lambda_k}^2}^{1\over2}}   
              \cr &=  
                  n^{1\over q} \paren{\tsuml_{k=1}^K \abs{\lambda_k}^2}^{1\over2}.}$$   

       Let $\tilde d_{k,\ell} = {1\over v_n} n^{-{1\over2}} \bar\lambda_k$ 
       for $1\le\ell\le n$ and $1\le k\le K$, and $\tilde d_{k,\ell} = 0$ otherwise,     
       where $\bar\lambda_k$ is the complex conjugate of $\lambda_k$.   
       Note that ${1\over v_n} n^{-{1\over2}} = n^{p-2\over2p} n^{-{1\over2}} = n^{-{1\over p}}$.   
       Hence  
       $$\eqalign{\paren{\tsuml_{k=1}^K \tsuml_{\ell=1}^n \abs{\tilde d_{k,\ell}}^p}^{1\over p}
                  &=     
                  n^{-{1\over p}} \paren{\tsuml_{k=1}^K \tsuml_{\ell=1}^n \abs{\lambda_k}^p }^{1\over p}
              \cr &=     
                  \paren{\tsuml_{k=1}^K \abs{\lambda_k}^p }^{1\over p}
              \cr &\le      
                  \paren{\tsuml_{k=1}^K \abs{\lambda_k}^2 }^{1\over 2} }$$
       and    
       $$\paren{\tsuml_{k=1}^K \tsuml_{\ell=1}^n \abs{\tilde d_{k,\ell} v_n}^2}^{1\over2} =        
         \paren{\tsuml_{k=1}^K \tsuml_{\ell=1}^n n^{-1} \abs{\lambda_k}^2 }^{1\over2} =        
         \paren{\tsuml_{k=1}^K \abs{\lambda_k}^2 }^{1\over2} .$$        
       Thus for $\paren{\sum_{k=1}^K \abs{\lambda_k}^2 }^{1\over2} \le 1$, $\set{\tilde d_{k,\ell}} \in {\cal B}_n$.     
       Moreover, for $\paren{\sum_{k=1}^K \abs{\lambda_k}^2 }^{1\over2} = 1$,    
       $$\abs{\tsuml_{k=1}^K \tsuml_{\ell=1}^n \lambda_k \tilde d_{k,\ell}} =  
         \tsuml_{k=1}^K \tsuml_{\ell=1}^n n^{-{1\over p}} \abs{\lambda_k}^2 = 
         n^{1\over q} \tsuml_{k=1}^K \abs{\lambda_k}^2 =
         n^{1\over q} \paren{\tsuml_{k=1}^K \abs{\lambda_k}^2}^{1\over2}.$$
       Hence  
       $$\sup_{\set{d_{k,\ell}}\in{\cal B}_n} \abs{\tsuml_{k=1}^K \tsuml_{\ell=1}^n \lambda_k d_{k,\ell}}
         \ge    
         n^{1\over q} \paren{\tsuml_{k=1}^K \abs{\lambda_k}^2}^{1\over2}.$$   
       It follows that  
       $$\sup_{\set{d_{k,\ell}}\in{\cal B}_n} \abs{\tsuml_{k=1}^K \tsuml_{\ell=1}^n \lambda_k d_{k,\ell}}
         = n^{1\over q} \paren{\tsuml_{k=1}^K \abs{\lambda_k}^2}^{1\over2}. \eqno{\TAG{2.20}{2.20}}$$   
       \QED

\xproclaim {Proposition 2.48}. Let $1<p<\infty$ where $p\ne2$. Then $B_p\not\cinjects\sumsum{\l{2}}{\l{p}}\oplus X_p$.  

\proof By duality, it suffices to show that $B_q\not\cinjects\sumsum{\l{2}}{\l{q}}\oplus X_q$ for $1<q<2$.   
       Let $1<q<2$ and suppose $B_q\cinjects\sumsum{\l{2}}{\l{q}}\oplus X_q$.
       Let $p$ be the conjugate index of $q$.  
       For each $n\in\NN$, let $v_n$ and ${\cal B}_n$ be as in Lemma 2.47.  
       Now $B_q\sim B_p^*\sim\sumsum{X_{p,v^{(n)}}^*}{\l{q}}$,     
       so $\sumsum{X_{p,v^{(n)}}^*}{\l{q}}\cinjects\sumsum{\l{2}}{\l{q}}\oplus X_q$.    
       Let $T\colon\sumsum{X_{p,v^{(n)}}^*}{\l{q}}\to\sumsum{\l{2}}{\l{q}}\oplus X_q$
       be an isomorphic imbedding with {complemented} range.      
       Let \hfil\break 
       $Q\colon\sumsum{\l{2}}{\l{q}}\oplus X_q\to\sumsum{\l{2}}{\l{q}}\oplus\set{0_{X_q}}$ be the obvious projection.   
       Then \hfil\break 
       $QT\colon\sumsum{X_{p,v^{(n)}}^*}{\l{q}}\to\sumsum{\l{2}}{\l{q}}\oplus\set{0_{X_q}}$ is a bounded linear operator.    

       We will show that there is a subspace $Y$ of $\sumsum{X_{p,v^{(n)}}^*}{\l{q}}$        
       isometric to \hfil\break 
       $\sumsum{\l{2}}{\l{q}}$ such that $\norm{Q|_{T(Y)}}<1$,    
       whence $(I-Q)|_{T(Y)}$ induces an isomorphic \hfil\break 
       imbedding of $\sumsum{\l{2}}{\l{q}}$ into $X_q$,    
       where $I$ is the formal identity mapping.     
       However \hfil\break 
       by \xciteplus{S}{Proposition 2}, presented below as Lemma 3.7,
       no such imbedding exists, and the proposition will follow.     


       Let $\set{e_{m,n}}$ be the standard basis of $\sumsum{X_{p,v^{(n)}}^*}{\l{q}}$,   
       where for each $n\in\NN$, $\setwlimits{e_{m,n}}{m=1}{\infty}$
       is isometrically equivalent to the standard basis of $X_{p,v^{(n)}}^*$ and \hbox{equivalent}
       to the standard basis of $\l{2}$.  
       Let $\set{\tilde e_{m,n}}$ be the standard basis of $\sumsum{X_{p,v^{(n)}}}{\l{p}}$, \hfil\break 
       where for each $n\in\NN$, $\setwlimits{\tilde e_{m,n}}{m=1}{\infty}$
       is isometrically equivalent to the standard basis of $X_{p,v^{(n)}}$.  

       For $K\in\NN$, let $\Gamma(K)$ denote a subset of $\NN$ having cardinality $K$.  
       Let $M\in\NN$.   
       Then for fixed $n\in\NN$,
       letting $\vector{\phantom0,\phantom0}$ denote the action of $X_{p,v^{(n)}}^*$ on $X_{p,v^{(n)}}$,  
       $$\eqalignno{\norm{\tsuml_{m\in\Gamma(M)} e_{m,n}}
         &= \sup_{\set{d_k}\in{\cal B}_n}
            \abs{\vector{\tsuml_{k=1}^{\infty} d_k \tilde e_{k,n}, \tsuml_{m\in\Gamma(M)} e_{m,n} }} \cr
         &= \sup_{\set{d_k}\in{\cal B}_n} \abs{\tsuml_{k\in\Gamma(M)} d_{k}} \cr
         &= \sup_{\set{d_k}\in{\cal B}_n} \abs{\tsuml_{k=1}^M d_k}. &\TAG{2.21}{2.21} \cr}$$
       Now for fixed $n\in\NN$,
       letting $M_n \le v_n^{-{2p\over p-2}} = n$ as in Lemma 2.47,   
       equations \TAG{2.21}{2.21} and \TAG{2.19}{2.19} yield   
       $$\norm{\tsuml_{m\in\Gamma(M_n)} e_{m,n}} = M_n^{{1\over q}},$$
       or upon normalization,   
       $$\norm{M_n^{-{1\over q}} \tsuml_{m\in\Gamma(M_n)} e_{m,n}} = 1.\eqno{\TAG{2.22}{2.22}}$$

       We now introduce a construction which will be used in two different settings.   
       Fix $n\in\NN$ and let $\tilde M_n = v_n^{-{2p\over p-2}} = n$.    
       Let $\setwlimits{E_k^{(n)}}{k=1}{\infty}$ be a sequence of disjoint subsets of $\NN$,    
       each of cardinality $\tilde M_n$.    
       Let $\set{\tau(m)}$ be an increasing sequence of positive integers.   
       For each $k\in\NN$, let $x_k^{(n)}=\tilde M_n^{-{1\over q}} \sum_{m\in E_k^{(n)}} e_{\tau(m),n}$.    
       Then each $x_k^{(n)}$ is of norm one by equation \TAG{2.22}{2.22}, and      
       $\setwlimits{x_k^{(n)}}{k=1}{\infty}$ is equivalent to the standard basis of $\l{2}$.  
       Recalling equation \TAG{2.20}{2.20} for the last step, for $K\in\NN$ and $\set{\lambda_k}\in\l{2}$,   
\vglue.3in 
\eject 
       $$\eqalignno{\norm{\tsuml_{k=1}^K \lambda_k x_k^{(n)}}
                 &= \tilde M_n^{-{1\over q}} \norm{\tsuml_{k=1}^K \lambda_k \tsuml_{m\in E_k^{(n)}} e_{\tau(m),n} }
             \cr &= n^{-{1\over q}} \sup_{\set{d_\ell}\in{\cal B}_n}
                    \abs{\vector{\tsuml_{\ell=1}^{\infty} d_{\ell} \tilde e_{\ell,n},
                                 \tsuml_{k=1}^K \lambda_k \tsuml_{m\in E_k^{(n)}} e_{\tau(m),n} }}
             \cr &= n^{-{1\over q}} \sup_{\set{d_\ell}\in{\cal B}_n}
                    \abs{\tsuml_{k=1}^K \lambda_k \tsuml_{\ell\in E_k^{(n)}} d_{\ell} }
             \cr &= \paren{\tsuml_{k=1}^K \abs{\lambda_k}^2}^{1\over2}. &\TAG{2.23}{2.23}}$$ 
       Hence $\set{x_k^{(n)}}_{k=1}^{\infty}$ is in fact isometrically equivalent to the standard basis of $\l{2}$.    

       We now distinguish two exhaustive but not mutually exclusive cases.  
       In the first case, there are infinitely many $n\in\NN$ such that    
       $\xlim_{m\to\infty} \norm{QT(e_{m,n})} = 0$.   
       In the second case, there are infinitely many $n\in\NN$ such that    
       $\xlim\xsup_{m\in\NN} \norm{QT(e_{m,n})} > 0$.   

       We will show that in either case,     
       there is an increasing sequence $\setwlimits{n(i)}{i=1}{\infty}$ of positive integers    
       and a sequence $\setwlimits{X_{n(i)}}{i=1}{\infty}$ of subspaces of $\sumsum{X_{p,v^{(n)}}^*}{\l{q}}$
       such that for each $i\in\NN$,  
       $X_{n(i)}$ is a subspace of $\span{e_{m,n(i)}\colon m\in\NN}$ isometric to $\l{2}$ with \hfil\break 
       $\big\|Q|_{T(X_{n(i)})}\big\| \le \norm{T^{-1}}\norm{QT|_{X_{n(i)}}} < {1\over2^i}$.   
       It will follow that there is a subspace \hfil\break 
       $Y=\sumsum{Y_n}{\l{q}}$ of $\sumsum{X_{p,v^{(n)}}^*}{\l{q}}$ 
       isometric to $\sumsum{\l{2}}{\l{q}}$ such that \hfil\break 
       $\big\|Q|_{T(Y)}\big\| \le \norm{T^{-1}}\big\|QT|_Y\big\| < 1$.  
       [$Y_{n(i)}=X_{n(i)}$ and $Y_k=\set{0}$ if $k\notin\set{n(i)}$.]   
       As noted before, the proposition will then follow.   

\noindent{\it{The first case\/.}}

       Fix $n\in\NN$ such that $\xlim_{m\to\infty} \norm{QT\paren{e_{m,n}}} = 0$.   
       Choose a subsequence \hfil\break 
       $\setwlimits{e_{\alpha(m),n}}{m=1}{\infty}$ of $\setwlimits{e_{m,n}}{m=1}{\infty}$
       such that for each $m\in\NN$, \hfil\break 
       $\norm{QT\paren{e_{\alpha(m),n}}} < \displaystyle{1\over2^{m+n}n^{1\over p}\norm{T^{-1}}}$.
       
       Let $\tilde M_n = v_n^{-{2p\over p-2}} = n$.  
       Let $\setwlimits{E_k^{(n)}}{k=1}{\infty}$ be a sequence of disjoint subsets of $\NN$,  
       each of cardinality $\tilde M_n$,  
       such that for each $k\in\NN$, $\inf E_k^{(n)}\ge k$.   
       Then for each $m\in E_k^{(n)}$, 
       $\norm{QT\paren{e_{\alpha(m),n}}} < \displaystyle{1\over2^{k+n}n^{1\over p}\norm{T^{-1}}}$.
       For each $k\in\NN$, let $x_k^{(n)}=\tilde M_n^{-{1\over q}} \sum_{m\in E_k^{(n)}} e_{\alpha(m),n}$.    
       Then each $x_k^{(n)}$ is of norm one by equation \TAG{2.22}{2.22},          
       $\setwlimits{x_k^{(n)}}{k=1}{\infty}$ is isometrically equiv\-alent to the standard basis of $\l{2}$
       as in equation \TAG{2.23}{2.23},   
       and for each $x_k^{(n)}$,   
       $$\eqalign{\norm{QT\paren{x_k^{(n)}}}
                = \tilde M_n^{-{1\over q}} \norm{\tsuml_{m\in E_k^{(n)}} QT\paren{e_{\alpha(m),n}}}   
            & \le n^{-{1\over q}} \tsuml_{m\in E_k^{(n)}} \norm{QT\paren{e_{\alpha(m),n}}}   
          \cr & < n^{-{1\over q}} n {1\over2^{k+n}n^{1\over p}\norm{T^{-1}}}
          \cr & = {1\over2^{k+n}\norm{T^{-1}}}.}$$
       Let $\set{\lambda_k}\in\l{2}$ be of norm one.  
       Then       
       $$\eqalign{\norm{QT\paren{\tsuml_{k=1}^{\infty} \lambda_k x_k^{(n)}}}
                = \norm{\tsuml_{k=1}^{\infty} \lambda_k QT\paren{x_k^{(n)}}}
            & \le \tsuml_{k=1}^{\infty} \norm{QT\paren{x_k^{(n)}}}
          \cr & < \tsuml_{k=1}^{\infty} {1\over2^{k+n}\norm{T^{-1}}}
          \cr & = {1\over2^n\norm{T^{-1}}}. }$$
       Letting $X_n=\span{x_k^{(n)}:k\in\NN}$, it follows that
       $\norm{Q|_{T(X_n)}} \le \norm{T^{-1}}\norm{QT|_{X_n}} < {1\over2^n}$.    

       Release $n\in\NN$ as a free variable.  
       Let $\setwlimits{n(i)}{i=1}{\infty}$ be an increasing sequence of positive integers such that
       for each $i\in\NN$, $\lim_{m\to\infty} \norm{QT(e_{m,n(i)})} = 0$.  
       Then for each $i\in\NN$,  
       there is a subspace $X_{n(i)}$ of $\sumsum{X_{p,v^{(n)}}^*}{\l{q}}$ isometric to $\l{2}$ such that    
       $X_{n(i)}$ is a subspace of $\span{e_{m,n(i)}\colon m\in\NN}$ with
       $\big\|Q|_{T(X_{n(i)})}\big\| \le \norm{T^{-1}}\norm{QT|_{X_{n(i)}}} < {1\over2^{n(i)}}\le{1\over2^i}$.   
       Thus the proposition follows in the first case.  

\noindent{\it{The second case\/.}}

       Fix $n\in\NN$ such that $c_n = \lim\sup_{m\in\NN} \norm{QT\paren{e_{m,n}}} > 0$.   
       Then $c_n\le\norm{QT}$.   
       Given $0<\epsilon<1$, we may choose a subsequence
       $\setwlimits{e_{\alpha(m),n}}{m=1}{\infty}$ of $\setwlimits{e_{m,n}}{m=1}{\infty}$ such that 
       $\xlim_{m\to\infty} \norm{QT\paren{e_{\alpha(m),n}}} = c_n$, with
       $\xsup_{m\in\NN} \abs{\norm{QT\paren{e_{\alpha(m),n}}} - c_n}<\epsilon c_n$, and such that    
       $\set{QT\paren{e_{\alpha(m),n}}}_{m=1}^{\infty}$
       is a basic sequence \xciteplus{B-P}{Theorem 3}, whence \hfil\break 
       $QT|_{\span{e_{\alpha(m),n}\colon m\in\NN}}$ is an isomorphic imbedding and  
       $\set{QT\paren{e_{\alpha(m),n}}}_{m=1}^{\infty}$ is equivalent to the standard basis of $\l{2}$.  
       Now by Proposition 2.19, given $0<\epsilon<1$ and such a \hfil\break 
       sequence $\setwlimits{e_{\alpha(m),n}}{m=1}{\infty}$,  
       we may choose a subsequence $\setwlimits{e_{\beta(m),n}}{m=1}{\infty}$ such that \hfil\break 
       $\set{QT\paren{e_{\beta(m),n}}}_{m=1}^{\infty}$ is
       $(1+\epsilon)$-equivalent to the standard basis of $\l{2}$.   
       
       Let $\tilde M_n = v_n^{-{2p\over p-2}} = n$.  
       Let $\setwlimits{E_k^{(n)}}{k=1}{\infty}$ be a sequence of disjoint subsets of $\NN$,  
       each of cardinality $\tilde M_n$.  
       Given $0<\epsilon<1$ and $\setwlimits{e_{\beta(m),n}}{m=1}{\infty}$ as above, for each $k\in\NN$ let 
       $x_k^{(n)}=\tilde M_n^{-{1\over q}} \sum_{m\in E_k^{(n)}} e_{\beta(m),n}$. 
       Then each $x_k^{(n)}$ is of norm one by equation \TAG{2.22}{2.22}, 
       $\setwlimits{x_k^{(n)}}{k=1}{\infty}$ is \hfil\break 
       isometrically equivalent
       to the standard basis of $\l{2}$ as in equation \TAG{2.23}{2.23},  
       and for each $x_k^{(n)}$,   
       $$\norm{QT\paren{x_k^{(n)}}} = \tilde M_n^{-{1\over q}} \norm{\tsuml_{m\in E_k^{(n)}} QT\paren{e_{\beta(m),n}}}   
                              \approx \tilde M_n^{-{1\over q}} \tilde M_n^{1\over\vphantom{q}2} c_n
                                    = \tilde M_n^{{1\over2}-{1\over q}} c_n, \eqno{\TAG{2.24}{2.24}}$$
       where the approximation can be improved to any degree by the choice of ($\epsilon$ and) \hfil\break 
       $\setwlimits{x_k^{(n)}}{k=1}{\infty}$.   
       
       Given $0<\epsilon<1$, we may choose a sequence $\setwlimits{x_k^{(n)}}{k=1}{\infty}$ as above such that 
       $\abs{\norm{QT\paren{x_k^{(n)}}} - \tilde M_n^{{1\over2}-{1\over q}} c_n}
        < \epsilon \tilde M_n^{{1\over2}-{1\over q}} c_n$, where  
       $QT|_{\span{x_k^{(n)}\colon k\in\NN}}$ is an isomorphic \hfil\break 
       imbedding and  
       $\set{QT\paren{x_k^{(n)}}}_{k=1}^{\infty}$ is equivalent to the standard basis of $\l{2}$.  
       Thus by \hfil\break 
       Proposition 2.19, given $0<\epsilon<1$ and such a sequence $\setwlimits{x_k^{(n)}}{k=1}{\infty}$,
       there is a {sub}\-{sequence} $\set{x_{\gamma(k)}^{(n)}}_{k=1}^{\infty}$ such that
       $\set{QT\paren{x_{\gamma(k)}^{(n)}}}_{k=1}^{\infty}$ 
       is $(1+\epsilon)$-equivalent to the standard basis of $\l{2}$.  
       Recalling \TAG{2.24}{2.24}, it follows that for $\set{\lambda_k}\in\l{2}$,     
       $$\norm{QT\paren{\tsuml_{k=1}^{\infty} \lambda_k x_{\gamma(k)}^{(n)}}}
               = \norm{\tsuml_{k=1}^{\infty} \lambda_k QT\paren{x_{\gamma(k)}^{(n)}}}
         \approx \paren{\tsuml_{k=1}^{\infty}
                        \abs{\lambda_k}^2}^{1\over2} \tilde M_n^{{1\over2}-{1\over q}} c_n, \eqno{\TAG{2.25}{2.25}}$$
       where the approximation can be improved to any degree by the choice of ($\epsilon$ and) \hfil\break 
       $\setwlimits{x_{\lambda(k)}^{(n)}}{k=1}{\infty}$.   

       Now $\set{x_k^{(n)}}_{k=1}^{\infty}$ is isometrically equivalent to the standard basis of $\l{2}$  
       as noted above,  
       and the same is true of $\set{x_{\gamma(k)}^{(n)}}_{k=1}^{\infty}$.   
       Let $X_n=\span{x_{\gamma(k)}^{(n)}:k\in\NN}$.  
       Then by \hfil\break 
       \TAG{2.25}{2.25}, it follows that   
       $$\norm{QT|_{X_n}} \approx \tilde M_n^{{1\over2}-{1\over q}} c_n \le n^{{1\over2}-{1\over q}} \norm{QT},
         \eqno{\TAG{2.26}{2.26}}$$  
       where the approximation can be improved to any degree as in \TAG{2.25}{2.25}.  

       Release $n$ as a free variable and note that $\xlim_{n\to\infty} n^{{1\over2}-{1\over q}}\norm{QT} = 0$.  
       Hence by the hypothesis of the second case and by \TAG{2.26}{2.26},    
       we may choose an increasing sequence $\setwlimits{n(i)}{i=1}{\infty}$ 
       of positive integers such that for each $i\in\NN$, \hfil\break 
       $c_{n(i)} = \lim\sup_{m\in\NN} \norm{QT(e_{m,n(i)})} > 0$  
       and
       there is a subspace $X_{n(i)}$ of $\sumsum{X_{p,v^{(n)}}^*}{\l{q}}$ isometric to $\l{2}$ such that 
       $X_{n(i)}$ is a subspace of $\span{e_{m,n(i)}\colon m\in\NN}$ with \hfil\break 
       $\big\|Q|_{T(X_{n(i)})}\big\| \le \norm{T^{-1}}\norm{QT|_{X_{n(i)}}} < 2^i$.   
       Thus the proposition follows in the second case, and in the general case.  
       \qquad\QED


\medskip
 Collecting our results and deducing simple consequences yields the following.  

\xproclaim {Proposition 2.49}. Let $1<p<\infty$ where $p\ne2$. Then  
                               \item{(a)} $B_p\not\cinjects\sumsum{\l{2}}{\l{p}}$,
                               \item{(b)} $\sumsum{\l{2}}{\l{p}}\oplus X_p\not\cinjects\sumsum{\l{2}}{\l{p}}$,  
                               \item{(c)} $B_p\not\cinjects\sumsum{\l{2}}{\l{p}}\oplus X_p$,  
                               \item{(d)} $\sumsum{\l{2}}{\l{p}}\oplus X_p\not\cinjects B_p$,  
                               \item{(e)} $B_p\oplus X_p\not\cinjects B_p$,  
                               \item{(f)} $B_p\oplus X_p\not\cinjects\sumsum{\l{2}}{\l{p}}\oplus X_p$, and   
                               \item{(g)} $\sumsum{X_p}{\l{p}}\not\cinjects B_p\oplus X_p$.   

\proof  
       \item{(a)} Part (a) is a restatement of part (d) of Proposition 2.37.

       \item{(b)} Part (b) follows from part (f) of Proposition 2.24:
                  $X_p\not\cinjects\sumsum{\l{2}}{\l{p}}$. 

       \item{(c)} Part (c) is a restatement of Proposition 2.48.  

       \item{(d)} Part (d) follows from part (g) of Proposition 2.37:
                  $X_p\not\cinjects B_p$.  

       \item{(e)} Part (e) follows from part (g) of Proposition 2.37:
                  $X_p\not\cinjects B_p$.  

       \item{(f)} Part (f) follows from part (c) above.  


       \item{(g)} Part (g) is a restatement of Proposition 2.46. 
       \qquad\QED
                               
\medskip
 Building on diagram \TAG{2.17}{2.17}, for $1<p<\infty$ where $p\ne2$, we have   
$$\matrix{&&&&B_p&&&&\L{p}\cr
          &&&\cusa&&\cdsa&&&\cua\cr
          \ltwo&&\sumsum{\ltwo}{\lp}&&&&B_p\oplus X_p&\cra&\sumsum{X_p}{\lp}.\cr 
          \cda&\cusa&&\cdsa&&\cusa&&&\cr
          \ltwo\oplus\lp&&&&\sumsum{\ltwo}{\lp}\oplus X_p&&&&\cr
          \cua&\cdsa&&\cusa&&&&&\cr
          \lp&&X_p&&&&&&\cr}\eqno{\TAG{2.27}{2.27}}$$ 

\preheadspace
\firsthead{Concluding Remarks}
\postheadspace
 Fix $1<p<\infty$ where $p\ne2$.  

 If $X$ and $Y$ are separable infinite-dimensional $\SL{p}$ spaces,   
 then $X \oplus Y$ is a \hfil\break 
 separable infinite-dimensional $\SL{p}$ space as well.   
 Suppose $X$ and $Y$ are as above and are isomorphic to their squares.   
 If $X$ and $Y$ are incomparable in the sense that \hfil\break 
 $X\not\cinjects Y$ and $Y\not\cinjects X$,  
 then $X\oplus Y$ is isomorphically distinct from both $X$ and $Y$, 
 while if $X\cinjects Y$, then $X\oplus Y\sim Y$.  

 From the list 
 $\l{p}$, $\l{2}\oplus\l{p}$, $\sumsum{\l{2}}{\l{p}}$, $X_p$, $B_p$, $\sumsum{X_p}{\l{p}}$, $\L{p}$  
 of seven spaces, the only incomparable pairs of spaces are   
 $\set{\sumsum{\l{2}}{\l{p}} , X_p}$ and $\set{B_p , X_p}$.   
 As has been shown, $\sumsum{\l{2}}{\l{p}}\oplus X_p$ and $B_p\oplus X_p$    
 are isomorphically distinct from
 each of the seven listed spaces and from each other.     
 Augmenting the list of seven spaces with the two new ones,  
 the only new incomparable pair of spaces is   
 $\set{B_p , \sumsum{\l{2}}{\l{p}}\oplus X_p}$.   
 However $\sumsum{\l{2}}{\l{p}}\cinjects B_p$, so
 $B_p\oplus\sumsum{\l{2}}{\l{p}} \sim B_p$, whence \hfil\break 
 $B_p\oplus\paren{\sumsum{\l{2}}{\l{p}}\oplus X_p}\sim     
  \paren{B_p\oplus\sumsum{\l{2}}{\l{p}}}\oplus X_p\sim B_p\oplus X_p$,
 which has already been included in the augmented list.   

 If $Z$ is a separable infinite-dimensional Banach space such that $Z\cinjects\L{p}$, then
 $\sumsum{Z}{\l{p}}$ is a separable infinite-dimensional $\SL{p}$ space.    
 However, from the augmented list of nine spaces above,  
 no space arises from this method of construction   
 which has not already been included in the list.   

%
%
%
\chapter{III}
\headspace
\chaptertitle{THE TENSOR PRODUCT CONSTRUCTION OF SCHECHTMAN}
\headspace
 Let $1<p<\infty$ where $p\ne2$.  
 Schechtman \xcite{S} constructed a sequence of isomorphically distinct separable infinite-dimensional $\SL{p}$ spaces 
 by iterating a certain tensor product of Rosenthal's space $X_p$ with itself.  
 Using $X_p^{\otimes n}$ to denote $X_p\otimes\cdots\otimes X_p$ \hfil\break 
 ($n$ factors),  
 the resulting sequence is $\set{X_p^{\otimes n}}_{n=1}^{\infty}$.  

 For closed subspaces $X$ and $Y$ of $\L{p}$, 
 $X \otimes Y$ is defined to be the closed linear span in $\L{p}([0,1]\times[0,1])$
 of products of the form $x(s)y(t)$ where $x\in X$ and $y\in Y$.
 It is a fairly routine matter to show that
 if $X$ and $Y$ are separable infinite-dimensional $\SL{p}$ spaces,    
 then $X\otimes Y$ is a separable infinite-dimensional $\SL{p}$ space.   
 More work is required to show that for $m \ne n$, 
 $X_p^{\otimes m} \not\sim X_p^{\otimes n}$.  
 
\preheadspace
\firsthead{The Tensor Product Construction}
\postheadspace
 
 We begin with some preliminary definitions and lemmas.   
 For each $k\in\NN$, let $I^k = [0,1]^k$.    
 Let $m,n\in\NN$. 

\definition Let $1\le p<\infty$ and let $X$ and $Y$ be closed subspaces of $\L{p}(I^m)$ and $\L{p}(I^n)$, 
            respectively. 
            Define the tensor product $X\otimes Y$ of $X$ and $Y$ by 
            $$X \otimes Y = \span{x(s)y(t): x\in X, y\in Y, s\in I^m, t\in I^n}_{\L{p}(I^{m+n})}.$$
            Denote the element $x(s)y(t)$ by $x\otimes y$. 


 Let $X$ and $Y$ be as above, and let $Z$ be a closed subspace of $\L{p}(I^k)$ for some \hfil\break 
 $k\in\NN$.  
 Then $X\otimes(Y\otimes Z) = (X\otimes Y)\otimes Z$.   
 Thus the expressions $X\otimes Y\otimes Z$ and \hfil\break 
 $\bigotimes_{i=1}^N X$ are unambiguous.   
 The tensor power $\bigotimes_{i=1}^N X$ will also be denoted $X^{\otimes N}$. 

 The following lemma will be used in the proof of the fact that
 the tensor product of complemented subspaces of $\L{p}$ is a complemented subspace of $\L{p}\paren{I^2}$. 

\xproclaim {Lemma 3.1}. Let $1\le p<\infty$. 
                        Then $\L{p}(I^m)\otimes\L{p}(I^n) = \L{p}(I^{m+n})$.

\proof Note that $\L{p}(I^m)\otimes\L{p}(I^n)$ is a closed subspace of $\L{p}(I^{m+n})$.
       Thus it will suffice to show that $\L{p}(I^m)\otimes\L{p}(I^n)$ is dense in $\L{p}(I^{m+n})$.
       Let $f\in\L{p}(I^{m+n})$ and let $\epsilon>0$.
       Choose $g\in C(I^{m+n})$ such that $\norm{f-g}_{\L{p}(I^{m+n})} < {\epsilon\over2}$. 
       By the Stone-Weierstrass theorem,
       choose $h\in\finitespan_{C(I^{m+n})} \set{h_1(s)h_2(t):h_1\in C(I^m),h_2\in C(I^n)}$ such that   
       $\norm{g-h}_{\L{p}(I^{m+n})} \le \norm{g-h}_{\L{\infty}(I^{m+n})} < {\epsilon\over2}$.  
       Then $\norm{f-h}_{\L{p}(I^{m+n})} < \epsilon$.   
       \qquad\QED

 The tensor product preserves the property of having an unconditional basis,
 as shown in the following lemma \xciteplus{S}{Lemma 3}.  

\xproclaim {Lemma 3.2}. Let $1\le p<\infty$ and let $X$ and $Y$ be as above. 
                        Suppose $\set{x_i}$ and $\set{y_j}$ are unconditional bases for $X$ and $Y$, respectively.  
                        Then $\set{x_i \otimes y_j}_{i,j\in\NN}$ is an \hfil\break 
                        unconditional basis for $X\otimes Y$.   

\proof Note that $\span{x_i\otimes y_j: i,j\in\NN} = X\otimes Y$.   
       Let $\set{r_k}$ be the sequence of Rademacher functions.  
       Then by the unconditionality of $\set{x_i(s)}$ for each $t$,
       Fubini's theorem, and a generalization of Khintchine's inequality, for scalars $a_{i,j}$ 
       
       \vglue1.5in 
       \eject 

       $$\eqalign{\norm{\tsuml_{i}\tsuml_{j} a_{i,j} (x_i \otimes y_j) }_{\L{p}(I^{m+n})}^p
         \! &= \dint\!\dint \abs{\tsuml_{i}\tsuml_{j} a_{i,j} x_i(s) y_j(t) }^p \,ds\,dt
     \cr &\approx \dint\!\dint\!\dint\!\dint \abs{\tsuml_{i}\tsuml_{j} a_{i,j} r_i(u) r_j(v) y_j(t) x_i(s) }^p \,ds\,du\,dv\,dt
     \cr &=       \dint\!\dint\!\dint\!\dint \abs{\tsuml_{i}\tsuml_{j} a_{i,j} x_i(s) y_j(t) r_i(u) r_j(v) }^p \,du\,dv\,ds\,dt
     \cr &\approx \dint\!\dint \paren{\tsuml_{i}\tsuml_{j} \abs{a_{i,j} x_i(s) y_j(t) }^2 }^{{p\over2}} \,ds\,dt .}$$

       If $\sum_{i}\sum_{j} a_{i,j} (x_i \otimes y_j) = 0$,
       then $\int\!\int \paren{\sum_{i}\sum_{j} \abs{a_{i,j} x_i(s) y_j(t) }^2 }^{{p\over2}}\,ds\,dt=0$ by the inequalities above, 
       and $a_{i,j}=0$ for all $i,j\in\NN$. 
       Hence $\set{x_i\otimes y_j}_{i,j\in\NN}$ is a basis for $X\otimes Y$. 
       The unconditionality of $\set{x_i \otimes y_j}_{i,j\in\NN}$ is similarly clear from the inequalities above. 
       \qquad\QED

\definition Let $1\le p<\infty$.  
            Let $X$ and $X'$ be closed subspaces of $\L{p}(I^m)$, and   
            let $Y$ and $Y'$ be closed subspaces of $\L{p}(I^n)$.  
            Suppose $S\colon X\to X'$ and $T\colon Y\to Y'$ are bounded linear operators.
            Define the tensor product $S\otimes T:X\otimes Y \to X'\otimes Y'$ of $S$ and $T$ by 
            $$(S\otimes T)\paren{\tsuml_{i} x_i(s) y_i(t)} = \tsuml_{i} S(x_i)(s)T(y_i)(t)$$
            for sequences $\set{x_i}$ in $X$ and $\set{y_i}$ in $Y$
            such that $\sum_i x_i(s) y_i(t) \in\L{p}\paren{I^{m+n}}$. 
 
 The tensor product of bounded linear operators is bounded and linear,
 as shown in the following lemma \xcite{S}.
 Moreover, the tensor product of projections is a projection,
 and the tensor product of isomorphisms is an isomorphism, 
 as shown in the subsequent lemma \xciteplus{S}{Lemmas 1 and 2}.   

\xproclaim {Lemma 3.3}. Let $1\le p<\infty$ and let $X$, $X'$, $Y$, $Y'$, $S$, and $T$ be as above.  
                        Then $S\otimes T$ is well-defined and linear, with $\norm{S\otimes T}\le\norm{S}\norm{T}$.  


\proof For $i\in\NN$, let $x_i\in X$ and $y_i\in Y$.
       Then $S\otimes T$ is formally linear by an easy computation.
       Suppose only finitely many elements of $\set{x_i}$ and $\set{y_i}$ are nonzero.
       Then by Fubini's theorem, 
       $$\eqalign{\norm{(S\otimes T) \paren{\tsuml_{i} x_i(s) y_i(t) } }_{\L{p}(I^{m+n})}^p 
               &= \dint\!\dint \abs{\tsuml_{i} S(x_i)(s) T(y_i)(t) }^p \,ds \,dt
           \cr &= \dint \norm{S\paren{\tsuml_{i} T(y_i)(t) x_i } }_{\L{p}(I^m)}^p \,dt
         \cr &\le \norm{S}^p \dint \norm{\tsuml_{i} T(y_i)(t) x_i }_{\L{p}(I^m)}^p \,dt       
           \cr &= \norm{S}^p \dint\!\dint \abs{\tsuml_{i} T(y_i)(t) x_i(s) }^p \,ds \,dt
           \cr &= \norm{S}^p \dint\!\dint \abs{\tsuml_{i} T(y_i)(t) x_i(s) }^p \,dt \,ds
           \cr &= \norm{S}^p \dint \norm{T\paren{\tsuml_{i} x_i(s) y_i } }_{\L{p}(I^n)}^p \,ds
         \cr &\le \norm{S}^p \norm{T}^p \dint \norm{\tsuml_{i} x_i(s) y_i }_{\L{p}(I^n)}^p \,ds 
           \cr &= \norm{S}^p \norm{T}^p \dint\!\dint \abs{\tsuml_{i} x_i(s) y_i(t) }^p \,dt \,ds 
           \cr &= \norm{S}^p \norm{T}^p \norm{\tsuml_{i} x_i(s) y_i(t) }_{\L{p}(I^{m+n})}^p .}$$
       If $z = \sum_{i} x_i(s) y_i(t) = 0$,
       then $(S\otimes T)(z) = 0$ by the inequality above,
       whence $(S\otimes T)(0)=0$ independently of the representation of $0$, and $S\otimes T$ is well-defined.  
       Moreover, $\norm{S\otimes T} \le \norm{S}\norm{T}$ by the inequality above.
       \qquad\QED

\xproclaim {Lemma 3.4}. Let $1\le p<\infty$ and let $X$, $X'$, $Y$, $Y'$, $S$, and $T$ be as above.  
                        \item{(a)} If $S$ and $T$ are projections, then $S\otimes T$ is a projection.  
                        \item{(b)} If $S$ and $T$ are isomorphisms, then $S\otimes T$ is an isomorphism.  

\proof \item{(a)} Suppose $S$ and $T$ are projections. 
                  Then \hfil\break 
                  $(S\otimes T)^2 = (S\otimes T)(S\otimes T) = S^2 \otimes T^2 = S\otimes T$. 
                  Hence $S\otimes T$ is a projection.  
       \item{(b)} Suppose $S$ and $T$ are isomorphisms. 
                  Then $S\otimes T$ and $S^{-1} \otimes T^{-1}$ are formal \hfil\break 
                  inverses,
                  and $\norm{S^{-1} \otimes T^{-1}} \le \norm{S^{-1}}\norm{T^{-1}}$ by Lemma 3.3.
                  Hence $S^{-1} \otimes T^{-1}$ \hfil\break 
                  is bounded and $S\otimes T$ is an isomorphism.  
                  \qquad\QED

 \remark Let $1\le p<\infty$.  
         Suppose $X\injects\L{p}(I^m)$ and $Y\injects\L{p}(I^n)$. 
         By part (b) above,
         $X\otimes Y$ is well-defined up to isomorphism if we identify $X\otimes Y$ with $X'\otimes Y'$ \hfil\break 
         for closed subspaces $X'$ and $Y'$ of $\L{p}(I^m)$ and $\L{p}(I^n)$
         isomorphic to $X$ and $Y$, \hfil\break 
         respectively. 

 The tensor product of complemented subspaces of $\L{p}$ is complemented, and
 the tensor product of $\SL{p}$ spaces is an $\SL{p}$ space, 
 as shown in the following proposition \hfil\break 
 \xciteplus{S}{Lemma 1}. 

\xproclaim {Proposition 3.5}. Let $1<p<\infty$ where $p\ne2$.  
                              Suppose $X$ and $Y$ are \hfil\break 
                              separable infinite-dimensional $\SL{p}$ spaces.  
                              Then $X\otimes Y$ is a separable infinite- \hfil\break 
                              dimensional $\SL{p}$ space.  

\proof It is clear that $X\otimes Y$ is separable and infinite-dimensional.
       Let $X'$ and $Y'$ be complemented subspaces of $\L{p}$ 
       isomorphic to $X$ and $Y$, respectively.  
       Then there are projections $P_{X'}\colon\L{p}\to X'$ and $P_{Y'}\colon\L{p}\to Y'$. 
       By part (a) of Lemma 3.4, \hfil\break 
       $P_{X'}\otimes P_{Y'}\colon \L{p}\otimes\L{p}\to X'\otimes Y'$ is a projection as well, 
       so $X'\otimes Y'$ is a complemented subspace of $\L{p}\otimes\L{p}$,
       which by Lemma 3.1 is equal to $\L{p}(I^2)$.  
       Hence \hfil\break 
       $X\otimes Y \sim X'\otimes Y' \cinjects \L{p}\otimes\L{p} = \L{p}(I^2) \sim \L{p}$.  

       It remains to show that $X\otimes Y \not\sim \l{2}$.  
       By \xciteplus{L-P}{Proposition 7.3}, $\l{p}\cinjects Z$ for every infinite-dimensional $\SL{p}$ space $Z$. 
       Now $\l{p}\cinjects X$ and $\span{y_0}\cinjects Y$ for $y_0\in Y\setminus\set{0}$, 
       whence $\l{p}\sim\l{p}\otimes\span{y_0}\cinjects X\otimes Y$. 
       It follows that $X\otimes Y \not\sim \l{2}$.  
       \qquad\QED

 Of course it follows that $X_p^{\otimes n}$ is an $\SL{p}$ space for $1<p<\infty$ with $p\ne2$. 


\preheadspace
\firsthead{The Isomorphic Distinctness of $X_p^{\otimes m}$ and $X_p^{\otimes n}$}
\postheadspace

 We now present results leading to the conclusion that the various tensor powers of $X_p$ are isomorphically distinct.   
 The main result is Theorem 3.10 below.  

 First we state some facts about stable random variables. 

 Let $1\le T\le2$.   
 Then there is a distribution $\mu$ such that $\int_{\RR} e^{i\alpha x} \,d\mu(\alpha) = e^{-\abs{x}^T}$   
 and a random variable $f\colon[0,1]\to\RR$ having distribution $\mu$.
 Such a random variable $f$ is said to be $T$-stable \xciteplus{W}{III.A.~13 and 14}.  

 If $f$ is a $T$-stable random variable,   
 then $f\in\L{t}$ for each $1\le t<T\le2$.   
 Let $\set{f_n}$ be a sequence of independent $T$-stable random variables.   
 Then for each $1\le t<T\le2$,   
 $\span{f_n}_{\L{t}}$ is isometric to $\l{T}$ \xciteplus{W}{III.A.~15 and 16}.   

 Let $1\le t<T\le2$,
 and let $\set{f_n}$ be a sequence of independent identically \hfil\break 
 distributed $T$-stable random variables normalized in $\L{t}$. 
 Then the sequence $\set{f_n}$ in $\L{t}$ is isometrically equivalent to the standard basis of $\l{t}$, 
 and equivalent to the standard basis of $\l{t'}$ for all $1\le t'<T\le2$ \xciteplus{RII}{Corollary 4.2}. 

 The following lemma is \xciteplus{S}{Proposition 1}.

\xproclaim {Lemma 3.6}. Let $1\le q<r<s\le2$. 
                        Let $X$ and $Y$ be closed subspaces of $\L{q}$ \hfil\break 
                        isomorphic to $\l{r}$ and $\l{s}$, respectively. 
                        Then $\l{r}\otimes\l{s} \sim X \otimes Y \sim \sumsum{\l{s}}{\l{r}}$ via \hfil\break 
                        equivalence of their standard bases. 

\proof Choose a sequence $\set{x_i}$ in $X$ of  
       independent identically distributed $r$-stable random variables normalized in $\L{q}$, 
       and a sequence $\set{y_j}$ in $Y$ of  
       independent identically distributed $s$-stable random variables normalized in $\L{r}$. 
       Then \hfil\break 
       $X\sim\l{r}\sim\span{x_i}_{\L{q}}$
       and 
       $Y\sim\l{s}\sim\span{y_j}_{\L{q}}$.

       For scalars $a_{i,j}$, by the $r$-stability and $q$-normalization of $\set{x_i}$ with $q<r$, we
       have 
       $$\eqalign{\norm{\tsuml_{i}\tsuml_{j} a_{i,j} \paren{x_i \otimes y_j} }_{\L{q}(I^2)}^q             
                &=
                  \dint\!\dint\abs{\tsuml_{i}\tsuml_{j} a_{i,j} x_i(u)y_j(v) }^q \,du\,dv  
      \cr &\approx
                  \dint\!\dint\abs{\tsuml_i
                  \paren{\tsuml_j a_{i,j} y_j(v) }x_i(u) }^q \,du\,dv  
            \cr &=
                  \dint\paren{\tsuml_i
                  \abs{\tsuml_j a_{i,j} y_j(v) }^r }^{{q\over r}} \,dv .}$$      
       Hence by the concavity of $(\phantom0)^{{q\over r}}$,
       and the $s$-stability and $r$-normalization of $\set{y_j}$ with $r<s$, we have 
       $$\eqalign{\norm{\tsuml_{i}\tsuml_{j} a_{i,j} \paren{x_i \otimes y_j} }_{\L{q}(I^2)}^q   
          &\approx
                  \dint\paren{\tsuml_i \abs{\tsuml_j a_{i,j} y_j(v) }^r }^{{q\over r}} \,dv 
          \cr &\le
                  \paren{\dint\tsuml_i \abs{\tsuml_j a_{i,j} y_j(v) }^r \,dv }^{{q\over r}}
            \cr &= 
                  \paren{\tsuml_i \dint\abs{\tsuml_j a_{i,j} y_j(v) }^r \,dv }^{{q\over r}}
            \cr &= 
                  \paren{\tsuml_i \paren{\tsuml_j \abs{a_{i,j}}^s }^{{r\over s}} }^{{q\over r}} .}$$ 
       Moreover, by the triangle inequality and the $s$-stability of $\set{y_j}$ with $q<s$, we have 
       $$\eqalign{\norm{\tsuml_{i}\tsuml_{j} a_{i,j} \paren{x_i \otimes y_j} }_{\L{q}(I^2)}^q 
          &\approx
                  \dint\paren{\tsuml_i \abs{\tsuml_j a_{i,j} y_j(v) }^r }^{{q\over r}} \,dv 
            \cr &=
                  \dint\norm{\set{\abs{\tsuml_j a_{i,j} y_j(v) }^q }_{i=1}^{\infty} }_{\l{{r\over q}}} \,dv
          \cr &\ge
                  \norm{\set{\dint\abs{\tsuml_j a_{i,j} y_j(v) }^q \,dv }_{i=1}^{\infty} }_{\l{{r\over q}}}
            \cr &\approx
                  \norm{\set{\paren{\tsuml_j \abs{a_{i,j}}^s }^{{q\over s}} }_{i=1}^{\infty} }_{\l{{r\over q}}} 
            \cr &=
                  \paren{\tsuml_i \paren{\tsuml_j \abs{a_{i,j}}^s }^{{r\over s}} }^{{q\over r}} .}$$ 
       Hence $\set{x_i \otimes y_j}$ is equivalent to the standard basis of $\sumsum{\l{s}}{\l{r}}$, 
       and \hfil\break 
       $\l{r}\otimes\l{s} \sim X \otimes Y \sim \sumsum{\l{s}}{\l{r}}$.
       \qquad\QED 


\medskip
 Let $1\le p<\infty$ and let $\set{x_i}$ be a sequence in $\L{p}$.
 Then $\set{x_i}$ is said to be \hfil\break 
 uniformly $p$-integrable if for each $\epsilon>0$, there is an $N\in\NN$ such that \hfil\break 
 $\int_{\set{t:\abs{x_i(t)}>N}} \abs{x_i(t)}^p \,dt < \epsilon^p$ for each $i\in\NN$.  

 A basis $\set{x_i}$ for a space $X$ is said to be symmetric if for all permutations $\tau$ of scalars $a_i$, 
 $\sum_i \tau\paren{a_i} x_i$ converges if and only if $\sum_i a_i x_i$ converges. 

 The following lemma is \xciteplus{S}{Proposition 2}.

\xproclaim {Lemma 3.7}. Let $1<q<r<s\le2$. 
                        Then there is no sequence $\set{x_{i,j}}_{i,j\in\NN}$ \hfil\break 
                        of independent random variables in $\L{q}$
                        equivalent to the standard basis of $\sumsum{\l{s}}{\l{r}}$. 

\proof Suppose $\set{x_{i,j}}_{i,j\in\NN}$ is a sequence of independent random variables in $\L{q}$ 
       equivalent to the standard basis of $\sumsum{\l{s}}{\l{r}}$, 
       where for each $j\in\NN$, $\set{x_{i,j}}_{i\in\NN}$ is equivalent to the standard basis of $\l{s}$. 
       Now $\l{q}\not\injects\sumsum{\l{s}}{\l{r}}$. 
       Hence $\set{x_{i,j}}_{i,j\in\NN}$ is \hfil\break 
       uniformly $q$-integrable \xciteplus{J-O}{third lemma}. 

       Let $\epsilon>0$, and choose $N\in\NN$ such that
       $\int_{\set{\abs{x_{i,j}} > N}} \abs{x_{i,j}}^q\,d\mu < \epsilon^q$
       for all $i,j\in\NN$. 
       Let $\delta={1\over D}$ for some $D\in\NN$, 
       and let $\set{I_k}_{k=1}^K$ be a partition of the interval $[-N,N]$
       into $K=D(2N+1)$ intervals of equal length $\abs{I_k}={2N\over K}={2N\over D(2N+1)}<\delta$.

       Let $\rho=\delta^{2q}$. 
       For each $j\in\NN$,
       choose a subsequence $\set{x_{i,j}}_{i\in M_j}$ of $\set{x_{i,j}}_{i\in\NN}$
       such that for each $i,i'\in M_j$ and $k\in\set{1,\ldots,K}$, 
       $$\abs{\mu\paren{\set{x_{i,j}\in I_k}}-\mu\paren{\set{x_{i',j}\in I_k}} } < {\rho\over3}.$$
       Then $\set{x_{i,j}}_{i\in M_j,j\in\NN}$ is still equivalent to the standard basis of $\sumsum{\l{s}}{\l{r}}$. 
       Without loss of generality, suppose $1\in M_j$ for each $j\in\NN$. 
       
       Choose a subsequence $\set{x_{1,j}}_{j\in L}$ of $\set{x_{1,j}}_{j\in\NN}$
       such that for each $j,j'\in L$ and $k\in\set{1,\ldots,K}$,
       $$\abs{\mu\paren{\set{x_{1,j}\in I_k}}-\mu\paren{\set{x_{1,j'}\in I_k}} } < {\rho\over3}.$$
       Then $\set{x_{i,j}}_{i\in M_j,j\in L}$ is still equivalent to the standard basis of $\sumsum{\l{s}}{\l{r}}$. 
       Without loss of generality, suppose $1\in L$. 
       Note that for each $j,j'\in L$, $i\in M_j$, $i'\in M_{j'}$, and $k\in\set{1,\ldots,K}$,
       $$\abs{\mu\paren{\set{x_{i,j}\in I_k}}-\mu\paren{\set{x_{i',j'}\in I_k}} } < \rho.$$

       For each $k\in\set{1,\ldots,K}$, let $c_k$ be the center of $I_k$. 
       Let $\set{z_{i,j}}_{i\in M_j,j\in L}$ be a \hfil\break 
       sequence of $\set{c_1,\ldots,c_K}$-valued
       independent random variables in $\L{q}$ such that 
       for each $j\in L$, $i\in M_j$, and $k\in\set{1,\ldots,K}$,
       $$\mu\paren{\set{z_{i,j}=c_k}} = \mu\paren{\set{x_{1,1}\in I_k}},$$
       and such that $\set{z_{i,j}=c_k}$ is chosen either as a subset of $\set{x_{i,j}\in I_k}$
       or as a superset of $\set{x_{i,j}\in I_k}$. 
       Then $\set{z_{i,j}}_{i\in M_j,j\in L}$ is identically distributed,  
       whence $\set{z_{i,j}}_{i\in M_j,j\in L}$ is a symmetric basis, 
       and for each $j\in L$, $i\in M_j$, and $k\in\set{1,\ldots,K}$,
       $$\abs{\mu\paren{\set{x_{i,j}\in I_k}} - \mu\paren{\set{z_{i,j}=c_k}} } < \rho.$$
       Hence for each $j\in L$, $i\in M_j$, and $k\in\set{1,\ldots,K}$,
       $$\mu\paren{\set{x_{i,j}\in I_k} \setminus \set{z_{i,j}=c_k} } < \rho.$$
       Now for each $j\in L$ and $i\in M_j$, 
       $$\eqalign{\norm{z_{i,j}-x_{i,j}}_q
         &\le
         \paren{\tint_{\set{\abs{x_{i,j}}>N}} \abs{z_{i,j}-x_{i,j}}^q }^{{1\over q}}
     \cr &\phantom{=} + 
         \paren{\tint_{\bigcup_{k=1}^{K}\paren{\set{x_{i,j}\in I_k} \cap \set{z_{i,j}=c_k}} } \abs{z_{i,j}-x_{i,j}}^q }^{{1\over q}}
     \cr &\phantom{=} + 
         \tsuml_{k=1}^K \paren{\tint_{\set{x_{i,j}\in I_k} \setminus \set{z_{i,j}=c_k}} \abs{z_{i,j}-x_{i,j}}^q }^{{1\over q}}
     \cr &<
         2\epsilon + {\delta\over2} + K\rho^{{1\over q}}(2N+1),}$$
       where $K\rho^{{1\over q}}(2N+1) = D(2N+1)\delta^2(2N+1) = \delta(2N+1)^2$.  


       Fix $J\in\NN$ and assume $\set{1,\ldots,J}$ is a subset of $L$ and each $M_j$. 
       Then 
       $$\eqalign{&\norm{\tsuml_{i=1}^J \tsuml_{j=1}^J a_{i,j} x_{i,j} - \tsuml_{i=1}^J \tsuml_{j=1}^J a_{i,j} z_{i,j} }_q 
          \cr &\le
                   \tsuml_{i=1}^J \tsuml_{j=1}^J \abs{a_{i,j}} \norm{x_{i,j}-z_{i,j}}_q 
          \cr &\le
                   \tsuml_{i=1}^J \tsuml_{j=1}^J \abs{a_{i,j}}
                   \max_{i,j\in\set{1,\ldots,J}} \norm{x_{i,j}-z_{i,j}}_q 
          \cr &\le
                   \paren{\tsuml_{i=1}^J \paren{\tsuml_{j=1}^J \abs{a_{i,j}}^s }^{{r\over s}} }^{{1\over r}}
                   \abs{J}^{\paren{1-{1\over r}} + \paren{1-{1\over s}}}
                   \paren{2\epsilon + {\delta\over2} + \delta(2N+1)^2}.}$$  
       For any $J\in\NN$ and $\gamma>0$, we can choose $\epsilon>0$ and $\delta>0$ such that \hfil\break 
       $\abs{J}^{\paren{1-{1\over r}} + \paren{1-{1\over s}}} \paren{2\epsilon + {\delta\over2} + \delta(2N+1)^2} < \gamma$. 
       Hence we can find a symmetric sequence equivalent to the standard basis of $\sumsum{\l{s}}{\l{r}}$, 
       contrary to fact. 
       \qquad\QED

\medskip
 A basis $\set{e_i}$ for a Banach space $E$ is said to be reproducible if    
 for each Banach space $X$ with basis $\set{x_i}$ such that $E\injects X$,   
 there is a block basic sequence $\set{z_i}$ with respect to $\set{x_i}$ equivalent to $\set{e_i}$.   
 For $r,s\in[1,\infty)$, the standard basis of $\sumsum{\l{s}}{\l{r}}$ \hfil\break 
 is reproducible \xciteplus{L-P 2}{Section 4}.   

 The following proposition has been extracted from the proof of \xciteplus{S}{Theorem}.
 The subsequent corollary is essentially \xciteplus{S}{Remark 1}.   

\xproclaim {Proposition 3.8}. Let $1<q<2$ and let $n\in\NN$.
                              Then $\bigotimes\limits_{i=1}^{2n} \l{r_i} \not\injects X_q^{\otimes n}$  
                              for \hfil\break 
                              $q <r_1 <r_2 <\cdots <r_{2n}\le2$.  

\proof
       Suppose $n=1$.  
       Let $q<r<s\le2$ and suppose $\l{r}\otimes\l{s} \injects X_q$.  
       Then by Lemma 3.6, $\sumsum{\l{s}}{\l{r}} \injects X_q$.  
       Now $X_q \sim \span{x_{i,j}}_{\L{q}}$ for some sequence $\set{x_{i,j}}$ of in\-\hbox{dependent} random variables in $\L{q}$.  
       By the reproducibility of the standard basis $\set{e_{i,j}}$ of $\sumsum{\l{s}}{\l{r}}$,   
       there is a block basic sequence $\set{z_{i,j}}$ with respect to $\set{x_{i,j}}$ equivalent to $\set{e_{i,j}}$.   
       However, $\set{z_{i,j}}$ is a sequence of independent random variables in $\L{q}$ equiv\-\hbox{alent} to $\set{e_{i,j}}$,   
       contrary to Lemma 3.7.  
       Hence the result holds for $n=1$.   


       Suppose the result is true for $n=k-1$,
       but there are \hfil\break 
       $q <r_1 <r_2 <\cdots <r_{2k}\le2$ such that   
       $\bigotimes\limits_{i=1}^{2k} \l{r_i} \injects X_q^{\otimes k}$ via a mapping $\tau$.  
       
       Let $\set{e_{j_1}\otimes e_{j_2}\otimes \cdots\otimes e_{j_{2k}}}_{j_1,j_2,\ldots,j_{2k} \in\NN}$ 
       be the standard basis of $\bigotimes\limits_{i=1}^{2k} \l{r_i}$, \hfil\break 
       and let $y_{j_1,j_2,\ldots,j_{2k}} = \tau\paren{e_{j_1}\otimes e_{j_2}\otimes \cdots\otimes e_{j_{2k}}}$
       for $j_1,j_2,\ldots,j_{2k} \in\NN$.   

       Let $\set{x_j}$ be a basis for $X_q$.  
       For each $m\in\NN$, 
       let $P_m$ be the obvious projection of $X_q^{\otimes k}$ onto  
       $\span{x_{j_1}\otimes x_{j_2}\otimes \cdots\otimes x_{j_k} : \max\set{j_1,j_2,\ldots,j_k} \le m}$,  
       and let $Q_m$ be the obvious projection of $X_q^{\otimes k}$ onto  
       $\span{x_{j_1}\otimes x_{j_2}\otimes \cdots\otimes x_{j_k} : \min\set{j_1,j_2,\ldots,j_k} > m}$.  

       Recalling that $X_q \sim X_q \oplus X_q$,  
       for each $s\in\NN$ 
       $$X_q^{\otimes s} \sim (X_q \oplus X_q) \otimes X_q^{\otimes(s-1)} \sim X_q^{\otimes s} \oplus X_q^{\otimes s}.$$ 
       Hence for each $s,t\in\NN$, 
       $$\tsuml_{i=1}^t \!{}^{\oplus} X_q^{\otimes s} \sim X_q^{\otimes s}.$$ 

       Note that for each $m\in\NN$, 
       $(I-Q_m)(X_q^{\otimes k}) \sim \sum\limits_{i=1}^t \!{}^{\oplus} X_q^{\otimes(k-1)}$ for some $t\in\NN$, 
       whence $(I-Q_m)(X_q^{\otimes k}) \sim X_q^{\otimes(k-1)}$.  

       Let $\set{e_{j_1}\otimes e_{j_2}}_{j_1,j_2 \in\NN}$ be the standard basis of $\l{r_1} \otimes \l{r_2}$  
       with order determined by a bijection $\phi\colon\NN\to\NN\times\NN$.  

       For each $j\in\NN$, 
       let $Y_j = \span{y_{\phi(j),j_3,j_4,\ldots,j_{2k}} : j_3,j_4,\ldots,j_{2k} \in\NN}$,  
       which is \hfil\break 
       isomorphic to $\bigotimes\limits_{i=1}^{2(k-1)} \l{r_{i+2}}$.     
       Then by the inductive hypothesis, for each $j,m\in\NN$ 
       $$Y_j \sim         \textstyle\bigotimes\limits_{i=1}^{2(k-1)} \l{r_{i+2}}
             \not\injects X_q^{\otimes(k-1)}
             \sim         (I-Q_m)(X_q^{\otimes k}),$$
       whence $(I-Q_m)|_{Y_j}$ is not an isomorphism.
       
       Let $\set{\epsilon_j}$ be a sequence of positive scalars.  
       Let $m_0=0$ and $Q_{m_0}=I$.
       Choose $z_1 \in Y_1$ with $\norm{z_1}=1$ and $m_1\in\NN$ such that
       $\norm{(I-Q_{m_0})(z_1)} < {\epsilon_1\over2}$ and \hfil\break 
       $\norm{(I-P_{m_1})(z_1)} < {\epsilon_1\over2}$.
       Choose $z_2 \in Y_2$ with $\norm{z_2}=1$ and a positive integer $m_2 > m_1$ such that
       $\norm{(I-Q_{m_1})(z_2)} < {\epsilon_2\over2}$ and 
       $\norm{(I-P_{m_2})(z_2)} < {\epsilon_2\over2}$.
       Continuing as above, we
       may inductively define a sequence $\set{z_j}$ and an increasing sequence $\set{m_j}$ of positive \hfil\break 
       integers 
       such that for each $j\in\NN$, 
       $z_j \in Y_j$ with $\norm{z_j}=1$,
       $\norm{(I-Q_{m_{j-1}})(z_j)} < {\epsilon_j\over2}$, and 
       $\norm{(I-P_{m_j})(z_j)} < {\epsilon_j\over2}$.
       Hence for each $j\in\NN$, 
       $\norm{(I-Q_{m_{j-1}} \circ P_{m_j})(z_j)} < \epsilon_j \norm{P_{m_j}}$. 
       Thus for an appropriate choice of $\set{\epsilon_j}$, 
       $\set{z_j}$ is equivalent to $\set{(Q_{m_{j-1}} \circ P_{m_j})(z_j)}$. 
       However,
       $\set{z_j}$ is equivalent to the standard basis 
       $\set{e_{j_1}\otimes e_{j_2}}_{j_1,j_2 \in\NN}$ of $\l{r_1} \otimes \l{r_2}$,    
       and 
       $\set{(Q_{m_{j-1}} \circ P_{m_j})(z_j)}$ is a sequence of independent random variables. 
       Hence there is a sequence of independent random variables equivalent to the standard basis of \hfil\break 
       $\l{r_1} \otimes \l{r_2}$,    
       contrary to Lemma 3.7. 
       \qquad\QED

\xproclaim {Corollary 3.9}. Let $1<q<2$. Then for each $n\in\NN$, $X_q^{\otimes(n+1)}\not\injects X_q^{\otimes n}$.   

\proof Let $n\in\NN$ and let $q<r_1<r_2<\cdots<r_{2n}\le2$.  
       Then for each $1\le i\le2n$, 
       $\l{r_i}\injects X_q$ by Lemma 2.35.  
       Hence $\bigotimes\limits_{i=1}^{2n} \l{r_i} \injects X_q^{\otimes2n}$.  
       However, $\bigotimes\limits_{i=1}^{2n} \l{r_i} \not\injects X_q^{\otimes n}$ by \hfil\break 
       Proposition 3.8.  
       It follows that $X_q^{\otimes2n}\not\injects X_q^{\otimes n}$.   

       Now suppose that $X_q^{\otimes(n+1)}\injects X_q^{\otimes n}$.   
       Then there is a chain
       $$\cdots\injects X_q^{\otimes(n+2)}\injects X_q^{\otimes(n+1)}\injects X_q^{\otimes n}.$$ 
       In particular, $X_q^{\otimes2n}\injects X_q^{\otimes n}$, contrary to fact.   
       It follows that $X_q^{\otimes(n+1)}\not\injects X_q^{\otimes n}$.   
       \qquad\QED
       \medskip

 Note that $X\cinjects X\otimes Y$ 
 [where $1\le p<\infty$, $X$ and $Y$ are isomorphic to closed subspaces of $\L{p}$, and $\xdim Y>0$],   
 since $X \sim X\otimes\span{y_0} \cinjects X\otimes Y$ for $y_0\in Y\setminus\set{0}$.   
 Hence for $n\in\NN$ and $1<p<\infty$ with $p\ne2$,   
 $X_p^{\otimes n} \cinjects X_p^{\otimes(n+1)}$.    

 For $1<q<2$, we have  
 $$X_q \ra X_q^{\otimes2} \ra X_q^{\otimes3} \ra \cdots \ra \L{q}. \eqno{\TAG{3.1}{3.1}}$$

 Note that $(X\otimes Y)^* \sim X^*\otimes Y^*$  
 [where $1<p<\infty$, and $X$ and $Y$ are isomorphic to closed subspaces of $\L{p}$].   
 Let $2<p<\infty$ with conjugate index $q$. 
 Then for each $k\in\NN$, 
 $\paren{X_p^{\otimes k}}^* \sim \paren{X_p^*}^{\otimes k} \sim X_q^{\otimes k}$.
 Let $n\in\NN$.
 Then the fact that $X_p^{\otimes(n+1)} \not\cinjects X_p^{\otimes n}$ follows from 
 $\paren{X_p^{\otimes(n+1)}}^* \sim X_q^{\otimes(n+1)} \not\cinjects X_q^{\otimes n} \sim \paren{X_p^{\otimes n}}^*$.

 For $1<p<\infty$ with $p\ne2$, we have 
 $$X_p \cra X_p^{\otimes2} \cra X_p^{\otimes3} \cra \cdots \cra \L{p}. \eqno{\TAG{3.2}{3.2}}$$

 Finally we have the main result \xciteplus{S}{Theorem}. 

\xproclaim {Theorem 3.10}. Let $1<p<\infty$ where $p\ne2$.  
                           Then $\set{X_p^{\otimes n}}_{n=1}^{\infty}$ is a sequence of mutually nonisomorphic $\SL{p}$ spaces.  

\proof Each $X_p^{\otimes n}$ is an $\SL{p}$ space by Proposition 3.5. 
       For $m\ne n$, the fact that $X_p^{\otimes m} \not\sim X_p^{\otimes n}$ follows from Corollary 3.9
       and the discussion leading to diagrams \TAG{3.1}{3.1} and \TAG{3.2}{3.2}. 
       In particular, if $X_p^{\otimes m} \sim X_p^{\otimes n}$ for $m<n$, then 
       $X_p^{\otimes(m+1)} \cinjects X_p^{\otimes n} \cinjects X_p^{\otimes m}$, 
       contrary to fact. 
       \qquad\QED

\preheadspace
\firsthead{The Sequence Space Realization of $X_p^{\otimes n}$}
\postheadspace

 For $n\in\NN$, $X_p^{\otimes n}$ has a realization as a sequence space, 
 as follows from Proposition 3.13 below.
 This proposition is essentially contained in \xciteplus{S}{Section 4},  
 although the presentation via Lemmas 3.11 and 3.12 owes more to Dale Alspach.  

\xproclaim {Lemma 3.11}. Let $2<p<\infty$ and $k\in\NN$.  
                         Let $\set{x_i}$ be a sequence of normalized independent mean zero random variables in $\L{p}$.  
                         Let $\set{y_j}$ be an unconditional basic sequence in $\L{p}(I^k)$ 
                         with closed linear span $Y = \span{y_j}_{\L{p}(I^k)}$. 
                         Let $\set{r_i}$ be the sequence of Rademacher functions. 
                         Then for scalars $a_{i,j}$ 
           $$\eqalign{&\norm{\tsuml_{i}\tsuml_{j} a_{i,j} (x_i \otimes y_j) }_{\L{p}(I^{k+1})}
                  \cr &\approx \max\set{
                       \paren{\tsuml_{i} \norm{\tsuml_{j} a_{i,j} y_j}_Y^p }^{{1\over p}},  
                       \paren{\int\norm{\tsuml_{j}\paren{\tsuml_{i}a_{i,j}\norm{x_i}_2 r_i(u) }y_j }_Y^p du }^{{1\over p}} }.}$$


\proof For each $i\in\NN$, let $f_i(t) = \sum_{j} a_{i,j} y_j(t)$. 
       Then for each $t\in[0,1]$,
       $\set{x_i(s) f_i(t)}_{i=1}^{\infty}$ is a sequence of independent mean zero random variables in $\L{p}$. 
       Thus by Theorem 2.2 [Rosenthal's inequality], for each $t\in[0,1]$
       $$\eqalign{&\paren{\dint\abs{\tsuml_{i} x_i(s) f_i(t) }^p \,ds }^{{1\over p}}
              \cr &\within{}{} \max\set{
                   \paren{\tsuml_{i} \dint\abs{x_i(s) f_i(t) }^p \,ds }^{{1\over p}} , 
                   \paren{\tsuml_{i} \dint\abs{x_i(s) f_i(t) }^2 \,ds }^{{1\over 2}} }.}$$ 
       Hence 
       $$\eqalign{&\paren{\dint\!\dint\abs{\tsuml_{i} x_i(s) f_i(t) }^p \,ds\,dt }^{{1\over p}}
              \cr &\within{}{} \max\set{
                   \paren{\tsuml_{i} \dint\!\dint\abs{x_i(s) f_i(t) }^p \,ds\,dt }^{{1\over p}} , 
                   \paren{\dint\paren{
                   \tsuml_{i} \dint\abs{x_i(s) f_i(t) }^2 \,ds }^{{p\over 2}} dt }^{{1\over p}} }.}$$ 
       Now
       $$\dint\!\int\abs{x_i(s) f_i(t) }^p \,ds\,dt = \norm{x_i}_p^p \norm{f_i}_{\L{p}(I^k)}^p
                                                             = \norm{f_i}_{\L{p}(I^k)}^p
                                                             = \norm{\tsuml_{j} a_{i,j} y_j }_Y^p$$  
       and 
       $$\eqalign{\dint\paren{\tsuml_{i} \dint\abs{x_i(s) f_i(t) }^2 \,ds }^{{p\over 2}} dt 
               &= \dint\paren{\tsuml_{i} \norm{x_i}_2^2 \abs{f_i(t) }^2 }^{{p\over 2}} dt  
     \cr &\approx \dint\!\dint\abs{\tsuml_{i} \norm{x_i}_2 f_i(t) r_i(u) }^p \,du\,dt  
     \cr       &= \dint\!\dint\abs{\tsuml_{i}
                                           \norm{x_i}_2 \tsuml_{j} a_{i,j} y_j(t) r_i(u) }^p \,dt\,du  
     \cr &\approx \dint\norm{\tsuml_{j} \paren{\tsuml_{i} 
                                      a_{i,j} \norm{x_i}_2 r_i(u) } y_j }_Y^p \,du .}$$  

       \vglue1.5in
       \eject 

       \noindent Hence           
       $$\eqalign{&\norm{\tsuml_{i}\tsuml_{j} a_{i,j} (x_i \otimes y_j) }_{\L{p}(I^{k+1})}
              \cr &= \paren{\dint\!\dint\abs{\tsuml_{i}\tsuml_{j}
                            a_{i,j} x_i(s) y_j(t) }^p \,ds\,dt }^{{1\over p}}
        \cr &\approx \paren{\dint\!\dint\abs{\tsuml_{i}
                            x_i(s) \tsuml_{j} a_{i,j} y_j(t) }^p \,ds\,dt }^{{1\over p}}
              \cr &= \paren{\dint\!\dint\abs{\tsuml_{i} x_i(s) f_i(t) }^p \,ds\,dt }^{{1\over p}}
        \cr &\approx \max\set{
                     \paren{\tsuml_{i} \dint\!\dint\abs{x_i(s) f_i(t) }^p \,ds\,dt }^{{1\over p}} , 
                     \paren{\dint\paren{\tsuml_{i} \dint
                            \abs{x_i(s) f_i(t) }^2 \,ds }^{{p\over 2}} dt }^{{1\over p}} } 
        \cr &\approx \max\set{
                     \paren{\tsuml_{i} \norm{\tsuml_{j} a_{i,j} y_j }_Y^p }^{{1\over p}} ,  
                     \paren{\dint\norm{\tsuml_{j}
                            \paren{\tsuml_{i} a_{i,j} \norm{x_i}_2 r_i(u) } y_j }_Y^p \,du }^{{1\over p}} }.}$$ 
       \QED

 Let $\set{r_j}$ be the sequence of Rademacher functions. 
 Kahane's inequality \hfil\break 
 \xciteplus{W}{Theorem III.A.18}
 states that for each $1\le p<\infty$, there is a constant $C_p$ such that for each Banach space $X$ 
 and for each finite sequence $\set{x_j}$ in $X$, \hfil\break 
 $\paren{\int\norm{\sum_j r_j(u) x_j }_X^p \,du }^{{1\over p}}
 \within{C_p}{1}
 \int\norm{\sum_j r_j(u) x_j }_X \,du$.

\xproclaim {Lemma 3.12}. Let $1\le p<\infty$ and let $\set{r_j}$ be the sequence of Rademacher \hfil\break 
                         functions. 
                         Then for scalars $a_{i,j}$ 
                         $$\dint\paren{\tsuml_i
                           \abs{\tsuml_j a_{i,j} r_j(u) }^2 }^{{p\over2}} \,du
                           \approx 
                           \paren{\tsuml_i \tsuml_j \abs{a_{i,j}}^2 }^{{p\over2}}.$$

\vglue1in
\eject 

\proof Let $\set{e_i}$ be the standard basis of $\l{2}$. 
       Then by Kahane's inequality, 
       $$\eqalign{\dint\paren{\tsuml_i
                  \abs{\tsuml_j a_{i,j} r_j(u) }^2 }^{{p\over2}} \,du
                &= 
                  \dint\norm{\tsuml_i
                  \paren{\tsuml_j a_{i,j} r_j(u) }e_i }_{\l{2}}^p \,du
            \cr &= 
                  \dint\norm{\tsuml_j r_j(u) \paren{\tsuml_i a_{i,j} e_i} }_{\l{2}}^p \,du
      \cr &\within{C_p^p}{1}
                  \paren{\dint\norm{\tsuml_j r_j(u)
                         \paren{\tsuml_i a_{i,j} e_i} }_{\l{2}} \,du}^p
      \cr &\within{1}{C_2^p}
                  \paren{\dint\norm{\tsuml_j r_j(u)
                         \paren{\tsuml_i a_{i,j} e_i} }_{\l{2}}^2 \,du}^{{p\over2}}
            \cr &= 
                  \paren{\dint\norm{\tsuml_i
                         \paren{\tsuml_j a_{i,j} r_j(u) }e_i }_{\l{2}}^2 \,du }^{{p\over2}}
            \cr &= 
                  \paren{\dint\tsuml_i
                         \abs{\tsuml_j a_{i,j} r_j(u) }^2 \,du }^{{p\over2}}
            \cr &= 
                  \paren{\tsuml_i \dint
                         \abs{\tsuml_j a_{i,j} r_j(u) }^2 \,du }^{{p\over2}}
            \cr &= 
                  \paren{\tsuml_i \tsuml_j \abs{a_{i,j}}^2 }^{{p\over2}} .}$$
       \QED

\xproclaim {Proposition 3.13}. Let $2<p<\infty$ and $n\in\NN$. 
                               Let $\set{x_i}$ be a sequence of normalized independent mean zero random variables in $\L{p}$.  
                               For each $i\in\NN$, let \hfil\break 
                               $w_i=\norm{x_i}_2$. 
                               Then for scalars $a_{i_1,\ldots,i_n}$
          $$\eqalign{&\norm{
                      \tsuml_{i_1,\ldots,i_n}
                      a_{i_1,\ldots,i_n} \paren{x_{i_1}\otimes\cdots\otimes x_{i_n}}}_{\L{p}(I^n)}
         \cr &\approx \max_{S_n}
                      \set{\paren{\tsuml_{i_k:\,k\in S_n}
                           \paren{\tsuml_{i_{\ell}:\,\ell\in S_n^{c}}
                           \abs{a_{i_1,\ldots,i_n} }^2
                           \textstyle\prod\limits_{\ell\in S_n^{c}} w_{i_{\ell}}^2
                                 }^{{p\over2}} }^{{1\over p}} } }$$
          where the $\max$ is taken over all subsets $S_n$ of $\set{1,\ldots,n}$, 
          and $S_n^{c}=\set{1,\ldots,n}\setminus S_n$.


\proof For $n=1$ [with $i_1=i$], the statement is 
       $$\eqalign{\norm{\tsuml_i a_i x_i}_p
                  &\approx 
                  \max\set{\paren{\paren{\tsuml_i \abs{a_i}^2 w_i^2 }^{{p\over2}} }^{{1\over p}},
                           \paren{\tsuml_i \paren{\abs{a_i}^2 }^{{p\over2}} }^{{1\over p}} } 
              \cr &=
                  \max\set{\paren{\tsuml_i \abs{a_i}^2 w_i^2 }^{{1\over2}},
                           \paren{\tsuml_i \abs{a_i}^p }^{{1\over p}} } ,}$$
       which is immediate from Corollary 2.3 [Rosenthal's inequality].  

       Assume the statement is true for $n=N$.
       We wish to prove the statement for $n=N+1$.   

       Let $\set{r_i}$ be the sequence of Rademacher functions.  
       By Lemma 3.11, 
       $$\norm{\tsuml_{i_1,\ldots,i_N} \tsuml_{i_{N+1}}
               a_{i_1,\ldots,i_{N+1}} \paren{x_{i_1}\otimes\cdots\otimes x_{i_N} }\otimes x_{i_{N+1}} }_{\L{p}(I^{N+1})}
         \approx \max\set{E_1,E_2} $$
       where 
       $$E_1 = \paren{\tsuml_{i_{N+1}} \norm{\tsuml_{i_1,\ldots,i_N}
               a_{i_1,\ldots,i_{N+1}} \paren{x_{i_1}\otimes\cdots\otimes x_{i_N} } }_{\L{p}(I^N)}^p }^{{1\over p}}$$
       and 
       $$E_2 = \paren{\dint\norm{\tsuml_{i_1,\ldots,i_N}
                      \paren{\tsuml_{i_{N+1}} a_{i_1,\ldots,i_{N+1}} \norm{x_{i_{N+1}}}_2 r_{i_{N+1}}(u) }
                      \paren{x_{i_1}\otimes\cdots\otimes x_{i_N} } }_{\L{p}(I^N)}^p \,du}^{{1\over p}}.$$
       Let
       $$A_{i_1,\ldots,i_N}(u)=\tsuml_{i_{N+1}} a_{i_1,\ldots,i_{N+1}} \norm{x_{i_{N+1}}}_2 r_{i_{N+1}}(u)$$ 
       and 
       $$B_{i_1,\ldots,i_{N+1}}^{(S_N^{c})} = a_{i_1,\ldots,i_{N+1}} \norm{x_{i_{N+1}}}_2
                                          \textstyle\prod\limits_{\ell\in S_N^{c}} w_{i_{\ell}}.$$
       By the inductive hypothesis, and then a rearrangement, we have 
       $$\eqalign{
       E_1 &\approx 
            \paren{\tsuml_{i_{N+1}} \displaystyle\max_{S_N}
                   \set{\tsuml_{i_k:\,k\in S_N}
                   \paren{\tsuml_{i_{\ell}:\,\ell\in S_N^{c}}
                   \abs{a_{i_1,\ldots,i_{N+1}} }^2
                   \textstyle\prod\limits_{\ell\in S_N^{c}} w_{i_{\ell}}^2
                   }^{{p\over2}} } }^{{1\over p}}
       \cr &\approx 
            \max_{S_N}\set{\paren{\tsuml_{i_k:\,k\in S_N \cup\set{N+1}}
                           \paren{\tsuml_{i_{\ell}:\,\ell\in S_N^{c}}
                           \abs{a_{i_1,\ldots,i_{N+1}} }^2
                           \textstyle\prod\limits_{\ell\in S_N^{c}} w_{i_{\ell}}^2
                           }^{{p\over2}} }^{{1\over p}} }
       \cr &= 
            \max_{\textstyle{{S_{N+1}: \atop N+1\in S_{N+1}}} } \set{\paren{\tsuml_{i_k:\,k\in S_{N+1}}
            \paren{\tsuml_{i_{\ell}:\,\ell\in S_{N+1}^{c}}
            \abs{a_{i_1,\ldots,i_{N+1}} }^2
            \textstyle\prod\limits_{\ell\in S_{N+1}^{c}} w_{i_{\ell}}^2
            }^{{p\over2}} }^{{1\over p}} }. }$$
       By the inductive hypothesis, a rearrangement, and Lemma 3.12, we have 
       $$\eqalign{E_2 &= \paren{\dint\norm{\tsuml_{i_1,\ldots,i_N}
                         A_{i_1,\ldots,i_N}(u) 
                         \paren{x_{i_1}\otimes\cdots\otimes x_{i_N} } }_{\L{p}(I^N)}^p \,du}^{{1\over p}} 
                  \cr &\approx 
                         \paren{\dint\displaystyle\max_{S_N}
                         \set{\tsuml_{i_k:\,k\in S_N}
                         \paren{\tsuml_{i_{\ell}:\,\ell\in S_N^{c}}
                         \abs{A_{i_1,\ldots,i_N}(u) }^2
                         \textstyle\prod\limits_{\ell\in S_N^{c}} w_{i_{\ell}}^2
                               }^{{p\over2}} } \,du }^{{1\over p}}
                  \cr &\approx 
                         \max_{S_N}
                         \set{\paren{\tsuml_{i_k:\,k\in S_N}
                         \dint
                         \paren{\tsuml_{i_{\ell}:\,\ell\in S_N^{c}}
                         \abs{A_{i_1,\ldots,i_N}(u) }^2
                         \textstyle\prod\limits_{\ell\in S_N^{c}} w_{i_{\ell}}^2
                               }^{{p\over2}} \,du }^{{1\over p}} }
                  \cr &=       
                         \max_{S_N}
                         \set{\paren{\tsuml_{i_k:\,k\in S_N}
                         \dint
                         \paren{\tsuml_{i_{\ell}:\,\ell\in S_N^{c}}
                         \abs{\tsuml_{i_{N+1}} B_{i_1,\ldots,i_{N+1}}^{(S_N^{c})} r_{i_{N+1}}(u) }^2
                               }^{{p\over2}} \,du }^{{1\over p}} }
                  \cr &\approx 
                         \max_{S_N}
                         \set{\paren{\tsuml_{i_k:\,k\in S_N}
                         \paren{\tsuml_{i_{\ell}:\,\ell\in S_N^{c} \cup\set{N+1}}
                         \abs{B_{i_1,\ldots,i_{N+1}}^{(S_N^{c})} }^2
                               }^{{p\over2}} }^{{1\over p}} }
                  \cr &=       
                         \max_{S_N}
                         \set{\paren{\tsuml_{i_k:\,k\in S_N}
                         \paren{\tsuml_{i_{\ell}:\,\ell\in S_N^{c} \cup\set{N+1}}
                         \abs{a_{i_1,\ldots,i_{N+1}} }^2
                         \textstyle\prod\limits_{\ell\in S_N^{c} \cup\set{N+1}} w_{i_{\ell}}^2
                               }^{{p\over2}} }^{{1\over p}} }
                  \cr &=       
                         \max_{\textstyle{{S_{N+1}: \atop N+1\notin S_{N+1}}} }
                         \set{\paren{\tsuml_{i_k:\,k\in S_{N+1}}
                         \paren{\tsuml_{i_{\ell}:\,\ell\in S_{N+1}^{c}}
                         \abs{a_{i_1,\ldots,i_{N+1}} }^2
                         \textstyle\prod\limits_{\ell\in S_{N+1}^{c}} w_{i_{\ell}}^2
                               }^{{p\over2}} }^{{1\over p}} }. }$$
       Hence  
       $$\eqalign{&\norm{
                   \tsuml_{i_1,\ldots,i_{N+1}}
                   a_{i_1,\ldots,i_{N+1}} \paren{x_{i_1}\otimes\cdots\otimes x_{i_{N+1}}}}_{\L{p}(I^{N+1})}
      \cr &\approx \max\set{E_1,E_2}
      \cr &\approx \max_{S_{N+1}}
                   \set{\paren{\tsuml_{i_k:\,k\in S_{N+1}}
                        \paren{\tsuml_{i_{\ell}:\,\ell\in S_{N+1}^{c}}
                        \abs{a_{i_1,\ldots,i_{N+1}} }^2
                        \textstyle\prod\limits_{\ell\in S_{N+1}^{c}} w_{i_{\ell}}^2
                              }^{{p\over2}} }^{{1\over p}} } .}$$
       \QED

 For $2<p<\infty$ and $n\in\NN$, 
 Proposition 3.13 yields a representation of $X_p^{\otimes n}$ as a sequence space, 
 taking $\set{x_i}$ to be a sequence of normalized independent mean zero random variables in $\L{p}$
 with $w=\set{w_i}=\set{\norm{x_i}_2}$ satisfying condition $(*)$ of Proposition 2.1. 


 In particular, for $n=2$ and $S_2\subset\set{i,j}$, for scalars $a_{i,j}$
 $$\norm{\sum_{i,j} a_{i,j} \paren{x_i\otimes y_j}}_{\L{p}(I^2)}
   \approx \max\set{{\cal N}_{[S_2=\emptyset]}, {\cal N}_{[S_2=\set{i}]}, {\cal N}_{[S_2=\set{j}]}, {\cal N}_{[S_2=\set{i,j}]}},$$
 where 
 $$\eqalign{\hbox to 0.65 in {${\cal N}_{[S_2=\emptyset]}$\hss}
                                     &= \paren{\paren{\tsuml_{i,j}\abs{a_{i,j}}^2 w_i^2 w_j^2}^{{p\over2}}}^{{1\over p}}
                                      = \paren{\tsuml_{i,j}\abs{a_{i,j}}^2 w_i^2 w_j^2}^{{1\over2}},\cr
            \hbox to 0.65 in {${\cal N}_{[S_2=\set{i}]}$\hss}
                                     &= \paren{\tsuml_i\paren{\tsuml_j\abs{a_{i,j}}^2 w_j^2}^{{p\over2}}}^{{1\over p}},\cr
            \hbox to 0.65 in {${\cal N}_{[S_2=\set{j}]}$\hss}
                                     &= \paren{\tsuml_j\paren{\tsuml_i\abs{a_{i,j}}^2 w_i^2}^{{p\over2}}}^{{1\over p}},\cr
            \hbox to 0.65 in {${\cal N}_{[S_2=\set{i,j}]}$\hss}
                                     &= \paren{\tsuml_{i,j}\paren{\abs{a_{i,j}}^2}^{{p\over2}}}^{{1\over p}}
                                      = \paren{\tsuml_{i,j}\abs{a_{i,j}}^p}^{{1\over p}}.}$$


%
%
%
\chapter{IV}
\headspace
\chaptertitle{THE INDEPENDENT SUM CONSTRUCTION OF ALSPACH}
\headspace
 Let $2<p<\infty$ and let $\Omega=\prod_{i=1}^{\infty}[0,1]$. 
 Alspach \xcite{A} developed a general method for constructing complemented subspaces of $\L{p}(\Omega)$, 
 given spaces $X_i$ of mean zero functions which are complemented in $\L{p}[0,1]$ in a special way. 
 The \hbox{construction} produces spaces $Z_i$ of mean zero functions which are similarly complemented in $\L{p}(\Omega)$,
 such that $Z_i$ is isometric to $X_i$, 
 each function in $Z_i$ depends only on component $i$ of $\Omega$, 
 there is a common supporting set $S_i$ for all functions in $Z_i$,
 and the measure of $S_i$ approaches zero slowly as $i$ increases.
 The independent sum of $\set{X_i}_{i=1}^{\infty}$ is then $\span{Z_i:i\in\NN}_{\L{p}(\Omega)}$.

 The rate at which the measure of $S_i$ approaches zero is controlled by a sequence $w$, 
 which plays a role similar to the role of $w$ in Rosenthal's space $\Xpw$. 
 Indeed, \hfil\break 
 Alspach's construction generalizes the construction of Rosenthal's space $\Xpw$.

 All of the $\SL{p}$ spaces of Chapter II can be constructed as independent sums in the above sense.
 The principal new separable infinite-dimensional $\SL{p}$ space constructed by Alspach as an independent sum is $D_p$, 
 which is the independent sum of copies of $\l{2}$, with $\l{2}$ realized as the span of the Rademachers in $\L{p}$. 
 Also new is $B_p \oplus D_p$. 
 The method of taking independent sums has the potential to generate a sequence of $\SL{p}$ spaces by iteration. 
 However, no general method has been developed for distinguishing the isomorphism types of the resulting spaces. 

\preheadspace
\firsthead{The Independent Sum $\sumsum{X_i}{I,w}$}
\postheadspace
 Fix $2<p<\infty$.  
 Let $\Omega=\prod_{i=1}^{\infty}[0,1]$. 
 For $t=\paren{t_1,t_2,\ldots}\in\Omega$ and $i\in\NN$, let $\pi_i:\Omega\to[0,1]$ be the projection $\pi_i(t)=t_i$.  
 Let $\Lz{p}[0,1]$ be the space of mean zero functions in $\L{p}[0,1]$. 
 For $0<k\le1$, identify $\L{p}[0,k]$ with the space of functions in $\L{p}[0,1]$ supported on $[0,k]$.  
 Let $\set{X_i}$ be a sequence of closed subspaces of $\Lz{p}[0,1]$.
 Let $w=\set{w_i}$ and $\set{k_i}$ be sequences of scalars from $(0,1]$ such that $k_i=w_i^{{2p\over p-2}}$. 
 Let $T_i:\L{p}[0,1]\to\L{p}[0,k_i]\subset\L{p}[0,1]$ be defined by
 $$T_i(f)(s) = \cases{k_i^{-{1\over p}} f\paren{\textstyle{s\over k_i}} & if $0\le s\le k_i$ \cr 0 & if $k_i<s\le1$ \cr}.$$
 Let $Y_i=T_i(X_i)$ and let $\tilde Y_i = \set{\tilde y_i = y_i\circ\pi_i : y_i \in Y_i} \subset \L{p}(\Omega)$.

\definition Let $p$, $\Omega$, $\pi_i$, $\set{X_i}$, $w=\set{w_i}$, $\set{k_i}$, $T_i$, $Y_i$, and $\tilde Y_i$ be as above. 
 \hfil\break 
 Suppose \hfil\break 
 (a) for each $i\in\NN$,
     the orthogonal projection of $\L{2}[0,1]$ onto $\overline{X_i}\subset\L{2}[0,1]$, when \hfil\break 
     restricted to $\L{p}[0,1]$,
     yields a bounded projection $P_i:\L{p}[0,1] \to X_i\subset\L{p}[0,1]$ onto $X_i$, and \hfil\break 
 (b) the sequence $\set{P_i}_{i=1}^{\infty}$ satisfies $\sup_{i\in\NN} \norm{P_i} < \infty$. \hfil\break 
 Define $\sumsum{X_i}{I,w}$, the independent sum of $\set{X_i}$ with respect to $w$, by 
 $$\sumsum{X_i}{I,w} = \span{\tilde Y_i:i\in\NN}_{\L{p}(\Omega)}.$$

\remark The mapping $T_i$ is an isometry, and the spaces $X_i$, $Y_i$, and $\tilde Y_i$ are isometric.
        If $\tilde y_i \in \tilde Y_i$ for each $i\in\NN$, 
        then $\set{\tilde y_i}_{i=1}^{\infty}$ is a sequence of independent mean zero random variables. 
        The sequence $w$ plays a role similar to the role of $w$ in Rosenthal's space $\Xpw$. 
        In particular, $w_i^{{2p\over p-2}}$ is related to the measure of the support of $\tilde y_i \in \tilde Y_i$. 

\xproclaim {Example 4.1}. Let $2<p<\infty$, 
                          let $r_1$ be the first Rademacher function \hfil\break 
                          $1_{[0,{1\over2})} - 1_{[{1\over2},1]}$,
                          let $X = \span{r_1}_{\L{p}[0,1]}$, and 
                          let $w=\set{w_i}$ be a sequence from $(0,1]$.
                          Then $\sumsum{X}{I,w}$ is isomorphic to $\l{2}$, $\l{p}$, $\l{2}\oplus\l{p}$, or $X_p$, 
                          where each can be realized by an appropriate choice of $w$ as in Proposition 2.1.

\proof Let $\set{k_i}$ and $\set{T_i}$ correspond with $w=\set{w_i}$ as above. 
       Let $y_i=T_i(r_1)$ and $\tilde y_i = y_i\circ\pi_i$.
       Then $\sumsum{X}{I,w} = \span{\tilde y_i:i\in\NN}_{\L{p}(\Omega)}$.
       Now $\set{\tilde y_i}_{i=1}^{\infty}$ is a sequence of
       independent symmetric three-valued random variables in $\L{p}(\Omega)$,  
       with $\tilde y_i$ supported on a set of measure $k_i = w_i^{{2p\over p-2}}$.  
       Moreover, $w_i = k_i^{{p-2\over2p}} = k_i^{{1\over2} - {1\over p}}
                      = \norm{\tilde y_i}_{\L{2}(\Omega)}\big/\norm{\tilde y_i}_{\L{p}(\Omega)}$.  
       Hence $\sumsum{X}{I,w} \sim \Xpw$ (essentially) by Corollary 2.3,
       so $\sumsum{X}{I,w}$ is \hfil\break 
       isomorphic to $\l{2}$, $\l{p}$, $\l{2}\oplus\l{p}$, or $X_p$,
       depending on $w$ as in Proposition 2.1 and the \hfil\break 
       definition of $X_p$.  
       \qquad\QED

\preheadspace
\firsthead{The Complementation of $\sumsum{X_i}{I,w}$ in $\L{p}(\Omega)$}
\postheadspace
 Fix $2<p<\infty$ and $0<k\le1$.  
 For $1\le r<\infty$, identify $\L{r}[0,k]$ with the space of functions in $\L{r}[0,1]$ supported on $[0,k]$,  
 and for a measure space $E$, let $\Lz{r}(E)$ be the space of mean zero functions in $\L{r}(E)$. 

 Let $T:\L{1}[0,1]\to\L{1}[0,k]\subset\L{1}[0,1]$ be defined by
 $$T(f)(s) = \cases{k^{-{1\over p}} f\paren{\textstyle{s\over k}} & if $0\le s\le k$ \cr 0 & if $k<s\le1$ \cr}.$$
 For $1\le r<\infty$, let $T_r=T|_{\L{r}[0,1]}$.

\xproclaim {Lemma 4.2}. Let $p$, $k$, and $T$ be as above.
                        For $1\le r<\infty$, let $f,g\in\L{r}[0,1]$. Then 
             \item{(a)} $\norm{T(f)}_r = k^{{p-r\over rp}}\norm{f}_r$, 
             \item{(b)} $T_r:\L{r}[0,1]\to\L{r}[0,k]\subset\L{r}[0,1]$, 
             \item{(c)} $T_r$ maps $\L{r}[0,1]$ onto $\L{r}[0,k]$,
             \item{(d)} $T_p$ is an isometry,
             \item{(e)} $T_p=T_2|_{\L{p}[0,1]}$,
             \item{(f)} $f$ has mean zero if and only if $T(f)$ has mean zero, and
             \item{(g)} $f$ and $g$ are orthogonal if and only if $T(f)$ and $T(g)$ are orthogonal.  

\proof Part (a) follows from the computation  \hfil\break 
       $\norm{T(f)}_r^r = \int_0^k \abs{T(f)(s)}^r \,ds
                        = \int_0^k \abs{k^{-{1\over p}} f\paren{{s\over k}} }^r \,ds
                        = k^{1-{r\over p}} \int_0^1 \abs{f(t)}^r \,dt
                        = k^{{p-r\over p}} \norm{f}_r^r$. \hfil\break 
       Part (b) follows from (a) and the definition of $T$.
       Considering $T_r$ as a mapping from $\L{r}[0,1]$ to $\L{r}[0,k]$,
       $T_r$ has inverse $T_r^{-1}:\L{r}[0,k]\to\L{r}[0,1]$ with $T_r^{-1}(h)(t)=k^{{1\over p}}h(kt)$, and (c) follows.
       Taking $r=p$, (d) follows from (a).
       Part (e) is clear.
       As in the \hfil\break 
       computation for (a), but taking $r=1$ and deleting the absolute values, \hfil\break 
       $\int_0^k T(f)(s) \,ds = k^{1-{1\over p}} \int_0^1 f(t) \,dt$, and (f) follows.
       Finally, 
       $\int_0^k T(f)(s) \cdot T(g)(s) \,ds = k^{-{2\over p}} \int_0^k f\paren{{s\over k}} \cdot g\paren{{s\over k}} \,ds
                                            = k^{1-{2\over p}} \int_0^1 f(t) \cdot g(t) \,dt$, 
       and (g) follows.
       \qquad\QED
  
 Let $R:\L{1}[0,1]\to\L{1}[0,k]$ be defined by $R(f) = 1_{[0,k]} \cdot f$.
 For $1\le r<\infty$, let $R_r=R|_{\L{r}[0,1]}$.

 Let $X$ be a closed subspace of $\Lz{p}[0,1]$ such that the orthogonal projection $P_2$
 of $\L{2}[0,1]$ onto $\overline{X}\subset\L{2}[0,1]$,
 when restricted to $\L{p}[0,1]$, yields a bounded projection 
 $P_p:\L{p}[0,1] \to X\subset\L{p}[0,1]$ onto $X$.
 Let $Y=T(X)$.

\xproclaim {Lemma 4.3}. Let $p$, $k$, $T$, $R$, $X$, $P_2$, $P_p$, and $Y$ be as above.
                        Let $1\le r<\infty$. Then 
             \item{(a)} $R_r:\L{r}[0,1]\to\L{r}[0,k]$ is a projection of $\L{r}[0,1]$ onto $\L{r}[0,k]$ with $\norm{R_r}=1$,  
             \item{(b)} $R_2$ is the orthogonal projection of $\L{2}[0,1]$ onto $\L{2}[0,k]$, 
             \item{(c)} $R_p=R_2|_{\L{p}[0,1]}$, 
             \item{(d)} $Y$ is a subspace of $\Lz{p}[0,k]$ isometric to $X$,
             \item{(e)} the closure of $X$ in $\L{2}[0,1]$ is contained in $\Lz{2}[0,1]$,
             \item{(f)} the closure of $Y$ in $\L{2}[0,k]$ is contained in $\Lz{2}[0,k]$,
             \item{(g)} $T_2\paren{\overline{X}}=\overline{Y}$, 
                        where $\overline{X}$ and $\overline{Y}$ are the closures of $X$ and $Y$ in $\L{2}[0,1]$,
             \item{(h)} $T_2 P_2 T_2^{-1}$ is the orthogonal projection of $\L{2}[0,k]$ onto $\overline{Y}\subset\L{2}[0,k]$,
             \item{(i)} $T_p P_p T_p^{-1}=\paren{T_2 P_2 T_2^{-1}}|_{\L{p}[0,k]}$, and
             \item{(j)} $T_p P_p T_p^{-1}$ maps $\L{p}[0,k]$ onto $Y$.
 
\proof Part (a) is clear.
       For $f,g\in\L{2}[0,1]$, $\paren{f-R_2(f)} \perp R_2(g)$, so \hfil\break 
       $\paren{f-R_2(f)}\in\paren{R_2\paren{\L{2}[0,1]}}^{\perp}$, and (b) follows. 
       Part (c) is clear.
       Part (d) follows from the fact that $T_p:\L{p}[0,1]\to\L{p}[0,k]$ is an isometry which preserves mean zero functions. 
       First noting that $X\subset\Lz{2}[0,1]$ and $Y\subset\Lz{2}[0,k]$, parts (e) and (f) are clear.
       Part (g) is clear.
       For $f,g\in\L{2}[0,k]$, $\paren{T_2^{-1}(f)-P_2\paren{T_2^{-1}(f)}} \perp P_2\paren{T_2^{-1}(g)}$,
       so $\paren{f-\paren{T_2 P_2 T_2^{-1}}(f)} \perp \paren{T_2 P_2 T_2^{-1}}(g)$,
       and (h) follows after noting (g).
       Parts (i) and (j) are clear.
       \qquad\QED

 For $r\in\set{2,p}$, let $Q_r=T_r P_r T_r^{-1} R_r$.

\xproclaim {Lemma 4.4}. Let $p$, $r$, $k$, $T$, $R$, $X$, $P_r$, $Y$, and $Q_r$ be as above. Then 
             \item{(a)} $Q_p:\L{p}[0,1]\to Y\subset\L{p}[0,1]$ maps $\L{p}[0,1]$ onto $Y$, 
             \item{(b)} $\norm{Q_p}=\norm{P_p}$, 
             \item{(c)} $Q_2$ is the orthogonal projection of $\L{2}[0,1]$ onto $\overline{Y}\subset\L{2}[0,1]$, 
             \item{(d)} $Q_p=Q_2|_{\L{p}[0,1]}$, and
             \item{(e)} $Q(1)=0$. 

\proof Note that $T_p^{-1} R_p:\L{p}[0,1] \to \L{p}[0,1]$ is surjective, with right inverse $T_p$.  
       Thus (a) follows, and $Q_p T_p = \paren{T_p P_p T_p^{-1} R_p} T_p = T_p P_p \paren{T_p^{-1} R_p T_p} = T_p P_p$.
       Since $T_p$ is an isometry, (b) follows.
       Part (c) follows from the fact that $R_2$ and $T_2 P_2 T_2^{-1}$ are orthogonal projections
       mapping $\L{2}[0,1]$ onto $\L{2}[0,k]$ and $\L{2}[0,k]$ onto $\overline{Y}\subset\L{2}[0,k]$, respectively.
       Part (d) follows from the fact that  
       $R_p=R_2|_{\L{p}[0,1]}$ and $T_p P_p T_p^{-1} = \paren{T_2 P_2 T_2^{-1}}|_{\L{p}[0,k]}$. 
       Noting that $1 \Tmapsto{R_p} 1_{[0,k]} \Tmapsto{T_p^{-1}} k^{{1\over p}}\cdot1_{[0,1]} \Tmapsto{P_p} 0 \Tmapsto{T_p} 0$,
       (e) follows. 
       \qquad\QED
   
 The relevant subspaces of $\L{p}[0,1]$ are related as in the diagram
   $$\matrix{\L{p}[0,1]            & \Tra{P_p}  & X        & \subset\Lz{p}[0,1]\subset\L{p}[0,1] \cr  
                                   & & & \cr 
             \TLda{R_p}\TRua{T_p^{-1}} & \Tdsa{Q_p} & \TRda{T_p} & \cr 
                                   & & & \cr 
             \L{p}[0,k]            & \Tra{T_pP_pT_p^{-1}}&Y& \subset\Lz{p}[0,k]\subset\L{p}[0,1] .}\eqno{\TAG{4.1}{4.1}}$$

 We now perform a similar construction for each $i\in\NN$. 

 Let $\set{k_i}$ be a sequence of scalars from $(0,1]$. 
 Then for $r\in\set{1,2,p}$, $\set{k_i}$ determines sequences $\set{T_{i,r}}$ and $\set{R_{i,r}}$ of mappings, 
 where $T_{i,r}$ and $R_{i,r}$ are simply $T_r$ and $R_r$, respectively, with $k_i$ replacing $k$. 
 Let $\set{X_i}$ be a sequence of closed subspaces of $\Lz{p}[0,1]$ 
 such that the orthogonal projection $P_{i,2}$ of $\L{2}[0,1]$ onto $\overline{X_i}\subset\L{2}[0,1]$,
 when restricted to $\L{p}[0,1]$,
 yields a bounded projection $P_{i,p}:\L{p}[0,1] \to X_i\subset\L{p}[0,1]$ onto $X_i$.
 Let $Y_i=T_{i,p}(X_i)$, and for $r\in\set{2,p}$, let $Q_{i,r} = T_{i,r} P_{i,r} T^{-1}_{i,r} R_{i,r}$. 
 Then $X_i$, $Y_i$, $P_{i,r}$, and $Q_{i,r}$ are simply $X$, $Y$, $P_r$, and $Q_r$, respectively, with $k_i$ replacing $k$. 
 Thus as in diagram \TAG{4.1}{4.1}, we have the diagram
 $$\matrix{\L{p}[0,1]                & \Tra{P_{i,p}} & X_i & \subset\Lz{p}[0,1]\subset\L{p}[0,1] \cr  
                                     & & & \cr 
           \TLda{R_{i,p}}\TRua{T_{i,p}^{-1}} & \Tdsa{Q_{i,p}} & \TRda{T_{i,p}} & \cr 
                                     & & & \cr 
           \L{p}[0,k_i]              & & Y_i       
                                       & \subset\Lz{p}[0,k_i]\subset\L{p}[0,1] ,}\eqno{\TAG{4.2}{4.2}}$$
 and Lemmas 4.2, 4.3, and 4.4 hold, with the obvious notational changes.

 Let $1\le r<\infty$ and let $i\in\NN$.
 Let $\Pi_{i,r}:\L{r}[0,1]\to\L{r}[\Omega]$ be the isometry \hfil\break 
 $\Pi_{i,r}(f)=f\circ\pi_i$.  
 Then for $f,g\in\L{r}[0,1]$, 
 $f$ has mean zero if and only if $\Pi_{i,r}(f)$ has mean zero, and
 $f$ and $g$ are orthogonal if and only if $\Pi_{i,r}(f)$ and $\Pi_{i,r}(g)$ are \hfil\break 
 orthogonal.  

 Given a closed subspace $Z_{i,r}$ of $\L{r}[0,1]$, let $\tilde Z_{i,r}=\Pi_{i,r}\paren{Z_{i,r}}\subset\L{r}(\Omega)$. 
 Let $\tilde L{}_i^r[0,1]=\Pi_{i,r}\paren{\L{r}[0,1]}$  
 and $\tilde L{}_{0,i}^r[0,1]=\Pi_{i,r}\paren{L_0^r[0,1]}$. 

 Given closed subspaces $Z_{i,r}$ and $Z'_{i,r}$ of $\L{r}[0,1]$ and a mapping \hfil\break 
 $L_{i,r}:Z_{i,r} \to Z'_{i,r}$, 
 let $\tilde L_{i,r} : \tilde Z_{i,r} \to \tilde Z{}'_{i,r}$ be the mapping $\tilde L_{i,r} = \Pi_{i,r} L_{i,r} \Pi_{i,r}^{-1}$.
 Then diagram \TAG{4.2}{4.2} induces the diagram
 $$\matrix{\tilde L{}_i^p[0,1] & \Tra{\tilde P_{i,p}}  & \tilde X_{i,p} 
                                   & \subset\tilde L{}_{0,i}^p[0,1]\subset\tilde L{}_i^p[0,1] \cr  
                                   & & & \cr 
           \TLda{\tilde R_{i,p}}\TRua{\tilde T{}_{i,p}^{-1}}
                                   & \Tdsa{\tilde Q_{i,p}} & \TRda{\tilde T_{i,p}} \cr 
                                   & & & \cr 
           \tilde L{}^p[0,k_i]     &                       & \tilde Y_{i,p}
                                   & \subset\tilde L{}_0^p[0,k_i]\subset\tilde L{}_i^p[0,1] ,}\eqno{\TAG{4.3}{4.3}}$$
 and results analogous to Lemmas 4.2, 4.3, and 4.4 hold.

 Let $E_i:\L{1}(\Omega)\to\tilde L{}_i^1[0,1]\subset\L{1}(\Omega)$ be the projection onto
 $\tilde L{}_i^1[0,1] = \Pi_{i,1}\paren{\L{1}[0,1]}$ of norm one defined by $E_i(f)=\E_{\S{B}_i}f$,
 where $\E_{\S{B}_i}$ is conditional expectation with \hfil\break 
 respect to the $\sigma$-algebra
 ${\cal B}_i = \set{\prod_{j=1}^{\infty}\!B_j:B_i\!\subset\![0,1] \hbox{ is measurable, } B_j\!=\![0,1] \hbox{ for } j\!\ne\!i}$. 
 For $1<r<\infty$, let $E_{i,r}=E_i|_{\L{r}(\Omega)}$. 
 [See Chapter V, The Complementation of $\Rpa$ \hfil\break 
 in $\L{p}$, Preliminaries, for properties of conditional expectation.]

\xproclaim {Lemma 4.5}. Let $p$, $\Pi_{i,r}$, $\tilde L{}_i^r[0,1]$, ${\cal B}_i$, and $E_i$ be as above 
                        for $1<r<\infty$ with conjugate index $s$,
                        and let $f\in\L{r}(\Omega)$.
                        Then 
             \item{(a)} $E_{i,r}:\L{r}(\Omega)\to\L{r}(\Omega)$ with $\norm{E_{i,r}} = 1$,
             \item{(b)} $E_{i,r}$ maps $\L{r}(\Omega)$ onto $\tilde L{}_i^r[0,1]=\Pi_{i,r}\paren{\L{r}[0,1]}$,
             \item{(c)} $f$ has mean zero if and only if $E_{i,r}(f)$ has mean zero, 
             \item{(d)} if $\set{f_i}_{i=1}^{\infty}$ is a sequence in $\L{r}(\Omega)$, 
                        then $\set{E_{i,r}(f_i)}_{i=1}^{\infty}$ is independent,  
             \item{(e)} $E_{i,r}^* = E_{i,s}$, 
             \item{(f)} $E_{i,2}$ is the orthogonal projection of $\L{2}(\Omega)$ onto $\tilde L{}_i^2[0,1]$, and  
             \item{(g)} $E_{i,p} = E_{i,2}|_{\L{p}(\Omega)}$. 

\proof By the convexity of $\abs{\pz}^r$, 
       $\int_{\Omega} \abs{E_i(f)}^r \le \int_{\Omega} E_i\paren{\abs{f}^r} = \int_{\Omega} \abs{f}^r$,
       and (a) follows.
       The fact that $E_{i,r}$ maps $\L{r}(\Omega)$ into $\tilde L{}_i^r[0,1]=\Pi_{i,r}\paren{\L{r}[0,1]}$
       follows from the choice of the $\sigma$-algebra ${\cal B}_i$.
       For $f\in\tilde L{}_i^r[0,1]=\Pi_{i,r}\paren{\L{r}[0,1]}$, $E_{i,r}(f)=f$, and (b) follows.
       Since $\int_{\Omega} E_i(f) = \int_{\Omega} f$, (c) follows.
       Part (d) follows from the choice of the $\sigma$-algebra ${\cal B}_i$.
       Noting that $\int_{\Omega} f \cdot E_{i,r}^*(g) = \int_{\Omega} E_{i,r}(f) \cdot g =
                    \int_{\Omega} E_i(f) \cdot g       = \int_{\Omega} E_i\paren{E_i(f) \cdot g} =
                    \int_{\Omega} E_i(f) \cdot E_i(g)  = \int_{\Omega} E_i\paren{E_i(g) \cdot f} =
                    \int_{\Omega} E_i(g) \cdot f       = \int_{\Omega} E_{i,s}(g) \cdot f$ 
       for $g\in\L{s}(\Omega)$, \hfil\break 
       (e) follows. 
       Part (g) is clear.

       Now $E_{i,2}:\L{2}(\Omega)\to\tilde L{}_i^2[0,1]\subset\L{2}(\Omega)$
       maps $\L{2}(\Omega)$ onto $\tilde L{}_i^2[0,1]=\Pi_{i,2}\paren{\L{2}[0,1]}$ by parts (a) and (b). 
       Let $f\in\L{2}(\Omega)$.  
       Then $\int_B \paren{f-E_i(f)} = 0$ for all $B\in{\cal B}_i$, 
       and $\int_{\Omega} \paren{f-E_i(f)} \cdot g = 0$ for all $g\in\tilde L{}_i^2[0,1]$. 
       Hence $f-E_i(f) \in \paren{\tilde L{}_i^2[0,1]}^{\perp}$, and (f) follows.
       \qquad\QED
   
 For $r\in\set{2,p}$, let $S_{i,r}=\tilde Q_{i,r} E_{i,r}$, where $\tilde Q_{i,r}$ and $E_{i,r}$ are as above.

\xproclaim {Lemma 4.6}. Let $p$, $r$, $P_i$, $\tilde Y_{i,p}$, and $S_{i,r}$ be as above.
                        Let $f\in\L{r}(\Omega)$ and \hfil\break 
                        $g\in\L{q}(\Omega)$,
                        where $q$ is the conjugate index of $p$. Then
             \item{(a)} $S_{i,p}:\L{p}(\Omega)\to\tilde Y_{i,p}\subset\L{p}(\Omega)$ maps $\L{p}(\Omega)$ onto $\tilde Y_{i,p}$,
             \item{(b)} $S_{i,2}$ is the orthogonal projection of $\L{2}(\Omega)$
                        onto $\overline{\tilde Y_{i,p}}\subset\L{2}(\Omega)$,
             \item{(c)} $S_{i,p}=S_{i,2}|_{\L{p}(\Omega)}$,
             \item{(d)} $\norm{S_{i,p}}\le\norm{P_i}$,   
             \item{(e)} $S_{i,p}(1)=0$, 
             \item{(f)} $\int S_{i,r}(f) = 0$, 
             \item{(g)} $\int S^*_{i,p}(g)  = 0$, 
             \item{(h)} $\set{S_{i,r}(f)}_{i=1}^{\infty}$ is independent, and 
             \item{(i)} if $\set{g_i}_{i=1}^{\infty}$ is a sequence in $\L{q}(\Omega)$, 
                        then $\set{S^*_{i,p}(g_i)}_{i=1}^{\infty}$ is independent.  

\proof Part (a) is clear.
       Since $S_{i,2}=\tilde Q_{i,2} E_{i,2}$ is the composition of orthogonal projections, 
       where $\L{2}(\Omega) \Tra{E_{i,2}} \tilde L{}_i^2[0,1]$ surjectively and 
       $\tilde L{}_i^2[0,1] \Tra{\tilde Q_{i,2}} \overline{\tilde Y_{i,p}}$ surjectively, (b) follows.
       Part (c) is clear.
       Noting that $\norm{S_{i,p}}\le\bigl\|\tilde Q_{i,p}\bigr\|\norm{E_{i,p}}=
                    \bigl\|\tilde Q_{i,p}\bigr\|=\norm{Q_i}=\norm{P_i}$, (d) follows.
       Since $E_{i,p}(1)=1$ and $\tilde Q_{i,p}(1)=0$, (e) follows.
       Since $\tilde Y_{i,p}\subset\tilde L{}_0^p[0,k_i]$
       and $\overline{\tilde Y_{i,p}}\subset\tilde L{}_0^2[0,k_i]$, (f) follows.
       Noting that $\int S^*_{i,p}(g) = \int g \cdot S_{i,p}(1) = \int g \cdot 0 = 0$, \hfil\break 
       (g) follows. 
       For reference, $S_{i,r} = \tilde Q_{i,r} E_{i,r}$ and $S^*_{i,p} = E^*_{i,p} \tilde Q{}^*_{i,p}$.
       Part (h) follows from an analogous property of $E_{i,r}$ which $\tilde Q_{i,r}$ preserves.
       Recalling that $E_{i,q}$ has an analogous property and $E^*_{i,p}=E_{i,q}$, (i) follows.
       \qquad\QED

 For $r\in\set{2,p}$, let $S_r=\sum_{i=1}^{\infty} S_{i,r}$. 
 We show below that the formal series defines a bounded linear operator on $\L{r}(\Omega)$. 

\xproclaim {Lemma 4.7}. Let $p$, $\tilde Y_{i,p}$, and $S_2$ be as above.
                        Then $S_2$ is the orthogonal \hfil\break 
                        projection of $\L{2}(\Omega)$ onto 
                        $\span{\overline{\tilde Y_{i,p}}\subset\L{2}(\Omega):i\in\NN}_{\L{2}(\Omega)}$. 

\proof For $f\in\L{2}(\Omega)$, $S_2(f)=\sum_{i=1}^{\infty}S_{i,2}(f)$, where \hfil\break 
       $S_{i,2}(f) \in \overline{\tilde Y_{i,p}} \subset \tilde L{}_0^2[0,k_i] \subset \L{2}(\Omega)$, 
       $S_{i,2}(f)$ is the orthogonal projection of $f$ onto the span of $S_{i,2}(f)$ in $\L{2}(\Omega)$,
       and $\set{S_{i,2}(f)}_{i=1}^{\infty}$ is an orthogonal sequence of random variables.   
       Hence $S_2:\L{2}(\Omega)\to\span{\overline{\tilde Y_{i,p}}\subset\L{2}(\Omega):i\in\NN}_{\L{2}(\Omega)}$
       is the orthogonal projection of $\L{2}(\Omega)$ onto
       $\span{\overline{\tilde Y_{i,p}}\subset\L{2}(\Omega):i\in\NN}_{\L{2}(\Omega)}$.
       \qquad\QED

\xproclaim {Theorem 4.8}. Let $2<p<\infty$ and let $w=\set{w_i}$ be a sequence of scalars from $(0,1]$.  
                          Let $\set{X_i}$ be a sequence of closed subspaces of $\Lz{p}[0,1]$ 
                          satisfying the hypotheses (a) and (b) in the definition of $\sumsum{X_i}{I,w}$.
                          Then $\sumsum{X_i}{I,w}$ is a complemented subspace of $\L{p}(\Omega)$ 
                          via the projection $S_p$.

\proof Let $f\in\L{p}(\Omega)$.
       Then $\set{S_{i,p}(f)}_{i=1}^{\infty}$ is a sequence of independent mean zero random variables in $\L{p}(\Omega)$. 
       Hence (essentially) by Theorem 2.2 [Rosenthal's \hfil\break 
       inequality],  
       $$\eqalign{\norm{S_p(f)}_{\L{p}(\Omega)} &= \norm{\tsuml_{i=1}^{\infty} S_{i,p}(f)}_{\L{p}(\Omega)}
         \cr &\within{K_p}{2} \max\set{\paren{\tsuml_{i=1}^{\infty} \norm{S_{i,p}(f)}_{\L{p}(\Omega)}^p }^{{1\over p}},
                                       \paren{\tsuml_{i=1}^{\infty} \norm{S_{i,p}(f)}_{\L{2}(\Omega)}^2 }^{{1\over 2}} } .}$$ 
       By the orthogonality of $\set{S_{i,p}(f)}_{i=1}^{\infty}$ 
       and the fact that $S_p=S_2|_{\L{p}(\Omega)}$ where $S_2$ is \hfil\break 
       orthogonal projection,
       $$\paren{\tsuml_{i=1}^{\infty} \norm{S_{i,p}(f)}_{\L{2}(\Omega)}^2 }^{{1\over 2}} 
                =\norm{\tsuml_{i=1}^{\infty} S_{i,p}(f)}_{\L{2}(\Omega)}
                =\norm{S_p(f)}_{\L{2}(\Omega)}
              \le\norm{f}_{\L{2}(\Omega)}
              \le\norm{f}_{\L{p}(\Omega)}.$$
       Let
   $G=\set{\set{g_i}_{i=1}^{\infty}:g_i\in\L{q}(\Omega),\paren{\sum_{i=1}^{\infty}\norm{g_i}_{\L{q}(\Omega)}^q}^{{1\over q}}\le1}$,
       where $q$ is the conjugate index of $p$. 
       Then for $g_i\in\L{q}(\Omega)$,
       $\set{S_{i,p}^*\paren{g_i}}_{i=1}^{\infty}$ is a sequence of independent mean zero random variables in $\L{q}(\Omega)$. 
       Hence by H\"older's inequality and (essentially) Lemma 2.4, 
       $$\eqalign{\paren{\tsuml_{i=1}^{\infty} \norm{S_{i,p}(f)}_{\L{p}(\Omega)}^p }^{{1\over p}} 
                 &=\sup_{\set{g_i}\in G} \abs{\tsuml_{i=1}^{\infty} \vector{S_{i,p}(f),g_i} }
             \cr &=\sup_{\set{g_i}\in G} \abs{\vector{f, \tsuml_{i=1}^{\infty} S_{i,p}^*(g_i)} }
         \cr &\le  \sup_{\set{g_i}\in G} \norm{\tsuml_{i=1}^{\infty} S_{i,p}^*(g_i)}_{\L{q}(\Omega)} \norm{f}_{\L{p}(\Omega)}
         \cr &\le 2\sup_{\set{g_i}\in G} \paren{\tsuml_{i=1}^{\infty} \norm{S_{i,p}^*(g_i)}_{\L{q}(\Omega)}^q }^{{1\over q}}
                         \norm{f}_{\L{p}(\Omega)}
         \cr &\le 2\sup_{i\in\NN} \norm{S_{i,p}^*}
                   \sup_{\set{g_i}\in G} \paren{\tsuml_{i=1}^{\infty} \norm{g_i}_{\L{q}(\Omega)}^q }^{{1\over q}}
                         \norm{f}_{\L{p}(\Omega)}
         \cr &\le 2\sup_{i\in\NN} \norm{P_i} \norm{f}_{\L{p}(\Omega)} .}$$
       It now follows that
       $\norm{S_p(f)}_{\L{p}(\Omega)} \le K_p \max\set{2\sup_{i\in\NN}\norm{P_i},1}\norm{f}_{\L{p}(\Omega)}$.
       Hence \hfil\break 
       $S_p:\L{p}(\Omega)\to\span{\tilde Y_{i,p}:i\in\NN}_{\L{p}(\Omega)}$
       maps $\L{p}(\Omega)$ onto $\span{\tilde Y_{i,p}:i\in\NN}_{\L{p}(\Omega)}$ with \hfil\break 
       $\norm{S_p}\le K_p \max\set{2\sup_{i\in\NN}\norm{P_i},1}$, and 
       $\sumsum{X_i}{I,w}=\span{\tilde Y_{i,p}:i\in\NN}_{\L{p}(\Omega)}$ \hfil\break 
       is complemented in $\L{p}(\Omega)$. 
       \qquad\QED

\preheadspace
\firsthead{Independent Sums with Basis}
\postheadspace

 Now suppose in addition to the hypotheses (a) and (b) in the definition of \hfil\break 
 $\sumsum{X_i}{I,w}$,
 the sequence $\set{X_i}$ of closed subspaces of $\Lz{p}[0,1]$ satisfies \hfil\break 
 (c) for each $i\in\NN$, $X_i$ has an unconditional orthogonal basis $\set{x_{i,n}}_{n=1}^{\infty}$. \hfil\break 
 Then of course $X_i=\span{x_{i,n}:n\in\NN}_{\L{p}[0,1]}$.

 Letting $Y_i=T_i\paren{X_i}$ as before,
 and letting $y_{i,n}=T_i\paren{x_{i,n}}$, we have \hfil\break 
 $Y_i=\span{y_{i,n}:n\in\NN}_{\L{p}[0,1]}$,
 and $\set{y_{i,n}}_{n=1}^{\infty}$ is an unconditional orthogonal basis for $Y_i$ 
 isometrically equivalent to $\set{x_{i,n}}_{n=1}^{\infty}$. 

 Letting $\tilde Y_i=\set{\tilde y_i=y_i\circ\pi_i:y_i\in Y_i}$ as before, 
 and letting $\tilde y_{i,n}=y_{i,n}\circ\pi_i$, \hfil\break 
 we have $\tilde Y_i=\span{\tilde y_{i,n}:n\in\NN}_{\L{p}(\Omega)}$,
 and $\set{\tilde y_{i,n}}_{n=1}^{\infty}$ is an unconditional orthogonal basis for $\tilde Y_i$ 
 isometrically equivalent to $\set{y_{i,n}}_{n=1}^{\infty}$ and $\set{x_{i,n}}_{n=1}^{\infty}$. 

 In this context, $\sumsum{X_i}{I,w}=\span{\tilde y_{i,n}:i,n\in\NN}_{\L{p}(\Omega)}$, 
 and $\set{\tilde y_{i,n}}_{i,n\in\NN}$ is an \hfil\break 
 unconditional orthogonal basis for $\sumsum{X_i}{I,w}$.

\remark Noting that $y_{i,n}=T_i\paren{x_{i,n}}$ and $k_i^{{p-2\over2p}}=w_i$, by part (a) of Lemma 4.2 we have
        $\norm{\tilde y_{i,n}}_{\L{2}(\Omega)}=\norm{y_{i,n}}_2=w_i\norm{x_{i,n}}_2$. 

\xproclaim {Proposition 4.9}. Let $2<p<\infty$ and let $w=\set{w_i}$ be a sequence of scalars from $(0,1]$.
                              Let $\set{X_i}$ be a sequence of closed subspaces of $\Lz{p}[0,1]$ 
                              such that each $X_i$ has an unconditional orthogonal basis $\set{x_{i,n}}_{n=1}^{\infty}$.
                              Let $\tilde y_{i,n}=\paren{T_i\paren{x_{i,n}}}\circ\pi_i\in\L{p}(\Omega)$, 
                              where $T_i$ and $\pi_i$ are as in the definition of $\sumsum{X_i}{I,w}$. 
                              Then for $K_p$ as in Theorem 2.2 and for scalars $a_{i,n}$, 
                      $$\norm{\tsuml_i \tsuml_n a_{i,n} \tilde y_{i,n} }_{\L{p}(\Omega)} \within{K_p}{2}
                        \max\set{\paren{\tsuml_i \norm{\tsuml_n a_{i,n} x_{i,n} }_p^p }^{{1\over p}}, 
                                 \paren{\tsuml_i w_i^2 \tsuml_n \abs{a_{i,n}}^2 \norm{x_{i,n}}_2^2 }^{{1\over2}} }.$$

\proof Let $z_i = \sum_n a_{i,n} \tilde y_{i,n}$.  
       Then $\set{z_i}$ is a sequence of independent mean zero random variables in $\L{p}(\Omega)$.  
       Hence (essentially) by Corollary 2.3 [Rosenthal's \hfil\break 
       inequality], 
       $$\norm{\tsuml_i z_i}_{\L{p}(\Omega)} \within{K_p}{2}
         \max\set{\paren{\tsuml_i \norm{z_i}_{\L{p}(\Omega)}^p }^{{1\over p}},
                  \paren{\tsuml_i \norm{z_i}_{\L{2}(\Omega)}^2 }^{{1\over2}} }.$$
       Note that
       $\norm{z_i}_{\L{p}(\Omega)}^p=\norm{\sum_n a_{i,n}\tilde y_{i,n}}_{\L{p}(\Omega)}^p=\norm{\sum_n a_{i,n}x_{i,n}}_p^p$.   
       Moreover, by the \hfil\break 
       orthogonality of $\set{\tilde y_{i,n}}_{n=1}^{\infty}$ and by the remark above, \hfil\break 
       $\norm{z_i}_{\L{2}(\Omega)}^2 = \norm{\sum_n a_{i,n} \tilde y_{i,n} }_{\L{2}(\Omega)}^2 
        = \sum_n \abs{a_{i,n}}^2 \norm{\tilde y_{i,n} }_{\L{2}(\Omega)}^2 
        = w_i^2 \sum_n \abs{a_{i,n}}^2 \norm{x_{i,n}}_2^2$. \hfil\break 
       The result now follows from the displayed inequality. 
       \qquad\QED 

\xproclaim {Corollary 4.10}. Let $2<p<\infty$ and let $w=\set{w_i}$ be a sequence of scalars from $(0,1]$.
                             Let $\set{X_i}$ be a sequence of closed subspaces of $\Lz{p}[0,1]$ 
                             satisfying the hypotheses (a) and (b) in the definition of $\sumsum{X_i}{I,w}$ 
                             such that each $X_i$ has an unconditional orthogonal basis $\set{x_{i,n}}_{n=1}^{\infty}$.
                             Suppose $\sum w_i^{{2p\over p-2}} < \infty$.
                             Then $\sumsum{X_i}{I,w} \sim \sumsum{X_i}{\l{p}}$.

\proof Let $\tilde y_{i,n}$ be as in Proposition 4.9.
       Let $K = \paren{\sum w_i^{{2p\over p-2}} }^{{p-2\over2p}}$. 
       By H\"older's inequality with conjugate indices $p'={p\over2}$ and $q'={p\over p-2}$, 
       and the orthogonality of $\set{x_{i,n}}_{n=1}^{\infty}$, for scalars $a_{i,n}$ we have 
       $$\eqalign{\paren{\tsuml_i w_i^2 \paren{\tsuml_n \abs{a_{i,n}}^2 \norm{x_{i,n}}_2^2 } }^{{1\over2}} 
         &\le     \paren{\paren{\tsuml_i w_i^{2{p\over p-2}} }^{{p-2\over p}} 
                  \paren{\tsuml_i \paren{\tsuml_n \abs{a_{i,n}}^2 \norm{x_{i,n}}_2^2 }^{{p\over2}} }^{{2\over p}} }^{{1\over2}} 
         \cr &=   \paren{\tsuml_i w_i^{{2p\over p-2}} }^{{p-2\over2p}} 
                  \paren{\tsuml_i \paren{\tsuml_n \norm{a_{i,n} x_{i,n}}_2^2 }^{{1\over2}p} }^{{1\over p}} 
         \cr &=   K \paren{\tsuml_i \norm{\tsuml_n a_{i,n} x_{i,n} }_2^p }^{{1\over p}} 
         \cr &\le K \paren{\tsuml_i \norm{\tsuml_n a_{i,n} x_{i,n} }_p^p }^{{1\over p}} .}$$
       Hence by Proposition 4.9 and the above bound, for $\tilde K = \max\set{1,K}$ we have 
       $$\eqalign{\norm{\tsuml_i \tsuml_n a_{i,n} \tilde y_{i,n} }_{\L{p}(\Omega)}
                  &\within{K_p}{2}
                  \max\set{\paren{\tsuml_i \norm{\tsuml_n a_{i,n} x_{i,n} }_p^p }^{{1\over p}}, 
                           \paren{\tsuml_i w_i^2 \tsuml_n \abs{a_{i,n}}^2 \norm{x_{i,n}}_2^2 }^{{1\over2}} }
              \cr &\within{\tilde K}{1} \paren{\tsuml_i \norm{\tsuml_n a_{i,n} x_{i,n} }_p^p }^{{1\over p}} .}$$
       It follows that $\sumsum{X_i}{I,w} \sim \sumsum{X_i}{\l{p}}$.
       \qquad\QED 

\xproclaim {Example 4.11}. Let $2<p<\infty$ and let $w=\set{w_i}$ be a sequence of scalars from $(0,1]$
                           such that $\sum w_i^{{2p\over p-2}}<\infty$.  
                           Then $\sumsum{\l{2}}{\l{p}}$, $\sumsum{X_p}{\l{p}}$, $B_p$,
                                $X_p\oplus\sumsum{\l{2}}{\l{p}}$, and $X_p \oplus B_p$ 
                                can be realized as $\sumsum{X_i}{I,w}$ for appropriately chosen $X_i$.

\proof Let $\set{x_n}$ be the sequence of Rademacher functions and let \hfil\break 
       $X=\span{x_n}_{\L{p}}\sim\l{2}$.
       Then $\sumsum{X}{I,w}\sim\sumsum{\l{2}}{\l{p}}$.  

       Let $\set{x_n}$ be a sequence of independent mean zero random variables in $\L{p}$  
       such that $v = \set{v_n} = \set{\norm{x_n}_2\big/\norm{x_n}_p}$ satisfies condition $(*)$ of Proposition 2.1, and   
       let $X=\span{x_n}_{\L{p}}\sim X_p$. 
       Then $\sumsum{X}{I,w}\sim\sumsum{X_p}{\l{p}}$.  

       For each $i\in\NN$, let $\set{x_{i,n}}_{n=1}^{\infty}$ be a sequence of
       independent mean zero random variables in $\L{p}$ such that
       $v^{(i)}\!=\!\set{v_{i,n}}_{n=1}^{\infty}\!=\!\set{\norm{x_{i,n}}_2\big/\norm{x_{i,n}}_p}_{n=1}^{\infty}$
       satisfies $v_{i,n}^{{2p\over p-2}} = {1\over i}$ for each $n\in\NN$.
       Let $X_i=\span{x_{i,n}:n\in\NN}_{\L{p}}\sim X_{p,v^{(i)}}$. 
       Then \hfil\break 
       $\sumsum{X_i}{I,w}\sim\sumsum{X_{p,v^{(i)}}}{\l{p}}\sim B_p$.  

       Let $\set{x_{1,n}}_{n=1}^{\infty}$ be a sequence of
       independent mean zero random variables in $\L{p}$ such that
       $v^{(1)}\!=\!\set{v_{1,n}}_{n=1}^{\infty}\!=\!\set{\norm{x_{1,n}}_2\big/\norm{x_{1,n}}_p}_{n=1}^{\infty}$
       satisfies condition $(*)$ of \hfil\break 
       Proposition 2.1, and let $X_1=\span{x_{1,n}:n\in\NN}_{\L{p}}\sim X_p$. 
       For each $i\in{\NN}\setminus\set{1}$,
       let $\set{x_{i,n}}_{n=1}^{\infty}$ be the sequence of Rademacher functions and let 
       $X_i=\span{x_{i,n}:n\in\NN}_{\L{p}}\sim\l{2}$.
       Then $\sumsum{X_i}{I,w} \sim \sumsum{X_i}{\l{p}}
                               \sim \paren{X_p\oplus\sum^{\oplus} \l{2} }_{\l{p}}
                               \sim X_p \oplus \sumsum{\l{2}}{\l{p}}$.

       Let $\set{x_{1,n}}_{n=1}^{\infty}$ be a sequence of
       independent mean zero random variables in $\L{p}$ such that
       $v^{(1)}\!=\!\set{v_{1,n}}_{n=1}^{\infty}\!=\!\set{\norm{x_{1,n}}_2\big/\norm{x_{1,n}}_p}_{n=1}^{\infty}$
       satisfies condition $(*)$ of \hfil\break 
       Proposition 2.1, and let $X_1=\span{x_{1,n}:n\in\NN}_{\L{p}}\sim X_p$. 
       For each $i\in{\NN}\setminus\set{1}$, let $\set{x_{i,n}}_{n=1}^{\infty}$ be a sequence of
       independent mean zero random variables in $\L{p}$ such that
       $v^{(i)}\!=\!\set{v_{i,n}}_{n=1}^{\infty}\!=\!\set{\norm{x_{i,n}}_2\big/\norm{x_{i,n}}_p}_{n=1}^{\infty}$
       satisfies $v_{i,n}^{{2p\over p-2}} = {1\over i}$ for each $n\in\NN$,
       and let $X_i=\span{x_{i,n}:n\in\NN}_{\L{p}}\sim X_{p,v^{(i)}}$.
       Then $\sumsum{X_i}{I,w} \sim \sumsum{X_i}{\l{p}} \sim \hfil\break 
                                    \paren{X_p\oplus\sum_{i\ge2}^{\oplus} X_{p,v^{(i)}} }_{\l{p}}
                               \sim X_p\oplus\paren{\sum_{i\ge2}^{\oplus} X_{p,v^{(i)}} }_{\l{p}}
                               \sim X_p \oplus B_p$.
       \qquad\QED

\preheadspace
\firsthead{The Independent Sum $\sumsum{X}{I}$}
\postheadspace

 Let $2<p<\infty$.
 Suppose $X$ is a closed subspace of $\Lz{p}[0,1]$ satisfying \hfil\break 
 (a${}'$) the orthogonal projection of $\L{2}[0,1]$ onto $\overline X\subset\L{2}[0,1]$, 
          when restricted to $\L{p}[0,1]$,
          yields a bounded projection $P:\L{p}[0,1]\to X\subset\L{p}[0,1]$ onto $X$, and \hfil\break 
 (c${}'$) $X$ has an unconditional orthogonal normalized basis $\set{x_n}$. \hfil\break 
 We adopt notation as before, with $X$ replacing $X_i$ and $x_n$ replacing $x_{i,n}$. 
 In particular, $\tilde y_{i,n} = \paren{T_i\paren{x_n}}\circ\pi_i\in\L{p}(\Omega)$, 
 where $T_i$ and $\pi_i$ are as in the definition of $\sumsum{X_i}{I,w}$. 

 For $2<p<\infty$, 
 we will show that for a fixed closed subspace $X$ of $\Lz{p}[0,1]$ \hfil\break 
 satisfying the hypotheses (a${}'$) and (c${}'$) above,  
 all spaces $\sumsum{X}{I,w}$ for sequences \hfil\break 
 $w=\set{w_i}$ from $(0,1]$ satisfying condition $(*)$ of Proposition 2.1 are mutually \hfil\break 
 isomorphic. 
 The following results follow the pattern of Propositions 2.7, 2.9, 2.10, 2.11, and Theorem 2.12, 
 where it is shown that the isomorphism type of $\Xpw$ does not \hfil\break 
 depend on $w$ as long as $w$ satisfies condition $(*)$.  

\xproclaim {Proposition 4.12}. Let $2<p<\infty$ and let $w=\set{w_i}$ be a sequence of scalars from $(0,1]$.
                               Let $X$ be a closed subspace of $\Lz{p}[0,1]$ satisfying the hypotheses (a${}'$) and (c${}'$) above.
                               Suppose $\set{E_j}$ is a sequence of disjoint nonempty finite subsets of $\NN$ such that
                               $\sum_{i\in E_j} w_i^{{2p\over p-2}} \le 1$ for each $j\in\NN$.  
                               Let $z_{j,n} = \sum_{i\in E_j} w_i^{{2\over p-2}} \tilde y_{i,n}$ and
                               let $\tilde z_{j,n}$ be the normalization of $z_{j,n}$ in $\L{p}(\Omega)$. 
                               Let $v_j = \paren{\sum_{i\in E_j} w_i^{{2p\over p-2}} }^{{p-2\over2p}}$ and $v=\set{v_j}$.
                               Then
       \item{(a)} $\set{\tilde z_{j,n}}$ is an unconditional basis for $\span{\tilde z_{j,n}:j,n\in\NN}_{\sumsum{X}{I,w}}$ 
                  which is equivalent to the standard basis of $\sumsum{X}{I,v}$, and   
       \item{(b)} there is a projection $P:\sumsum{X}{I,w}\to\span{\tilde z_{j,n}:j,n\in\NN}_{\sumsum{X}{I,w}}$.

\proof First we establish some notation.
       Let $Y_{p,\set{x_n}}$ be the Banach space of all sums of the form $y=\sum_i\sum_n a_{i,n} \tilde y_{i,n}$ 
       (for scalars $a_{i,n}$) such that \hfil\break 
       $\norm{y}_{Y_{p,\set{x_n}} } = \paren{\sum_i \norm{\sum_n a_{i,n} \tilde y_{i,n} }_{\L{p}(\Omega)}^p }^{{1\over p}}
                                    = \paren{\sum_i \norm{\sum_n a_{i,n} x_n}_p^p }^{{1\over p}} < \infty$.
       Let $Y_{2,w,\set{x_n}}$ be the Hilbert space of all sums of the form $y=\sum_i\sum_n a_{i,n} \tilde y_{i,n}$ 
       (for scalars $a_{i,n}$) such that 
       $\norm{y}_{Y_{2,w,\set{x_n}} } = \paren{\sum_i \norm{\sum_n a_{i,n} \tilde y_{i,n} }_{\L{2}(\Omega)}^2 }^{{1\over2}}
                                      = \paren{\sum_i w_i^2 \sum_n \abs{a_{i,n}}^2 \norm{x_n}_2^2 }^{{1\over2}} < \infty$,
       where the inner product in $Y_{2,w,\set{x_n}}$ is defined by \hfil\break 
       $\vector{y_a,y_b} = \sum_i \int \paren{\sum_n a_{i,n} \tilde y_{i,n} }
                                       \overline{\paren{\sum_n b_{i,n} \tilde y_{i,n} } }
                         = \sum_i w_i^2 \sum_n a_{i,n} \bar b_{i,n} \norm{x_n}_2^2$ \hfil\break 
       (where $y_a=\sum_i\sum_n a_{i,n} \tilde y_{i,n}$, $y_b=\sum_i\sum_n b_{i,n} \tilde y_{i,n}$, 
       and bar is complex conjugation). 

       Let $\xnorm{\phantom0}$ be the norm on $\sumsum{X}{I,w}$ defined by \hfil\break 
       $\xnorm{y} = \max\set{\norm{y}_{Y_{p,\set{x_n}}}, \norm{y}_{Y_{2,w,\set{x_n}}} }$.
       By Proposition 4.9, $\xnorm{\phantom0}$ is equivalent to the standard norm on $\sumsum{X}{I,w}$. 
       Without loss of generality,
       we will proceed in the \hfil\break 
       context of $\sumsum{X}{I,w}$ endowed with the norm $\xnorm{\phantom0}$. 

       We now find the normalizing factor for $z_{j,n}$. 
       Let $\sigma_j = \sum_{i\in E_j} w_i^{{2p\over p-2}}$. 
       Noting that $2+{4\over p-2}={2p\over p-2}$, 
       $1 = \norm{x_n}_p \ge \norm{x_n}_2$, and $\sigma_j^{{1\over p}} \ge \sigma_j^{{1\over2}}$, we have 
       $$\eqalign{\xnorm{z_{j,n}} =
                  \xnorm{\tsuml_{i\in E_j} w_i^{{2\over p-2}} \tilde y_{i,n}}
                  &=
                  \max\set{\paren{\tsuml_{i\in E_j} \norm{w_i^{{2\over p-2}} x_n }_p^p }^{{1\over p}}, 
                           \paren{\tsuml_{i\in E_j} w_i^2 w_i^{{4\over p-2}} \norm{x_n}_2^2 }^{{1\over2}} }
              \cr &=
                  \max\set{\paren{\tsuml_{i\in E_j} w_i^{{2p\over p-2}} \norm{x_n}_p^p }^{{1\over p}}, 
                           \paren{\tsuml_{i\in E_j} w_i^{{2p\over p-2}} \norm{x_n}_2^2 }^{{1\over2}} }
              \cr &=
                  \max\set{\sigma_j^{{1\over p}} \norm{x_n}_p, 
                           \sigma_j^{{1\over2}} \norm{x_n}_2 }
              \cr &= \sigma_j^{{1\over p}} .}$$
       Hence $\tilde z_{j,n} = \sigma_j^{-{1\over p}} z_{j,n} 
                             = \sigma_j^{-{1\over p}} \sum_{i\in E_j} w_i^{{2\over p-2}} \tilde y_{i,n}$. 

\item{(a)} 
       The unconditionality of $\set{\tilde z_{j,n}}$ follows from
       the unconditionality of $\set{\tilde y_{i,n}}$ in \hfil\break 
       $\sumsum{X}{I,w}$.
       We now examine the equivalence of the bases.  
       For scalars $a_{j,n}$, we have  
       $$\eqalignno{\norm{\tsuml_j\tsuml_n a_{j,n} \tilde z_{j,n}}^p_{Y_{p,\set{x_n}}}
                 &= \norm{\tsuml_j\tsuml_n a_{j,n} \sigma_j^{-{1\over p}} \tsuml_{i\in E_j} w_i^{{2\over p-2}} \tilde y_{i,n} }
                         ^p_{Y_{p,\set{x_n}}}
             \cr &= \norm{\tsuml_j\tsuml_{i\in E_j}\tsuml_n \sigma_j^{-{1\over p}} w_i^{{2\over p-2}} a_{j,n} \tilde y_{i,n} }
                         ^p_{Y_{p,\set{x_n}}}
             \cr &= \tsuml_j\tsuml_{i\in E_j} \norm{\tsuml_n \sigma_j^{-{1\over p}} w_i^{{2\over p-2}} a_{j,n} x_n }_p^p
             \cr &= \tsuml_j \sigma_j^{-1} \tsuml_{i\in E_j} w_i^{{2p\over p-2}} \norm{\tsuml_n a_{j,n} x_n }_p^p
             \cr &= \tsuml_j \norm{\tsuml_n a_{j,n} x_n }_p^p &\TAG{4.4}{4.4}
             \cr &= \norm{\tsuml_j\tsuml_n a_{j,n} \tilde y{}_{j,n}^{(v)} }^p_{Y_{p,\set{x_n}}},}$$
       and noting that $2+{4\over p-2}={2p\over p-2}$ and $1-{2\over p}={p-2\over2p}2$, 
       $$\eqalignno{\norm{\tsuml_j\tsuml_n a_{j,n} \tilde z_{j,n}}^2_{Y_{2,w,\set{x_n}}}
                 &= \norm{\tsuml_j\tsuml_n a_{j,n} \sigma_j^{-{1\over p}} \tsuml_{i\in E_j} w_i^{{2\over p-2}} \tilde y_{i,n} }
                         ^2_{Y_{2,w,\set{x_n}}}
             \cr &= \norm{\tsuml_j\tsuml_{i\in E_j}\tsuml_n \sigma_j^{-{1\over p}} w_i^{{2\over p-2}} a_{j,n} \tilde y_{i,n} }
                         ^2_{Y_{2,w,\set{x_n}}}
       \cr &= \tsuml_j\tsuml_{i\in E_j} w_i^2 \tsuml_n \abs{\sigma_j^{-{1\over p}} w_i^{{2\over p-2}} a_{j,n} }^2 \norm{x_n}_2^2
       \cr &= \tsuml_j \sigma_j^{-{2\over p}} \tsuml_{i\in E_j} w_i^{{2p\over p-2}} \tsuml_n \abs{a_{j,n} }^2 \norm{x_n}_2^2
       \cr &= \tsuml_j \paren{\sigma_j^{{p-2\over2p}}}^2 \tsuml_n \abs{a_{j,n} }^2 \norm{x_n}_2^2
       \cr &= \tsuml_j v_j^2 \tsuml_n \abs{a_{j,n} }^2 \norm{x_n}_2^2 &\TAG{4.5}{4.5} 
       \cr &= \norm{\tsuml_j\tsuml_n a_{j,n} \tilde y{}_{j,n}^{(v)} }^2_{Y_{2,v,\set{x_n}}},}$$ 
       where $\tilde y{}_{j,n}^{(v)}$ is analogous to $\tilde y_{j,n}$ with $v$ replacing $w$. 
       Hence 
       $$\eqalign{\xnorm{\tsuml_j\tsuml_n a_{j,n} \tilde z_{i,n} }
               &= \max\set{\norm{\tsuml_j\tsuml_n a_{j,n} \tilde z_{j,n}}_{Y_{p,\set{x_n}}},
                           \norm{\tsuml_j\tsuml_n a_{j,n} \tilde z_{j,n}}_{Y_{2,w,\set{x_n}}} }
           \cr &= \max\set{\norm{\tsuml_j\tsuml_n a_{j,n} \tilde y{}_{j,n}^{(v)} }_{Y_{p,\set{x_n}}},
                           \norm{\tsuml_j\tsuml_n a_{j,n} \tilde y{}_{j,n}^{(v)} }_{Y_{2,v,\set{x_n}}} }
           \cr &= \xnorm{\tsuml_j\tsuml_n a_{j,n} \tilde y{}_{j,n}^{(v)} }_v,}$$ 
       where $\xnorm{\phantom0}_v$ is analogous to $\xnorm{\phantom0}$ with $v$ replacing $w$. 
       Hence $\set{\tilde z_{j,n}}$ is equivalent to the standard basis $\set{\tilde y{}_{j,n}^{(v)} }$ of $\sumsum{X}{I,v}$.  

\item{(b)}
       Let $\pi:Y_{2,w,\set{x_n}}\to\span{z_{j,n}:j,n\in\NN}_{Y_{2,w,\set{x_n}}}$ be
       the orthogonal projection
       onto $\span{z_{j,n}:j,n\in\NN}_{Y_{2,w,\set{x_n}}}$ defined by 
       $$\pi(y)=\tsuml_j\tsuml_n \displaystyle{\vector{y,z_{j,n}}\over\vector{z_{j,n},z_{j,n}}} z_{j,n}.$$
       Let $y\in\sumsum{X}{I,w} = Y_{p,\set{x_n}} \cap Y_{2,w,\set{x_n}}$.
       Then $\norm{\pi(y)}_{Y_{2,w,\set{x_n}}}\le\norm{y}_{Y_{2,w,\set{x_n}}}$. 
       We will show that $\norm{\pi(y)}_{Y_{p,\set{x_n}}}\le\norm{y}_{Y_{p,\set{x_n}}}$ as well, whence 
       $$\eqalign{\xnorm{\pi(y)}&=\max\set{\norm{\pi(y)}_{Y_{p,\set{x_n}}}, \norm{\pi(y)}_{Y_{2,w,\set{x_n}}} }
                          \cr &\le\max\set{\norm{y}_{Y_{p,\set{x_n}}}, \norm{y}_{Y_{2,w,\set{x_n}}} } = \xnorm{y} .}$$ 
       Thus letting $P:\sumsum{X}{I,w}\to\span{z_{j,n}:j,n\in\NN}_{\sumsum{X}{I,w}}$
       be the restriction of $\pi$ to $\sumsum{X}{I,w}$, 
       $P$ will satisfy our requirements.

\item{}
       Fix $y = \sum_i\sum_n a_{i,n} \tilde y_{i,n} \in \sumsum{X}{I,w}$. 
       Let $\lambda_{j,n}=\vector{y,z_{j,n}}\big/\vector{z_{j,n},z_{j,n}}$, so that $\pi(y) = \sum_j\sum_n \lambda_{j,n} z_{j,n}$.
       Noting that $2+{2\over p-2}={2(p-1)\over p-2}$ and $2+{4\over p-2}={2p\over p-2}$, we have 
       $$\eqalign{\lambda_{j,n}&=\vector{y,z_{j,n}}\big/\vector{z_{j,n},z_{j,n}}
              \cr &=\vector{\tsuml_i\tsuml_n a_{i,n} \tilde y_{i,n}, \tsuml_{i\in E_j} w_i^{{2\over p-2}} \tilde y_{i,n} }
                    \Bigg/
                    \vector{\tsuml_{i\in E_j}w_i^{{2\over p-2}}\tilde y_{i,n},\tsuml_{i\in E_j}w_i^{{2\over p-2}}\tilde y_{i,n} } 
              \cr &=\paren{\tsuml_{i\in E_j} w_i^2 a_{i,n} w_i^{{2\over p-2}} \norm{x_n}_2^2 }
                    \Bigg/
                    \paren{\tsuml_{i\in E_j} w_i^2 w_i^{{4\over p-2}} \norm{x_n}_2^2 }
              \cr &=\paren{\tsuml_{i\in E_j} w_i^{{2(p-1)\over p-2}} a_{i,n} } 
                    \Bigg/
                    \paren{\tsuml_{i\in E_j} w_i^{{2p\over p-2}} }
              \cr &=\sigma_j^{-1} \tsuml_{i\in E_j} w_i^{{2(p-1)\over p-2}} a_{i,n} .}$$ 
       Thus we have
       $$\eqalign{\norm{\pi(y)}_{Y_{p,\set{x_n}}} &= \norm{\tsuml_j\tsuml_n \lambda_{j,n} z_{j,n} }_{Y_{p,\set{x_n}}}
                  \cr &= \norm{\tsuml_j\tsuml_n \lambda_{j,n}\tsuml_{i\in E_j} w_i^{{2\over p-2}} \tilde y_{i,n} }_{Y_{p,\set{x_n}}}
                  \cr &= \norm{\tsuml_j\tsuml_{i\in E_j}\tsuml_n \lambda_{j,n} w_i^{{2\over p-2}} \tilde y_{i,n} }_{Y_{p,\set{x_n}}}
                  \cr &= \paren{\tsuml_j\tsuml_{i\in E_j} \norm{\tsuml_n \lambda_{j,n} w_i^{{2\over p-2}} x_n}_p^p }^{{1\over p}}
                  \cr &= \paren{\tsuml_j\tsuml_{i\in E_j} w_i^{{2p\over p-2}} \norm{\tsuml_n \lambda_{j,n} x_n}_p^p }^{{1\over p}}
                  \cr &= \paren{\tsuml_j \sigma_j
                  \norm{\tsuml_n \sigma_j^{-1} \tsuml_{i\in E_j} w_i^{{2(p-1)\over p-2}} a_{i,n} x_n}_p^p }^{{1\over p}}
                  \cr &= \paren{\tsuml_j \sigma_j^{1-p}
                         \norm{\tsuml_n \tsuml_{i\in E_j} w_i^{{2(p-1)\over p-2}} a_{i,n} x_n}_p^p }^{{1\over p}} ,}$$
       where by H\"older's inequality, letting $q$ be the conjugate index of $p$ and noting that $(p-1)q=p$ and ${p\over q}=p-1$, 
       $$\eqalign{\norm{\tsuml_n \tsuml_{i\in E_j} w_i^{{2(p-1)\over p-2}} a_{i,n} x_n}_p^p
                      &= \int \abs{\tsuml_{i\in E_j} \paren{w_i^{{2(p-1)\over p-2}} } \paren{\tsuml_n a_{i,n} x_n} }^p
                \cr &\le \int \abs{\paren{\tsuml_{i\in E_j} \paren{w_i^{{2(p-1)\over p-2}} }^q }^{{1\over q}}
                                   \paren{\tsuml_{i\in E_j} \abs{\tsuml_n a_{i,n} x_n }^p }^{{1\over p}} }^p
                  \cr &= \paren{\tsuml_{i\in E_j} w_i^{{2p\over p-2}} }^{{p\over q}}
                                \tsuml_{i\in E_j} \displaystyle\int \abs{\tsuml_n a_{i,n} x_n }^p  
                  \cr &= \sigma_j^{p-1} \tsuml_{i\in E_j} \norm{\tsuml_n a_{i,n} x_n }_p^p ,}$$
       whence
       $$\norm{\pi(y)}_{Y_{p,\set{x_n}}} \le \paren{\tsuml_j \tsuml_{i\in E_j} \norm{\tsuml_n a_{i,n} x_n }_p^p }^{{1\over p}} 
                                         \le \paren{\tsuml_i \norm{\tsuml_n a_{i,n} x_n }_p^p }^{{1\over p}} 
                                           = \norm{y}_{Y_{p,\set{x_n}}}.$$

       \QED\medskip

\remark We have actually shown that for $\sumsum{X}{I,w}$ and $\sumsum{X}{I,v}$
        endowed with the norms $\xnorm{\pz}$ and $\xnorm{\pz}_v$, respectively,  
        $\set{\tilde z_{j,n}}$ is isometrically equivalent to the standard basis of $\sumsum{X}{I,v}$, and   
        there is a projection \hfil\break 
        $P:\sumsum{X}{I,w}\to\span{\tilde z_{j,n}:j,n\in\NN}_{\sumsum{X}{I,w}}$ with $\norm{P}=1$.

\xproclaim {Proposition 4.13}. Let $2<p<\infty$ and let $X$ be a closed subspace of $\Lz{p}[0,1]$ \hfil\break 
                               satisfying the hypotheses (a${}'$) and (c${}'$) above.
                               Let $w=\set{w_i}$ and $w'=\set{w'_i}$ be \hfil\break 
                               sequences of scalars from $(0,1]$ satisfying condition $(*)$ of Proposition 2.1.
                               Then \hfil\break 
                               $\sumsum{X}{I,w'}\cinjects\sumsum{X}{I,w}$. 

\proof By condition $(*)$, we may choose a sequence $\set{E_j}$ of disjoint nonempty finite subsets of $\NN$
       such that for each $j\in\NN$,
       $\paren{{w'_j\over2}}^{{2p\over p-2}} \le \sum_{i\in E_j} w_i^{{2p\over p-2}} \le \paren{w'_j}^{{2p\over p-2}}$.
       Then for $v_j = \paren{\sum_{i\in E_j} w_i^{{2p\over p-2}} }^{{p-2\over2p}}$, ${w'_j\over2}\le v_j\le w'_j$.
       Let $v=\set{v_j}$ and let $y\in\sumsum{X}{I,w'}$.
       Then ${1\over2}\norm{y}_{\sumsum{X}{I,w'}}\le\norm{y}_{\sumsum{X}{I,v}}\le\norm{y}_{\sumsum{X}{I,w'}}$.
       Hence \hfil\break 
       $\sumsum{X}{I,w'}\sim\sumsum{X}{I,v}$.
       However, $\sumsum{X}{I,v}\cinjects\sumsum{X}{I,w}$ by Proposition 4.12.
       It follows that $\sumsum{X}{I,w'}\cinjects\sumsum{X}{I,w}$.
       \qquad\QED

 Let $2<p<\infty$ and let $X$ be a closed subspace of $\Lz{p}[0,1]$ satisfying the \hfil\break 
 hypotheses (a${}'$) and (c${}'$) above.
 For each sequence $v=\set{v_i}$ from $(0,1]$,
 define spaces $Y_{p,\set{x_n}}$ and $Y_{2,v,\set{x_n}}$ as in the proof of Proposition 4.12. 
 For each $k\in\NN$, let \hfil\break 
 $v^{(k)}=\set{v_i^{(k)}}_{i=1}^{\infty}$ be a sequence from $(0,1]$, 
 and let $Y_k$ be a closed subspace of \hfil\break 
 $\sumsum{X}{I,v^{(k)}}$.
 Let $\paren{Y_1 \oplus Y_2 \oplus\cdots}_{p,2,\set{v^{(k)}}}$
 be the Banach space of all sequences $\set{y_k}$ with $y_k\in Y_k$ such that
 $\norm{\set{y_k}}=\max\set{\paren{\sum\norm{y_k}^p_{Y_{p,\set{x_n}}} }^{{1\over p}}\!,
                            \paren{\sum\norm{y_k}^2_{Y_{2,v^{(k)},\set{x_n}}} }^{{1\over 2}}\!}<\infty$.

 For each sequence $v=\set{v_i}$ from $(0,1]$,
 let $S(X,v)$ denote $\sumsum{X}{I,v}$, and 
 let $\tilde S\paren{X,v}$ denote $\seqsum{S(X,v)}{p,2,\set{v}}$,
 where $\set{v}$ is the sequence $\set{v,v,\ldots}$. 

\xproclaim {Proposition 4.14}. Let $2<p<\infty$ and let $X$ be a closed subspace of $\Lz{p}[0,1]$
                               satisfying the hypotheses (a${}'$) and (c${}'$) above.
                               Let $w=\set{w_i}$ be a sequence of scalars from $(0,1]$
                               satisfying condition $(*)$ of Proposition 2.1.
                               Let $S(X,w)$ and $\tilde S(X,w)$ be as above.
                               Then $\tilde S(X,w) \cinjects S(X,w)$.

\proof By condition $(*)$, we may choose a sequence $\set{N_k}$ of disjoint infinite \hfil\break 
       subsets of $\NN$ such that for each $\epsilon>0$ and for each $k$, 
       $$\sum_{{\scriptstyle w_i<\epsilon \atop \scriptstyle i\in N_k}} {w_i}^{{2p\over p-2}} = \infty.$$
       Hence for each $k$, we may choose a sequence 
       $\setwlimits{E_j^{(k)}}{j=1}{\infty}$ of disjoint nonempty finite subsets of $N_k$ such that for each $j$, 
       $$\paren{{w_j\over2}}^{{2p\over p-2}} \le \sum_{i\in E_j^{(k)}} {w_i}^{{2p\over p-2}} \le w_j^{{2p\over p-2}}.$$  
       Then for $v_j^{(k)} = \paren{\sum_{i\in E_j^{(k)}} {w_i}^{{2p\over p-2}}}^{{p-2\over 2p}}$,   
       ${w_j\over2} \le v_j^{(k)} \le w_j$.  
       Hence for $v^{(k)}=\setwlimits{v_j^{(k)}}{j=1}{\infty}$ and $y_k\in S\paren{X,w}$, 
       ${1\over2}\norm{y_k}_{Y_{2,w,\set{x_n}}} \le\norm{y_k}_{Y_{2,v^{(k)},\set{x_n}}} \le\norm{y_k}_{Y_{2,w,\set{x_n}}}$.
       Hence
       $$\tilde S(X,w) = \seqsum{S(X,w)}{p,2,\set{w}}
         \sim
         \paren{S\!\paren{X,v^{(1)}} \oplus S\!\paren{X,v^{(2)}} \oplus\cdots}_{\!p,2,\set{v^{(k)}}} \eqno{\TAG{4.6}{4.6}}$$   
       via the formal identity mapping. 

       Let $z_{j,n}^{(k)}=\sum_{i\in E_j^{(k)}} {w_i}^{{2\over p-2}} \tilde y_{i,n}$ 
       and let $\tilde z{}_{j,n}^{(k)}$ be the normalization of $z_{j,n}^{(k)}$ in $\L{p}(\Omega)$.
       Then by part (a) of Proposition 4.12, for each $k$ there is an isomorphism \hfil\break 
       $J_k\colon S\paren{X,v^{(k)}} \to \spansub{\tilde z{}_{j,n}^{(k)}:j,n\in\NN}{S(X,w)}$.
       Moreover, for $y_k\in S\paren{X,v^{(k)}}$, \hfil\break 
       $\normsub{J_k(y_k)}{Y_{p,\set{x_n}}} = \normsub{y_k}{Y_{p,\set{x_n}}}$ and 
       $\normsub{J_k(y_k)}{Y_{2,w,\set{x_n}}} = \normsub{y_k}{Y_{2,v^{(k)},\set{x_n}}}$ 
       by equations \TAG{4.4}{4.4} and \TAG{4.5}{4.5}, respectively. Hence  
       $$\paren{S\paren{X,v^{(1)}} \oplus S\paren{X,v^{(2)}} \oplus \cdots}_{p,2,\set{v^{(k)}}} 
         \sim   
         \paren{\spansub{\tilde z{}_{j,n}^{(1)}}{S(X,w)}\oplus\spansub{\tilde z{}_{j,n}^{(2)}}{S(X,w)}\oplus\cdots}
         _{p,2,\set{w}}  \eqno{\TAG{4.7}{4.7}}$$   
       via the isometry $\set{y_k}\mapsto\set{J_k(y_k)}$. 

       The direct sum on the right side of \TAG{4.7}{4.7} should be thought of as an internal \hfil\break 
       direct sum of subspaces of $S(X,w)$.
       We next show that
       $$\paren{\spansub{\tilde z{}_{j,n}^{(1)}}{S(X,w)}\oplus\spansub{\tilde z{}_{j,n}^{(2)}}{S(X,w)}\oplus\cdots}_{p,2,\set{w}} 
         \sim \spansub{\tilde z{}_{j,n}^{(k)}:j,n,k\in\NN}{S(X,w)} \eqno{\TAG{4.8}{4.8}}$$   
       via the mapping $\set{s_k}\mapsto\sum s_k$.
       For each $k$ and for scalars $a_{j,n}^{(k)}$, let \hfil\break 
       $s_k=\sum_j\sum_n a_{j,n}^{(k)} \tilde z{}_{j,n}^{(k)}
            \in \spansub{\tilde z{}_{j,n}^{(k)}:j,n\in\NN}{S(X,w)}.$
       Then by equations \TAG{4.4}{4.4} and \TAG{4.5}{4.5}, 

       \vglue1.3in 
       \eject      
       $$\eqalign{
         \norm{\set{s_k}}
            &=\max\set{\paren{\tsuml\normsubpower{s_k}{Y_{p,\set{x_n}}}{p\vphantom{2}}}^{{1\over p}},
                       \paren{\tsuml\normsubpower{s_k}{Y_{2,w,\set{x_n}}}{2}}^{{1\over 2}}} 
         \cr&=\max\set{\paren{\tsuml_k
              \normsubpower{\tsuml_j\tsuml_n a_{j,n}^{(k)} \tilde z{}_{j,n}^{(k)}}{Y_{p,\set{x_n}}}{p\vphantom{2}}}^{{1\over p}},
                       \paren{\tsuml_k
              \normsubpower{\tsuml_j\tsuml_n a_{j,n}^{(k)} \tilde z{}_{j,n}^{(k)}}{Y_{2,w,\set{x_n}}}{2}}^{{1\over 2}} }
         \cr&=\max\set{\paren{\tsuml_k \tsuml_j \norm{\tsuml_n a_{j,n}^{(k)} x_n}_p^p }^{{1\over p}},
                       \paren{\tsuml_k \tsuml_j \paren{v_j^{(k)}}^2 \tsuml_n \abs{a_{j,n}^{(k)}}^2 \norm{x_n}_2^2 }^{{1\over 2}} }
         \cr&=\max\set{\paren{\norm{\tsuml_k\tsuml_j\tsuml_n
                       a_{j,n}^{(k)} \tilde z{}_{j,n}^{(k)} }_{Y_{p,\set{x_n}}}^p }^{{1\over p}},
                       \paren{\norm{\tsuml_k\tsuml_j\tsuml_n
                       a_{j,n}^{(k)} \tilde z{}_{j,n}^{(k)} }_{Y_{2,w,\set{x_n}}}^2 }^{{1\over 2}} }
         \cr&=\xnorm{\tsuml_k\tsuml_j\tsuml_n a_{j,n}^{(k)} \tilde z{}_{j,n}^{(k)}}
       \approx\normsub{\tsuml_k\tsuml_j\tsuml_n a_{j,n}^{(k)} \tilde z{}_{j,n}^{(k)}}{S(X,w)}
             =\normsub{\tsuml s_k}{S(X,w)},}$$
       where $\xnorm{\pz}$ is as in the proof of Proposition 4.12. 
       Hence the mapping $\set{s_k}\mapsto\sum s_k$ is an isomorphism.

       By part (b) of Proposition 4.12, we have 
       $$\spansub{\tilde z{}_{j,n}^{(k)}:j,n,k\in\NN}{S(X,w)} \cinjects \sumsum{X}{I,w} = S(X,w). \eqno{\TAG{4.9}{4.9}}$$   
       Combining \TAG{4.6}{4.6}, \TAG{4.7}{4.7}, \TAG{4.8}{4.8}, and \TAG{4.9}{4.9} yields 
       $\tilde S(X,w) \cinjects S(X,w)$.
       \qquad\QED

\xproclaim {Proposition 4.15}. Let $2<p<\infty$ and let $X$ be a closed subspace of $\Lz{p}[0,1]$ \hfil\break 
                               satisfying the hypotheses (a${}'$) and (c${}'$) above.
                               Let $w=\set{w_i}$ be a sequence of scalars from $(0,1]$
                               satisfying condition $(*)$ of Proposition 2.1. Then \hfil\break 
                               $\sumsum{X}{I,w}\sim\sumsum{X}{I,w}\oplus\sumsum{X}{I,w}$. 

\proof Let $S(X,w)$ and $\tilde S(X,w)$ be as in Proposition 4.14. Then \hfil\break 
       $\tilde S(X,w) \cinjects S(X,w)$.
       Let $Y$ be a closed subspace of $S(X,w)$ such that \hfil\break 
       $S(X,w) \sim \tilde S(X,w) \oplus Y$. 
       Note that $ \tilde S(X,w) \sim S(X,w) \oplus \tilde S(X,w)$. Hence \hfil\break 
       $S(X,w) \oplus S(X,w) \sim S(X,w) \oplus \tilde S(X,w) \oplus Y
                                       \sim \tilde S(X,w) \oplus Y
                                       \sim S(X,w)$.
       \qquad\QED

\xproclaim {Theorem 4.16}. Let $2<p<\infty$ and let $X$ be a closed subspace of $\Lz{p}[0,1]$ \hfil\break 
                           satisfying the hypotheses (a${}'$) and (c${}'$) above.
                           Let $w=\set{w_i}$ and $w'=\set{w'_i}$ be \hfil\break 
                           sequences of scalars from $(0,1]$ satisfying condition $(*)$ of Proposition 2.1. Then \hfil\break 
                           $\sumsum{X}{I,w}\sim\sumsum{X}{I,w'}$. 

\proof The spaces $\sumsum{X}{I,w}$ and $\sumsum{X}{I,w'}$ satisfy the hypotheses of \hfil\break 
       Lemma 2.8. \qquad\QED

\definition Let $2<p<\infty$.
            Let $X$ be a closed subspace of $\Lz{p}[0,1]$ satisfying \hfil\break 
   (a${}'$) the orthogonal projection of $\L{2}[0,1]$ onto $\overline X\subset\L{2}[0,1]$, 
            when restricted to $\L{p}[0,1]$,
            yields a bounded projection $P:\L{p}[0,1]\to X\subset\L{p}[0,1]$ onto $X$, and \hfil\break 
   (c${}'$) $X$ has an unconditional orthogonal normalized basis $\set{x_n}$. \hfil\break 
            Define $\sumsum{X}{I}$, the independent sum of $X$, to be (the isomorphism type of) \hfil\break 
            $\sumsum{X}{I,w}$ for any sequence $w=\set{w_i}$ of scalars from $(0,1]$ satisfying condition $(*)$ 
            of Proposition 2.1.

 By Theorem 4.16, $\sumsum{X}{I}$ is well-defined.


\preheadspace
\firsthead{The Space $D_p$}
\postheadspace

\definition Let $2<p<\infty$, let $\set{x_n}$ be the sequence of Rademacher functions, and let $X=\span{x_n}_{\L{p}}\sim\l{2}$. 
            Define $D_p$ to be $\sumsum{X}{I}$.
            For the conjugate index $q$, define $D_q$ to be $D_p^*$.  

\xproclaim {Proposition 4.17}. Let $1<p<\infty$ where $p\ne2$. Then
                               \item{(a)} $X_p \cinjects D_p$, 
                               \item{(b)} $\sumsum{\l{2}}{\l{p}} \cinjects D_p$, and
                               \item{(c)} $\sumsum{\l{2}}{\l{p}}\oplus X_p \cinjects D_p$. 

\proof It suffices to prove the result for $2<p<\infty$, since the result for \hfil\break 
       $1<p<2$ will then follow by duality.

       Suppose $2<p<\infty$.
       Realize $D_p$ as $\sumsum{X}{I,w}$,  
       where $X$ and $\set{x_n}$ are as in the definition of $D_p$,
       and $w=\set{w_i}$ is a sequence of scalars from $(0,1]$ satisfying condition $(*)$ of Proposition 2.1.
       Then $D_p = \span{\tilde y_{i,n}:i,n\in\NN}_{\L{p}(\Omega)}$, where \hfil\break 
       $\tilde y_{i,n} = \paren{T_i\paren{x_n}}\circ\pi_i\in\L{p}(\Omega)$,
       and $T_i$ and $\pi_i$ are as in the definition of $\sumsum{X_i}{I,w}$. 

  \item{(a)}
       Let $D_p^{(1)} = \span{\tilde y_{i,1}:i\in\NN}_{\L{p}(\Omega)}$.  
       Then $D_p^{(1)}$ is a complemented subspace of $D_p$ by the unconditionality of $\set{\tilde y_{i,n}}$,   
       and $D_p^{(1)} = \sumsum{X^{(1)}}{I,w}$
       where $X^{(1)} = \span{x_1}_{\L{p}}$ and $x_1=1_{[0,{1\over2})}-1_{[{1\over2},1]}$.
       As noted in Example 4.1, $\sumsum{X^{(1)}}{I,w} \sim X_p$.  
       Hence $X_p \sim D_p^{(1)} \cinjects D_p$. 

  \item{(b)}
       Choose an increasing sequence $\set{i_k}$ of positive integers such that $\sum w_{i_k}^{{2p\over p-2}}<\infty$,
       and let $w'=\set{w_{i_k}}$.
       Let $D'_p = \span{\tilde y_{i_k,n}:k,n\in\NN}_{\L{p}(\Omega)}$.  
       Then $D'_p$ is a complemented subspace of $D_p$ by the unconditionality of $\set{\tilde y_{i,n}}$, and \hfil\break 
       $D'_p = \sumsum{X}{I,w'} \sim \sumsum{X}{\l{p}} \sim \sumsum{\l{2}}{\l{p}}$ by Corollary 4.10. Hence \hfil\break 
       $\sumsum{\l{2}}{\l{p}} \sim D'_p \cinjects D_p$. 

  \item{(c)}
       By Proposition 4.15 and parts (a) and (b) above, \hfil\break 
       $\sumsum{\l{2}}{\l{p}}\oplus X_p \cinjects D_p\oplus D_p \sim D_p$. 

       \QED

 For $2<p<\infty$, it is clear that $D_p\not\cinjects B_p$, 
 since otherwise $X_p\cinjects D_p\cinjects B_p$ by part (a) of Proposition 4.17, 
 so $X_p\cinjects B_p$, contrary to part (g) of Proposition 2.37. 

 We now present results leading to the conclusion that $B_p\not\cinjects D_p$ \xcite{A}. 
 We begin with a definition and some preliminary observations used in the proof of the \hfil\break 
 subsequent lemma.  

 Let $2<p<\infty$ and let $\set{r_n}$ be the sequence of Rademacher functions. 
 Given a sequence $w=\set{w_i}$ of positive scalars, 
 let $\tilde y_{i,n} = T_i(r_n)\circ\pi_i$,
 where $T_i$ and $\pi_i$ are as in the definition of $\sumsum{X_i}{I,w}$.
 Let $P_0:D_p\to D_p$ be the zero mapping. 
 For each \hfil\break 
 $m\in\NN$, let $P_m:D_p\to D_p$ be the natural projection
 of $D_p$ onto \hfil\break 
 $\span{\tilde y_{i,n}: i\in\set{1,\ldots,m},\>n\in\NN}_{D_p}$. 
 A sequence $\set{z_k}$ in $D_p$ will be said to be strip \hfil\break 
 disjoint if there is an increasing sequence $\set{m_k}$ in $\NN$ such that \hfil\break 
 $\norm{\paren{P_{m_k}-P_{m_{k-1}}}\paren{z_k} }_{D_p} \ge \paren{1-{1\over2^k}}\norm{z_k}_{D_p}$
 for all $k\in\NN$. 

 Let $2<p<\infty$, let $w$ be a positive scalar, and let $\set{w}=\set{w,w,\ldots}$. 
 Let $\set{e_n}$ be the standard basis for $X_{p,\set{w}}$. 
 Let $T:X_{p,\set{w}} \to D_p$ be an isomorphic imbedding.
 Suppose $\epsilon>0$ is such that for each $m\in\NN$,
 $\norm{P_m\paren{T\paren{e_n}}}_{D_p}<\epsilon$ for infinitely many $n\in\NN$.

 Then we may choose increasing sequences $\set{\gamma(n)}$ and $\set{m(n)}$ in $\NN$ such that \hfil\break 
 $T\paren{e_{\gamma(n)}}=x_n+y_n$, where $x_n=P_{m(n)}\paren{T\paren{e_{\gamma(n)}}}$, $\norm{x_n}_{D_p}<\epsilon$, 
 $\set{y_n}$ is strip disjoint, and
 $\set{x_n}$ and $\set{y_n}$ are block basic sequences with respect to the standard basis of $D_p$. 

 There are constants $K$ and $C$ such that for each finite $F\subset\NN$,
 $$\norm{T^{-1}}^{-1} \norm{\tsuml_{n\in F} e_{\gamma(n)} }_{X_{p,\set{w}}}
   \le 
   \norm{\tsuml_{n\in F} T\paren{e_{\gamma(n)}} }_{D_p}
   \le 
   \norm{\tsuml_{n\in F} x_n }_{D_p} + \norm{\tsuml_{n\in F} y_n }_{D_p},$$
 where [letting $\abs{F}$ denote the cardinality of $F$] 
 $$\norm{\tsuml_{n\in F} e_{\gamma(n)} }_{X_{p,\set{w}}} = \max\set{\abs{F}^{{1\over p}}, \abs{F}^{{1\over2}} w },$$ 
 $$\norm{\tsuml_{n\in F} x_n }_{D_p} \le K \paren{\tsuml_{n\in F} \norm{x_n}_{D_p}^2 }^{{1\over2}} 
                                     \le K \paren{\tsuml_{n\in F} \epsilon^2 }^{{1\over2}}  
                                      =  K \abs{F}^{{1\over2}} \epsilon,$$
 and 
 $$\eqalign{\norm{\tsuml_{n\in F} y_n }_{D_p}
           &\le 
            C \max\set{\paren{\tsuml_{n\in F} \norm{y_n}_p^p }^{{1\over p}},
                       \paren{\tsuml_{n\in F} \norm{y_n}_2^2 }^{{1\over2}} } 
       \cr &\le 
            C \max\set{\abs{F}^{{1\over p}} \norm{T},
                       \abs{F}^{{1\over2}} \max_{n\in F} \norm{y_n}_2 } .}$$
 Thus for $F$ such that
 $\abs{F}^{{1\over2}} w > \abs{F}^{{1\over p}}$ and 
 $\abs{F}^{{1\over2}} \max_{n\in F} \norm{y_n}_2 > \abs{F}^{{1\over p}} \norm{T}$,  
 $$\norm{T^{-1}}^{-1} \abs{F}^{{1\over2}} w
   \le  
   K \abs{F}^{{1\over2}} \epsilon + C \abs{F}^{{1\over2}} \max_{n\in F} \norm{y_n}_2,$$
 so
 $$\max_{n\in F} \norm{y_n}_2 \ge {\norm{T^{-1}}^{-1} w - K \epsilon \over C}.$$ 
 Hence we may choose an increasing sequence $\set{\beta(n)}$ in $\NN$
 such that for all $n\in\NN$ 
 $$\norm{y_{\beta(n)}}_2 \ge {\norm{T^{-1}}^{-1} w - K \epsilon \over C}.$$

\xproclaim {Lemma 4.18}. Let $2<p<\infty$.
                         Let $\set{e_{i,n}}$ be the standard basis for $B_p$ and let \hfil\break 
                         $w_i=\paren{{1\over i}}^{{p-2\over2p}}$.  
                         Suppose $T:B_p\to D_p$ is an isomorphic imbedding. 
                         Then there is an $\epsilon>0$ such that for all but a finite number of $i\in\NN$,
                         there is an $m_i\in\NN$ and an infinite $K_i\subset\NN$ 
                         such that $\norm{P_{m_i}\paren{T\paren{e_{i,n}}} }_{D_p} \ge w_i \epsilon$ 
                         for all $n\in K_i$. 

\proof Suppose the conclusion is false. 
       Then for each $\epsilon>0$, there is an infinite ${\NN}_{\epsilon}\subset\NN$ such that
       for all $i\in{\NN}_{\epsilon}$, all $m\in\NN$, and all infinite $K\subset\NN$, 
       there is an $n\in K$ such that $\norm{P_m\paren{T\paren{e_{i,n}}} }_{D_p} < w_i \epsilon$. 

       Fix $\epsilon>0$ and let $\epsilon_i={\epsilon\over2^i}$. 
       For $i\in\NN$, choose $\alpha(i)\in{\NN}_{\epsilon_i}$ such that 
       $\set{\alpha(i)}$ is an increasing sequence in $\NN$. 
       Let $i\in\NN$. 
       Then for each $m\in\NN$, \hfil\break 
       $\norm{P_m\paren{T\paren{e_{\alpha(i),n}}} }_{D_p} < w_{\alpha(i)} \epsilon_i = {w_{\alpha(i)}\over2^i} \epsilon$  
       for infinitely many $n\in\NN$. 

       We may choose increasing sequences $\set{\gamma_i(n)}$ and $\set{m_i(n)}$ in $\NN$ such that \hfil\break 
       $T\paren{e_{\alpha(i),\gamma_i(n)}} = x_{i,n}+y_{i,n}$, where
       $x_{i,n} = P_{m_i(n)}\paren{T\paren{e_{\alpha(i),\gamma_i(n)} }}$, 
       $\norm{x_{i,n}}_{D_p} < {w_{\alpha(i)}\over2^i} \epsilon$, 
       $\set{y_{i,n}}_{i,n\in\NN}$ is strip disjoint, and 
       $\set{x_{i,n}}_{i,n\in\NN}$ and $\set{y_{i,n}}_{i,n\in\NN}$
       are block basic sequences with respect to the standard basis of $D_p$. 

       There are constants $K$ and $C$,
       and there is an increasing sequence $\set{\beta_i(n)}$ in $\NN$, 
       such that for all $n\in\NN$  
       $$\norm{y_{i,\beta_i(n)}}_2 \ge {\norm{T^{-1}}^{-1} w_{\alpha(i)} - K {w_{\alpha(i)}\over2^i} \epsilon \over C}.$$

       By the fact that $\L{p}$ is of type $2$ \xciteplus{W}{III.A.17,23}, 
       and by H\"older's inequality for conjugate indices $p'={p\over2}$ and $q'={p\over p-2}$,
       there is a constant $K$ such that for scalars $a_{i,n}$ 
       $$\eqalign{\norm{\tsuml_i\tsuml_n a_{i,n} x_{i,n} }_{D_p}
                 &\le K
                  \paren{\tsuml_i\tsuml_n \abs{a_{i,n}}^2 \norm{x_{i,n}}_{D_p}^2 }^{{1\over2}}
             \cr &\le K
                  \paren{\tsuml_i\tsuml_n \abs{a_{i,n}}^2 \paren{{w_{\alpha(i)}\over2^i}\epsilon}^2 }^{{1\over2}}
             \cr &= K \epsilon
                  \paren{\tsuml_i\paren{\tsuml_n \abs{a_{i,n}}^2 w_{\alpha(i)}^2 } \paren{{1\over2^i}}^2 }^{{1\over2}}
             \cr &\le K \epsilon
                  \paren{
                  \paren{\tsuml_i\paren{\tsuml_n \abs{a_{i,n}}^2 w_{\alpha(i)}^2 }^{{p\over2}} }^{{2\over p}}
                  \paren{\tsuml_i \paren{{1\over2^i}}^{2{p\over p-2}} }^{{p-2\over p}}
                        }^{{1\over2}}
             \cr &= K \epsilon
                  \paren{\tsuml_i\paren{\tsuml_n \abs{a_{i,n}}^2 w_{\alpha(i)}^2 }^{{1\over2}p} }^{{1\over p}} 
                  \paren{\tsuml_i \paren{{1\over2^i}}^{{2p\over p-2}} }^{{p-2\over2p}}
             \cr &\le K \epsilon 
                  \norm{\tsuml_i\tsuml_n a_{i,n} e_{\alpha(i),n} }_{B_p}
                  \paren{\tsuml_i \paren{{1\over2^i}}^{{2p\over p-2}} }^{{p-2\over2p}}
             \cr &= K \epsilon 
                  \norm{\tsuml_i\tsuml_n a_{i,n} e_{\alpha(i),\gamma_i(n)} }_{B_p}
                  \paren{\tsuml_i \paren{{1\over2^i}}^{{2p\over p-2}} }^{{p-2\over2p}} .}$$

       Thus given $\delta>0$,
       $\norm{\sum_i\sum_n a_{i,n} x_{i,n} }_{D_p} \le\delta \norm{\sum_i\sum_n a_{i,n} e_{\alpha(i),\gamma_i(n)} }_{B_p}$
       for $\epsilon$ \hfil\break 
       sufficiently small.
       Define $S: \span{e_{\alpha(i),\gamma_i(n)} : i,n\in{\NN}}_{B_p} \to D_p$ by \hfil\break 
       $S\paren{\sum_i\sum_n a_{i,n} e_{\alpha(i),\gamma_i(n)} } = \sum_i\sum_n a_{i,n} y_{i,n}$.  
       Then for $\epsilon$ sufficiently small, $S$ is an \hfil\break 
       isomorphic imbedding. 
       Since $\set{y_{i,n}}_{i,n\in\NN}$ is strip disjoint,
       $\span{y_{i,n} : i,n\in{\NN}}_{D_p} \sim X_{p,v}$ for some $v$. 
       However, $\span{e_{\alpha(i),\gamma_i(n)} : i,n\in{\NN}}_{B_p} \sim B_p$.
       Since $X_{p,v}\injects\l{2}\oplus\l{p}$ by Proposition 2.1, Theorem 2.12, and part (a) of Proposition 2.24, 
       $B_p \sim  \span{e_{\alpha(i),\gamma_i(n)} : i,n\in{\NN}}_{B_p} 
            \injects \span{y_{i,n} : i,n\in{\NN}}_{D_p} \sim X_{p,v}
            \injects \l{2}\oplus\l{p}$, 
       so $B_p\injects\l{2}\oplus\l{p}$, contrary to Lemma 2.23 and part (a) of Proposition 2.37. 
       \qquad\QED

\xproclaim {Lemma 4.19}. Let $2<p<\infty$.
                         Let $w=\set{w_i}$ where $w_i=\paren{1\over i}^{{p-2\over2p}}$, 
                         and let $\tilde y_{i,n}$ be as above. 
                         Let $\set{E_{\ell}}$ be a sequence of disjoint nonempty finite subsets of $\NN$.  
                         Let $\set{z_{k,\ell}}$ be a sequence in $D_p$ which is 
                         normalized with respect to $\xnorm{\pz}_{D_p}$
                         such that for each $\ell\in\NN$,
                         $z_{k,\ell}\in\span{\tilde y_{i,n}:i\in E_{\ell},n\in{\NN}}_{D_p}$ for all $k\in\NN$ and 
                         $\set{z_{k,\ell}}_{k\in\NN}$ is equivalent to the standard basis of $\l{2}$. 
                         Then there is an infinite $L\subset\NN$,
                         and for each $\ell\in L$ there is an infinite $K_{\ell}\subset\NN$, 
                         such that $\set{z_{k,\ell}}_{k\in K_{\ell},\>\ell\in L}$
                         is equivalent to either the standard basis of $\l{2}$
                         or the standard basis of $\sumsum{\l{2}}{\l{p}}$. 

\proof By passing to a subsequence,
       we may assume that $\set{z_{k,\ell}}$ is a block basic sequence with respect to the standard basis of $D_p$.  

       Let $z_{k,\ell}=\sum_{i\in E_{\ell}} v_{i,k,\ell}$ where $v_{i,k,\ell}=\sum_{n\in N_{i,k,\ell}} b_{i,n} \tilde y_{i,n}$
       for $N_{i,k,\ell}\subset\NN$ and scalars $b_{i,n}$.
       Let $\lambda_{i,k,\ell}=\paren{\sum_{n\in N_{i,k,\ell}} \abs{b_{i,n}}^2 }^{{1\over2}}$.
       Then for scalars $a_{k,\ell}$ 
       $$\eqalign{&\xnorm{\tsuml_{\ell} \tsuml_k a_{k,\ell} z_{k,\ell}}_{D_p} = 
                   \xnorm{\tsuml_{\ell} \tsuml_k a_{k,\ell} \tsuml_{i\in E_{\ell}} v_{i,k,\ell} }_{D_p} = 
                   \xnorm{\tsuml_{\ell} \tsuml_k a_{k,\ell} \tsuml_{i\in E_{\ell}}
                          \tsuml_{n\in N_{i,k,\ell}} b_{i,n} \tilde y_{i,n} }_{D_p}
              \cr &= 
                   \xnorm{\tsuml_{\ell} \tsuml_{i\in E_{\ell}} \tsuml_k
                          \tsuml_{n\in N_{i,k,\ell}} a_{k,\ell} b_{i,n} \tilde y_{i,n} }_{D_p}
           \cr &= \max\set{
               \paren{\tsuml_{\ell} \tsuml_{i\in E_{\ell}}
               \paren{\tsuml_k \tsuml_{n\in N_{i,k,\ell}} \abs{a_{k,\ell} b_{i,n}}^2 }^{{1\over2}p} }^{{1\over p}},
               \paren{\tsuml_{\ell} \tsuml_{i\in E_{\ell}} w_i^2 \tsuml_k
                      \tsuml_{n\in N_{i,k,\ell}} \abs{a_{k,\ell} b_{i,n}}^2 }^{{1\over2}} } 
           \cr &= \max\set{
               \paren{\tsuml_{\ell} \tsuml_{i\in E_{\ell}}
               \paren{\tsuml_k \abs{a_{k,\ell}}^2 \lambda_{i,k,\ell}^2 }^{{p\over2}} }^{{1\over p}},
               \paren{\tsuml_{\ell} \tsuml_{i\in E_{\ell}} w_i^2
                      \tsuml_k \abs{a_{k,\ell}}^2 \lambda_{i,k,\ell}^2 }^{{1\over2}} }.}$$ 
       As a special case of the above, 
       $$1=\xnorm{z_{k,\ell}}_{D_p} 
          =\max\set{
          \paren{\tsuml_{i\in E_{\ell}} \lambda_{i,k,\ell}^p}^{{1\over p}},
          \paren{\tsuml_{i\in E_{\ell}} w_i^2 \lambda_{i,k,\ell}^2}^{{1\over2}} }
          \ge\paren{\sum_{i\in E_{\ell}} \lambda_{i,k,\ell}^p}^{{1\over p}},$$
       whence $\lambda_{i,k,\ell}\le1$ for $k,\ell\in\NN$ and $i\in E_{\ell}$. 
       Let $\set{\epsilon_k}$ be a sequence of positive scalars with limit zero.
       For each $\ell\in\NN$, choose an increasing sequence $\set{\alpha_{\ell}(k)}$ in $\NN$ 
       and scalars $\Lambda_i$ for $i\in E_{\ell}$ such that 
       $\abs{\lambda_{i,\alpha_{\ell}(k),\ell}-\Lambda_i}<\epsilon_k$ for $k\in\NN$ and $i\in E_{\ell}$. 
       Then
       $$\eqalign{&\xnorm{\tsuml_{\ell} \tsuml_k a_{k,\ell} z_{\alpha_{\ell}(k),\ell}}_{D_p}  
           \cr &= \max\set{
                  \paren{\tsuml_{\ell} \tsuml_{i\in E_{\ell}}
                  \paren{\tsuml_k \abs{a_{k,\ell}}^2 \lambda_{i,\alpha_{\ell}(k),\ell}^2 }^{{p\over2}} }^{{1\over p}},
                  \paren{\tsuml_{\ell} \tsuml_{i\in E_{\ell}} w_i^2
                         \tsuml_k \abs{a_{k,\ell}}^2 \lambda_{i,\alpha_{\ell}(k),\ell}^2 }^{{1\over2}} } 
           \cr &\approx \max\set{
                  \paren{\tsuml_{\ell} \tsuml_{i\in E_{\ell}} \Lambda_i^p
                  \paren{\tsuml_k \abs{a_{k,\ell}}^2 }^{{1\over2}p} }^{{1\over p}},
                  \paren{\tsuml_{\ell} \tsuml_{i\in E_{\ell}} w_i^2 \Lambda_i^2
                  \paren{\tsuml_k \abs{a_{k,\ell}}^2 }^{{1\over2}2} }^{{1\over2}} }  
           \cr &= \max\set{
                  \paren{\tsuml_{\ell} \tsuml_{i\in E_{\ell}} \Lambda_i^p
                         \norm{\set{a_{k,\ell}}_k}_{\l{2}}^p }^{{1\over p}},
                  \paren{\tsuml_{\ell} \tsuml_{i\in E_{\ell}} w_i^2 \Lambda_i^2
                         \norm{\set{a_{k,\ell}}_k}_{\l{2}}^2 }^{{1\over2}} } ,}$$
       where the approximation can be improved to any degree by the choice of
       $\set{\epsilon_k}$ and $\set{\alpha_{\ell}(k)}$. 
       As a special case of the above, \hfil\break 
       $1 = \xnorm{z_{\alpha_{\ell}(k),\ell}}_{D_p} \approx
            \max\set{\paren{\sum_{i\in E_{\ell}} \Lambda_i^p}^{{1\over p}},
                     \paren{\sum_{i\in E_{\ell}} w_i^2 \Lambda_i^2}^{{1\over2}} }$,
       where the approximation can be improved to any degree by the choice of
       $\set{\epsilon_k}$ and $\set{\alpha_{\ell}(k)}$. 
       Hence $\set{z_{\alpha_{\ell}(k),\ell}}$
       can be chosen to be equivalent to the standard basis of $\sumsum{\l{2}}{I,W}$ 
       where $W=\set{W_{\ell}}$ and 
       $$W_{\ell} =
         {\paren{\sum_{i\in E_{\ell}} w_i^2 \Lambda_i^2 }^{{1\over2}}   
          \over
          \paren{\sum_{i\in E_{\ell}} \Lambda_i^p }^{{1\over p}} }.$$
       If $\xinf_{\ell\in\NN} W_{\ell}>0$,
       then $\set{z_{\alpha_{\ell}(k),\ell}}$ is equivalent to the standard basis of $\l{2}$. \hfil\break 
       If $\xinf_{\ell\in\NN} W_{\ell}=0$,
       then $\set{z_{\alpha_{\ell}(k),\ell}}$ is equivalent to the standard basis of $\sumsum{\l{2}}{\l{p}}$. 
       \qquad\QED

\remark As a special case of the first display in the above proof, \hfil\break 
        $\xnorm{v_{i,k,\ell}}_{D_p}=\max\set{\lambda_{i,k,\ell},w_i\lambda_{i,k,\ell}}=\lambda_{i,k,\ell}$.

\xproclaim {Lemma 4.20}. Let $2<p<\infty$. 
                         Suppose $T:B_p\to D_p$ is an isomorphic imbedding. 
                         Then $B_p$ has a complemented subspace $X$ isomorphic to $B_p$,  
                         and $D_p$ has a closed subspace $Y$ isomorphic to
                         $\l{2}\oplus X_{p,v}$ or $\sumsum{\l{2}}{\l{p}} \oplus X_{p,v}$ for some $v$,   
                         such that $T(X) \subset Y$.

\proof Choose (as we may by Lemma 4.18) $\epsilon>0$ and ${\NN}'\subset\NN$ with finite \hfil\break 
       complement
       such that for each $i\in{\NN}'$, there is an $m_i\in\NN$ and an infinite $K_i\subset\NN$
       such that $\norm{P_{m_i}\paren{T\paren{e_{i,n}}} }_{D_p} \ge w_i \epsilon$ for all $n\in K_i$. 

       For each $i\in{\NN}'$ and $n\in K_i$, 
       let $T\paren{e_{i,n}}=x_{i,n}+y_{i,n}$, where $x_{i,n}=P_{m_i}\paren{T\paren{e_{i,n}}}$.  
       For each $i\in{\NN}'$, choose an infinite $H_i\subset K_i$ such that $y_{i,n}=r_{i,n}+s_{i,n}$ for $n\in H_i$, 
       where $\norm{r_{i,n}}_{D_p} < {w_i\over2^i} \epsilon$ for $n\in H_i$,
       and $\set{s_{i,n}}_{n\in H_i}$ is strip disjoint.
       Choose infinite $G_i\subset H_i$ for $i\in{\NN}'$ such that 
       $\set{s_{i,n}}_{i\in{\NN}',n\in G_i}$ is strip disjoint.

       Now for $i\in{\NN}'$ and $n\in G_i$, $T\paren{e_{i,n}}=x_{i,n}+r_{i,n}+s_{i,n}$, where \hfil\break 
       $x_{i,n}=P_{m_i}\paren{T\paren{e_{i,n}}}$, $\norm{r_{i,n}}_{D_p} < {w_i\over2^i} \epsilon$,
       and $\set{s_{i,n}}_{i\in{\NN}',n\in G_i}$ is strip disjoint.

       For each $i\in{\NN}'$, choose an infinite $F_i\subset G_i$ such that
       $\set{x_{i,n}\big/\norm{x_{i,n}}_{D_p} }_{n\in F_i}$
       is $\paren{1+{1\over2^i}}$-equivalent to the standard basis of $\l{2}$.
       Choose (as we may by Lemma 4.19) an infinite ${\NN}''\subset{\NN}'$,
       and for each $i\in{\NN}''$ choose an infinite $E_i\subset F_i$, 
       such that $\span{x_{i,n} : i\in{\NN}'',n\in E_i }_{D_p}$ is isomorphic to $\l{2}$ or $\sumsum{\l{2}}{\l{p}}$.
       Now \hfil\break 
       $\span{x_{i,n} : i\in{\NN}'',n\in E_i }_{D_p} \sim \span{x_{i,n}+r_{i,n} : i\in{\NN}'',n\in E_i }_{D_p}$,
       since \hfil\break 
       $\norm{r_{i,n}}_{D_p} < {w_i\over2^i} \epsilon \le {\norm{x_{i,n}}_{D_p}\over2^i}$
       for $i\in{\NN}''$ and $n\in E_i$, and $\set{r_{i,n}}_{n}$ has an upper $\l{2}$ \hfil\break 
       estimate.

       Let $X=\span{e_{i,n} : i\in{\NN}'',n\in E_i }_{B_p} \sim B_p$, and let \hfil\break 
       $Y=\span{x_{i,n}+r_{i,n} : i\in{\NN}'',n\in E_i }_{D_p} \oplus \span{s_{i,n} : i\in{\NN}'',n\in E_i }_{D_p}$.
       Then \hfil\break 
       $T(X)=\span{x_{i,n}+r_{i,n}+s_{i,n} : i\in{\NN}'',n\in E_i }_{D_p}\subset Y$, and \hfil\break 
       $Y\sim\span{x_{i,n} : i\in{\NN}'',n\in E_i }_{D_p} \oplus \span{s_{i,n} : i\in{\NN}'',n\in E_i }_{D_p}$ 
       is isomorphic to $\l{2}\oplus X_{p,v}$ or $\sumsum{\l{2}}{\l{p}} \oplus X_{p,v}$ for some $v$.
       \qquad\QED

\xproclaim {Proposition 4.21}. Let $1<p<\infty$ where $p\ne2$. Then $B_p \not\cinjects D_p$.

\proof Suppose $2<p<\infty$ and $B_p \cinjects D_p$. Then \hfil\break 
       $B_p \cinjects \sumsum{\l{2}}{\l{p}} \oplus X_{p,v}$ for some $v$ by Lemma 4.20,  
       but $X_{p,v} \cinjects X_p$ for all $v$ by \hfil\break 
       Proposition 2.1, Theorem 2.12, and part (b) of Proposition 2.24.
       Hence \hfil\break 
       $B_p \cinjects \sumsum{\l{2}}{\l{p}} \oplus X_p$,  
       contrary to Proposition 2.48. 
       The result for $1<p<2$ now follows by duality. 
       \qquad\QED

\preheadspace
\firsthead{Sums Involving $D_p$}
\postheadspace

 A few more $\SL{p}$ spaces can be constructed by forming sums involving $D_p$. 
 The resulting spaces are $B_p\oplus D_p$ and $\sumsum{X_p}{\l{p}}\oplus D_p$. 

 We first present results leading to the conclusion that $D_p\not\cinjects\sumsum{X_p}{\l{p}}$ \xcite{A}. 

 Given $E\subset\NN$, let $P_E:\sumsum{\l{2}}{\l{r}}\to\sumsum{\l{2}}{\l{r}}$ 
 be the natural projection onto the subspace $\sumsum{X_i}{\l{r}}$  
 with $X_i=\l{2}$ if $i\in E$ and $X_i=\set{0}$ otherwise.   
 Given $M\in\NN$, let $P_M=P_{\set{1,\ldots,M}}$.  

 Given $F\subset\NN$, let $P'_F:\sumsum{X_q}{\l{q}}\to\sumsum{X_q}{\l{q}}$  
 be the natural projection onto the subspace $\sumsum{Y_i}{\l{q}}$ 
 with $Y_i=X_q$ if $i\in F$ and $Y_i=\set{0}$ otherwise.   
 Given $N\in\NN$, let $P'_N=P'_{\set{1,\ldots,N}}$.  

\xproclaim {Lemma 4.22}. Let $1<q<r<2$. 
                         Then $\sumsum{\l{2}}{\l{r}} \not\injects \sumsum{X_q}{\l{q}}$.  

\proof Suppose $\sumsum{\l{2}}{\l{r}} \injects \sumsum{X_q}{\l{q}}$.  
       Let $T:\sumsum{\l{2}}{\l{r}} \to \sumsum{X_q}{\l{q}}$ be an isomorphic imbedding.  
       Then given $n\in\NN$,
       $P'_n \circ T : \sumsum{\l{2}}{\l{r}} \to \sumsum{X_q}{\l{q}}$ is not an isomorphic imbedding,  
       essentially by Lemma 3.7. 
       Thus given $\epsilon>0$ and $m\in\NN$,
       there is an $x\in\sumsum{\l{2}}{\l{r}}$ with $P_m(x)=0$ such that 
       $\norm{P'_n\paren{T\paren{x}}}_{} < {\epsilon\over2\norm{T^{-1}}}\norm{x}_{}$. 
       Hence there is an $M\in\NN$ with $m<M$ such that \hfil\break 
       $\norm{P'_n(T(P_M(x)))}_{} <  {\epsilon\over2\norm{T^{-1}}} \norm{P_M(x)}_{}
                                  \le{\epsilon\over2} \norm{T(P_M(x))}_{}$. 
       Letting $y = P_M\paren{x}$ and \hfil\break 
       $E=\set{m+1,\ldots,M}$, 
       $P_E(y)=y$ and 
       $\norm{P'_n(T(y))}_{} < {\epsilon\over2} \norm{T(y)}_{}$. 
       Now there is an $N\in\NN$ with $n<N$ such that 
       $\norm{P'_N(T(y))}_{} > \paren{1-{\epsilon\over2}} \norm{T(y)}_{}$. 
       Letting $F=\set{n+1,\ldots,N}$, 
       $(1-\epsilon) \norm{T(y)}_{} < \norm{P'_F(T(y))}_{} \le \norm{T(y)}_{}$. 
       
       Given $\epsilon_1,\epsilon_2,\ldots>0$,
       we will inductively find disjoint nonempty finite sets \hfil\break 
       $E_1,E_2,\ldots\subset\NN$ with $\max E_i < \min E_{i'}$ for $i<i'$,
       $y_1,y_2,\ldots\in\sumsum{\l{2}}{\l{r}}$ with \hfil\break 
       $P_{E_i}(y_i)=y_i$, and
       disjoint nonempty finite sets $F_1,F_2,\ldots\subset\NN$ with $\max F_i < \min F_{i'}$ for $i<i'$, 
       such that $(1-\epsilon_i) \norm{T(y_i)} < \norm{P'_{F_i}(T(y_i))} \le \norm{T(y_i)}$ 
       for each $i\in\NN$. 

       Given $\epsilon_1>0$, the argument above with $n=1$ and $m=1$ shows how to find 
       a finite $E_1\subset\NN$ and $y_1\in\sumsum{\l{2}}{\l{r}}$ with $P_{E_1}(y_1)=y_1$,
       and a finite $F_1\subset\NN$, such that \hfil\break 
       $(1-\epsilon_1) \norm{T(y_1)} < \norm{P'_{F_1}(T(y_1))} \le \norm{T(y_1)}$. 

       Let $\set{\epsilon_i}$ be a sequence of positive scalars and let $k\in\NN$.
       Suppose $E_1,\ldots,E_k$, $y_1,\ldots,y_k$, and $F_1,\ldots,F_k$ 
       satisfying our requirements for all $i\in\set{1,\ldots,k}$ have been found.  
       The argument above with $n>\max F_k$ and $m>\max E_k$
       shows how to find a finite $E_{k+1}\subset\NN$ and $y_{k+1}\in\sumsum{\l{2}}{\l{r}}$
       with $\max E_k<\min E_{k+1}$ and \hfil\break 
       $P_{E_{k+1}}\paren{y_{k+1}}=y_{k+1}$, 
       and a finite $F_{k+1}\subset\NN$ with $\max F_k<\min F_{k+1}$,
       such that \hfil\break 
       $\paren{1-\epsilon_{k+1}} \norm{T\paren{y_{k+1}}}    
        <   \norm{P'_{F_{k+1}}\paren{T\paren{y_{k+1}}}} 
        \le \norm{T\paren{y_{k+1}}}$. 
       Thus $\set{E_i}$, $\set{y_i}$, and $\set{F_i}$ can be found as claimed. 
       
       For $\set{\epsilon_i}$ approaching zero rapidly and $\set{y_i}$ normalized,
       $\span{y_i}\sim\l{r}$, but \hfil\break 
       $\span{T\paren{y_i}} \sim \span{P'_{F_i}\paren{T\paren{y_i}}} \sim \l{q}$. 
       Hence $\l{r}\injects\l{q}$, contrary to fact. 
       It follows that no such isomorphic imbedding $T$ exists. 
       \qquad\QED

\xproclaim {Lemma 4.23}. Let $1\le q<\infty$
                         and let $\set{x_i}$ be unconditional in $\L{q}$.
                         Let $C$ be the sign-unconditional constant for $\set{x_i}$
                         and let $K_q$ be Khintchine's constant for $\L{q}$. 
                         Then for scalars $d_i$, 
                         $$\norm{\tsuml_i d_i x_i}_q^q
                           \within{C^q K_q^q}{C^q K_q^q}
                           \int\paren{\tsuml_i \abs{d_i x_i(s)}^2 }^{{1\over2}q} \,ds
                           =
                           \norm{\tsuml_i \abs{d_i x_i}^2 }_{{q\over2}}^{{q\over2}}.$$

\proof Let $\set{r_i}$ be the sequence of Rademacher functions. 
       Then by the \hfil\break 
       unconditionality of $\set{x_i}$, Fubini's theorem, and Khintchine's inequality, we have 
       $$\eqalign{\norm{\tsuml_i d_i x_i}_q^q
         & \within{C^q}{C^q}
         \int\paren{\int \abs{\tsuml_i d_i r_i(t) x_i(s)}^q \,ds }\,dt
     \cr & =
         \int\paren{\int \abs{\tsuml_i d_i x_i(s) r_i(t)}^q \,dt }\,ds
     \cr & \within{K_q^q}{K_q^q}
         \int\paren{\tsuml_i \abs{d_i x_i(s)}^2 }^{{1\over2}q} \,ds
     \cr & =
         \norm{\tsuml_i \abs{d_i x_i}^2 }_{{q\over2}}^{{q\over2}}.}$$
     \QED

\xproclaim {Lemma 4.24}. Let $1<q<r<2$. Then $\sumsum{\l{2}}{\l{r}} \injects D_q$. 

\proof Let $p$ be the conjugate index of $q$,
       let $\set{r_n}$ be the sequence of \hfil\break 
       Rademacher functions, 
       let $\Omega=\prod_{i=1}^{\infty} [0,1]$, and
       let $\set{N_i}$ be a sequence of disjoint \hfil\break 
       infinite subsets of $\NN$
       with ${\NN} = \bigcup_{i\in{\NN}} N_i$. 
       For each $i\in\NN$,
       let $\set{r_{i,n}}_{n\in{\NN}}=\set{r_n}_{n\in N_i}$, and 
       let \hbox{$z_i:[0,1]\to{\RR}$} be the normalization in $\L{q}$ of $1_{[0,k_i]}$, 
       where $k_i=w_i^{{2p\over p-2}}$ and
       $\set{w_i}$ is a sequence of positive scalars satisfying condition $(*)$ of Proposition 2.1.

       Let $u=\paren{u_1,u_2,\ldots}$ and $v=\paren{v_1,v_2,\ldots}$. 
       Now $\set{z_i(u_i)r_i(v_i)}_{i\in\NN}$,
       being a \hfil\break 
       sequence of independent symmetric three-valued random variables,
       and is equivalent to the standard basis of $X_{q,\set{w_i}}$.
       Thus by \xciteplus{RII}{Corollary 4.2},
       we may choose a \hfil\break 
       sequence $\set{a_i}$ of scalars and a sequence $\set{F_j}$
       of nonempty finite intervals in $\NN$ with \hfil\break 
       ${\NN} = \bigcup_{j\in{\NN}} F_j$ and 
       $1+\max F_j=\min F_{j+1}$, 
       such that for $y_j(u,v) = \sum_{i\in F_j} a_i z_i(u_i) r_i(v_i)$,
       \kern3pt
       $\set{y_j(u,v)}$ is a (perturbation of)
       a sequence of independent $r$-stable normalized random variables in $\L{q}(\Omega^2)$. 
       Then for scalars $b_{j,n}$,
       letting $c_j=\paren{\sum_n \abs{b_{j,n}}^2}^{{1\over2}}$,
       by Khintchine's inequality, Lemma 4.23, and the $r$-stability of $\set{y_j(u,v)}$,
       for $t=(t_{i,n})_{i\in\NN,\>n\in N_i}$ we have 

       $$\eqalign{&\norm{\tsuml_j\tsuml_n b_{j,n}
                        \tsuml_{i\in F_j} a_i z_i(u_i) r_i(v_i)
                        r_{i,n}(t_{i,n}) }_{\L{q}(\Omega^3)}^q 
              \cr &=
                  \int_{\Omega} \int_{\Omega} \paren{\int_{\Omega}
                  \abs{\tsuml_j\tsuml_n b_{j,n} \tsuml_{i\in F_j}
                       a_i z_i(u_i) r_i(v_i) r_{i,n}(t_{i,n}) }^q \,dt }\,du\,dv
              \cr &\approx
                  \int_{\Omega} \int_{\Omega} 
                  \paren{\tsuml_j\tsuml_n \abs{b_{j,n}}^2 \tsuml_{i\in F_j} \abs{a_i z_i(u_i) r_i(v_i)}^2 }^{{1\over2}q} \,du\,dv
              \cr &=
                  \int_{\Omega} \int_{\Omega} 
                  \paren{\tsuml_j \tsuml_{i\in F_j} \abs{c_j a_i z_i(u_i) r_i(v_i) }^2 }^{{1\over2}q} \,du\,dv
              \cr &\approx
                  \norm{\tsuml_j \tsuml_{i\in F_j} c_j a_i z_i(u_i) r_i(v_i) }_{\L{q}(\Omega^2)}^q
              \cr &=
                  \norm{\tsuml_j c_j y_j(u,v)}_{\L{q}(\Omega^2)}^q
              \approx
                  \paren{\tsuml_j \abs{c_j}^r }^{{1\over r}q}
              =
                  \paren{\tsuml_j \paren{\tsuml_n
                  \abs{b_{j,n}}^2 }^{{1\over2}r} }^{{1\over r}q}.}$$
       Hence 
       $$\span{\tsuml_{i\in F_j} a_i z_i(u_i) r_i(v_i) r_{i,n}(t_{i,n}):j,n\in{\NN} }_{\L{q}(\Omega^3)}
         \sim \sumsum{\l{2}}{\l{r}}.$$
       Moreover, by the choice of $\set{z_i}$, 
       $$\span{\tsuml_{i\in F_j} a_i z_i(u_i) r_i(v_i) r_{i,n}(t_{i,n}):j,n\in{\NN} }_{\L{q}(\Omega^3)} \injects D_q.$$
       It follows that $\sumsum{\l{2}}{\l{r}} \injects D_q$. 
       \qquad\QED

\xproclaim {Proposition 4.25}. Let $1<p<\infty$ where $p\ne2$. 
                               Then $D_p\not\cinjects\sumsum{X_p}{\l{p}}$.

\proof Suppose $1<q<2$ and $D_q\cinjects\sumsum{X_q}{\l{q}}$.
       Then for $1<q<r<2$, $\sumsum{\l{2}}{\l{r}} \injects D_q \cinjects \sumsum{X_q}{\l{q}}$ by Lemma 4.24, 
       so $\sumsum{\l{2}}{\l{r}} \injects \sumsum{X_q}{\l{q}}$, contrary to Lemma 4.22. 
       Hence $D_q\not\cinjects\sumsum{X_q}{\l{q}}$,
       and the result for $2<p<\infty$ holds by duality. 
       \qquad\QED

 Next we present results leading to the conclusion that $\sumsum{X_p}{\l{p}}\not\cinjects B_p\oplus D_p$ \xcite{A}. 

 Let $1<q<r<2$, and let $p$ be the conjugate index of $q$.  
 Let $\set{e_i}$ be the \hfil\break 
 standard basis of $\l{r}$. 
 Let $\set{z_{i,j}}$ be the standard basis of $D_p$,
 and let $\set{z_{i,j}^*}$ be the \hfil\break 
 corresponding dual basis of $D_q$,
 where for each $j\in\NN$, $\span{z_{i,j}:i\in\NN}_{D_p}\sim\l{2}$.

 Given $E\subset\NN$,
 let $P_E:\l{r}\to\l{r}$ be the natural projection onto the subspace \hfil\break 
 $\ell^r_{(E)} = \span{e_i:i\in E}_{\l{r}}$. 
 Given $M\in\NN$, let $P_M=P_{\set{1,\ldots,M}}$.  

 Given $F\subset\NN$, let $P'_F:D_q\to D_q$  
 be the natural projection onto the subspace $D_q^{(F)} = \span{z_{i,j}^*:i\in\NN,\>j\in F}_{D_q}$. 
 Given $N\in\NN$, let $P'_N=P'_{\set{1,\ldots,N}}$ and let \hfil\break 
 $D_q^{(N)}=D_q^{\set{1,\ldots,N}}$.  

\xproclaim {Lemma 4.26}. Let $1<q<r<2$. 
                         Suppose $T:\l{r}\to D_q$ is an \hbox{isomorphic} imbedding. 
                         Then for each sequence $\set{\epsilon_i}$ of positive scalars,
                         there is a normalized block basic sequence $\set{y_i}$ in $\l{r}$ 
                         and a sequence $\set{F_i}$ of disjoint nonempty finite subsets of $\NN$
                         with $\max F_i<\min F_{i'}$ for $i<i'$, such that
                         $\l{r} \sim \span{y_i}_{\l{r}}
                                \sim \span{T(y_i)}_{D_q}
                                \sim \span{P'_{F_i}\paren{T(y_i)}}_{D_q}$
                         via equivalence of natural bases, 
                         with 
                         $\paren{1-\epsilon_i} \norm{T\paren{y_i}} < \norm{P'_{F_i}\paren{T\paren{y_i}}}
                                                                 \le \norm{T\paren{y_i}}$ 
                         for each $i\in\NN$. 

\proof Given $n\in\NN$, $D_q^{(n)}\sim\l{2}$, so
       $P'_n \circ T : \l{r} \to D_q^{(n)}$ is not an isomorphic imbedding.  
       Thus given $\epsilon>0$ and $m\in\NN$,
       there is an $x\in\l{r}$ with $P_m(x)=0$ such that \hfil\break 
       $\norm{P'_n\paren{T\paren{x}}}_{} < {\epsilon\over2\norm{T^{-1}}}\norm{x}_{}$. 
       Hence there is an $M\in\NN$ with $m<M$ such that \hfil\break 
       $\norm{P'_n(T(P_M(x)))}_{} <  {\epsilon\over2\norm{T^{-1}}} \norm{P_M(x)}_{}
                                  \le{\epsilon\over2} \norm{T(P_M(x))}_{}$. 
       Letting $y = P_M\paren{x}$ and \hfil\break 
       $E=\set{m+1,\ldots,M}$, 
       $P_E(y)=y$ and 
       $\norm{P'_n(T(y))}_{} < {\epsilon\over2} \norm{T(y)}_{}$. 
       Now there is an $N\in\NN$ with $n<N$ such that 
       $\norm{P'_N(T(y))}_{} > \paren{1-{\epsilon\over2}} \norm{T(y)}_{}$. 
       Letting $F=\set{n+1,\ldots,N}$, 
       $(1-\epsilon) \norm{T(y)}_{} < \norm{P'_F(T(y))}_{} \le \norm{T(y)}_{}$. 
       
       Given $\epsilon_1,\epsilon_2,\ldots>0$,  
       we will inductively find disjoint nonempty finite sets \hfil\break 
       $E_1,E_2,\ldots\subset\NN$ with $\max E_i < \min E_{i'}$ for $i<i'$,
       $y_1,y_2,\ldots\in\l{r}$ with $P_{E_i}(y_i)=y_i$, and
       disjoint nonempty finite sets $F_1,F_2,\ldots\subset\NN$ with $\max F_i < \min F_{i'}$ for $i<i'$, 
       such that \hfil\break 
       $(1-\epsilon_i) \norm{T(y_i)} < \norm{P'_{F_i}(T(y_i))} \le \norm{T(y_i)}$ 
       for each $i\in\NN$. 

       Given $\epsilon_1>0$, the argument above with $n=1$ and $m=1$ shows how to find 
       a finite $E_1\subset\NN$ and $y_1\in\l{r}$ with $P_{E_1}(y_1)=y_1$,
       and a finite $F_1\subset\NN$, such that \hfil\break 
       $(1-\epsilon_1) \norm{T(y_1)} < \norm{P'_{F_1}(T(y_1))} \le \norm{T(y_1)}$. 

       Let $\set{\epsilon_i}$ be a sequence of positive scalars and let $k\in\NN$.
       Suppose $E_1,\ldots,E_k$, $y_1,\ldots,y_k$, and $F_1,\ldots,F_k$ 
       satisfying our requirements for all $i\in\set{1,\ldots,k}$ have been found.  
       The argument above with $n>\max F_k$ and $m>\max E_k$
       shows how to find a finite $E_{k+1}\subset\NN$ and $y_{k+1}\in\l{r}$ with 
       $\max E_k<\min E_{k+1}$ and $P_{E_{k+1}}\paren{y_{k+1}}=y_{k+1}$, 
       and a finite $F_{k+1}\subset\NN$ with $\max F_k<\min F_{k+1}$,
       such that \hfil\break 
       $\paren{1-\epsilon_{k+1}} \norm{T\paren{y_{k+1}}}    
        <   \norm{P'_{F_{k+1}}\paren{T\paren{y_{k+1}}}} 
        \le \norm{T\paren{y_{k+1}}}$. 
       Thus $\set{E_i}$, $\set{y_i}$, and $\set{F_i}$ can be found as claimed. 
       
       For $\set{\epsilon_i}$ approaching zero rapidly and $\set{y_i}$ normalized, \hfil\break 
       $\l{r}\sim\span{y_i}_{\l{r}}\sim\span{T\paren{y_i}}_{D_q} \sim \span{P'_{F_i}\paren{T\paren{y_i}}}_{D_q}$ 
       via equivalence of natural bases. 
       \qquad\QED

\xproclaim {Lemma 4.27}. Let $1<q<r<2$. 
                         Then $\sumsum{\l{r}}{\l{q}} \not\injects D_q$. 

\proof Suppose $\sumsum{\l{r}}{\l{q}} \injects D_q$. 
       Let $T:\sumsum{\l{r}}{\l{q}} \to D_q$ be an \hbox{isomorphic} imbedding. 
       Let $\set{e_{i,j}}$ be the standard basis of $\sumsum{\l{r}}{\l{q}}$, 
       where for each $j\in\NN$,
       $\set{e_{i,j}}_{i\in\NN}$ is isometrically equivalent to the standard basis of $\l{r}$. 
       For each $j\in\NN$, \hfil\break 
       let $\ell^r_{(j)}=\span{e_{i,j}:i\in\NN}$, 
       and for a sequence $\set{\epsilon_i^{(j)}}_{i\in\NN}$ of positive scalars,
       choose \hfil\break 
       (as we may by Lemma 4.26) a normalized block basic sequence $\set{y_i^{(j)}}_{i\in\NN}$ in $\ell^r_{(j)}$ 
       \hfil\break 
       and disjoint nonempty finite subsets $F_1^{(j)},F_2^{(j)},\ldots$ of $\NN$ with 
       $\max F_i^{(j)}<\min F_{i'}^{(j)}$ \hfil\break 
       for $i<i'$, such that \hfil\break 
       $\l{r} \sim \ell^r_{(j)}
              \sim \span{y_i^{(j)}:i\in\NN}_{\ell^r_{(j)}}
              \sim \span{T\paren{y_i^{(j)}}:i\in\NN}_{D_q}
              \sim \span{P'_{F_i^{(j)}}\paren{T\paren{y_i^{(j)}}}:i\in\NN}_{D_q}$ 
       via equivalence of natural bases,
       with \hfil\break 
       $\paren{1-\epsilon_i^{(j)}} \norm{T\paren{y_i^{(j)}}} <  \norm{P'_{F_i^{(j)}}\paren{T\paren{y_i^{(j)}}}}
                                                            \le \norm{T\paren{y_i^{(j)}}}$ 
       for each $i\in\NN$. 

       For $\epsilon_i^{(j)}$ approaching zero rapidly
       and for infinite subsets $M_1,M_2,\ldots$ of $\NN$
       chosen so that $\set{F_i^{(j)}}_{i\in M_j,\>j\in\NN}$ is disjoint, \hfil\break 
       $\sumsum{\l{r}}{\l{q}} \sim \span{T\paren{y_i^{(j)}}:i\in M_j,\>j\in\NN}_{D_q}
                              \sim \span{P'_{F_i^{(j)}} \paren{T \paren{y_i^{(j)}}}:i\in M_j,\>j\in\NN}_{D_q}$ 
       via \hfil\break 
       equivalence of natural bases. 
       Hence the standard basis of $\sumsum{\l{r}}{\l{q}}$ is equivalent to
       the span in $\L{q}$ of a sequence of independent random variables, 
       contrary to Lemma 3.7. 
       It follows that $\sumsum{\l{r}}{\l{q}} \not\injects D_q$. 
       \qquad\QED

\xproclaim {Lemma 4.28}. Let $1<q<r<2$. 
                         Then $\sumsum{\l{r}}{\l{q}} \not\injects B_q\oplus D_q$. 

\proof Suppose $\sumsum{\l{r}}{\l{q}} \injects B_q\oplus D_q$. 
       Let $T:\sumsum{\l{r}}{\l{q}} \to B_q\oplus D_q$ be an isomorphic imbedding. 
       Let $Q:B_q\oplus D_q\to B_q\oplus\set{0_{D_q}}$ be the obvious \hbox{projection.} 
       Then $QT:\sumsum{\l{r}}{\l{q}} \to B_q\oplus\set{0_{D_q}}$ is a bounded linear operator. 
       As in the proof of Lemma 2.45,
       there is a subspace $X$ of $\sumsum{\l{r}}{\l{q}}$,
       isometric to $\sumsum{\l{r}}{\l{q}}$, 
       such that $\norm{Q|_{T(X)}}<1$,
       whence $(I-Q)|_{T(X)}$ induces an isomorphic imbedding of $\sumsum{\l{r}}{\l{q}}$ into $D_q$.
       However by Lemma 4.27, no such imbedding exists.
       It follows that \hfil\break 
       $\sumsum{\l{r}}{\l{q}} \not\injects B_q\oplus D_q$. 
       \qquad\QED

\xproclaim {Proposition 4.29}. Let $1<p<\infty$ where $p\ne2$. 
                               Then $\sumsum{X_p}{\l{p}}\not\cinjects B_p\oplus D_p$. 

\proof First let $1<q<2$ and suppose $\sumsum{X_q}{\l{q}}\cinjects B_q\oplus D_q$. For \hfil\break 
       $1<q<r<2$, $\l{r}\injects X_q$ by Lemma 2.35, 
       so \hbox{$\sumsum{\l{r}}{\l{q}} \injects \sumsum{X_q}{\l{q}} \cinjects B_q\oplus D_q$.}
       Hence $\sumsum{\l{r}}{\l{q}} \injects B_q\oplus D_q$, contrary to Lemma 4.28. 
       It follows that \hfil\break 
       $\sumsum{X_q}{\l{q}}\not\cinjects B_q\oplus D_q$. 
       The result now holds for $2<p<\infty$ by duality. 
       \qquad\QED

 Finally, we distinguish $D_p$, $B_p\oplus D_p$, and $\sumsum{X_p}{\l{p}}\oplus D_p$ 
 from each other and from the $\SL{p}$ spaces of Rosenthal.  

\xproclaim {Proposition 4.30}. Let $1<p<\infty$ where $p\ne2$. Then 
           \item{(a)} $D_p \not\cinjects B_p$,
           \item{(b)} $B_p \not\cinjects D_p$, 
           \item{(c)} $B_p \oplus X_p \not\cinjects D_p$, 
           \item{(d)} $B_p \oplus D_p \not\cinjects D_p$, 
           \item{(e)} $\sumsum{X_p}{\l{p}} \not\cinjects D_p$,
           \item{(f)} $D_p \not\cinjects \sumsum{X_p}{\l{p}}$, 
           \item{(g)} $B_p \oplus D_p \not\cinjects \sumsum{X_p}{\l{p}}$,
           \item{(h)} $\sumsum{X_p}{\l{p}} \oplus D_p \not\cinjects \sumsum{X_p}{\l{p}}$, 
           \item{(i)} $D_p \not\cinjects B_p \oplus X_p$, 
           \item{(j)} $B_p \oplus D_p \not\cinjects B_p \oplus X_p$, 
           \item{(k)} $D_p \not\cinjects \lpsumltwo \oplus X_p$,
           \item{(l)} $\sumsum{X_p}{\l{p}} \not\cinjects B_p \oplus D_p$, and
           \item{(m)} $\sumsum{X_p}{\l{p}} \oplus D_p \not\cinjects B_p \oplus D_p$. 

\proof Suppose $2<p<\infty$. 
       \item{(a)} Suppose $D_p \cinjects B_p$.
                  Then $X_p \cinjects D_p \cinjects B_p$ by part (a) of Proposition 4.17,
                  so $X_p \cinjects B_p$, contrary to part (g) of Proposition 2.37. 
       \item{(b)} Part (b) is a restatement of Proposition 4.21. 
       \item{(c)} Part (c) is immediate from part (b). 
       \item{(d)} Part (d) is immediate from part (b). 
       \item{(e)} Suppose $\sumsum{X_p}{\l{p}} \cinjects D_p$.
                  Then $B_p \cinjects \sumsum{X_p}{\l{p}} \cinjects D_p$ by Proposition 2.27,
                  so $B_p \cinjects D_p$, contrary to part (b) above. 
       \item{(f)} Part (f) is a restatement of Proposition 4.25. 
       \item{(g)} Part (g) is immediate from part (f). 
       \item{(h)} Part (h) is immediate from part (f). 
       \item{(i)} Suppose $D_p \cinjects B_p \oplus X_p$. 
                  Then $D_p \cinjects B_p \oplus X_p \cinjects \sumsum{X_p}{\l{p}}$
                  by part (a) of Proposition 2.43,
                  so $D_p \cinjects \sumsum{X_p}{\l{p}}$, contrary to part (f) above. 
       \item{(j)} Part (j) is immediate from part (i). 
       \item{(k)} Suppose $D_p \cinjects \lpsumltwo \oplus X_p$.
                  Then $D_p \cinjects \lpsumltwo \oplus X_p \cinjects B_p \oplus X_p$ by Proposition 2.32,
                  so $D_p \cinjects B_p \oplus X_p$, contrary to part (i) above. 
       \item{(l)} Part (l) is a restatement of Proposition 4.29. 
       \item{(m)} Part (m) is immediate from part (l).

       The result for $1<p<2$ follows by duality. 
       \qquad\QED

\medskip
 Building on diagram \TAG{2.27}{2.27}, for $1<p<\infty$ where $p\ne2$, we have

$$\matrix{B_p&&&&\sumsum{X_p}{\lp}&&&&\cr
          &\cdsa&&\cusa&&\cdsa&&&\cr
          &&B_p \oplus X_p&&&&\sumsum{X_p}{\lp} \oplus D_p&\cra&\L{p}.\cr
          &\cusa&&\cdsa&&\cusa&&&\cr
          \lpsumltwo \oplus X_p&&&&B_p \oplus D_p&&&&\cr
          &\cdsa&&\cusa&&&&&\cr
          &&D_p&&&&&&\cr}\eqno{\TAG{4.10}{4.10}}$$

%
%
%
\chapter{V}
\headspace
\chaptertitle{THE CONSTRUCTION AND ORDINAL INDEX OF BOURGAIN,} 
\chaptertitle{ROSENTHAL, AND SCHECHTMAN} 
\headspace
 Let $1<p<\infty$ and let $B$ and $B_1,B_2,\ldots$ be separable Banach spaces with \hfil\break 
 $B\injects\L{p}$ and $B_i\injects\L{p}$. 
 Bourgain, Rosenthal, and Schechtman \xcite{B-R-S} iterate and \hfil\break 
 intertwine two constructions, 
 a disjoint sum of $B$ with itself and an independent sum of $B_1,B_2,\ldots$, 
 to produce a chain $\set{\Rpa}_{\alpha<\omega_1}$ of separable $\SL{p}$ spaces.
 An ordinal \hbox{index} is introduced which assigns to each separable Banach space $B$ an ordinal number \hfil\break 
 $h_p(B)$. 
 The index $h_p(\phantom0)$ proves to be an isomorphic invariant, 
 and is used to select a subchain $\set{R^p_{\tau(\alpha)}}_{\alpha<\omega_1}$
 of [infinite-dimensional] isomorphically distinct spaces. 
 Thus Bourgain, Rosenthal, and Schechtman show that there are uncountably many \hfil\break 
 separable infinite-dimensional $\SL{p}$ spaces [up to isomorphism]. 
  
\preheadspace
\firsthead{Preliminaries}
\postheadspace

 We let $\omega_1$ denote the first uncountable ordinal,
 and we let $\omega$ denote the first \hfil\break 
 infinite ordinal [except in some contexts where $\omega$ will denote an element of a space $\Omega$].

 A strict partial order on a nonempty set $X$ is a relation $\prec$ on $X$ which is \hfil\break 
 transitive and anti-reflexive. 

 A tree is a nonempty set $T$ with a strict partial order $\prec$ such that for each $x\in T$,
 $\set{y\in T:y\prec x}$ is well-ordered by $\prec$. 
 We say that a tree $\paren{T,\prec}$ is a CFRE (countable finite-ranked elements) tree
 if $T$ is finite or countable,
 and for each $x\in T$, \hfil\break 
 $\set{y\in T:y\prec x}$ is finite. 


 Let $\paren{T,\prec}$ be a tree. 
 A subtree of $T$ is a nonempty subset $S$ of $T$ with partial order $\prec$
 [suitably restricted] such that for each $x\in S$,
 the set $\set{y\in T:y\prec x}$ is \hfil\break 
 contained in $S$. 

 Let $\paren{T,\prec}$ be a tree. 
 A branch of $T$ is a maximal totally ordered subset of $T$. 
 Suppose $\paren{T,\prec}$ is a CFRE tree. 
 We say that $B$ is a finite branch of $T$ if $B$ is of the form
 $\set{y\in T:y\preceq x}$ for some $x\in T$. 
 We call $\set{y\in T:y\preceq x}$ the finite branch of $T$ \hfil\break 
 generated by $x$. 
 Note that a finite branch of $T$ need not be a branch of $T$, 
 although a finite branch of $T$ is a branch of some subtree of $T$. 

\medskip

 Let $\triless$ be a relation on a nonempty set $X$. 

 An infinite $\triless$-chain $x_1\triless x_2\triless\cdots$ in $X$ is a sequence $\set{x_n}_{n\in\NN}$ in $X$ 
 such that $x_n\triless x_{n+1}$ for all $n\in\NN$. 
 A finite $\triless$-chain $x_1\triless\cdots\triless x_N$ in $X$ is a sequence $\set{x_n}_{n=1}^N$ in $X$ 
 such that $x_n\triless x_{n+1}$ for all $1\le n<N$. 
 An $x\in X$ is $\triless$-terminal in $X$ if there is no $y\in X$ such that $x\triless y$. 

 The relation $\triless$ is well-founded in $X$
 if there is no infinite $\triless$-chain $x_1\triless x_2\triless\cdots$ in $X$.  
 Note that if $\triless$ is well-founded, 
 then $\triless$ must be anti-reflexive and there can be no finite $\triless$-chain 
 $x_1\triless\cdots\triless x_N$ with $x_1=x_N$. 

\medskip

 For $n\in\NN$,
 an $n$-string is an $n$-tuple which is not delimited by punctuation. 
 We will identify a $0$-string with the empty set.
 For $n\in\NN\cup\set{0}$,
 let $D_n$ be the set of all $n$-strings of 0's and 1's. 
 Then $D_n=\set{t_1\cdots t_n:t_i\in\set{0,1} \hbox{ for all } 1\le i\le n}$ for $n\in\NN$, 
 and $D_0=\set{\emptyset}$. 
 Fix $n\in\NN\cup\set{0}$.
 Then $D_n$ has cardinality $2^n$. 
 There is a natural \hfil\break 
 identification of $D_n$ with $S_n=\set{0,\ldots,2^n-1}$, 
 namely $t_1\cdots t_n\mapsto\sum_{i=1}^n t_i 2^{n-i}$ for $n\in\NN$,  
 and $\set{\emptyset}\mapsto 0$. 
 Thus for $n\in\NN$,
 $t_1\cdots t_n\in D_n$ is the $n$-place binary expansion [possibly with leading 0's] of some $r\in S_n$.

 Let $n,m\in\NN\cup\set{0}$. 
 Given $t\in D_n$ and $s\in D_m$, 
 let $t\con s$ be the element of $D_{n+m}$
 formed by the concatenation of $t$ and $s$. 

\medskip

 Let $\paren{\Omega,\S{M},\mu}$ and $\paren{\Omega',\S{M}',\mu'}$ be probability spaces, 
 and let $X$ and $X'$ be spaces of measurable functions on $\Omega$ and $\Omega'$, respectively. 
 We say that $X$ and $X'$ are \hfil\break 
 distributionally isomorphic, denoted $X\distiso X'$,
 if there is a linear bijection $T:X\to X'$ such that $\dist(Tx)=\dist(x)$ for all $x\in X$. 

\preheadspace
\firsthead{The Ordinal Index}
\postheadspace

 Before introducing the ordinal index $h_p$, 
 we introduce a general ordinal index $h$ based on essentially the same concept, 
 but applicable to a simpler class of spaces. 

\preheadspace
\secondhead{A General Ordinal Index $h$}
\postheadspace

 Let $\triless$ be a relation on a nonempty set $X$. 

 For each ordinal $\alpha$, we define a subset $H_{\alpha}(\triangleleft)$ of $X$. 
 Let $H_0(\triangleleft)=X$. 
 If $\alpha=\beta+1$ and $H_{\beta}(\triangleleft)$ has been defined, 
 let $H_{\alpha}(\triangleleft) = \set{x\in H_{\beta}(\triangleleft) :
      x\triangleleft y {\rm\ for\ some\ } y\in H_{\beta}(\triangleleft)}$. \hfil\break 
 If $\alpha$ is a limit ordinal and $H_{\beta}(\triangleleft)$ has been defined for all $\beta<\alpha$, let \hfil\break 
 $H_{\alpha}(\triangleleft) = \bigcap_{\beta<\alpha} H_{\beta}(\triangleleft)$. 

 If $\beta<\alpha$, then $H_{\beta}(\triangleleft) \supset H_{\alpha}(\triangleleft)$. 
 The members of the nonincreasing family $\paren{H_{\alpha}(\triangleleft)}$ cannot all be distinct. 
 For suppose the members are distinct. 
 Then there is a family $(x_{\alpha})$ of distinct elements of $X$, 
 with $x_{\alpha} \in H_{\alpha}(\triangleleft) \setminus H_{\alpha+1}(\triangleleft)$.  
 Thus for a sufficiently large ordinal $\Gamma$, 
 $\set{x_{\alpha} : \alpha<\Gamma}$ has cardinality larger than the cardinality of $X$, 
 contrary to $\set{x_{\alpha} : \alpha<\Gamma} \subset X$. 
 Hence there is a least ordinal $\gamma$ such that 
 $H_{\gamma}(\triangleleft) = H_{\gamma+1}(\triangleleft)$.
 Let $h(\triangleleft)$ denote this least ordinal $\gamma$,
 and let $S(\triangleleft)$ denote the stable set $H_{\gamma}(\triangleleft)$. 
 Then the cardinality of $h(\triangleleft)$ is bounded by the cardinality of $X$. 
 Note that if $H_{\gamma}(\triangleleft) = H_{\gamma+1}(\triangleleft)$,
 then $H_{\gamma}(\triangleleft) = H_{\gamma'}(\triangleleft)$ for all $\gamma'>\gamma$.


 Suppose $\triless$ is not well-founded. 
 Then there is an infinite $\triless$-chain $x_1\triless x_2\triless\cdots$ in $X$.
 For such a chain,
 $\set{x_1,x_2,\ldots}\subset H_{\alpha}(\triless)$ for all $\alpha$.
 Thus $\set{x_1,x_2,\ldots}\subset S(\triless)$ and $S(\triless)\ne\emptyset$.  
 For the converse, suppose $S(\triless)\ne\emptyset$. 
 Let $x\in S(\triless)$. 
 Then $x$ is not $\triless$-terminal in $S(\triless)$, 
 so there is some $y\in S(\triless)$ with $x\triless y$. 
 By induction, there is an infinite $\triless$-chain 
 $x_1\triless x_2\triless\cdots$ in $S(\triless)\subset X$. 
 Thus $\triless$ is not well-founded. 
 It follows that $\triless$ is well-founded if and only if $S(\triless)=\emptyset$.  

 Let $\triangleleft$ and $\triangleleft'$ be relations on nonempty sets $X$ and $X'$, respectively. 
 A function $\tau:(X,\triangleleft)\to(X',\triangleleft')$ preserves relations if
 $\tau x\,\triangleleft'\,\tau y$ whenever $x \triangleleft y$. 

 The following lemma \xciteplus{B-R-S}{Lemma 2.4} establishes a property of
 the ordinal index $h$ with respect to relation-preserving maps.

\xproclaim {Lemma 5.1}. Let $\triangleleft$ and $\triangleleft'$ be relations
                        on nonempty sets $X$ and $X'$, respectively. 
                        Suppose $\tau:(X,\triangleleft)\to(X',\triangleleft')$ preserves relations.
                        Then 
                        $\tau\paren{H_{\alpha}(\triangleleft)} \subset H_{\alpha}(\triangleleft')$ for all ordinals $\alpha$. 
                        If in addition $\triangleleft'$ is well-founded, 
                        then $h(\triangleleft) \le h(\triangleleft')$. 

\proof Clearly $\tau\paren{H_0(\triangleleft)} = \tau(X) \subset X' = H_0(\triangleleft')$. 
       Suppose $\alpha=\beta+1$ and $\tau\paren{H_{\beta}(\triangleleft)} \subset H_{\beta}(\triangleleft')$. 
       Then $\tau : H_{\beta}(\triangleleft) \to H_{\beta}(\triangleleft')$ [suitably restricted]. 
       Since $\tau$ preserves relations,
       if $x$ is not $\triangleleft$-terminal in $H_{\beta}(\triangleleft)$,
       then $\tau(x)$ is not $\triangleleft'$-terminal in $H_{\beta}(\triangleleft')$. 
       Hence $\tau\paren{H_{\alpha}(\triangleleft)} \subset H_{\alpha}(\triangleleft')$. 
       Suppose $\alpha$ is a limit ordinal and
       $\tau\paren{H_{\beta}(\triangleleft)} \subset H_{\beta}(\triangleleft')$ for all $\beta<\alpha$. 
       Then 
       $\tau\paren{H_{\alpha}(\triangleleft)} = \tau\paren{\bigcap_{\beta<\alpha} H_{\beta}(\triangleleft)}
        \subset \bigcap_{\beta<\alpha} \tau\paren{H_{\beta}(\triangleleft)}
        \subset \bigcap_{\beta<\alpha} H_{\beta}(\triangleleft')
              = H_{\alpha}(\triangleleft')$. 

       Suppose $\triangleleft'$ is well-founded.
       Let $\gamma=h(\triangleleft)$ and $\gamma'=h(\triangleleft')$. Then \hfil\break 
       $\tau\paren{H_{\gamma'}(\triangleleft)} \subset H_{\gamma'}(\triangleleft') = \emptyset$. 
       Thus $H_{\gamma'}(\triangleleft) = \emptyset$ as well. 
       Hence $\gamma\le\gamma'$ and $h(\triangleleft) \le h(\triangleleft')$. 
       \qquad\QED
       
\preheadspace
\secondhead{Motivation from $\L{p}$}
\postheadspace

 Let $1\le p<\infty$. 
 Let $\set{g_n}_{n\in\NN}$ be the sequence of normalized functions in $\L{p}$ given by 
 $g_1=1_{[0,1]}$, $g_2=2^{{1\over p}} 1_{[0,1/2]}$, $g_3=2^{{1\over p}} 1_{[1/2,1]}$, $\ldots$ , 
 $g_n=2^{{k\over p}} 1_{[r/2^k,(r+1)/2^k]}$, 
 $\ldots$ ,
 where $n=2^k+r$ such that $k\in\NN\cup\set{0}$ and $0\le r< 2^k$. 
 For $n$, $k$, and $r$ as above, 
 $2n=2^{k+1}+2r$ where $0\le 2r< 2^{k+1}$, and 
 $2n+1=2^{k+1}+(2r+1)$ where $0< 2r+1< 2^{k+1}$.
 Thus 
 $g_{2n}=2^{{(k+1)\over p}} 1_{[2r/2^{k+1},(2r+1)/2^{k+1}]}  
        =2^{{(k+1)\over p}} 1_{[r/2^k,(r+1/2)/2^k]}$, \hfil\break 
 $g_{2n+1}=2^{{(k+1)\over p}} 1_{[(2r+1)/2^{k+1},(2r+2)/2^{k+1}]}  
          =2^{{(k+1)\over p}} 1_{[(r+1/2)/2^k,(r+1)/2^k]}$, and \hfil\break 
 $g_n=2^{-{1\over p}} \paren{g_{2n}+g_{2n+1}}$. 
 This reflects the fact that $\xsupp g_n = \xsupp g_{2n} \cup \xsupp g_{2n+1}$ 
 [with the union being essentially disjoint].
 The coefficient $2^{-{1\over p}}$ is simply a \hfil\break 
 normalization factor. 
 Thus the functions $g_1,g_2,\ldots$ can be arranged in a binary tree
 $$\matrix{\hbox{[level 0:]}&   &   &   &g_1   &   &   &   \cr
           \hbox{[level 1:]}&   &g_2&   &      &   &g_3&   \cr
           \hbox{[level 2:]}&g_4&   &g_5&      &g_6&   &g_7\cr
           \vdots                 &   &   &   &\vdots&   &   &   \cr}$$
 according to their supports, 
 where the functions at level $k$ are of the form $g_{2^k+r}$ with $0\le r< 2^k$. 

 Indexing by binary expansions, 
 $g_t=2^{-{1\over p}} \paren{g_{t\con 0}+g_{t\con 1}}$, 
 where $t$ is the binary \hfil\break 
 expansion of $n\in\NN$, 
 and $t\con 0$ and $t\con 1$ are the binary expansions of $2n$ and $2n+1$, \hfil\break 
 respectively. 
 The corresponding tree is 
 $$\matrix{\hbox{[level 0:]}&       &      &       &g_1   &       &      &       \cr
           \hbox{[level 1:]}&       &g_{10}&       &      &       &g_{11}&       \cr
           \hbox{[level 2:]}&g_{100}&      &g_{101}&      &g_{110}&      &g_{111}\cr
           \vdots                 &       &      &       &\vdots&       &      &      ,\cr}$$
 where the functions at level $k$ are of the form $g_{1\con s}$ where $s$ is the $k$-place binary \hfil\break 
 expansion of $r$ for $0\le r< 2^k$. 

 Dropping the superfluous leading 1's and indexing by strings of 0's and 1's, \hfil\break 
 $g_s=2^{-{1\over p}} \paren{g_{s\con 0}+g_{s\con 1}}$, 
 where $s$ is a string of 0's and 1's. 
 The corresponding tree is 
 $$\matrix{\hbox{[level 0:]}&      &     &      &g_{\emptyset}&      &     &      \cr
           \hbox{[level 1:]}&      &g_{0}&      &             &      &g_{1}&      \cr
           \hbox{[level 2:]}&g_{00}&     &g_{01}&             &g_{10}&     &g_{11}\cr
           \vdots                 &      &     &      &\vdots&      &      &     &     ,\cr}$$
 where the functions at level $k$ are indexed by $k$-strings of 0's and 1's.


 Level $k$ itself can be thought of as a $2^k$-tuple of elements of $\L{p}$. 
 Recalling that $D_k$ is the set of all $k$-strings of 0's and 1's, 
 the cardinality of $D_k$ is $2^k$. 
 Thus level $k$ can be thought of as a function from $D_k$ to $\L{p}$, 
 or an element of $\paren{\L{p}}^{D_k}$.
 Letting $u_k$ denote level $k$, 
 $$\matrix{u_0&=&\big(&      &     &      &g_{\emptyset}&      &     &      &\big)\cr
           u_1&=&\big(&      &g_{0}&      &,            &      &g_{1}&      &\big)\cr
           u_2&=&\big(&g_{00}&,    &g_{01}&,            &g_{10}&,    &g_{11}&\big)\cr
           \vdots & &     &      &     &      &\vdots       &      &     &      &    ,\cr}\eqno{\TAG{5.1}{5.1}}$$
 where for each $s\in D_k$,
 $u_k(s) = 2^{-{1\over p}} \paren{u_{k+1}(s\con 0) + u_{k+1}(s\con 1)}$. 
 Moreover, for each $s\in D_k$ and each $d\in\NN$,
 $$u_k(s) = 2^{-{d\over p}} \tsuml_{r\in D_d} u_{k+d}(s\con r).\eqno{\TAG{5.2}{5.2}}$$ 
 Furthermore, for each $k\in\NN\cup\set{0}$ and each $c\in\RR^{D_k}$,
 $$\int\abs{\tsuml_{s\in D_k} c(s) u_k(s) }^p = 
   \tsuml_{s\in D_k} \abs{c(s)}^p \displaystyle\int\abs{u_k(s)}^p = 
   \tsuml_{s\in D_k} \abs{c(s)}^p.\eqno{\TAG{5.3}{5.3}}$$ 

\preheadspace
\secondhead{The Space $\paren{\overline{B}^{\delta},\prec}$}
\postheadspace
 
 For $n\in\NN\cup\set{0}$, 
 recall that $D_n$ is the set of all $n$-strings of 0's and 1's, 
 and there is a natural identification of $D_n$ with $\set{0,\ldots,2^n-1}$,
 namely $t_1\cdots t_n\mapsto\sum_{i=1}^n t_i 2^{n-i}$ for $n\in\NN$, and $\set{\emptyset}\mapsto0$.  
 For a vector space $B$, 
 $B^{D_n}$ is the set of all functions from $D_n$ to $B$, 
 which can be identified with the set of all $2^n$-tuples $\paren{b_0,\ldots,b_{2^n-1}}$ of elements of $B$.
 We identify $B^{D_0}$ with $B$. 

 We do not assign an independent meaning to $\S{D}$,  
 but given a vector space $B$, 
 we let $B^{\SD}$ denote $\bigcup_{n=0}^{\infty} B^{D_n}$. 

 Let $B$ be a vector space. 
 If $u\in B^{\SD}$, then $u\in B^{D_n}$ for a unique $n\in{\Bbb N}\cup\set{0}$, denoted $\abs{u}$. 
 Define $\prec$ on $B^{\SD}$ by $u\prec v$ if $\abs{u}<\abs{v}$ and for $k=\abs{v}-\abs{u}$, \hfil\break 
 $u(t) = 2^{-{k\over p}} \sum_{s\in D_k} v(t\con s)$ for all $t\in D^{\abs{u}}$. 
 Then $\prec$ is a strict partial order. 


\definition Suppose $B$ is a separable Banach space, $1\le p<\infty$, and $0<\delta\le1$. 
            Let $\overline{B}^{\delta}$ be the set of all $u\in B^{\SD}$ such that 
            $$\delta\paren{\tsuml_{t\in D_{\abs{u}}} \abs{c(t)}^p}^{{1\over p}}
              \le \norm{\tsuml_{t\in D_{\abs{u}}} c(t)u(t)}_B
              \le \paren{\tsuml_{t\in D_{\abs{u}}} \abs{c(t)}^p}^{{1\over p}}$$
            for all $c\in\RR^{D_{\abs{u}}}$. 
            Let $\prec$ on $\overline{B}^{\delta}$ be
            the strict partial order $\prec$ on $B^{\SD}$ [suitably restricted]. 

\remark For $B=\L{p}$, equation \TAG{5.2}{5.2} implies that $u_0\prec u_1\prec\cdots$,
        and equation \TAG{5.3}{5.3} implies that 
        $u_k\in\overline{\L{p}}^1$ for all $k\in\NN\cup\set{0}$, 
        whence $\paren{\overline{\L{p}}^1,\prec}$ is not well-founded. 

\preheadspace
\secondhead{A Characterization of $\L{p}\injects B$}
\postheadspace

 The following proposition \xciteplus{B-R-S}{Proposition 2.2} characterizes those spaces $B$ for which $\L{p}\injects B$. 
 Essentially, the issue is whether or not $B$ contains a sequence which simulates the behavior of the sequence 
 $\set{u_k(t)}_{k\ge0,t\in D_k}$ in $\L{p}$. 

\xproclaim {Proposition 5.2}. Let $B$ be a separable Banach space and let $1\le p<\infty$. 
                              Then $\L{p}\injects B$ if and only if there is a $0<\delta\le1$ 
                              such that $\paren{\overline{B}^{\delta},\prec}$ is not well-founded. 

\proof Suppose $\L{p}\injects B$.
       Let $T:\L{p} \to B$ be an isomorphic imbedding with $\norm{T}\le1$,
       and let $0<\delta\le1$ be such that
       $\delta\norm{x}_p\le\norm{T(x)}_B\le\norm{x}_p$ for all $x\in \L{p}$. 
       Let $\tau : \paren{\L{p}}^{\SD} \to B^{\SD}$ be defined by $(\tau u)(t) = T(u(t))$
       for $u\in\paren{\L{p}}^{\SD}$ and $t\in D_{\abs{u}}$. 
       Then $\tau$ preserves order by the linearity of $T$. 

       Let $u\in\overline{\L{p}}^{1}$.
       Then for all $c\in\RR^{D_{\abs{u}}}$, 
       $$\eqalign{&\delta\paren{\tsuml_{t\in D_{\abs{u}}} \abs{c(t)}^p}^{{1\over p}} 
          \cr &= \delta\norm{\tsuml_{t\in D_{\abs{u}}} c(t)u(t)}_p
             \le \norm{T\paren{\tsuml_{t\in D_{\abs{u}}} c(t)u(t)}}_B
             \le \norm{\tsuml_{t\in D_{\abs{u}}} c(t)u(t)}_p
          \cr &= \paren{\tsuml_{t\in D_{\abs{u}}} \abs{c(t)}^p}^{{1\over p}}.}$$
       Since $\norm{\sum_{t\in D_{\abs{u}}} c(t)(\tau u)(t)}_B = \norm{T\paren{\sum_{t\in D_{\abs{u}}} c(t)u(t)}}_B$, 
       it follows that $\tau u\in\overline{B}^{\delta}$. 
       Hence $\tau : \overline{\L{p}}^{1} \to \overline{B}^{\delta}$ [suitably restricted]. 

       As noted in the remark above, 
       there is a sequence $\set{u_k}$ in $\overline{\L{p}}^1$ with $u_0\prec u_1\prec\cdots$. 
       Since $\tau : \overline{\L{p}}^{1} \to \overline{B}^{\delta}$ preserves order, 
       $\tau u_0\prec \tau u_1\prec\cdots$ in $\overline{B}^{\delta}$. 
       Hence $\paren{\overline{B}^{\delta},\prec}$ is not well-founded. 

       For the converse, suppose there is a $0<\delta\le1$
       such that $\paren{\overline{B}^{\delta},\prec}$ is not well-founded. 
       Then there is a sequence $\set{v_k}$ in $\overline{B}^{\delta}$ with $v_0\prec v_1\prec\cdots$. 
       Let $\set{r(k)}$ be the increasing sequence in $\NN\cup\set{0}$ with $r(k)=\abs{v_k}$ for all $k$. 
       For $\set{u_k}$ as in \TAG{5.1}{5.1}, \hfil\break 
       let $\set{\tilde u_k}$ be the subsequence of $\set{u_k}$ such that $\abs{\tilde u_k}=r(k)=\abs{v_k}$ for all $k$. 
       For \hfil\break 
       $k\in\NN\cup\set{0}$, 
       let $X_k = \span{\tilde u_k(t) : t\in D_{r(k)}}_{\L{p}}$, 
       let $B_k = \span{v_k(t) : t\in D_{r(k)}}_B$, and 
       let $T_k: X_k \to B_k$ be defined by
       $$T_k\paren{\tsuml_{t\in D_{r(k)}} c(t) \tilde u_k(t) } = \tsuml_{t\in D_{r(k)}} c(t) v_k(t)$$ 
       for $c\in\RR^{D_{r(k)}}$. 
       Then $T_k$ is well-defined and linear,
       and $T_i = T_j |_{X_i}$ for $i<j$. 
       Since
       $\norm{\sum_{t\in D_{r(k)}} c(t) \tilde u_k(t) }_p = \paren{\sum_{t\in D_{r(k)}} \abs{c(t)}^p }^{{1\over p}}$
       by equation \TAG{5.3}{5.3}, and \hfil\break 
       $\delta\paren{\sum_{t\in D_{r(k)}} \abs{c(t)}^p }^{{1\over p}}
        \le\norm{\sum_{t\in D_{r(k)}} c(t) v_k(t) }_{B}
        \le\paren{\sum_{t\in D_{r(k)}} \abs{c(t)}^p }^{{1\over p}}$, 
       we have
       $$\delta\norm{\tsuml_{t\in D_{r(k)}} c(t) \tilde u_k(t) }_p 
         \le\norm{T_k\paren{\tsuml_{t\in D_{r(k)}} c(t) \tilde u_k(t) } }_{B}
         \le\norm{\tsuml_{t\in D_{r(k)}} c(t) \tilde u_k(t) }_p,$$ 
       whence $\delta\norm{x}_p \le \norm{T_k(x)}_B \le \norm{x}_p$
       for $k\in\NN\cup\set{0}$ and $x\in X_k$. 

       Given $x\in\bigcup_{k=0}^{\infty} X_k$, $x\in X_k$ for some $k\in\NN\cup\set{0}$. 
       Let $\tilde T:\bigcup_{k=0}^{\infty} X_k \to \bigcup_{k=0}^{\infty} B_k$
       be defined by $\tilde T(x) = T_k(x)$ for $x\in X_k$. 
       Then $\delta\norm{x}_p \le \norm{\tilde T(x)}_B \le \norm{x}_p$ for all $x\in\bigcup_{k=0}^{\infty} X_k$.
       Since $\bigcup_{k=0}^{\infty} X_k$ is dense in $\L{p}$, 
       $\tilde T$ extends to an isomorphic imbedding of $\L{p}$ into $B$. 
       \qquad\QED

\preheadspace
\secondhead{The Ordinal Index $h_p(\delta,\pz)$}
\postheadspace

 The ordinal index $h(\triangleleft)$ serves as a model for the ordinal index $h_p(\delta,B)$, 
 for which
 the underlying set is $\overline{B}^{\delta}$. 
 The ordinal index $h_p(B)$ is then derived from the indices $h_p(\delta,B)$. 

\definition Suppose $B$ is a separable Banach space, $1\le p<\infty$, and $0<\delta\le1$.
            Let $H_0^{\delta}(B) = \overline{B}^{\delta}$. 
            If $\alpha=\beta+1$ and $H_{\beta}^{\delta}(B)$ has been defined, let \hfil\break 
            $H_{\alpha}^{\delta}(B) = \set{u\in H_{\beta}^{\delta}(B):u\prec v \hbox{ for some } v\in H_{\beta}^{\delta}(B)}$. 
            If $\alpha$ is a limit ordinal and $H_{\beta}^{\delta}(B)$ has been defined for all $\beta<\alpha$, 
            let $\Had(B) = \bigcap_{\beta<\alpha} H_{\beta}^{\delta}(B)$. 

\definition Suppose $B$ is a separable Banach space, $1\le p<\infty$, and $0<\delta\le1$.
            Let $h_p(\delta,B)$ be the least ordinal $\alpha$ such that $\Had(B)=H_{\alpha+1}^{\delta}(B)$. 

 The following proposition \xciteplus{B-R-S}{Proposition 2.3} leads to one half of the \hfil\break 
 characterization contained in Theorem 5.5.

\xproclaim {Proposition 5.3}. Let $B$ be a separable Banach space.
                              Let $1\le p<\infty$ and $0<\delta\le1$. 
                              If $\L{p}\not\injects B$, then $h_p(\delta,B) < \omega_1$.  

\proof Suppose $\L{p} \not\injects B$. 
       Let $B_{\omega}$ be a countable dense subset of $B$.  
       Let $\overline{B_{\omega}}^{\delta,2}$ be the countable set of all $u \in B_{\omega}^{\SD}$ such that  
       $${\delta\over2} \paren{\tsuml_{t\in D_{\abs{u}}} \abs{c(t)}^p}^{{1\over p}}
         \le \norm{\tsuml_{t\in D_{\abs{u}}} c(t)u(t)}_B
         \le 2 \paren{\tsuml_{t\in D_{\abs{u}}} \abs{c(t)}^p}^{{1\over p}}$$
       for all $c\in\RR^{D_{\abs{u}}}$. 
       Let $\triangleleft$ be the relation on $\overline{B_{\omega}}^{\delta,2}$ defined by $u \triangleleft v$ if
       (a) $\abs{u}<\abs{v}$ and
       (b) for $k=\abs{v}-\abs{u}$ and for $\delta_{\ell} = \delta 4^{-(\ell+1)}$,
       $\norm{u(t) - 2^{-{k\over p}} \sum_{s\in D_k} v(t\con s) }_B \le \delta_{\abs{u}}$ 
       for all $t\in D_{\abs{u}}$. 

       We will show that $\triangleleft$ is well-founded and there is a relation-preserving map \hfil\break 
       $\tau : \paren{\overline{B}^{\delta},\prec} \to \paren{\overline{B_{\omega}}^{\delta,2},\triangleleft}$. 
       It will follow by Lemma 5.1 that $h_p(\delta,B) \le h(\triangleleft) < \omega_1$. 
 
       First we show that $\triangleleft$ is well-founded. 
       Suppose $\triangleleft$ is not well-founded. Let \hfil\break 
       $u_1\triangleleft u_2\triangleleft\cdots$ be an infinite $\triangleleft$-chain in $\overline{B_{\omega}}^{\delta,2}$. 
       We will show that there is a corresponding infinite $\prec$-chain
       $\bar u_1\prec \bar u_2\prec\cdots$ in $\overline{B}^{\delta}$, 
       whence $\L{p} \injects B$ by Proposition 5.2, contrary to hypothesis.
       It will follow that $\triangleleft$ is well-founded. 


       Given $i,j\in\NN$ with $i<j$, let $\Delta(i,j) = \abs{u_j}-\abs{u_i}$.
       Fix $i\in\NN$.
       For $i<j\in\NN$ and $t\in D_{\abs{u_i}}$,
       let $\tilde u{}^{(i)}_j (t) = 2^{-{\Delta(i,j)\over p}} \sum_{s\in D_{\Delta(i,j)}} u_j (t\con s)$. 
       Then $\tilde u{}^{(i)}_j \prec u_j$.
       For $t\in D_{\abs{u_i}}$,
       $$\eqalign{&
         \norm{\tilde u{}^{(i)}_j (t) - \tilde u{}^{(i)}_{j+1} (t)}_B
       \cr &=
         \norm{2^{-{\Delta(i,j)\over p}} \tsuml_{s\in D_{\Delta(i,j)}} u_j (t\con s) 
               \,-\,
               2^{-{\Delta(i,j+1)\over p}} \tsuml_{x\in D_{\Delta(i,j+1)}} u_{j+1} (t\con x) }_B
       \cr &=
         \norm{2^{-{\Delta(i,j)\over p}} \tsuml_{s\in D_{\Delta(i,j)}} u_j (t\con s) 
               \,-\,
               2^{-{\Delta(i,j)+\Delta(j,j+1)\over p}}
               \tsuml_{s\in D_{\Delta(i,j)}} \, \tsuml_{r\in D_{\Delta(j,j+1)}}
               u_{j+1} (t\con s\con r) }_B
       \cr &\le
               2^{-{\Delta(i,j)\over p}} \tsuml_{s\in D_{\Delta(i,j)}}
         \norm{u_j (t\con s) \,-\, 2^{-{\Delta(j,j+1)\over p}} \tsuml_{r\in D_{\Delta(j,j+1)}} u_{j+1} (t\con s\con r) }_B
       \cr &\le
               2^{-{\Delta(i,j)\over p}} \cdot 2^{\Delta(i,j)} \cdot \delta_{\abs{u_j}}
       \cr &=
               2^{\Delta(i,j){p-1\over p}} \cdot \delta_{\abs{u_j}}
       \cr &<
               2^{\abs{u_j}} \cdot \delta_{\abs{u_j}}.}$$
       Hence for $i<j<k\in\NN$ and $t\in D_{\abs{u_i}}$,
       $$\eqalign{\norm{\tilde u{}^{(i)}_j (t) - \tilde u{}^{(i)}_{j+k} (t)}_B
                  &\le
                  \tsuml_{n=j}^{j+k-1} \norm{\tilde u{}^{(i)}_n (t) - \tilde u{}^{(i)}_{n+1} (t)}_B
              \cr &<
                  \tsuml_{n=j}^{j+k-1} 2^{\abs{u_n}} \cdot \delta_{\abs{u_n}}
              \cr &<
                  \tsuml_{n=j}^{\infty} 2^{\abs{u_n}+1} \cdot \delta 4^{-(\abs{u_n}+1)}
              \cr &=
                  \delta \tsuml_{n=j}^{\infty} 2^{-(\abs{u_n}+1)} 
              \cr &\le
                  \delta \tsuml_{n=j}^{\infty} 2^{-n}
              \cr &=
                  \delta 2^{1-j}.}$$
       Now $\xlim_{j\to\infty} \delta 2^{1-j} = 0$,
       so $\set{\tilde u{}^{(i)}_j(t)}_{j=i+1}^{\infty}$ is Cauchy. 
       Let $\overline{u}_i(t) = \lim_{j\to\infty} \tilde u{}^{(i)}_j(t)$. \hfil\break 
       Releasing $i$ as a free variable,
       $\overline{u}_i(t)$ is defined for all $i\in\NN$ and all $t\in D_{\abs{u_i}}$. 


       Fix $i,j\in\NN$ with $i<j$. 
       Then for $t\in D_{\abs{u_i}}$,
       $$\eqalign{\overline{u}_i(t) = \lim_{k\to\infty} \tilde u{}^{(i)}_k(t)
                  &=
                  \lim_{k\to\infty} 2^{-{\Delta(i,k)\over p}} \tsuml_{x\in D_{\Delta(i,k)}} u_k (t\con x) 
              \cr &=
                  \lim_{k\to\infty} 2^{-{\Delta(i,j)+\Delta(j,k)\over p}}
                  \tsuml_{s\in D_{\Delta(i,j)}} \, \tsuml_{r\in D_{\Delta(j,k)}}
                  u_k (t\con s\con r) 
              \cr &=
                  2^{-{\Delta(i,j)\over p}} \tsuml_{s\in D_{\Delta(i,j)}} 
                  \lim_{k\to\infty} 2^{-{\Delta(j,k)\over p}} \tsuml_{r\in D_{\Delta(j,k)}} u_k (t\con s\con r) 
              \cr &=
                  2^{-{\Delta(i,j)\over p}} \tsuml_{s\in D_{\Delta(i,j)}} 
                  \lim_{k\to\infty} \tilde u{}^{(j)}_k (t\con s) 
              \cr &=
                  2^{-{\Delta(i,j)\over p}} \tsuml_{s\in D_{\Delta(i,j)}} \overline{u}_j (t\con s).}$$
       Hence $\overline{u}_i \prec \overline{u}_j$. 
       More generally, $\overline{u}_1\prec\overline{u}_2\prec\cdots$. 
       As noted previously, it follows that $\L{p} \injects B$, 
       contrary to hypothesis, so $\triangleleft$ is well-founded. 

       We next show that there is a relation-preserving map 
       $\tau : \paren{\overline{B}^{\delta},\prec} \to \paren{\overline{B_{\omega}}^{\delta,2},\triangleleft}$. 
       Let $u\in\overline{B}^{\delta}$. 
       For each $t\in D_{\abs{u}}$, choose $v(t)\in B_{\omega}$ such that $\norm{u(t)-v(t)}_B\le\epsilon_{\abs{u}}$, 
       where $\epsilon_{\ell} = \delta 8^{-(\ell+1)}$ for $\ell\in\NN$.
       Let $\tau u = v$. 

       First we show that $\tau u \in \overline{B_{\omega}}^{\delta,2}$. 
       Note that $2^{\ell} \cdot \epsilon_{\ell} = 2^{\ell} 8^{-(\ell+1)} \delta < {\delta\over2} < 1$. 
       Thus for $t\in D_{\abs{u}}$ and $c\in\RR^{D_{\abs{u}}}$, 
       $$\eqalign{\norm{\tsuml_{t\in D_{\abs{u}}} c(t)v(t) }_B
                 &=
                  \norm{\tsuml_{t\in D_{\abs{u}}} c(t)u(t) + \tsuml_{t\in D_{\abs{u}}} c(t) \paren{v(t)-u(t)} }_B
             \cr &\le
                  \norm{\tsuml_{t\in D_{\abs{u}}} c(t)u(t) }_B + 
                  \tsuml_{t\in D_{\abs{u}}} \abs{c(t)} \cdot \epsilon_{\abs{u}}
             \cr &\le
                  \paren{\tsuml_{t\in D_{\abs{u}}} \abs{c(t)}^p }^{{1\over p}} + 
                  2^{\abs{u}} \cdot \epsilon_{\abs{u}} \cdot \paren{\tsuml_{t\in D_{\abs{u}}} \abs{c(t)}^p }^{{1\over p}} 
             \cr &=
                  \paren{\tsuml_{t\in D_{\abs{u}}} \abs{c(t)}^p }^{{1\over p}}  
                  \paren{1 + 2^{\abs{u}} \cdot \epsilon_{\abs{u}} }
             \cr &\le
                  2 \paren{\tsuml_{t\in D_{\abs{u}}} \abs{c(t)}^p }^{{1\over p}} }$$  
       and 
       $$\eqalign{\norm{\tsuml_{t\in D_{\abs{u}}} c(t)v(t) }_B
                 &=
                  \norm{\tsuml_{t\in D_{\abs{u}}} c(t)u(t) - \tsuml_{t\in D_{\abs{u}}} c(t) \paren{u(t)-v(t)} }_B
             \cr &\ge
                  \norm{\tsuml_{t\in D_{\abs{u}}} c(t)u(t) }_B - 
                  \tsuml_{t\in D_{\abs{u}}} \abs{c(t)} \cdot \epsilon_{\abs{u}}
             \cr &\ge
                  \delta \paren{\tsuml_{t\in D_{\abs{u}}} \abs{c(t)}^p }^{{1\over p}} - 
                  2^{\abs{u}} \cdot \epsilon_{\abs{u}} \cdot \paren{\tsuml_{t\in D_{\abs{u}}} \abs{c(t)}^p }^{{1\over p}} 
             \cr &=
                  \paren{\tsuml_{t\in D_{\abs{u}}} \abs{c(t)}^p }^{{1\over p}}  
                  \paren{\delta - 2^{\abs{u}} \cdot \epsilon_{\abs{u}} }
             \cr &\ge
                  {\delta\over2} \paren{\tsuml_{t\in D_{\abs{u}}} \abs{c(t)}^p }^{{1\over p}}.}$$  
       Hence $\tau u = v \in \overline{B_{\omega}}^{\delta,2}$. 

       We next show that $\tau$ preserves relations. 
       Suppose $u,v\in\overline B^{\delta}$ with $u\prec v$.
       Let $k=\abs{v}-\abs{u}$. 
       Then for all $t\in D_{\abs{u}}$, 
       $$\eqalign{&\norm{\tau u(t) - 2^{-{k\over p}} \tsuml_{s\in D_k} \tau v(t\con s) }_B
              \cr &\le 
                   \norm{\tau u(t) - u(t)}_B
                 + \norm{u(t) - 2^{-{k\over p}} \tsuml_{s\in D_k} v(t\con s) }_B
                 + 2^{-{k\over p}} \tsuml_{s\in D_k} \norm{v(t\con s) - \tau v(t\con s) }_B
               \cr &\le 
                   \epsilon_{\abs{u}} + 0 + 2^{-{k\over p}} \cdot 2^k \cdot \epsilon_{\abs{v}} 
               \cr &< 
                   \epsilon_{\abs{u}} + 2^k \cdot \epsilon_{\abs{u}+k} 
               \cr &= 
                   \delta \paren{{1\over8^{\abs{u}+1}} + {2^k\over8^{\abs{u}+k+1}} }
                 < \delta {2\over8^{\abs{u}+1}}
                 < {\delta\over4^{\abs{u}+1}}
                 = \delta_{\abs{u}}
                 = \delta_{\abs{\tau u}}.}$$
       Hence $\tau u \triangleleft \tau v$ and $\tau$ preserves relations. 
       As noted previously, since $\triangleleft$ is well-founded, \hfil\break 
       it follows that $h_p(\delta,B) \le h(\triangleleft) < \omega_1$. 
       \qquad\QED

       \medskip

 The following lemma \xcite{B-R-S} provides useful information
 about the behavior of $h_p(\delta,B)$ as a function of $\delta$. 

\xproclaim {Lemma 5.4}. Let $B$ be a separable Banach space and let $1\le p<\infty$. 
                        Suppose $0<\delta_1<\delta_2\le1$.
                        Then $H_{\alpha}^{\delta_1}(B) \supset H_{\alpha}^{\delta_2}(B)$ for each ordinal $\alpha$.
                        If in addition $\L{p} \not\injects B$,
                        then $h_p(\delta_1,B) \ge h_p(\delta_2,B)$,
                        whence $h_p(\delta,B)$ is a nonincreasing function of $\delta$. 


\proof Let $0<\delta_1<\delta_2\le1$. 
       Then $H_0^{\delta_1}(B) = \overline{B}^{\delta_1} \supset \overline{B}^{\delta_2} = H_0^{\delta_2}(B)$. 
       Suppose $\alpha = \beta+1$ and $H_{\beta}^{\delta_1}(B) \supset H_{\beta}^{\delta_2}(B)$. 
       If $x\in H_{\alpha}^{\delta_2}(B)$,
       then $x$ is nonmaximal in $H_{\beta}^{\delta_2}(B)$, 
       so $x$ is nonmaximal in $H_{\beta}^{\delta_1}(B)$, 
       whence $x\in H_{\alpha}^{\delta_1}(B)$.
       Hence $H_{\alpha}^{\delta_1}(B) \supset H_{\alpha}^{\delta_2}(B)$. \hfil\break 
       Suppose $\alpha$ is a limit ordinal and
       $H_{\beta}^{\delta_1}(B) \supset H_{\beta}^{\delta_2}(B)$ for all $\beta<\alpha$. 
       Then \hfil\break 
            $H_{\alpha}^{\delta_1}(B) = \bigcap_{\beta<\alpha} H_{\beta}^{\delta_1}(B)
             \supset
             \bigcap_{\beta<\alpha} H_{\beta}^{\delta_2}(B) = H_{\alpha}^{\delta_2}(B)$. 
       It follows that for each ordinal $\alpha$,
       $H_{\alpha}^{\delta_1}(B) \supset H_{\alpha}^{\delta_2}(B)$. 

       Suppose $\L{p} \not\injects B$.
       Then by Proposition 5.2,
       $\paren{\overline{B}^{\delta},\prec}$ is well-founded for all \hfil\break 
       $0<\delta\le1$,
       so $H_{\gamma_i}^{\delta_i}(B) = \emptyset$ for $\gamma_i = h_p\paren{\delta_i,B}$.
       Thus $H_{\gamma_2}^{\delta_1}(B) \supset H_{\gamma_2}^{\delta_2}(B) = \emptyset$, 
       so \hbox{$\gamma_1 \ge \gamma_2$} and $h_p(\delta_1,B) \ge h_p(\delta_2,B)$.  
       Hence $h_p(\delta,B)$ is a nonincreasing function of $\delta$.
       \qquad\QED

\preheadspace
\secondhead{The Ordinal Index $h_p$}
\postheadspace

 Finally we define the ordinal index $h_p$. 

\definition Suppose $B$ is a separable Banach space and $1\le p<\infty$.
            If $\L{p}\not\injects B$, let $h_p(B) = \sup_{0<\delta\le1} h_p(\delta,B)$. 
            If $\L{p}\injects B$, let $h_p(B) = \omega_1$. 

 We presently show that if $\L{p}\not\injects B$,
 then $\set{h_p(\delta,B):0<\delta\le1}$ is bounded,
 whence $h_p(B)$ is well-defined. 
 Note that the hypothesis $\L{p}\not\injects B$ is equivalent to \hfil\break 
 asserting that for each $0<\delta\le1$,
 there is an ordinal $\alpha$ such that $\Had(B)=\emptyset$. 

 The following two results \xciteplus{B-R-S}{Theorem 2.1} 
 establish a countability criterion for $h_p$ and the monotonicity of $h_p$. 

\xproclaim {Theorem 5.5}. Let $B$ be a separable Banach space and let $1\le p<\infty$.
                          Then $h_p(B)\le\omega_1$, with $h_p(B)<\omega_1$ if and only if $\L{p}\not\injects B$. 

\proof If $\L{p}\injects B$, then $h_p(B) = \omega_1$. 
       Henceforth suppose $\L{p}\not\injects B$.
       Now $h_p(\delta,B)$ is a nonincreasing function of $\delta$ by Lemma 5.4,
       and $h_p(\delta,B)<\omega_1$ for all $0<\delta\le1$ by Proposition 5.3. 
       Hence $h_p(B) = \sup_{0<\delta\le1} h_p(\delta,B) = \sup_{n\in\NN} h_p\paren{{1\over n},B} <\omega_1$.
       \qquad\QED


\xproclaim {Theorem 5.6}. Let $X$ and $Y$ be separable Banach spaces and let $1\le p<\infty$.
                          If $X \injects Y$, then $h_p(X) \le h_p(Y)$. 

\proof Suppose $X\injects Y$. 
       If $\L{p}\injects Y$, then $h_p(X) \le \omega_1 = h_p(Y)$ by Theorem 5.5. 
       Henceforth suppose $\L{p}\not\injects Y$, whence $\L{p}\not\injects X$. 
       Then by Proposition 5.2, $\paren{\overline{Y}^{\gamma},\prec}$ is well-founded for each $0<\gamma\le1$.

       Let $T:X \to Y$ be an isomorphic imbedding with $\norm{T}\le1$,
       and let $0<\eta\le1$ be such that for each $x\in X$, 
       $\eta\norm{x}_X\le\norm{T(x)}_Y\le\norm{x}_X$. 
       Let $\tau : X^{\SD} \to Y^{\SD}$ be defined by $(\tau u)(t) = T(u(t))$ for $u\in X^{\SD}$ and $t\in D_{\abs{u}}$. 
       Then $\tau$ preserves order by the linearity of $T$. 

       Fix $0<\delta\le1$ and let $u\in\overline{X}^{\delta}$.
       Then for all $c\in\RR^{D_{\abs{u}}}$, 
       $$\eqalign{&\eta\delta\paren{\tsuml_{t\in D_{\abs{u}}} \abs{c(t)}^p}^{{1\over p}}
             \cr &\le\eta\norm{\tsuml_{t\in D_{\abs{u}}} c(t)u(t)}_X
                  \le\norm{T\paren{\tsuml_{t\in D_{\abs{u}}} c(t)u(t)}}_Y
                  \le\norm{\tsuml_{t\in D_{\abs{u}}} c(t)u(t)}_X
             \cr &\le \paren{\tsuml_{t\in D_{\abs{u}}} \abs{c(t)}^p}^{{1\over p}}.}$$
       Since $\norm{\sum_{t\in D_{\abs{u}}} c(t)(\tau u)(t)}_Y = \norm{T\paren{\sum_{t\in D_{\abs{u}}} c(t)u(t)}}_Y$, 
       it follows that $\tau u\in\overline{Y}^{\eta\delta}$. 
       Hence $\tau : \overline{X}^{\delta} \to \overline{Y}^{\eta\delta}$ [suitably restricted]. 
       Since $\tau$ preserves order and $\paren{\overline{Y}^{\eta\delta},\prec}$ is well-founded,
       $h_p(\delta,X) \le h_p(\eta\delta,Y)$ by Lemma 5.1. 
       Releasing $\delta$ as a free variable,
       \hfil\break 
       $h_p(X) = \sup_{0<\delta\le1} h_p(\delta,X) \le \sup_{0<\delta\le1} h_p(\eta\delta,Y)
               = \sup_{0<\gamma\le\eta} h_p(\gamma,Y) = h_p(Y)$, 
       since $h_p(\gamma,Y)$ is a nonincreasing function of $\gamma$ by Lemma 5.4.
       \qquad\QED

 \remark It follows that $h_p(\phantom0)$ is an isomorphic invariant. 

\preheadspace
\firsthead{The Disjoint and Independent Sum Constructions}
\postheadspace

 Let $(\Omega,\mu)$ be a probability space,  
 let $\paren{\Omega^{\NN},\mu^{\NN}}$ be the corresponding product \hfil\break 
 space, and 
 let $(\set{0,1},m)$ be the probability space with $m(0)={1\over2}=m(1)$. Suppose \hfil\break 
 $1\le p<\infty$, and let $B$ and $B_1,B_2,\ldots$ be closed subspaces of $\L{p}(\Omega)$.

 Given $b_0,b_1\in B$, let $b(\omega,\epsilon)$ be the element of $\L{p}(\Omega\times\set{0,1})$ such that \hfil\break 
 $b(\omega,0)=2^{{1\over p}}b_0(\omega)$ and $b(\omega,1)=2^{{1\over p}}b_1(\omega)$
 for all $\omega\in\Omega$.
 Let $b_0\oplus b_1$ denote the element $b(\omega,\epsilon)$ of 
 $\L{p}(\Omega\times\set{0,1})$ corresponding to $b_0,b_1\in B$. 

\definition Let $1\le p<\infty$ and let $B$ be a closed subspace of $\L{p}(\Omega)$.
            Define the $\L{p}$-disjoint sum $\psum{B}$ to be any space of random variables distributionally \hfil\break 
            isomorphic to the subspace $\tilde B$ of $\L{p}(\Omega\times\set{0,1})$ defined by 
            $$\tilde B = \set{b(\omega,\epsilon) \in \L{p}(\Omega\times\set{0,1}) :
                              b(\omega,\epsilon) = b_0\oplus b_1 \hbox{ for some } b_0,b_1 \in B}.$$

 Note that $1_{\Omega}\oplus1_{\Omega}=2^{{1\over p}}\cdot1_{\Omega\times\set{0,1}}$,
 and if $b(\omega,\epsilon)=b_0\oplus b_1$, then
 $$\eqalign{\norm{b_0\oplus b_1}_{\oplus}^p = \norm{b(\omega,\epsilon)}_{\tilde B}^p
                                           &= \int_{\Omega\times\set{0,1}} \abs{b(\omega,\epsilon)}^p
                                       \cr &= \int_{\Omega\times\set{0}} \abs{b(\omega,\epsilon)}^p
                                            + \int_{\Omega\times\set{1}} \abs{b(\omega,\epsilon)}^p
                                       \cr &= {1\over2} \int_{\Omega} 2\abs{b_0(\omega)}^p
                                            + {1\over2} \int_{\Omega} 2\abs{b_1(\omega)}^p
                                       \cr &= \norm{b_0}_B^p + \norm{b_1}_B^p.}$$
 Hence for $b\in B$, $\norm{b\oplus0}_{\oplus} = \norm{b}_B = \norm{0\oplus b}_{\oplus}$.

 Given $i\in\NN$ and $b_i\in B_i$, let $\tilde b_i$ be the element of $\L{p}\paren{\Omega^{\NN}}$ such that \hfil\break 
 $\tilde b_i(\omega_1,\omega_2,\ldots)=b_i(\omega_i)$ for all $\omega_1,\omega_2,\ldots\in\Omega$.

\definition Let $1\le p<\infty$ and let $B_1,B_2,\ldots$ be closed subspaces of $\L{p}(\Omega)$.
            For each $i\in\NN$, let
            $$\tilde B_i = \set{b\in\L{p}\paren{\Omega^{\NN}} : b=\tilde b_i \hbox{ for some } b_i \in B_i}.$$
            Define the $\L{p}$-independent sum $\Ipsum{B_i}{}$ to be any space of random variables \hfil\break 
            distributionally isomorphic to $\span{\tilde B_i : i\in\NN}_{\L{p}\paren{\Omega^{\NN}}}$. 

 Finally, the spaces $\Rpa$ for $0<\alpha<\omega_1$ are defined as disjoint or independent sums, 
 depending on whether $\alpha$ is a successor or limit ordinal, respectively. 


\definition Let $1\le p<\infty$. 
            Let $R^p_0=\span{1}_{\L{p}}$. 
            Suppose $0<\alpha<\omega_1$. 
            If $\alpha=\beta+1$ and $\Rpb$ has been defined, let $\Rpa = \psum{\Rpb}$. 
            If $\alpha$ is a limit ordinal and $\Rpb$ has been defined for all $\beta<\alpha$,
            let $\Rpa = \Ipsum{\Rpb}{\beta<\alpha}$. 

\xremark{1} It is shown in \xciteplus{B-R-S}{Proposition 2.8} that for $1<p<\infty$ and \hfil\break 
            $\alpha<\omega_1$, $\Rpa$ has an unconditional basis. 

\xremark{2} Technically, $\Rpa = \Ipsum{R_{\beta_i}^p}{\beta_i<\alpha}$
            for an enumeration $\set{\beta_i}$ of the ordinals less than $\alpha$, 
            but it is clear that the definition of $\Rpa$ does not depend on the order. 

 The following two results serve as lemmas for the subsequent theorem \xciteplus{B-R-S}{Proposition 2.7},
 which distinguishes $\Rpa$ from $\L{p}$ isomorphically. 
 Proposition 5.7 is \hfil\break 
 a corollary of \xciteplus{J-M-S-T}{Theorem 9.1}. 
 Proposition 5.8 is \xciteplus{B-R-S}{Theorem 1.1}. 

\xproclaim {Proposition 5.7}. Let $1<p<\infty$. 
                              Suppose $X$ is a closed subspace of $\L{p}$ such that $\L{p}\injects X$. 
                              Then $\L{p}\cinjects X$. 

\proof Let $Y$ be a closed subspace of $X$ such that $\L{p}\sim Y\subset\L{p}$. 
       By \xciteplus{J-M-S-T}{Theorem 9.1},
       choose a closed subspace $Z$ of $Y$ such that $\L{p}\sim Z$ where $Z$ is \hfil\break 
       complemented in $\L{p}$. 
       Let $P$ be a projection from $\L{p}$ onto $Z$. 
       Since $P(Z)=Z$ and \hfil\break 
       $Z\subset X\subset\L{p}$,
       the restriction of $P$ to $X$ is a projection from $X$ onto $Z$. Hence \hfil\break 
       $\L{p}\sim Z\cinjects X$. 
       \qquad\QED

\xproclaim {Proposition 5.8}. Let $1<p<\infty$.
                              Let $X$ be a Banach space with an \hfil\break 
                              unconditional Schauder decomposition
                              $\set{X_i}$ such that $\L{p} \cinjects X$. 
                              Then either $\L{p} \cinjects X_i$ for some $i$, 
                              or there is a block basic sequence with respect to $\set{X_i}$
                              equivalent to the Haar basis of $\L{p}$, 
                              with closed linear span complemented in $X$. 

 The proof of Proposition 5.8 consumes \xciteplus{B-R-S}{Section 1}, and will not be \hfil\break 
 presented here.

\xproclaim {Theorem 5.9}. Let $1<p<\infty$ where $p\ne2$, and let $\alpha<\omega_1$. 
                          Then $\L{p}\not\injects\Rpa$. 

\proof Clearly $\L{p} \not\injects \span{1}_{\L{p}} = R_0^p$.  


       Suppose $\alpha=\beta+1$ and $\L{p}\not\injects\Rpb$. 
       Suppose for the moment that $\L{p}\injects\Rpa$. 
       Then $\L{p}\injects\tilde R^p_{\alpha}\subset\L{p}$
       for some $\tilde R^p_{\alpha}$ distributionally isomorphic to $\Rpa$.
       Hence $\L{p}\cinjects \Rpa$ by Proposition 5.7.
       Now $\Rpa = \psum{\Rpb}$,
       so $\L{p} \cinjects\psum{\Rpb}$,
       whence $\L{p}\cinjects\Rpb$ by Proposition 5.8,
       contrary to the inductive hypothesis. 
       Hence $\L{p}\not\injects\Rpa$. 

       Suppose $\alpha$ is a limit ordinal and $\L{p}\not\injects\Rpb$ for all $\beta<\alpha$. 
       Suppose for the moment that $\L{p}\injects\Rpa$. 
       Then $\L{p} \cinjects \Rpa$ as above. 
       Let $\set{\beta_i}_{i=0}^{\infty}$ be an enumeration of the \hfil\break 
       ordinals less than $\alpha$, with $\beta_0=0$. 
       Let $X_0=R^p_{\beta_0}=R^p_0=\span{1}_{\L{p}}$, and for $i\ge1$, 
       let $X_i = \paren{R^p_{\beta_i}}_0$, the space of mean zero functions in $R^p_{\beta_i}$. 
       Now $\L{p} \cinjects \Ipsum{X_i}{i\ge0}$, 
       since $\Rpa = \Ipsum{\Rpb}{\beta<\alpha} = \Ipsum{X_i}{i\ge0}$, 
       but $\L{p} \not\injects X_i$ for $i\ge0$. Let \hfil\break 
       $\tilde X_i = \set{x\in\L{p}\paren{[0,1]^{\NN}} :
                          x=\tilde x_i \hbox{ for some } x_i\in X_i}$,
       with notation as in the definition of $\Ipsum{B_i}{}$. 
       Then by Proposition 5.8,
       there is a block basic sequence $\set{z_i}_{i\ge0}$ with respect to $\set{\tilde X_i}_{i\ge0}$
       [with at most $z_0$ not mean zero] equivalent to the Haar basis of $\L{p}$.
       Hence $\L{p} \sim \span{z_i:i\ge0}_{\L{p}\paren{[0,1]^{\NN}}} \sim \span{z_i:i\ge1}_{\L{p}\paren{[0,1]^{\NN}}}$.
       Since $\set{z_i}_{i\ge1}$ is a sequence of 
       independent mean zero random variables in $\L{p}\paren{[0,1]^{\NN}}$, 
       $\span{z_i:i\ge1}_{\L{p}\paren{[0,1]^{\NN}}} \injects X_p$ 
       [by Corollary 2.3, Proposition 2.1, Theorem 2.12, and part (b) of Proposition 2.24 for $2<p<\infty$, and
        by \xciteplus{RII}{Corollary 4.3} for $1<p<2$].

       Hence $\L{p} \injects X_p$,
       directly contrary to part (g) of Proposition 2.24 for $2<p<\infty$, 
       and indirectly contrary to the same result for $1<p<2$ as we presently show.
       Thus it will follow that $\L{p}\not\injects\Rpa$. 

       Suppose $\L{s}\injects X_s$ for $1<s<2$, and let $r$ be the conjugate index of $s$. 
       Then $\L{s}\injects X_s\subset\L{s}$, whence $\L{s}\cinjects X_s$ by Proposition 5.7.
       Hence $\L{r}\cinjects X_r$, contrary to part (g) of Proposition 2.24. 
       \qquad\QED

\remark As shown in \xcite{B-R-S}, Theorem 5.9 is true for $p=1$ as well, but the proof is not identical.


\preheadspace
\secondhead{The Interaction of the Constructions and the Ordinal Index}
\postheadspace

 The disjoint and independent sum constructions are designed to force the ordinal index $h_p\paren{\Rpa}$ to increase 
 [not necessarily strictly, but in the sense that the set \hfil\break 
 $\set{h_p\paren{\Rpa}:\alpha<\omega_1}$ has no maximum]. 
 The first results in this direction are the \hfil\break 
 following proposition \xciteplus{B-R-S}{Lemma 2.5} and corollary \xcite{B-R-S}.

\xproclaim {Proposition 5.10}. Let $1\le p<\infty$, $0<\delta\le1$, and $\alpha<\omega_1$. 
                               Suppose $B$ is a closed subspace of $\L{p}$. 
                               Then for each $e\in\Had(B)$, there is some $\bar e \in H_{\alpha+1}^{\delta} \psum{B}$. 

\proof Suppose $e=x_0\in B^{D_0}$.
       Let $\tau e = \paren{x_0 \oplus 0, 0 \oplus x_0} \in (B \oplus B)_p^{D_1}$. 
       Then $\tau e (0) = x_0\oplus0 \in \psum{B}$ and
                            $\tau e (1) = 0\oplus x_0 \in \psum{B}$.
       Let
       $$\bar e = {x_0 \oplus x_0 \over 2^{{1\over p}} }.\eqno{\TAG{5.4}{5.4}}$$  
       Then $\bar e \in \paren{B\oplus B}_p^{D_0}$ 
       and $\bar e = 2^{-{1\over p}} \paren{\tau e (0) + \tau e (1)}$.
       Hence $\bar e \prec \tau e$. 

       Let $k\in\Bbb N$ and suppose $e=\paren{x_0,\ldots,x_{2^k-1}} \in B^{D_k}$.
       Then $e(t)\in B$ for $t\in D_k$. \hfil\break 
       Let $\tau e=\paren{x_0\oplus0,\ldots,x_{2^k-1}\oplus0,0\oplus x_0,\ldots,0\oplus x_{2^k-1}}\in\paren{B\oplus B}_p^{D_{k+1}}$.
       Then for \hfil\break 
       $t\in D_k$, $\tau e (0\con t) = e(t)\oplus0 \in \psum{B}$ and
                   $\tau e (1\con t) = 0\oplus e(t) \in \psum{B}$.
       Let
       $$\bar e= \paren{ {x_0+x_1\over2^{{1\over p}}}\oplus0, \ldots, {x_{2^k-2}+x_{2^k-1}\over2^{{1\over p}}}\oplus0,
                         0\oplus{x_0+x_1\over2^{{1\over p}}}, \ldots, 0\oplus{x_{2^k-2}+x_{2^k-1}\over2^{{1\over p}}} }.$$
       Then $\bar e \in \paren{B\oplus B}_p^{D_k}$ 
       and $\bar e (t) =  2^{-{1\over p}} \paren{\tau e (t\con 0) + \tau e (t\con 1)}$ for $t\in D_k$.
       Hence $\bar e \prec \tau e$. 

       We will show that if $e\in\Had(B)$, then $\tau e \in \Had\psum{B}$. 
       Since $\bar e \prec \tau e$, it will follow that $\bar e$ is a nonmaximal element of $\Had\psum{B}$, 
       so $\bar e \in H_{\alpha+1}^{\delta} \psum{B}$. 

       First we show that $\tau$ preserves order.
       Suppose $e \prec d$. Without loss of generality suppose $\abs{d}-\abs{e}=1$. 
       Then for $t\in D_{\abs{e}}$, $e(t)=2^{-{1\over p}}(d(t\con0)+d(t\con1))$. 
       Thus for $t\in D_{\abs{e}}$ 
       $$\tau e (0\con t) = e(t)\oplus 0 =
         {(d(t\con 0)\oplus 0) + (d(t\con 1)\oplus 0) \over 2^{{1\over p}} } =  
         {\tau d(0\con t\con 0) + \tau d(0\con t\con 1) \over 2^{{1\over p}} }$$   
       and
       $$\tau e (1\con t) = 0\oplus e(t) =
         {(0\oplus d(t\con 0)) + (0\oplus d(t\con 1)) \over 2^{{1\over p}} } =  
         {\tau d(1\con t\con 0) + \tau d(1\con t\con 1) \over 2^{{1\over p}} }.$$   
       Hence for $s=(0\con t)$ or $s=(1\con t)$, 
       $\tau e (s) = 2^{-{1\over p}} \paren{\tau d(s\con 0) + \tau d(s\con 1)}$, 
       so $\tau e \prec \tau d$ \hfil\break 
       and $\tau$ preserves order.   

       We now show by induction on $\alpha$ that if $e\in\Had(B)$, then $\tau e \in \Had\psum{B}$. 

       Suppose $\alpha=0$ and let $e\in H_0^{\delta}(B)=\overline{B}^{\delta}$. 
       Then for $k=\abs{e}$ and $c\in{\Bbb R}^{D_{k+1}}$, 
       $$\eqalign{\norm{\tsuml_{\textstyle{t\in D_k \atop b\in\set{0,1}}} c(b\con t) \tau e (b\con t) }^p_{\oplus}
                 &=
                  \norm{\paren{\tsuml_{t\in D_k} c(0\con t) \tau e (0\con t) } + 
                        \paren{\tsuml_{t\in D_k} c(1\con t) \tau e (1\con t) } }^p_{\oplus}
             \cr &=
                  \norm{\paren{\tsuml_{t\in D_k} c(0\con t) (e(t) \oplus 0) } + 
                        \paren{\tsuml_{t\in D_k} c(1\con t) (0 \oplus e(t)) } }^p_{\oplus}
             \cr &=
                  \norm{\paren{\tsuml_{t\in D_k} c(0\con t) e(t)} \oplus
                        \paren{\tsuml_{t\in D_k} c(1\con t) e(t)} }^p_{\oplus}
             \cr &=
                  \norm{\tsuml_{t\in D_k} c(0\con t) e(t)}_B^p +
                  \norm{\tsuml_{t\in D_k} c(1\con t) e(t)}_B^p
             \cr &\within{1}{\delta^{-p}}
                  \tsuml_{t\in D_k} \abs{c(0\con t)}^p +
                  \tsuml_{t\in D_k} \abs{c(1\con t)}^p
             \cr &=
                  \tsuml_{\textstyle{t\in D_k \atop b\in\set{0,1}}} \abs{c(b\con t)}^p.}$$
       Hence $\tau e \in \overline{\psum{B}}^{\delta} = H_0^{\delta}\psum{B}$. 
  
       Suppose $\alpha=\beta+1$,
       where if $d\in H_{\beta}^{\delta}(B)$, then $\tau d\in H_{\beta}^{\delta}\psum{B}$. 
       Let $e\in\Had(B)$. 
       Then $e\in H_{\beta}^{\delta}(B)$,
       there is some $d\in H_{\beta}^{\delta}(B)$ such that $e \prec d$,
       and $\tau d \in H_{\beta}^{\delta}\psum{B}$. 
       Since $\tau$ preserves order, $\tau e \prec \tau d$. 
       Thus $\tau e$ is a nonmaximal element of $H_{\beta}^{\delta}\psum{B}$,
       whence $\tau e \in \Had\psum{B}$. 

       Suppose $\alpha$ is a limit ordinal, where for each $\beta<\alpha$,
       if $d\in H_{\beta}^{\delta}(B)$, then \hfil\break 
       $\tau d\in H_{\beta}^{\delta}\psum{B}$. 
       Let $e\in\Had(B)$. 
       Then $e\in H_{\beta}^{\delta}(B)$ for all $\beta<\alpha$, and \hfil\break 
       $\tau e\in H_{\beta}^{\delta}\psum{B}$ for all $\beta<\alpha$, 
       whence $\tau e \in \Had\psum{B}$. 


       Hence if $e\in\Had(B)$, then $\tau e \in \Had\psum{B}$. 
       Now as previously noted, \hfil\break 
       if $e\in\Had(B)$, then $\bar e \prec \tau e \in \Had\psum{B}$, 
       so $\bar e \in H_{\alpha+1}^{\delta} \psum{B}$. 
       \qquad\QED

\xproclaim {Corollary 5.11}. Let $1\le p<\infty$ and $\alpha<\omega_1$. 
                             Suppose $B$ is a closed subspace of $\L{p}$ such that $\L{p}\not\injects B$. 
                             If $h_p(B) > \alpha$, then $h_p \psum{B} > \alpha+1$.  

\proof Suppose $h_p(B) > \alpha$. 
       Then $h_p(\delta,B) > \alpha$ for some $0<\delta\le1$.  
       Thus $\Had(B)\ne\emptyset$,  
       so $H_{\alpha+1}^{\delta}\psum{B}\ne\emptyset$ by Proposition 5.10. Hence \hfil\break 
       $h_p\paren{\delta,\psum{B}} > \alpha+1$, 
       so $h_p \psum{B} > \alpha+1$.  
       \qquad\QED

\remark It follows that if $h_p(B)$ is a successor ordinal, then \hfil\break 
        $h_p(B) < h_p\psum{B}$, 
        while if $h_p(B)$ is a limit ordinal, then $h_p(B) \le h_p\psum{B}$. 
        Thus this result is not sufficient to force $h_p\paren{\Rpa}$ to increase. 

 For each ordinal $\alpha<\omega_1$, we define a probability space $\Omega_{\alpha}$. 
 Let $\Omega_0=[0,1]$. 
 If $\alpha=\beta+1$ and $\Omega_{\beta}$ has been defined, let $\Omega_{\alpha}=\Omega_{\beta}\times\set{0,1}$. 
 If $\alpha$ is a limit ordinal and $\Omega_{\beta}$ has been defined for all $\beta<\alpha$,
 let $\Omega_{\alpha}=\prod_{\beta<\alpha}\Omega_{\beta}$. 

 The following theorem \xciteplus{B-R-S}{Theorem 2.6} leads almost immediately to the subsequent corollary
 \xciteplus{B-R-S}{Theorem B(2)}, 
 which is the key to forcing $h_p\paren{\Rpa}$ to increase in the sense mentioned previously. 

\xproclaim {Theorem 5.12}. Let $1\le p<\infty$ and $\alpha<\omega_1$. 
                           Then $1_{\Omega_{\alpha}} \in H^1_{\alpha}\paren{\Rpa}$. 

\proof First we show that $1_{\Omega_{\alpha}} \in \Rpa$. 
       Clearly $1_{\Omega_0} \in\span{1}_{\L{p}}=R^p_0$. 
       Suppose $\alpha=\beta+1$ and $1_{\Omega_{\beta}}\in\Rpb$. 
       Then $1_{\Omega_{\alpha}} = 2^{-{1\over p}}(1_{\Omega_{\beta}}\oplus1_{\Omega_{\beta}}) \in \psum{\Rpb} = \Rpa$. 
       Suppose $\alpha$ is a limit ordinal and $1_{\Omega_{\beta}}\in\Rpb$ for all $\beta<\alpha$. 
       Fix $\beta<\alpha$, so $1_{\Omega_{\beta}}\in\Rpb$. 
       Now $\Rpb$ is distributionally isomorphic to some closed subspace $\tilde R^p_{\beta}$ of $\Rpa$. 
       Let $T : \Rpb \to \tilde R^p_{\beta} \subset \Rpa$ be the distributional isomorphism. 
       Then $T(1_{\Omega_{\beta}}) = 1_{\Omega_{\alpha}} \in \tilde R^p_{\beta} \subset \Rpa$. 
       Hence $1_{\Omega_{\alpha}}\in\Rpa$. 

       We now show that $1_{\Omega_{\alpha}} \in H^1_{\alpha}\paren{\Rpa}$. 
       Clearly $1_{\Omega_0} \in \overline{\span{1}}^1_{\L{p}}= H_0^1\paren{\span{1}_{\L{p}}} = H_0^1\paren{R_0^p}$. 
       Suppose $\alpha=\beta+1$ and $1_{\Omega_{\beta}}\in H^1_{\beta}\paren{\Rpb}$. 
       Then $1_{\Omega_{\beta}}\in\Rpb$,
       so $\bar 1_{\Omega_{\beta}} = 2^{-{1\over p}}(1_{\Omega_{\beta}}\oplus1_{\Omega_{\beta}})$ for
       $\bar 1_{\Omega_{\beta}}$ as in equation \TAG{5.4}{5.4}.
       Hence by Proposition 5.10, 
       $1_{\Omega_{\alpha}} = 2^{-{1\over p}} \paren{1_{\Omega_{\beta}}\oplus1_{\Omega_{\beta}}}
                            = \bar 1_{\Omega_{\beta}}
                          \in H^1_{\alpha} \psum{\Rpb} = H^1_{\alpha} \paren{\Rpa}$. 
       Suppose $\alpha$ is a limit ordinal and $1_{\Omega_{\beta}}\in H^1_{\beta} \paren{\Rpb}$ for all $\beta<\alpha$. 
       Fix $\beta<\alpha$,
       so $1_{\Omega_{\beta}}\in H^1_{\beta} \paren{\Rpb}$.
       Let $T : \Rpb \to \tilde R^p_{\beta} \subset \Rpa$ be as above. Let \hfil\break 
       $\tau : \paren{\Rpb}^{\SD} \to \paren{\Rpa}^{\SD}$ be defined by 
       $(\tau u)(t) = T(u(t))$ for $u\in\paren{\Rpb}^{\SD}$ and $t\in D_{\abs{u}}$. 
       Since $T$ is an isometry, 
       $\tau$ maps $\overline{\Rpb}^1$ into $\overline{\Rpa}^1$. 
       Hence $\tau : \overline{\Rpb}^1 \to \overline{\Rpa}^1$ [suitably \hfil\break 
       restricted]. 
       Since $1_{\Omega_{\beta}}\in\paren{\Rpb}^{D_0}$, 
       $\tau 1_{\Omega_{\beta}} = T(1_{\Omega_{\beta}}) = 1_{\Omega_{\alpha}}$.  
       Since $T$ is linear, $\tau$ preserves order. 
       Thus by Lemma 5.1, 
       $\tau\paren{H^1_{\beta}\paren{\Rpb}} \subset H^1_{\beta}\paren{\Rpa}$. 
       Hence $1_{\Omega_{\alpha}}=\tau1_{\Omega_{\beta}} \in H^1_{\beta}\paren{\Rpa}$.  
       Now $1_{\Omega_{\alpha}} \in H^1_{\beta}\paren{\Rpa}$ for all $\beta<\alpha$.  
       Hence $1_{\Omega_{\alpha}} \in \bigcap_{\beta<\alpha} H^1_{\beta}\paren{\Rpa} = H^1_{\alpha}\paren{\Rpa}$.  
       \qquad\QED

\xproclaim {Corollary 5.13}. Let $1<p<\infty$ where $p\ne2$, and let $\alpha<\omega_1$. Then \hfil\break 
                             $h_p\paren{\Rpa} \ge \alpha+1$. 

\proof By Theorem 5.9, $\L{p}\not\injects\Rpa$, 
       and $H^1_{\alpha}\paren{\Rpa}\ne\emptyset$ by Theorem 5.12.
       Thus $h_p\paren{1,\Rpa} > \alpha$,
       whence $h_p\paren{\Rpa} \ge h_p\paren{1,\Rpa} \ge \alpha+1$. 
       \qquad\QED

 We collect our main results concerning the ordinal index $h_p$, the spaces $\Rpa$, and their interaction. 
 The proof of the subsequent theorem \xciteplus{B-R-S}{Theorem A} will make implicit use of these results. 

\xproclaim {Proposition 5.14}. Let $1<p<\infty$ where $p\ne2$.
                               Let $B$, $X$, and $Y$ be separable Banach spaces.
                               Let $\alpha,\beta<\omega_1$.
                               Then
             \item{(a)} $\L{p}\not\injects B$ if and only if $h_p(B)<\omega_1$, 
             \item{(b)} if $X\injects Y$, then $h_p(X)\le h_p(Y)$, 
             \item{(c)} $\L{p}\not\injects\Rpa$, 
             \item{(d)} if $\alpha<\beta$, then $\Rpa\cinjects\Rpb$, 
             \item{(e)} $h_p\paren{\Rpa}<\omega_1$, and
             \item{(f)} $h_p\paren{\Rpa}\ge\alpha+1$. 

\proof Parts (a), (b), (c), and (f) are restatements of Theorem 5.5, Theorem 5.6,
       Theorem 5.9, and Corollary 5.13, respectively. 
       Part (d) is clear from definitions.
       Part (e) is clear from parts (c) and (a). 
       \qquad\QED

\xproclaim {Theorem 5.15}. Let $1<p<\infty$ where $p\ne2$. 
                           There is a strictly increasing function $\tau:\omega_1\to\omega_1$
                           such that for $\gamma,\delta<\omega_1$, 
                \item{(a)} if $\gamma<\delta$, 
                           then $R_{\tau(\gamma)}^p \cinjects R_{\tau(\delta)}^p$ but
                           $R_{\tau(\delta)}^p \not\injects R_{\tau(\gamma)}^p$, and 
                \item{(b)} if $Y$ is a separable Banach space such that 
                           $R_{\tau(\alpha)}^p \injects Y$ for all $\alpha<\omega_1$, \hfil\break 
                           \indent then $\L{p} \injects Y$. 

\proof Let $\tau(0)=\omega<\omega_1$ [so $R^p_{\tau(0)}$ is infinite-dimensional].
       If $\tau(\beta)$ has been defined with $\tau(\beta)<\omega_1$,
       let $\tau(\beta+1) = h_p\paren{R^p_{\tau(\beta)}} < \omega_1$. Then \hfil\break 
       $h_p\paren{R^p_{\tau(\beta+1)}} \ge \tau(\beta+1)+1 > \tau(\beta+1) = h_p\paren{R^p_{\tau(\beta)}}$. 
       More generally,
       if $0<\alpha<\omega_1$ and $\tau(\beta)$ has been defined with $\tau(\beta)<\omega_1$ for all $\beta<\alpha$,
       let $\tau(\alpha) = \sup_{\beta<\alpha} h_p\paren{R^p_{\tau(\beta)}} < \omega_1$ 
       [each $h_p\paren{R^p_{\tau(\beta)}} < \omega_1$ and $\set{\beta:\beta<\alpha}$ is countable]. 
       Then \hfil\break 
       $h_p\paren{R^p_{\tau(\alpha)}} \ge \tau(\alpha)+1 > \tau(\alpha) = \sup_{\beta<\alpha} h_p\paren{R^p_{\tau(\beta)}}$, 
       so $h_p\paren{R^p_{\tau(\alpha)}} > h_p\paren{R^p_{\tau(\beta)}}$ for all $\beta<\alpha$. 
       Thus $R^p_{\tau(\alpha)}\not\injects R^p_{\tau(\beta)}$ for all $\beta<\alpha$,
       so $\tau(\alpha)>\tau(\beta)$ for all $\beta<\alpha$, 
       and $\tau$ is strictly increasing. 

       \item{(a)} Suppose $\gamma<\delta<\omega_1$. 
                  Then $\tau(\gamma)<\tau(\delta)$
                  and $R_{\tau(\gamma)}^p \cinjects R_{\tau(\delta)}^p$,
                  but $R_{\tau(\delta)}^p \not\injects R_{\tau(\gamma)}^p$ \hfil\break 
                  as shown above. 

       \item{(b)} Let $Y$ be a separable Banach space such that $R_{\tau(\alpha)}^p \injects Y$ for all $\alpha<\omega_1$. 
                  Then $\alpha < \tau(\alpha)+1 \le h_p\paren{R^p_{\tau(\alpha)}} \le h_p(Y) \le \omega_1$
                  for all $\alpha<\omega_1$. 
                  Thus $h_p(Y)=\omega_1$,
                  whence $\L{p}\injects Y$. 

                  \noindent\QED

\remark Let $1<p<\infty$ where $p\ne2$. 
        We will show that $\Rpa\cinjects\L{p}$ for all \hfil\break 
        $\alpha<\omega_1$. 
        Thus part (a) will yield uncountably many isomorphically distinct $\SL{p}$ spaces 
        [at most one $\Rpa\sim\l{2}$].
        By \xciteplus{J-M-S-T}{Corollary 9.2}, if $\L{p}\injects Y\cinjects\L{p}$, 
        then $Y\sim\L{p}$.  
        Thus part (b) will imply that there is no separable $\SL{p}$ space $Y$, other than $\L{p}$ itself, 
        such that $R^p_{\tau(\alpha)}\injects Y$ for all $\alpha<\omega_1$. 

\preheadspace
\firsthead{The Complementation of $\Rpa$ in $\L{p}$}
\postheadspace

 This section is devoted to the proof that $\Rpa\cinjects\L{p}$ for $1<p<\infty$ and $\alpha<\omega_1$. 
 We proceed by showing that $\Rpa\sim Z_{T_{\alpha}}^p\cinjects Z_{\NN}^p\sim\L{p}$
 for spaces $Z_{T_{\alpha}}^p$ and $Z_{\NN}^p$ to be defined. 
 The major components of the proof are Theorem 5.22, Proposition 5.25, and Proposition 5.26. 

\preheadspace
\secondhead{Preliminaries}
\postheadspace

 Let $\TT$ be a countable set, and let $\set{0,1}^{\TT}$ be the standard product space. 

 We say that a measurable function $f$ on $\set{0,1}^{\TT}$ depends on $E\subset\TT$
 if $f(x)=f(y)$ for all $x,y\in\set{0,1}^{\TT}$ such that $x|_E = y|_E$. 
 We say that a measurable set $S\subset\set{0,1}^{\TT}$ depends on $E\subset\TT$
 if the indicator function $1_S$ depends on $E$. 
 Thus $S\subset\set{0,1}^{\TT}$ \hfil\break 
 depends on $E\subset\TT$
 if $1_S(x)=1_S(y)$ for all $x,y\in\set{0,1}^{\TT}$ such that $x|_E = y|_E$. 

 It is easy to check that given $E\subset\TT$, 
 the set $\S{A}$ of all measurable $S\subset\set{0,1}^{\TT}$ which depend on $E$
 is a $\sigma$-algebra,
 which we call the $\sigma$-algebra corresponding to $E$. 
 Given $E\subset\TT$, 
 let $\S{A}_E$ be the $\sigma$-algebra corresponding to $E$.
 It is easy to check that 
 \item{(a)} if $A\subset B\subset\TT$, then $\S{A}_A\subset\S{A}_B$, and 
 \item{(b)} if $A,B\subset\TT$, then $\S{A}_{A\cap B} = \S{A}_A \cap \S{A}_B$. 

 \noindent
 Let $f$ be a measurable function on $\set{0,1}^{\TT}$ and let $E\subset\TT$. 
 It is easy to check that 
 \item{(c)} $f$ is $\S{A}_E$-measurable if and only if $f$ depends on $E$. 

\medskip

 Let $\paren{\Omega,\S{M},\mu}$ be a probability space. 
 Given a sub $\sigma$-algebra $\S{A}$ of $\S{M}$, 
 let $\E_{\S{A}}$ be the conditional expectation operator with respect to $\S{A}$. 

 Let $\S{A}$ be a sub $\sigma$-algebra of $\S{M}$. 
 Then for each integrable function $f$ on $\Omega$, 
 \item{(a)} $\E_{\S{A}} f$ is $\S{A}$-measurable, and 
 \item{(b)} $\int_S \E_{\S{A}} f = \int_S f$ for all $S\in\S{A}$. 

 \noindent
 Moreover, $\E_{\S{A}} f$ is essentially defined by these two conditions. 

 Let $\S{A}$ and $\S{B}$ be sub $\sigma$-algebras of $\S{M}$, 
 let $f$ and $g$ be integrable functions on $\Omega$, and
 let $1\le p<\infty$. 
 Conditional expectation has the following properties (\xcite{Ch}, \xcite{Db}, and \xcite{Stn}): 
 \item{(c)} if $f$ is $\S{A}$-measurable, then $\E_{\S{A}} f = f$,                                      
 \item{(d)} $\E_{\S{A}} \E_{\S{A}} f = \E_{\S{A}} f$,                                                   
 \item{(e)} if $f\in\L{p}(\Omega)$,
            then $\E_{\S{A}}f\in\L{p}(\Omega)$, with $\norm{\E_{\S{A}} f}_p \le \norm{f}_p$,            
 \item{(f)} if $f,g\in\L{2}(\Omega)$, then $\int g\E_{\S{A}} f = \int f\E_{\S{A}} g$,                   
 \item{(g)} if $f\in\L{2}(\Omega)$, then $f = \E_{\S{A}} f + f'$, 
            where $f'\in\L{2}(\Omega)$ such that $\int f'h=0$
            for all $\S{A}$-measurable $h\in\L{2}(\Omega)$,                                             
 \item{(h)} if $\S{A}\subset\S{B}$,
            then $\E_{\S{A}} f = \E_{\S{B}} f$ if and only if $\E_{\S{B}} f$ is $\S{A}$-measurable, and 
 \item{(i)} if $\S{A}\subset\S{B}$,
            then $\E_{\S{A}} \E_{\S{B}} f = \E_{\S{A}} f = \E_{\S{B}} \E_{\S{A}} f$.                    

 \noindent
 Suppose $\E_{\S{A}}$ and $\E_{\S{B}}$ commute. 
 Then $\E_{\S{A}} \E_{\S{B}} f$, which is equal to $\E_{\S{B}} \E_{\S{A}} f$, 
 is in turn $\S{A}$-measurable and $\S{B}$-measurable,
 whence $\S{A}\cap\S{B}$-measurable.
 Now $F=\E_{\S{A}} f$ is integrable on $\Omega$,
 $\S{A}\cap\S{B}\subset\S{B}$, 
 and $\E_{\S{B}} F = \E_{\S{B}} \E_{\S{A}} f$ is $\S{A}\cap\S{B}$-measurable.
 Thus \hfil\break 
 $\E_{\S{A}\cap\S{B}} f = \E_{\S{A}\cap\S{B}} \E_{\S{A}} f 
 =\E_{\S{A}\cap\S{B}} F = \E_{\S{B}} F = \E_{\S{B}} \E_{\S{A}} f$. 
 Hence
 \item{(j)} if $\E_{\S{A}} \E_{\S{B}} = \E_{\S{B}} \E_{\S{A}}$,
            then $\E_{\S{A}} \E_{\S{B}} = \E_{\S{A}\cap\S{B}} = \E_{\S{B}} \E_{\S{A}}$.  

\medskip

 Let $\paren{\set{0,1}^{\NN},\S{M},\mu}$ be the standard product space. 
 Let $A$ and $B$ be subsets of $\NN$, with corresponding $\sigma$-algebras $\S{A}$ and $\S{B}$, respectively. 
 Let $f$ be an integrable \hfil\break 
 function on $\set{0,1}^{\NN}$. 
 Consider $f$ as a function of $t=\paren{t_1,t_2,\ldots}$ where $t_i\in\set{0,1}$. 
 Then $\E_{\S{A}} f$ is given by integration with respect to those $t_i$ such that $i\in\NN\setminus A$. 
 Hence
 \item{(a)} $\E_{\S{A}} \E_{\S{B}} f = \E_{\S{B}} \E_{\S{A}} f$, and
 \item{(b)} $\E_{\S{A}} \E_{\S{B}} f = \E_{\S{A}\cap\S{B}} f = \E_{\S{B}} \E_{\S{A}} f$. 

\preheadspace
\secondhead{The Isomorphism of $Z_{\NN}^p$ and $\L{p}$}
\postheadspace

 Let $\set{A_n}$ be a sequence of sets. 
 We say that $\set{A_n}$ is monotonic if it is either nondecreasing or nonincreasing, 
 and $\set{A_n}$ is compatible if there is a permutation $\tau$ such that $\set{A_{\tau(n)}}$ is monotonic. 

 The following result \xciteplus{Stn}{Theorem 8} substitutes for \xciteplus{B-R-S}{Lemma 3.2}. 
 We do not present the proof,
 but apply the result in the proof of the subsequent corollary,
 which substitutes for \xciteplus{B-R-S}{Lemma 3.3}. 
 This alternative approach was suggested in a remark of \xcite{B-R-S}. 

\xproclaim {Proposition 5.16}. Let $1<p<\infty$, let $\paren{\Omega,\S{M},\mu}$ be a probability space, 
                               and let $\set{f_n}$ be a sequence of integrable functions on $\Omega$. 
                               Suppose $\set{\S{A}_n}$ is a compatible sequence of sub $\sigma$-algebras of $\S{M}$. 
                               Then there is a constant $A_p$, depending only on $p$, such that 
                               $$\norm{\paren{\tsuml_n \abs{\S{E}_{\S{A}_n} f_n}^2 }^{{1\over2}} }_p 
                                 \le A_p
                                 \norm{\paren{\tsuml_n \abs{f_n}^2 }^{{1\over2}} }_p.$$

\xproclaim {Corollary 5.17}. Let $1<p<\infty$, let $\paren{\Omega,\S{M},\mu}$ be a probability space, 
                             let $\set{f_n}$ be a sequence of integrable functions on $\Omega$, 
                             and let $\set{\S{B}_n}$ be a sequence of sub $\sigma$-algebras of $\S{M}$. 
                             Suppose $\set{\S{L}_n}$, $\set{\S{R}_n}$, and $\set{\S{T}_n}$
                             are sequences of sub $\sigma$-algebras \hfil\break 
                             of $\S{M}$ such that 
                  \item{(a)} each of $\set{\S{L}_n}$, $\set{\S{R}_n}$, and $\set{\S{T}_n}$ is compatible, 
                  \item{(b)} for each $n$, $\S{E}_{\S{L}_n}$, $\S{E}_{\S{R}_n}$, and $\S{E}_{\S{T}_n}$ commute, and
                  \item{(c)} for each $n$, $\S{B}_n = \S{L}_n \cap \S{R}_n \cap \S{T}_n$. 

\vglue -2 pt
                  \noindent{\sl Then for $A_p$ as above,} 
                             $$\norm{\paren{\tsuml_n \abs{\S{E}_{\S{B}_n} f_n}^2 }^{{1\over2}} }_p 
                               \le A_p^3
                               \norm{\paren{\tsuml_n \abs{f_n}^2 }^{{1\over2}} }_p.$$ 


\vglue 2 pt
\proof By part (c), 
       $\S{E}_{\S{B}_n} = \S{E}_{\S{L}_n \cap \S{R}_n \cap \S{T}_n}$.  
       By part (b), 
       $\S{E}_{\S{L}_n \cap \S{R}_n \cap \S{T}_n} = \S{E}_{\S{L}_n} \S{E}_{\S{R}_n} \S{E}_{\S{T}_n}$.  
       Thus $\S{E}_{\S{B}_n} = \S{E}_{\S{L}_n} \S{E}_{\S{R}_n} \S{E}_{\S{T}_n}$.  
       Hence by Proposition 5.16 (applied three times), we have
       $$\norm{\paren{\tsuml_n \abs{\S{E}_{\S{B}_n} f_n}^2 }^{{1\over2}} }_p 
         =
         \norm{\paren{\tsuml_n \abs{
               \S{E}_{\S{L}_n} \paren{\S{E}_{\S{R}_n} \paren{\S{E}_{\S{T}_n} f_n}}
               }^2 }^{{1\over2}} }_p 
         \le A_p^3
         \norm{\paren{\tsuml_n \abs{f_n}^2 }^{{1\over2}} }_p.$$
       \noindent\QED

 Let $n\in\NN$. 
 Then $n$ has a unique expression as $n=2^k+r$ for $k\in\NN\cup\set{0}$ and $0\le r<2^k$. 
 For $n=2^k+r$ as above, let $\lambda(n)=k$. 

 Let $D'_0=\set{1}$.
 For $k\in\NN$,
 let $D'_k = \set{t_0\cdots t_k: t_0\!=\!1 \hbox{ and } t_i\!\in\!\set{0,1} \hbox{ for } 1\le i\le k}$. 
 Let $\SD' = \bigcup_{k=0}^{\infty} D'_k$.  

 Now $\SD'$ has a natural strict partial order $\prec$ defined by $s_0\cdots s_{k_1} \prec t_0\cdots t_{k_2}$
 if $k_1<k_2$ and $s_i=t_i$ for all $0\le i\le k_1$. 

 Let $\gamma:\paren{\NN,<}\to\paren{\SD',\prec}$ be defined by
 $\gamma(n)=t_0\cdots t_k\in D'_k$ for $k=\lambda(n)$, 
 where $t_0\cdots t_k$ is the binary expansion of $n$. 
 Then $\gamma$ is a bijection, and $\gamma^{-1}$ preserves order.
 Let $\dotprec$ be the strict partial order on $\NN$ induced by $\prec$ via $\gamma$
 [$m\dotprec n \iff \gamma(m)\prec\gamma(n)$]. 
 Then $<$ extends $\dotprec$. 

 The following application of Corollary 5.17 substitutes for \xciteplus{B-R-S}{Scholium 3.4}.
 The result serves as a lemma for Theorem 5.22.

\xproclaim {Proposition 5.18}. Let $1<p<\infty$,
                               let $\paren{\set{0,1}^{\NN},\S{M},\mu}$ be the standard product space, 
                               and let $\set{f_n}$ be a sequence of integrable functions on $\set{0,1}^{\NN}$. 
                               Given $n\in\NN$, let $B_n=\set{m\in\NN:m\dotpreceq n}$, 
                               and let $\B_n$ be the corresponding sub $\sigma$-algebra of $\S{M}$. 
                               Then for $A_p$ as above and $N\in\NN$, 
                               $$\norm{\paren{\tsuml_{n=1}^N \abs{\S{E}_{\S{B}_n} f_n}^2 }^{{1\over2}} }_p 
                                 \le A_p^3
                                 \norm{\paren{\tsuml_{n=1}^N \abs{f_n}^2 }^{{1\over2}} }_p.$$ 

\proof Given $k\in\NN\cup\set{0}$,
       let $\Lambda_k = \set{m\in\NN : \lambda(m)=k}$, and let \hfil\break 
       $T_{[k]} = \set{m\in\NN : \lambda(m)\le k}$. 
       Given $n\in\NN$,
       let $\Lambda(n) = \set{m\in\NN : \lambda(m)=\lambda(n)}$, and let 
       $$T_n = \set{m\in\NN : m\le n},$$ 
       $$B_n = \set{m\in\NN : m\dotpreceq n}$$ 
       as above, which is the branch of $\paren{T_n,\dotprec}$ generated by $n$,
       $$L_n = \set{m\in\NN : m\dotpreceq n' \hbox{ for some } n'\in\Lambda(n) \hbox{ with } n'\le n},$$
       the union of the branches $B_{n'}$ for $n'\in\Lambda(n)$ with $n'\le n$, and
       $$R_n = \set{m\in\NN : m\dotpreceq n' \hbox{ for some } n'\in\Lambda(n) \hbox{ with } n'\ge n},$$
       the union of the branches $B_{n'}$ for $n'\in\Lambda(n)$ with $n'\ge n$. 

       Fix $K\in\NN\cup\set{0}$.
       For each $n\in T_{[K]}$, 
       choose $N(n)\in\Lambda_K$ such that $n\dotpreceq N(n)$.
       Then given $n\in T_{[K]}$, 
       $B_{N(n)}$ is an extension of $B_n$ to a branch of $T_{[K]}$, and
       $$B_n = B_{N(n)} \cap T_n = L_{N(n)} \cap R_{N(n)} \cap T_n.$$ 
       Note that $\set{L_N}_{N\in\Lambda_K}$, 
                 $\set{R_N}_{N\in\Lambda_K}$, and
                 $\set{T_n}_{n\in T_{[K]}}$ 
       are each monotonic. Hence \hfil\break 
             $\set{L_{N(n)}}_{n\in T_{[K]}}$, 
             $\set{R_{N(n)}}_{n\in T_{[K]}}$, and
             $\set{T_n}_{n\in T_{[K]}}$
       are each compatible. 

       For $n\in T_{[K]}$, 
       let $\S{B}_n$, $\S{L}_n$, $\S{R}_n$, and $\S{T}_n$ be the $\sigma$-algebras 
       corresponding to $B_n$, $L_n$, $R_n$, and $T_n$, respectively. 
       Then $\set{\S{L}_{N(n)}}_{n\in T_{[K]}}$, 
            $\set{\S{R}_{N(n)}}_{n\in T_{[K]}}$, and
            $\set{\S{T}_n}_{n\in T_{[K]}}$
       are each compatible. 
       Moreover, for $n\in T_{[K]}$, 
       $\S{B}_n = \S{L}_{N(n)} \cap \S{R}_{N(n)} \cap \S{T}_n$, 
       and $\S{E}_{\S{L}_{N(n)}}$, $\S{E}_{\S{R}_{N(n)}}$, and $\S{E}_{\S{T}_n}$ commute. 
 
       Hence for $N=2^{K+1}-1\in T_{[K]}$ and $f_1,\ldots,f_N$ integrable on $\set{0,1}^{\NN}$, \hfil\break 
       $\norm{\paren{\sum_{n=1}^N \abs{\S{E}_{\S{B}_n} f_n}^2 }^{{1\over2}} }_p
       \le A_p^3
       \norm{\paren{\sum_{n=1}^N \abs{f_n}^2 }^{{1\over2}} }_p$
       by Corollary 5.17. Releasing \hfil\break 
       $K\in\NN\cup\set{0}$ as a free variable,
       we have the same result for arbitrary $N\in\NN$. 
       \qquad\QED

 The following square function inequality \xciteplus{Burk}{Theorem 9}
 is quoted in \xciteplus{B-R-S}{Scholium 3.5}. 
 We do not present the proof, but apply the result in the proof of \hfil\break 
 Theorem 5.22. 


\xproclaim {Proposition 5.19}. Let $1<p<\infty$,
                               let $\paren{\Omega,\S{M},\mu}$ be a probability space, and 
                               let $\set{\S{T}_n}_{n=0}^{\infty}$ be a nondecreasing sequence of sub $\sigma$-algebras of $\S{M}$. 
                               Suppose $\set{g_n}_{n=0}^{\infty}$ is a sequence in $\L{p}(\Omega)$
                               such that $g_n$ is $\S{T}_n$-measurable for all $n\in\NN\cup\set{0}$ 
                               and $\E_{\S{T}_{n-1}} g_n = 0$ for all $n\in\NN$. 
                               Then there is a constant $K_p$, depending only on $p$, such that 
                               $${1\over K_p}
                                 \norm{\paren{\tsuml_n \abs{g_n}^2 }^{{1\over2}} }_p 
                                 \le
                                 \norm{\tsuml_n g_n}_p
                                 \le K_p
                                 \norm{\paren{\tsuml_n \abs{g_n}^2 }^{{1\over2}} }_p.$$

 For $n\in\NN$, 
 let $B_n$, $T_n$, $\B_n$, and $\T_n$ be as above. 
 Then for $n\in\NN$, 
 $T_n$ is the subtree $\set{1,\ldots,n}$ of $\paren{\NN,\dotprec}$, 
 $B_n$ is the branch of $T_n$ generated by $n$, and 
 $\T_n$ and $\B_n$ are the $\sigma$-algebras corresponding to $T_n$ and $B_n$, respectively. 
 Let $T_0=B_0=\emptyset$ and let $\T_0$ and $\B_0$ be the trivial algebras. 
 Let
 $$\eqalign{Z_{\NN}^p
     &= \span{f: f \hbox{ is } \B_n\hbox{-measurable for some } n\in\NN}_{\L{p}\paren{\set{0,1}^{\NN}}}
 \cr &= \span{f: f \hbox{ is measurable and depends on } B_n \hbox{ for some } n\in\NN}_{\L{p}\paren{\set{0,1}^{\NN}}}. }$$
 Let $\Delta_0=\Gamma_0=\set{\hbox{constant functions on }\set{0,1}^{\NN}}$. 
 For $n\in\NN$, let
 $$\Delta_n = \set{f \hbox{ on } \set{0,1}^{\NN} : f \hbox{ is } \T_n\hbox{-measurable and } \E_{\T_{n-1}} f = 0}$$
 and
 $$\Gamma_n = \set{f\in\Delta_n : f \hbox{ is } \B_n\hbox{-measurable} \vphantom{0^N} }.$$

 Suppose $f$ is measurable and $n\in\NN$. 
 Then $\paren{\E_{\T_n} - \E_{\T_{n-1}} } f$ is $\T_n$-measurable, 
 and $\E_{\T_{n-1}} \paren{\E_{\T_n} - \E_{\T_{n-1}} } f = \E_{\T_{n-1}} f - \E_{\T_{n-1}} f = 0$, 
 whence $\paren{\E_{\T_n} - \E_{\T_{n-1}} } f \in\Delta_n$. 
 Note that if $f\in\Delta_n$, then $f = \paren{\E_{\T_n} - \E_{\T_{n-1}} } f$. 
 Hence for $n\in\NN$, 
 $$\Delta_n = \set{\paren{\E_{\T_n} - \E_{\T_{n-1}} } f : f \hbox{ on } \set{0,1}^{\NN} \hbox{ is measurable}}.$$

 The following lemmas for Theorem 5.22 have been extracted from the proof of \hfil\break 
 \xciteplus{B-R-S}{Theorem 3.1}.


\xproclaim {Lemma 5.20}. Let $1\le p<\infty$, and let $Z_{\NN}^p$ and $\Gamma_n$ be as above. 
                         Then \hfil\break 
                         $Z_{\NN}^p = \span{\Gamma_n:n\ge0}_{\L{p}\paren{\set{0,1}^{\NN}}}$. 

 \proof Note that $\Gamma_n\subset Z_{\NN}^p$ for $n\in\NN\cup\set{0}$,
        whence $\span{\Gamma_n:n\ge0}_{\L{p}\paren{\set{0,1}^{\NN}}} \subset Z_{\NN}^p$. 
        We now show that $Z_{\NN}^p \subset \span{\Gamma_n:n\ge0}_{\L{p}\paren{\set{0,1}^{\NN}}}$, 
        whence $Z_{\NN}^p = \span{\Gamma_n:n\ge0}_{\L{p}\paren{\set{0,1}^{\NN}}}$. 

        Let $n\in\NN$ and let $f$ be $\B_n$-measurable.
        Now $B_n \subset T_n$, so $\B_n \subset \T_n$, whence $f$ is $\T_n$-measurable and $\E_{\T_n} f = f$. 
        Moreover, 
        $\E_{\T_0} f$ is $\T_0$-measurable, whence $\E_{\T_0} f$ is constant, 
        and $\int\E_{\T_0} f = \int f$, whence $\E_{\T_0} f = \int\E_{\T_0} f = \int f$. 
        Thus 
        $$f = \tint f - \E_{\T_0} f + \E_{\T_n} f 
            = \tint f + \tsuml_{i=1}^n \paren{\E_{\T_i} - \E_{\T_{i-1}} } f.$$
        Let $1\le i\le n$.
        Then $\paren{\E_{\T_i} - \E_{\T_{i-1}} } f \in\Delta_i$. 
        We now show that $\paren{\E_{\T_i} - \E_{\T_{i-1}} } f$ is $\B_i$-measurable, 
        whence it will follow that $\paren{\E_{\T_i} - \E_{\T_{i-1}} } f \in \Gamma_i$. 

        Note that $f = \E_{\B_n} f$, whence 
        $$\paren{\E_{\T_i} - \E_{\T_{i-1}} } f  
        = \paren{\E_{\T_i} - \E_{\T_{i-1}} } \E_{\B_n} f
        = \E_{\T_i} \E_{\B_n} f - \E_{\T_{i-1}} \E_{\B_n} f
        = \E_{\T_i \cap \B_n} f - \E_{\T_{i-1} \cap \B_n} f.$$
        Suppose first that $i\not\in B_n$. 
        Then $T_i \cap B_n = T_{i-1} \cap B_n$, so $\T_i \cap \B_n = \T_{i-1} \cap \B_n$, 
        whence 
        $$\paren{\E_{\T_i} - \E_{\T_{i-1}} } f  
        = \E_{\T_i \cap \B_n} f - \E_{\T_{i-1} \cap \B_n} f
        = 0,$$
        which is $\B_i$-measurable. 
        Next suppose that $i\in B_n$. 
        Then $T_i \cap B_n = B_i$, so $\T_i \cap \B_n = \B_i$, 
        and $T_{i-1} \cap B_n \subset B_i$, so $\T_{i-1} \cap \B_n \subset \B_i$, 
        whence
        $$\paren{\E_{\T_i} - \E_{\T_{i-1}} } f  
        = \E_{\T_i \cap \B_n} f - \E_{\T_{i-1} \cap \B_n} f
        = \E_{\B_i} f - \E_{\B'_i} f$$
        for some $\B'_i\subset\B_i$. 
        Now $\E_{\B_i} f$ is $\B_i$-measurable,
        and $\E_{\B'_i} f$ is $\B'_i$-measurable, whence \hfil\break 
        $\B_i$-measurable. 
        Thus $\paren{\E_{\T_i} - \E_{\T_{i-1}} } f$ is $\B_i$-measurable
        [now in both cases]. 
        As noted above,
        it follows that $\paren{\E_{\T_i} - \E_{\T_{i-1}} } f \in \Gamma_i$. 

        We now have 
        $$f = \tint f + \tsuml_{i=1}^n \paren{\E_{\T_i} - \E_{\T_{i-1}} } f \in
              \span{\Gamma_i:0\le i\le n}_{\L{p}\paren{\set{0,1}^{\NN}}}.$$ 
        Thus $f\in\span{\Gamma_n:n\ge0}_{\L{p}\paren{\set{0,1}^{\NN}}}$. 
        It follows that $Z_{\NN}^p \subset \span{\Gamma_n:n\ge0}_{\L{p}\paren{\set{0,1}^{\NN}}}$, 
        whence $Z_{\NN}^p = \span{\Gamma_n:n\ge0}_{\L{p}\paren{\set{0,1}^{\NN}}}$. 
        \qquad\QED 

\medskip

\xproclaim {Lemma 5.21}. Let $2\le p<\infty$, and let $\Delta_i$ be as above. 
                         Then $\set{\Delta_i}_{i\ge0}$ is an unconditional Schauder decomposition of 
                         $\L{p}\paren{\set{0,1}^{\NN}}$. 

\proof Suppose $f,g\in\L{2}\paren{\set{0,1}^{\NN}}$, and let $i\in\NN$. 
       If $f\in\Delta_i$ and $g\in\Delta_j$ for $i<j\in\NN$,
       then $\E_{\T_{j-1}}g=0$ and $f$ is $\T_{j-1}$-measurable, so
       $$\tint fg = \tint f \paren{g - \E_{\T_{j-1}}g}  
                  = \tint fg - \tint f\E_{\T_{j-1}}g
                  = \tint fg - \tint g\E_{\T_{j-1}}f
                  = \tint fg - \tint gf
                  = 0,$$
       whence $f$ and $g$ are orthogonal. 
       If $f\in\Delta_i$ and $g\in\Delta_0$,
       then $g$ is constant, and \hfil\break 
       $\int f = \int\E_{\T_{i-1}}f$, but $\E_{\T_{i-1}}f = 0$, so 
       $$\tint fg = g\tint f = g\tint \E_{\T_{i-1}}f = 0,$$
       whence $f$ and $g$ are orthogonal. 
       Hence $\set{\Delta_i}_{i\ge0}$ is orthogonal. 

       Suppose $f\in\L{2}\paren{\set{0,1}^{\NN}}$. 
       Let $f_0=\E_{\T_0}f\in\Delta_0$,
       and for $i\in\NN$, let \hfil\break 
       $f_i=\paren{\E_{\T_i}-\E_{\T_{i-1}}}f\in\Delta_i$. 
       Then for $n\in\NN$, 
       $$\tsuml_{i=0}^n f_i 
       = \E_{\T_0}f + \tsuml_{i=1}^n \paren{\E_{\T_i}-\E_{\T_{i-1}}}f
       = \E_{\T_n}f.$$
       Note that $\L{p}\paren{\set{0,1}^{\NN}} \subset \L{2}\paren{\set{0,1}^{\NN}}$. 
       If $f\in\L{p}\paren{\set{0,1}^{\NN}}$, then \hfil\break 
       $\lim_{n\to\infty} \norm{f-\E_{\T_n}f}_p=0$, 
       whence $f=\sum_{i=0}^{\infty} f_i$ in $\L{p}\paren{\set{0,1}^{\NN}}$. 
       By the orthogonality of $\set{\Delta_i}_{i\ge0}$, 
       the representation $f=\sum_{i=0}^{\infty} f'_i$ with $f'_i\in\Delta_i$ is unique. 
       By Proposition 5.19, the convergence is unconditional.  
       Hence $\set{\Delta_i}_{i\ge0}$ is an unconditional Schauder decomposition of 
       $\L{p}\paren{\set{0,1}^{\NN}}$. 
       \qquad\QED

\remark The above result can be viewed as a consequence of the \hfil\break 
        unconditionality of the Haar system. 

 We are now prepared to prove the following theorem \xciteplus{B-R-S}{Theorem 3.1}, \hfil\break 
 which is a major component of the proof that $\Rpa\cinjects\L{p}$. 


\xproclaim {Theorem 5.22}. Let $1<p<\infty$, and let $Z_{\NN}^p$ be as above. 
                           Then \hfil\break 
                           $Z_{\NN}^p \cinjects \L{p}\paren{\set{0,1}^{\NN}}$. 

\proof First suppose $2\le p<\infty$, whence $\L{p}\paren{\set{0,1}^{\NN}} \subset \L{2}\paren{\set{0,1}^{\NN}}$. 
       Fix \hfil\break 
       $i\in\NN\cup\set{0}$.
       Let $f\in\Delta_i$ and let $g=\E_{\B_i}f$. 
       If $i=0$, then $\Gamma_i=\Delta_i$, $\E_{\B_i}f=f$, and $\E_{\B_i}|_{\Delta_i}$ is the identity mapping. 
       Suppose $i\in\NN$. 
       Then $g$ is $\B_i$-measurable. 
       Now $B_i\subset T_i$, so $\B_i\subset\T_i$, 
       whence $g$ is $\T_i$-measurable. 
       Moreover, $\E_{\T_{i-1}}g = \E_{\T_{i-1}} \E_{\B_i} f = \E_{\B_i} \E_{\T_{i-1}} f = 0$. 
       Thus $g$ is a $\B_i$-measurable element of $\Delta_i$, 
       whence $g\in\Gamma_i$. 
       If $f\in\Gamma_i$, then $\E_{\B_i}f=f$. 
       Hence for $i\in\NN\cup\set{0}$, 
       $\E_{\B_i}|_{\Delta_i}$ is the orthogonal projection of $\Delta_i$ onto $\Gamma_i$. 

       By Lemma 5.21, 
       $\set{\Delta_i}_{i\ge0}$ is an unconditional Schauder decomposition of \hfil\break 
       $\L{2}\paren{\set{0,1}^{\NN}}$. 
       For $f\in\L{2}\paren{\set{0,1}^{\NN}}$, 
       let $\set{f_i}$ be the unique sequence with $f_i\in\Delta_i$ such that $f=\sum_{i=0}^{\infty}f_i$. 
       Let $\pi : \L{2}\paren{\set{0,1}^{\NN}} \to \L{2}\paren{\set{0,1}^{\NN}}$ be defined by  
       $$\pi f = \tsuml_{i=0}^{\infty} \E_{\B_i} f_i.$$
       Then $\pi$ is the orthogonal projection of $\L{2}\paren{\set{0,1}^{\NN}}$
       onto $\span{\Gamma_i:i\ge0}_{\L{2}\paren{\set{0,1}^{\NN}}}$, 
       where $\span{\Gamma_i:i\ge0}_{\L{2}\paren{\set{0,1}^{\NN}}} = Z_{\NN}^2$ by Lemma 5.20. 

       Let $P$ be the restriction of $\pi$ to $\L{p}\paren{\set{0,1}^{\NN}}$, 
       let $f\in\L{p}\paren{\set{0,1}^{\NN}}$, and let $\set{f_i}$ be as above. 
       Then by Proposition 5.19, Proposition 5.18, and Proposition 5.19 again, for $n\in\NN$ we have 
               $$\norm{\tsuml_{i=0}^n \E_{\B_i} f_i}_p
         \le K_p \norm{\paren{\tsuml_{i=0}^n \abs{\E_{\B_i} f_i}^2 }^{{1\over2}} }_p
   \le K_p A_p^3 \norm{\paren{\tsuml_{i=0}^n \abs{f_i}^2 }^{{1\over2}} }_p
 \le K_p^2 A_p^3 \norm{\tsuml_{i=0}^n f_i}_p,$$
       where the constants $K_p$ and $A_p$ are as in the cited propositions. Hence \hfil\break 
       $\norm{Pf}_p \le K_p^2 A_p^3 \norm{f}_p$,
       and $P : \L{p}\paren{\set{0,1}^{\NN}} \to \L{p}\paren{\set{0,1}^{\NN}}$ is bounded. 
       Of course \hfil\break 
       $P$ is a projection,
       and $P$ maps $\L{p}\paren{\set{0,1}^{\NN}}$
       onto $\span{\Gamma_i:i\ge0}_{\L{p}\paren{\set{0,1}^{\NN}}}$, where \hfil\break 
       $\span{\Gamma_i:i\ge0}_{\L{p}\paren{\set{0,1}^{\NN}}} = Z_{\NN}^p$ by Lemma 5.20. 

       For $2<p<\infty$ with conjugate index $q$,
       the adjoint of $P$ induces a bounded projection of $\L{q}\paren{\set{0,1}^{\NN}}$ onto $Z_{\NN}^q$. 
       \qquad\QED


\remark While $Z_{\NN}^p \cinjects \L{p}$ is our major concern, 
        in fact $Z_{\NN}^p\sim\L{p}$.

\preheadspace
\secondhead{The Complementation of $\Rpa$ in $Z_{\NN}^p$} 
\postheadspace

 Recall that a tree $\paren{T,\prec}$ is a CFRE tree if $T$ is finite or countable, 
 and for each $x\in T$, $\set{y\in T:y\prec x}$ is finite. 
 Let $\paren{T,\prec}$ be a CFRE tree. 
 For $t\in T$, let $B_t$ be the finite branch of $T$ generated by $t$.
 For $1\le p<\infty$, let
 $$Z_T^p = \span{f: f \hbox{ is measurable and depends on } B_t \hbox{ for some } t\in T}_{\L{p}\paren{\set{0,1}^T}}.$$
 The space $Z_T^p$ is similar to the previously defined space $Z_{\NN}^p$. 

 Let $S$ be a nonempty subset of $\NN$. 
 Then $\paren{S,\dotprec}$ is a CFRE tree, where $\dotprec$ is the \hfil\break 
 previously introduced partial order on $\NN$ [suitably restricted].
 The finite branches of $S$ are precisely those sets of the form $B_n\cap S$ for $n\in S$, 
 where $B_n$ is the finite branch of $\paren{\NN,\dotprec}$ generated by $n$. 
 For $1\le p<\infty$, 
 $\L{p}\paren{\set{0,1}^S}$ is isomorphic to the subspace of $\L{p}\paren{\set{0,1}^{\NN}}$
 consisting of those functions which depend on $S$, 
 and $Z_S^p$ is isomorphic to the space
 $$\tilde Z_S^p =
   \span{f: f \hbox{ is measurable and depends on } B_n\cap S \hbox{ for some } n\in S}_{\L{p}\paren{\set{0,1}^{\NN}}}.$$

 The following lemmas \xciteplus{B-R-S}{Lemmas 3.6 and 3.7} lead to the subsequent \hfil\break 
 proposition \xciteplus{B-R-S}{Theorem 3.8}, 
 which is a component of the proof that $\Rpa\cinjects\L{p}$. 

\xproclaim {Lemma 5.23}. Let $1\le p<\infty$ and let $\emptyset\ne S\subset\NN$. 
                         Then $Z_S^p \cinjects Z_{\NN}^p$. 

\proof Let $\S{S}$ be the $\sigma$-algebra corresponding to $S$,
       and let $P:Z_{\NN}^p\to Z_{\NN}^p$ be \hfil\break 
       defined by $Pf=\E_{\S{S}}f$. 
       Note that $\tilde Z_S^p\subset Z_{\NN}^p$. 
       If $f\in Z_{\NN}^p$ depends on $B_n$, then $Pf$ depends on $B_n\cap S$,
       which is either the empty set or a finite branch of $S$ of the form
       $B_m\cap S$ for some $m\in S$, 
       whence $P$ maps $Z_{\NN}^p$ into $\tilde Z_S^p$. 
       Now $Pf=f$ for $f\in\tilde Z_S^p$.
       Hence $P$ maps $Z_{\NN}^p$ onto $\tilde Z_S^p$, and $P^2=P$.
       Finally, $\norm{Pf}_p=\norm{\E_{\S{S}}f}_p\le\norm{f}_p$, whence $\norm{P}=1$.  
       Hence $Z_S^p \sim \tilde Z_S^p \cinjects Z_{\NN}^p$. 
       \qquad\QED

 For $n\in\NN$, let $N_n=\set{t_1\cdots t_n:t_i\in\NN \hbox{ for all } 1\le i\le n}$. 
 Let $\S{N}=\bigcup_{n=1}^{\infty}N_n$, 
 and define a strict partial order $\prec$ on $\S{N}$ by
 $s_1\cdots s_n \prec t_1\cdots t_m$ if $n<m$ and $s_i=t_i$ for all $1\le i\le n$. 
       
\xproclaim {Lemma 5.24}. Let $\paren{T,\prec}$ be a CFRE tree.
                         Then $\paren{T,\prec}$ is order-isomorphic to a subset of $\paren{\NN,\dotprec}$. 

\proof Clearly $T$ is order-isomorphic to a subset of $\S{N}$. 
       We will show that $\S{N}$ is order-isomorphic to a subset of $\S{D}'$. 
       The result will then follow upon noting that 
       $\S{D}'$ is order-isomorphic to $\NN$ endowed with $\dotprec$. 

       We describe a subset $\S{S}$ of $\S{D}'$ such that $\S{N}$ is order-isomorphic to $\S{S}$. 
       Given \hfil\break 
       $t\in\S{D}'$, let $S(t)=\set{t\con1,t\con01,t\con001,\ldots}$. 
       Then $S(t)$ is a countable set of distinct and \hfil\break 
       mutually incomparable successors of $t$. 
       Moreover, if $s$ and $t$ are distinct and \hfil\break 
       incomparable elements of $\S{D}'$, 
       then $S(s)$ and $S(t)$ are disjoint, 
       and the elements of $S(s)\cup S(t)$ are mutually incomparable elements of $\S{D}'$. 
       For $A\subset\S{D}'$, let \hfil\break 
       $S(A)=\bigcup_{a\in A}S(a)$.   
       Finally, let $\S{S}=S(1)\cup S(S(1))\cup\cdots$. 
       Then $\S{N}$ is order-isomorphic to $\S{S}\subset\S{D}'$, 
       and the result follows as noted above. 
       \qquad\QED

\xproclaim {Proposition 5.25}. Let $1\le p<\infty$ and let $T$ be a CFRE tree.
                               Then $Z_T^p\cinjects Z_{\NN}^p$. 

\proof If trees $T$ and $T'$ are order-isomorphic, then $Z^p_T\sim Z^p_{T'}$.  
       Thus by Lemma 5.24, we may choose $T'\subset\NN$ such that $Z_T^p\sim Z_{T'}^p$. 
       Now $Z_{T'}^p\cinjects Z_{\NN}^p$ by Lemma 5.23. 
       Hence $Z_T^p\cinjects Z_{\NN}^p$. 
       \qquad\QED

\remark By Proposition 5.25 and Theorem 5.22,
        for $1<p<\infty$ and $T$ a CFRE tree, 
        $Z_T^p\cinjects\L{p}\paren{\set{0,1}^{\NN}}$,  
        whence $Z_T^p\cinjects\L{p}\paren{\set{0,1}^T}$. 


 The following proposition \xciteplus{B-R-S}{Lemma 3.9}
 is the final component of the proof that $\Rpa\cinjects\L{p}$. 

\xproclaim {Proposition 5.26}. Let $1\le p<\infty$ and $\alpha<\omega_1$. 
                               Then there is a well-founded CFRE tree $T_{\alpha}$
                               such that $\Rpa$ is distributionally isomorphic to $Z_{T_{\alpha}}^p$. 

\proof Clearly $R_0^p=\span{1}_{\L{p}}$ is distributionally isomorphic to $Z_{T_0}^p$ where $T_0=\emptyset$. 
       Moreover, 
       $R_1^p=\psum{R_0^p}$ is distributionally isomorphic to $Z_{T_1}^p$ where $T_1=\set{1}$.

       Suppose $\alpha=\beta+1>1$, 
       where $\Rpb$ is distributionally isomorphic to $Z_{T_{\beta}}^p$
       for some well-founded CFRE tree $\paren{T_{\beta},\prec_{\beta}}$. 
       Without loss of generality, suppose $\Rpb=Z_{T_{\beta}}^p$. 
       Choose $\theta\not\in T_{\beta}$. 
       Let $T_{\alpha}=T_{\beta}\cup\set{\theta}$, 
       and let $\prec_{\alpha}$ extend $\prec_{\beta}$ by declaring $\theta\prec_{\alpha}\tau$ for all $\tau\in T_{\beta}$. 
       Then $\paren{T_{\alpha},\prec_{\alpha}}$ is a well-founded CFRE tree. 
       For the case $\alpha=\beta+1>1$, 
       it remains to show that $\Rpa$ is distributionally isomorphic to $Z_{T_{\alpha}}^p$. 

       Let $\bar0,\bar1\in\set{0,1}^{\set{\theta}}$ be defined by
       $\bar0(\theta)=0$ and $\bar1(\theta)=1$, so that $\bar\jmath(\theta)=j$.
       Note that $\set{0,1}^{\set{\theta}}=\set{\bar0,\bar1}$. 
       Let $e_0,e_1:\set{0,1}^{\set{\theta}}\to\set{0,1}$ be defined by
       $e_0(t)=1-t(\theta)$ and $e_1(t)=t(\theta)$. 
       Then $e_i(\bar\jmath)=1$ if $i=j$ and $e_i(\bar\jmath)=0$ if $i\ne j$. 

       Given $s\in\set{0,1}^{T_{\beta}}$ and $t\in\set{0,1}^{\set{\theta}}$,
       we associate $(s,t) \in \set{0,1}^{T_{\beta}}\times\set{0,1}^{\set{\theta}}  
                             = \set{0,1}^{T_{\beta}}\times\set{\bar0,\bar1}$
       with the element $[s,t]\in\set{0,1}^{T_{\alpha}}$ which extends both $s$ and $t$. 
       Thus there is an association $J:\L{p}\paren{\set{0,1}^{T_{\beta}}\times\set{\bar0,\bar1}}
                                    \to\L{p}\paren{\set{0,1}^{T_{\alpha}}}$.  
       Let $\psum{Z_{T_{\beta}}^p}$ be identified with the subspace of 
       $\L{p}\paren{\set{0,1}^{T_{\beta}}\times\set{\bar0,\bar1}}$ 
       which is related to $Z_{T_{\beta}}^p$ as in the definition of $\psum{B}$.
       Let $\brackpsum{Z_{T_{\beta}}^p} = J \psum{Z_{T_{\beta}}^p}$. Then \hfil\break 
       $\brackpsum{Z_{T_{\beta}}^p} \distiso \psum{Z_{T_{\beta}}^p}$.

       Let $b_0,b_1\in Z_{T_{\beta}}^p$. 
       Then $b_i\otimes e_i\in Z_{T_{\alpha}}^p$, 
       where $\paren{b_i\otimes e_i}[s,t]=2^{{1\over p}}b_i(s)e_i(t)$ for \hfil\break 
       $s\in\set{0,1}^{T_{\beta}}$ and $t\in\set{0,1}^{\set{\theta}}=\set{\bar0,\bar1}$.
       If $b=b_0\otimes e_0+b_1\otimes e_1$, then \hfil\break 
       $b[s,\bar\jmath]=2^{{1\over p}}b_0(s)e_0(\bar\jmath)+2^{{1\over p}}b_1(s)e_1(\bar\jmath)$,
       so $b[s,\bar0]=2^{{1\over p}}b_0(s)$ and $b[s,\bar1]=2^{{1\over p}}b_1(s)$,  
       whence $b\in\brackpsum{Z_{T_{\beta}}^p}$. 
       Conversely, if $b\in\brackpsum{Z_{T_{\beta}}^p}$,
       then $b=b_0\otimes e_0+b_1\otimes e_1$ for \hfil\break 
       $b_0(s)=2^{-{1\over p}}b[s,\bar0]$ and $b_1(s)=2^{-{1\over p}}b[s,\bar1]$. Hence \hfil\break 
       $\brackpsum{Z_{T_{\beta}}^p} = \set{b_0\otimes e_0+b_1\otimes e_1 : b_0,b_1\in Z_{T_{\beta}}^p}
                                      \subset Z_{T_{\alpha}}^p$.

       Let $f\in Z_{T_{\alpha}}^p$.
       For $s\in\set{0,1}^{T_{\beta}}$, 
       let $b_0(s)=2^{-{1\over p}}f[s,\bar0]$ and $b_1(s)=2^{-{1\over p}}f[s,\bar1]$.  
       Then $b_i\in Z_{T_{\beta}}^p$, and $f=b_0\otimes e_0 + b_1\otimes e_1$, 
                                      so $f\in\brackpsum{Z_{T_{\beta}}^p}$.
       Thus \hfil\break 
       $Z_{T_{\alpha}}^p \subset \brackpsum{Z_{T_{\beta}}^p}$,
       whence
       $Z_{T_{\alpha}}^p = \brackpsum{Z_{T_{\beta}}^p}$.
       For the case $\alpha=\beta+1>1$, 
       it now follows that $\Rpa = \psum{\Rpb} = \psum{Z_{T_{\beta}}^p}
                                   \distiso \brackpsum{Z_{T_{\beta}}^p} = Z_{T_{\alpha}}^p$. 

       Suppose $\alpha$ is a limit ordinal, where for each $\beta<\alpha$, 
       $\Rpb$ is distributionally \hfil\break 
       isomorphic to $Z_{T_{\beta}}^p$ for some well-founded CFRE tree $\paren{T_{\beta},\prec_{\beta}}$.
       Without loss of \hfil\break 
       generality, suppose $\Rpb=Z_{T_{\beta}}^p$ for all $\beta<\alpha$, 
       and suppose $T_{\gamma}\cap T_{\beta}=\emptyset$ for all $\gamma\ne\beta$ with $\gamma,\beta<\alpha$. 
       Let $T_{\alpha}=\bigcup_{\beta<\alpha}T_{\beta}$,
       and let $\sigma\prec_{\alpha}\tau$ if there is some $\beta<\alpha$
       such that $\sigma,\tau\in T_{\beta}$ with $\sigma\prec_{\beta}\tau$.
       Then $\paren{T_{\alpha},\prec_{\alpha}}$ is a well-founded CFRE tree. 

       Note that $B$ is a finite branch of $T_{\alpha}$ if and only if
       $B$ is a finite branch of $T_{\beta}$ for some $\beta<\alpha$. 
       Thus $f$ depends on a finite branch $B$ of $T_{\alpha}$ if and only if 
       $f$ depends on a finite branch $B$ of $T_{\beta}$ for some $\beta<\alpha$, 
       so $Z_{T_{\alpha}}^p = \span{Z_{T_{\beta}}^p : \beta<\alpha}_{\L{p}\paren{\set{0,1}^{T_{\alpha}}}}$.
       Since $\set{T_{\beta}}_{\beta<\alpha}$ is disjoint, 
       $\span{Z_{T_{\beta}}^p : \beta<\alpha}_{\L{p}\paren{\set{0,1}^{T_{\alpha}}}}
       \distiso \Ipsum{Z_{T_{\beta}}^p}{\beta<\alpha}$. Hence \hfil\break 
       $Z_{T_{\alpha}}^p = \span{Z_{T_{\beta}}^p : \beta<\alpha}_{\L{p}\paren{\set{0,1}^{T_{\alpha}}}}
                  \distiso \Ipsum{Z_{T_{\beta}}^p}{\beta<\alpha} 
                         = \Ipsum{\Rpb}{\beta<\alpha} 
                         = \Rpa$.
       \qquad\QED

\medskip
 The following theorem \xciteplus{B-R-S}{Theorem B(3)} is now almost immediate. 

\xproclaim {Theorem 5.27}. Let $1<p<\infty$ and $\alpha<\omega_1$. 
                           Then $\Rpa\cinjects\L{p}$. 

\proof By Proposition 5.26, we may choose a well-founded CFRE tree $T_{\alpha}$ such that $\Rpa \sim Z_{T_{\alpha}}^p$.
       Then $Z_{T_{\alpha}}^p \cinjects Z_{\NN}^p$ by Proposition 5.25, 
       and $Z_{\NN}^p \cinjects \L{p}\paren{\set{0,1}^{\NN}}$ by Theorem 5.22.
       Hence $\Rpa \cinjects \L{p}\paren{\set{0,1}^{\NN}} \sim \L{p}$. 
       \qquad\QED

\preheadspace
\firsthead{Concluding Remarks}
\postheadspace

 Let $1<p<\infty$ where $p\ne2$. 


 Conceivably $R_{\tau(\alpha)}^p\sim\l{2}$ for some $\alpha<\omega_1$, 
 but in light of part (a) of Theorem 5.15,
 this is possible only for $\alpha=0$. 
 Thus as in the remark following Theorem 5.15, 
 $\set{R_{\tau(\alpha)}^p}_{0<\alpha<\omega_1}$ is an uncountable chain of isomorphically distinct $\SL{p}$ spaces, 
 and there is no separable $\SL{p}$ space $Y$, other than $\L{p}$ itself, 
 such that $R_{\tau(\alpha)}^p\injects Y$ for all $\alpha<\omega_1$. 
 By Theorem 5.27 and part (a) of Theorem 5.15, for $\gamma<\delta<\omega_1$ we have 
 $$R_{\tau(\gamma)}^p \cra R_{\tau(\delta)}^p \cra \L{p}.\eqno{\TAG{5.5}{5.5}}$$

 The isomorphism type of $\Rpa$ for $\omega<\alpha<\omega_1$ is not well understood. 
 Recent work by Dale Alspach indicates that $R^p_{\omega}\sim X_p$. 

 We know that $\set{h_p\paren{\Rpa}}_{\alpha<\omega_1}$ is a nondecreasing chain of ordinals such that \hfil\break 
 $\set{h_p\paren{\Rpa}:\alpha<\omega_1}$ has no maximum, 
 but little is known about the specific values of $h_p\paren{\Rpa}$ for $\omega\le\alpha<\omega_1$,
 or precisely where the increases occur. 

 Part (b) of Theorem 5.15 reflects one way in which $\set{\Rpa}_{\alpha<\omega_1}$ reaches toward $\L{p}$. 
 However, it is not known whether
 for each separable $\SL{p}$ space $Y\not\sim\L{p}$,
 there is an $\alpha<\omega_1$ such that $Y\injects\Rpa$, 
 nor whether there is an $\alpha<\omega_1$ such that $Y\injects\Rpa$ for uncountably many $\SL{p}$ spaces $Y$.


\vfill\eject
\endgroup 
%
%
\begingroup               
\frenchspacing            
\parindent=5em            
\baselineskip=12pt        
\def\bcite#1{\hbox to 5em {[{\bf#1}]\hfil}} 
\def\papercite#1#2#3#4#5#6#7{\item{\bcite{#1}} 
                             {#2}, {\it{#3}}, {#4} {\bf#5} {(#6)}, {#7}. \par\smallskip} 
\def\othercite#1#2{\item{\bcite{#1}} {#2} \par\smallskip} 
\topglue 1 true in        
\pageno=160 
\mark{}\mark{2}
\centerline{{\bf BIBLIOGRAPHY}}
\vglue 27pt               
\othercite{A}{{D.~Alspach, }{{\it Another method of construction of ${\cal L}_p$ spaces},}
              {unpublished \break manuscript, }{1974.}}
\papercite{B-P}{C.~Bessaga and A.~Pe\char32lczy\'nski}
               {On bases and unconditional convergence of series in Banach spaces}
               {Studia Math.}{17}{1958}{151--164}
\papercite{B-R-S}{J.~Bourgain, H.~P.~Rosenthal, and G.~Schechtman}
                 {An ordinal $L^p$-index for Banach spaces,  
                  with application to complemented subspaces of $L^p$}
                 {Ann.~of Math.~(2)}{114}{1981}{193--228} 
\papercite{Burk}{D.~L.~Burkholder}{Martingale transforms}
                {Ann. Math. Statist.}{37}{1966}{1494--1504}
\othercite{Ch}{{K.~L.~Chung, }{{\it A course in probability theory},}
               {2nd ed., }{Academic Press, }{San Diego, }{1974.}}
\othercite{Db}{{J.~L.~Doob, }{{\it Measure theory},}
               {Graduate Texts in Mathematics, Vol.~143,}\break{Springer-Verlag, }{New York, }{1994.}} 
\papercite{J-M-S-T}{W.~B.~Johnson, B.~Maurey, G.~Schechtman, and L.~Tzafriri}
                   {Symmetric\break structures in Banach spaces}
                   {Mem.~Amer.~Math.~Soc.}{217}{1979}{1--298} 
\papercite{J-O}{W.~B.~Johnson and E.~Odell} 
               {Subspaces of $L_p$ which embed into $\ell_p$}
               {Compositio Math.}{28}{1974}{37--49}
\papercite{J-S-Z}{W.~B.~Johnson, G.~Schechtman, and J.~Zinn}
                 {Best constants in moment\break inequalities for linear combinations 
                  of independent and exchangeable random variables}
                 {Ann.~Probab.}{13}{1985}{234--253} 
\papercite{L-P}{J.~Lindenstrauss and A.~Pe\char32lczy\'nski}
               {Absolutely summing operators in ${\cal L}_p$-spaces and their applications}
               {Studia Math.}{29}{1968}{275--326}
\papercite{L-P 2}{J.~Lindenstrauss and A.~Pe\char32lczy\'nski}
                 {Contributions to the theory of the\break classical Banach spaces}
                 {J.~Funct.~Anal.}{8}{1971}{225--249}
\papercite{L-R}{J.~Lindenstrauss and H.~P.~Rosenthal}{The ${\cal L}_p$ spaces}
               {Israel J.~Math.}{7}{1969}{325--349}
\othercite{L-T}{{J.~Lindenstrauss and L.~Tzafriri, }{{\it Classical Banach spaces},}
                {Lecture Notes in Mathematics, Vol.~338, }{Springer-Verlag, }{New York, }{1973.}} 
\papercite{P}{A.~Pe\char32lczy\'nski}{On the isomorphism of the spaces $m$ and $M$}
             {Bull.~Acad.~Polon. Sci.~S\'er.~Sci.~Math.~Astr.~Phys.}
             {6}{1958}{695--696}
\papercite{R}{H.~P.~Rosenthal}
             {On quasi-complemented subspaces of Banach spaces, 
              with an appendix on compactness of operators from $L^p(\mu)$ to $L^r(\nu)$}
             {J.~Funct.~Anal.}{4}{1969}{176--214}   
\papercite{RI}{H.~P.~Rosenthal}
              {On the subspaces of $L^p$ $(p>2)$  
               spanned by sequences of independent random variables}
              {Israel J.~Math.}{8}{1970}{273--303} 
\vfill\eject 
\mark{}\mark{2}
\othercite{RII}{{H.~P.~Rosenthal,}
                {{\it On the span in $L^p$ of sequences of independent random\break variables} (II),}
                {Proceedings of the Sixth Berkeley Symposium on Mathematical Statistics and Probability, Vol.~2,}
                {University of California Press, }{Berkeley and Los Angeles,}
                {1972, }{pp.~149--167.}} 
\papercite{S}{G.~Schechtman}
             {Examples of ${\cal L}_p$ spaces $(1<p\neq2<\infty)$}
             {Israel J.~Math.}{22}{1975}{138--147}
\papercite{Stn}{E.~M.~Stein}
               {Topics in harmonic analysis related to the Littlewood-Paley\break theory}
               {Ann. of Math. Stud.}{63}{1970}{1--146} 
\othercite{W}{{P.~Wojtaszczyk,}
              {{\it Banach spaces for analysts\/},}
              {Cambridge University Press,}
              {New York, }{1991.}} 
\vfill
\eject
\endgroup
\end